\documentclass[12pt,a4paper,twoside]{report}
\usepackage{mathrsfs}
\usepackage{amsfonts}

\usepackage{amsthm,amssymb,amsmath,amscd,mathrsfs}
\usepackage[all]{xy}
\usepackage{graphicx}
\usepackage{mathptm,pslatex}
\usepackage{fancyhdr}

\begin{document}
\evensidemargin=0.80cm \oddsidemargin=2cm


\thispagestyle{empty} {\huge \bf
\begin{center} Growth Estimates and Integral Representations of
Harmonic and Subharmonic Functions
\end{center}}

\bigskip
\bigskip
\centerline{{\Large by}} {\LARGE\begin{center}{\sc Guoshuang Pan}
\end{center}}

\bigskip
\begin{center}{\Large written under the supervision of}\vspace{0.5cm}
{\Large \\Professor {\sc Guantie Deng}}
\end{center}

 \bigskip
\bigskip
\bigskip

\bigskip
{\Large \rm
\begin{center} A dissertation submitted to the \\ Graduate School in fulfillment of
\\ the requirements for the degree of \\{\bf Doctor of Philosophy}\vspace{0.7cm}
\\ School of Mathematical Sciences \\ Beijing Normal University \\ Beijing, People's Republic
of China
\end{center}}

\bigskip
\bigskip
{\Large
\begin{center} {\rm April\ \ 2009}
\end{center}}


\newpage
\thispagestyle{empty} \quad


\newpage\addcontentsline{toc}{chapter}{\numberline{}{\bf Abstract}}
\begin{center}
  {\bf\Huge Abstract}
\end{center}

    There are ten chapters in this dissertation, which focuses on nine
contents: growth estimates for a class of subharmonic functions in
the half plane; growth estimates for a class of subharmonic
functions in the half space; a generalization of harmonic majorants;
properties of limit for Poisson integral; a lower bound for a class
of harmonic functions in the half space; the Carleman formula of
subharmonic functions in the half space; a generalization of the
Nevanlinna formula for analytic functions in the right half plane;
integral representations of harmonic functions in the half plane;
integral representations of harmonic functions in the half space.

  The outline of the paper is arranged as follows:

  Chapter 1 presents the background, basic notations, some basic definitions,
lemmas, theorems and propositions of the research;

  In Chapter 2, we prove that a class of subharmonic  functions
represented by the modified  kernels have the  growth estimates
 at infinity in the upper half plane  ${\bf C}_{+}$, which
 generalizes the growth properties of analytic
functions and harmonic functions;

   In Chapter 3, a class of subharmonic  functions represented by
 the modified  kernels are proved to have the  growth estimates
  at infinity in the upper half space of ${\bf R}^{n}$, which
  generalizes the growth properties
 of analytic functions and harmonic functions;

  In Chapter 4, we extend the harmonic majorant of a nonnegative and subharmonic function in
${\bf C}_+$ to the harmonic majorant represented by the modified
Poisson kernel and  to the upper half space;

  In Chapter 5, we extend the properties of limit for Poisson integral
in the upper half plane to the properties of limit for Poisson
integral represented by the modified Poisson kernel and  to the
upper half space;

  In Chapter 6, we derive a lower bound for a class of harmonic
functions in the upper half space of ${\bf R}^{n}$ from the upper
bound by using the generalization of the Carleman formula
 for harmonic functions in the upper half
space and the generalization of the Nevanlinna formula for harmonic
functions in the upper half ball;

  In Chapter 7, the object of this chapter is to generalize the
Carleman formula for harmonic functions in the upper half plane to
subharmonic functions in the upper half space;

  In Chapter 8, we generalize the
Nevanlinna formula for analytic functions to the right half plane;

  In Chapter 9, using a modified Poisson kernel in the upper half
plane, we prove that a harmonic  function $u(z)$ in the upper half
plane with its positive part $u^{+}(z)=\max\{u(z),0\}$ satisfying a
slowly growing  condition can be represented by  its integral in the
boundary of the upper half plane, the integral representation is
unique up to the addition of a harmonic polynomial, vanishing in the
boundary of the upper half plane and that its negative part
$u^{-}(z)=\max\{-u(z),0\}$ can be  dominated by a similar slowly
growing  condition, this  improves some classical results about
harmonic functions in the upper half plane;

 In Chapter 10, using a modified Poisson kernel in the upper
half space, we prove that a harmonic  function $u(x)$ in the upper
half space with its positive part $u^{+}(x)=\max\{u(x),0\}$
satisfying a slowly  growing  condition can be represented by  its
integral in the boundary of the upper half space, the integral
representation is unique up to the addition of a harmonic
polynomial, vanishing in the boundary of the upper half space and
that its negative part $u^{-}(x)=\max\{-u(x),0\}$ can be  dominated
by a similar slowly growing  condition, this  improves some
classical results about harmonic functions in the upper half space.\\



\renewcommand{\headrulewidth}{0.7pt}
\tableofcontents

\newpage
\ \thispagestyle{empty}
\newpage

\chapter{Introduction }
\pagenumbering{arabic}

  The present chapter consists of three sections with the first
providing the background for the research project; the second
presenting the basic notations; the third section providing us some
basic definitions, lemmas, theorems and propositions.

\section{Background}

   A complex-valued function $h$ on an open subset $\Omega$ of the
complex plane ${\bf C}$ is called harmonic on $\Omega$ if  $h \in
  C^2(\Omega)$ and
$$
\triangle h\equiv 0
$$
on $\Omega$. Here
$$
\triangle h =\frac{\partial^2 h}{\partial x^2} + \frac{\partial^2
h}{\partial y^2}
$$
is the Laplacian of $h$. We often assume that $\Omega$ is a region
(that is, an open and connected set) even when connectivity is not
needed, and we are mainly interested in the case in which $\Omega$
is a disk or half plane.

  Harmonic functions arise in the study of analytic functions
(we use the terms analytic and holomorphic synonymously). If $f$ is
analytic on a region $\Omega$, then by the Cauchy-Riemann equations,
each of the functions $f$, $\overline{f}$, $\Re f$ is harmonic on
$\Omega$. The theory of harmonic functions is needed in the study of
analytic functions on a disk or half plane.

  Harmonic functions-the solutions of Laplace's equation-play a
crucial role in many areas of mathematics, physics, and engineering.
So it is necessary to extend harmonic functions to ${\bf R}^n$,
where $n$ denotes a fixed positive integer greater than $1$. Let
$\Omega$ be an open, nonempty subset of ${\bf R}^n$. A twice
continuously differentiable, complex-valued function $u$ defined on
$\Omega$ is harmonic on $\Omega$ if
$$
\triangle u\equiv 0,
$$
where $\triangle  =D_1^2+\cdots+D_n^2$ and $D_j^2$ denotes the
second partial derivative with respect to the $j^{th}$ coordinate
variable. The operator $\triangle$ is called the Laplacian, and the
equation $ \triangle u\equiv 0$ is called Laplace's equation.

  We let $x=(x_1,\cdots,x_n)$ denote a typical point in ${\bf R}^n$
and let $|x|=(x_1^2+\cdots+x_n^2)^{1/2}$ denote the Euclidean norm
of $x$.

  The simplest nonconstant harmonic functions are the coordinate
functions; for example, $u(x)=x_1$. A slightly more complex example
is the function on ${\bf R}^3$ defined by
$$
u(x)=x_1^2+x_2^2-2x_3^2+ix_2.
$$
As we will see later, the function
$$
u(x)=|x|^{2-n}
$$
is vital to harmonic function theory when $n>2$; it is obvious that
this function is harmonic on ${\bf R}^n-\{0\}$.

  We can obtain additional examples of harmonic functions by
differentiation, noting that for smooth functions the Laplacian
commutes with any partial derivative. In particular, differentiating
the last example with respect to $x_1$ shows that $x_1|x|^{-n}$ is
harmonic on ${\bf R}^n-\{0\}$ when $n>2$.

  The function $x_1|x|^{-n}$ is harmonic on ${\bf R}^n-\{0\}$ even when $n=2$.
This can be verified directly or by noting that $x_1|x|^{-2}$ is a
partial derivative of $\log|x|$, a harmonic function on ${\bf
R}^2-\{0\}$. The function $\log|x|$ plays the same role when $n=2$
that $|x|^{2-n}$ plays when $n>2$. Notice that $\lim_{x\rightarrow
\infty}\log|x|=\infty$, but $\lim_{x\rightarrow \infty}|x|^{2-n}=0$;
note also $\log|x|$ is neither bounded above nor below, but
$|x|^{2-n}$ is always positive. These facts hint at the contrast
between harmonic function theory in the plane and in higher
dimensions. Another key difference arises from the close connection
between holomorphic and harmonic functions in the plane-a
real-valued function on $\Omega \subset {\bf R}^2$ is harmonic if
and only if it is locally the real part of a holomorphic function.
No comparable result exists in highter dimensions.

   Let $\Omega$ be a region in the complex plane.
A real-valued function $u$ on an open subset $\Omega$ of the
  complex plane ${\bf C}$ is defined to be subharmonic if  $u \in
  C^2(\Omega)$ and
$$
\triangle u\geq 0
$$
on $\Omega$. A broader definition that relaxes the smoothness
assumption and permits $u$ to take the value $-\infty$. Examples of
subharmonic functions include $\log|f|$,
$\log^{+}|f|=\max(\log|f|,0)$ and $|f|^p(0<p<\infty)$, where $f$ is
any analytic function on $\Omega$.

  Elementary properties of subharmonic functions are often one-sided
versions of properties of harmonic functions. For example, a
subharmonic function $u$ on $\Omega$ has a sub-mean value property:
$$
u(a)\leq \frac{1}{2\pi}\int_0^{2\pi}u(a+Re^{i\theta})d\theta.
$$
This property characterizes subharmonic functions.

  One of the most fundamental results in the theory of subharmonic
functions is due to F. Riesz and states that any such function
$u(x)$ can be locally written as the sum of a potential plus a
harmonic function, i.e.
$$
u(x)=p(x)+h(x).
$$
In other words, if $u(x)$ is subharmonic in a domain $D$ in ${\bf
R}^m$, there exists a positive measure $d\mu$, finite on compact
subsets of $D$, and uniquely determined by $u(x)$, such that if $E$
is a compact subset of $D$ and
$$
 p(x)=\left\{\begin{array}{ll}
 \int_E \log|x-\xi| d\mu e_{\xi},  &   \mbox{if }   m=2 ,\\
 -\int_E |x-\xi|^{2-m} d\mu e_{\xi} &  \mbox{if}\   m>2,
 \end{array}\right.
$$
then
$$
h(x)=u(x)-p(x).
$$
is harmonic in the interior of $E$.

  By means of this theorem many of the local properties of
  subharmonic functions can be deduced from those of potentials such
  as $p(x)$. The mass distribution $d\mu$ also plays a fundamental
  role in more delicate questions concerning $u$. Thus for instance
  if $m=2$ and $u(z)=\log|f(z)|$, where $f$ is a regular function of
  the complex variable $z$, then $\mu(E)$ reduces to the number of
  zeros of $f(z)$ on the set $E$. From this point of view the main
  difference between this case and that of a general subharmonic
  function is that in the latter case the "zeros" can have an
  arbitrary mass distribution instead of occurring in units of one.

    In higher dimension we may regard $d\mu$ as the gravitational or
    electric charge, giving rise to the potential $p(x)$. For this
    reason the theory of subharmonic functions is frequently called
    potential theory.

    We now come to a famous problem in harmonic function theory: given
 a continuous function $f$ on $S$, does there exist a continuous
 function function $u$ on $\overline{B}$, with u harmonic on $B$,
 such that $u=f$ on $S$? If so, how do we find $u$? This is Dirichlet
 problem for the ball.

  The  Dirichlet problem of the upper half plane is to find a function
 $u$ satisfying
$$
 u\in C^2({\bf C}_{+}),
$$
$$
 \Delta u=0,   z\in {\bf C}_{+},
$$
$$
 \lim_{z\rightarrow x}u(z)=f(x)\ {\rm nontangentially  \  a.e.}x\in \partial {\bf C}_{+},
$$
where $f$ is a measurable function of ${\bf R}$. The Poisson
integral of the upper half plane is defined by
$$
u(z)=P[f](z)=\int_{{\bf R}}P(z,\xi)f(\xi)d\xi.\eqno{(1.1.1)}
$$
 As we all know, the Poisson integral $P[f]$ exists if
$$
\int_{{\bf R}}\frac{|f(\xi)|}{1+|\xi|^2} d\xi<\infty.
$$
We will generalize these results from harmonic functions to
subharmonic functions.

 Write the subharmonic function
$$
u(z)= v(z)+h(z), \quad z\in {\bf C}_{+},
$$
where $v(z)$ is the harmonic function defined by (1.1.1), $h(z)$ is
defined by
$$
h(z)= \int_{{\bf C}_{+}} G(z,\zeta)d\mu(\zeta)
$$
and $G(z,\zeta)$ is called Green function.

 Hayman \cite{LLST} has proved that the asymptotic behaviour of
subharmonic functions
$$
u(z)= o(|z|),  \quad  {\rm as}  \ |z|\rightarrow\infty
$$
holds everywhere in the upper half plane outside some exceptional
set of disks under the following conditions:
$$
\int_{{\bf R}}\frac{|f(\xi)|}{1+|\xi|^2} d\xi<\infty
$$
and
$$
\int_{{\bf C}_{+}}\frac{\eta}{1+|\zeta|^{2}} d\mu(\zeta)<\infty,
$$
where  $\mu$ is a positive Borel measure and $\zeta=\xi+i\eta.$

  The first aim in this dissertation is to extend the classic results
to the modified Poisson kernel $P_m(z,\xi)$ and the modified Green
function $G_m(z,\zeta)$. That is to say, if
$$
v(z)= \int_{{\bf R}}P_m(z,\xi)f(\xi)d\xi ,
$$
$$
h(z)= \int_{{\bf C}_{+}} G_m(z,\zeta)d\mu(\zeta),
$$
 we will prove that the asymptotic behaviour of subharmonic functions
$$
v(z)= o(y^{1-\alpha}|z|^{m+\alpha}),  \quad  {\rm as}  \
|z|\rightarrow\infty
$$
holds everywhere in the upper half plane outside some exceptional
set of disks under the following conditions:
$$
\int_{{\bf R}}\frac{|f(\xi)|}{1+|\xi|^{2+m}} d\xi<\infty
$$
and
$$
\int_{{\bf C}_{+}}\frac{\eta}{1+|\zeta|^{2+m}} d\mu(\zeta)<\infty.
$$

  Next, we can also conclude that the asymptotic behaviour of subharmonic functions
$$
u(z)= o\big(y^{1-\frac{\alpha}{p}}(\log|z|)^{\frac{1}{q}}
|z|^{\frac{\gamma}{p}+\frac{1}{q}-2+\frac{\alpha}{p}}\big), \quad
{\rm as}  \ |z|\rightarrow\infty
$$
holds everywhere in the upper half plane outside some exceptional
set of disks by replacing the two conditions above into
$$
\int_{{\bf R}}\frac{|f(\xi)|^p}{(1+|\xi|)^\gamma} d\xi<\infty
$$
and
$$
\int_{{\bf C}_{+}}\frac{\eta^p}{(1+|\zeta|)^\gamma}
d\mu(\zeta)<\infty,
$$
where $1< p<\infty,\ \frac{1}{p}+\frac{1}{q}=1$ and $1-p<\gamma
<1+p$.

  On the other hand, we will generalize these results from the upper half plane to
the upper half space.

  The  Dirichlet problem of the upper half space is to find a function
 $u$ satisfying
$$
 u\in C^2(H),
$$
$$
 \Delta u=0,   x\in H,
$$
$$
 \lim_{x\rightarrow x'}u(x)=f(x')\ {\rm nontangentially  \  a.e.}x'\in \partial H,
$$
where $f$ is a measurable function of ${\bf R}^{n-1} $. The Poisson
integral of the upper half space is defined by
$$
u(x)=P[f](x)=\int_{{\bf R}^{n-1}}P(x,y')f(y')dy'.
$$
 As we all know, the Poisson integral $P[f]$ exists if
$$
\int_{{\bf R}^{n-1}}\frac{|f(y')|}{1+|y'|^n} dy'<\infty.
$$

  Write the harmonic function
$$
v(x)= \int_{{\bf R}^{n-1}}P_m(x,y')f(y')dy', \quad x\in H,
$$
Siegel-Talvila \cite{STS} have proved that the asymptotic behaviour
of
$$
v(x)= o(x_n^{1-n}|x|^{m+n}),  \quad  {\rm as}  \
|x|\rightarrow\infty
$$
holds everywhere in the upper half space outside some exceptional
set of balls under the following condition:
$$
\int_{{\bf R}^{n-1}}\frac{|f(y')|}{1+|y'|^{n+m}} dy'<\infty.
$$

  We will generalize these results from harmonic functions to
subharmonic functions, then we will obtain some further results.

 In addition, we also discuss some other problems about harmonic and
subharmonic functions, such as the generalization of harmonic
majorants, properties of limit for Poisson integral, the Carleman
formula and Nevanlinna formula and integral representations.

\section{Basic Notations}

  Let ${\bf C}$  denote the complex plane
with points $z=x+iy$, where $x, y \in {\bf R}$. The boundary and
closure of an open $\Omega$ of ${\bf C}$ are denoted by
$\partial{\Omega}$
 and $\overline{\Omega}$ respectively.
 The upper half plane is the set
 ${\bf C}_+=\{z=x+iy \in {\bf C}:\; y>0\}$, whose boundary is
 $\partial{{\bf C}_+}.$
 We identify ${\bf C}$ with ${\bf R} \times {\bf R}$ and
${\bf R} $ with $ {\bf R} \times \{0\}$, with this convention we
then have $ \partial {{\bf C}_+}={\bf R}$.

   A twice continuously differentiable function $u(z)$
 defined on an open set
 $\Omega $ is harmonic if $\triangle u\equiv 0,$ where
 $\triangle=\frac{\partial^2}{\partial x^2}+\frac{\partial^2}{\partial y^2}$
 is Laplace operator in $z$.
    We  write $B_R$ and $\partial B_R $ for the open ball and the
    circle
 of radius $R$ in ${\bf C}$ centered at the origin
  and  $B_R^{+}=B_R\bigcap {\bf C}_+$ and $\partial B_R^{+}$ for the open upper half
  ball and the upper half circle
 of radius $R$ in ${\bf C}$ centered at the origin.

   Similarly, let ${\bf R}^{n} (n\geq3)$  denote the  $n$-dimensional
Euclidean space with points
$x=(x_1,x_2,\cdots,x_{n-1},x_{n})=(x',x_n)$, where $x' \in {\bf
R}^{n-1}$ and $x_{n} \in {\bf R}$.  The boundary and closure of an
open  $\Omega$ of ${\bf R}^{n}$ are denoted by $\partial{\Omega}$
 and $\overline{\Omega}$ respectively.
 The upper half space is the set
 $H=\{x=(x',x_n)\in {\bf R}^{n}:\; x_n>0\}$, whose boundary is
 $\partial{H}$ .
 We identify ${\bf R}^{n}$ with ${\bf R}^{n-1}\times {\bf R}$ and
${\bf R}^{n-1}$ with $ {\bf R}^{n-1}\times \{0\}$,
  with this convention we then have $ \partial {H}={\bf R}^{n-1}$,
  writing typical points $x,\ y \in {\bf R}^{n}$ as $x=(x',x_n),\
y=(y',y_n),$  where $x'=(x_1,x_2,\cdots,x_{n-1}),\
y'=(y_1,y_2,\cdots y_{n-1}) \in {\bf R}^{n-1}$ and putting
$$
x\cdot y=\sum_{j=1}^{n}x_jy_j=x'\cdot y'+x_ny_n,\ \ |x|=\sqrt{x\cdot
x},\ \ |x'|=\sqrt{x'\cdot x'}.
$$
where $|x|$ is the Euclidean norm.

  A twice continuously differentiable function $u(x)$ defined on an open set
 $\Omega $ is harmonic if $\triangle_x u\equiv 0,$ where
 $\triangle_x=\frac{\partial^2}{\partial x_1^2}+\frac{\partial^2}{\partial x_2^2}
 +\cdots+\frac{\partial^2}{\partial x_n^2}$ is Laplace operator in $x$.
 The upper half space $H$ is the set
 $H=\{x=(x',x_n)\in {\bf R}^{n}:\; x_n>0,\}.$
    We  write $B_R$ and $\partial B_R $ for the open ball and the sphere
 of radius $R$ in ${\bf R}^{n}$ centered at the origin
  and  $B_R^{+}=B_R\bigcap H$ and $\partial B_R^{+}$ for the open upper half
  ball and the upper half sphere
 of radius $R$ in ${\bf R}^{n}$ centered at the origin.
In the sense of Lebesgue measure $dx'=dx_1\cdots dx_{n-1},\;
 dx=dx'dx_n$ and  let $\sigma$ denote $(n-1)$-dimensional
 surface-area measure.

   Throughout the dissertation, let $A$
denote various positive constants independent of the variables in
question.

\section{Preliminary  Results}

  In this section, we will introduce some definitions, lemmas, theorems and
propositions that will be used in the following chapters.

\vspace{0.2cm}
 \noindent
{\bf Definition $A_1$\cite{RR} } {\it  Let $X$ be a metric space. A
function $u: X\rightarrow [-\infty,\infty)$ is said to be upper
semicontinuous or usc if
$$
\{x:\ u(x)<a\}
$$
is an open set in $X$ for each real number $a$, or, equivalently, if
for every $x \in X$,
$$
\limsup_{y\rightarrow x}u(y)\leq u(x).
$$
}

\vspace{0.2cm}
 \noindent
{\bf Definition $A_2$\cite{RR} } {\it  Let $\Omega$ be an open set
in the complex plane. We say that a function $u: \Omega\rightarrow
[-\infty,\infty)$ is subharmonic on $\Omega$ if

\noindent {\rm (1)}\
 $u$ is usc on $\Omega$;

 \noindent{\rm (2)}\
 for every open set $A$ with compact closure $\overline{A} \subseteq \Omega$
and every continuous function $h: \overline{A}\rightarrow
(-\infty,\infty)$ whose restriction to $A$ is harmonic, if $u\leq h$
on $\partial A$, then $u\leq h$ on $A$.}

\vspace{0.2cm}
 \noindent
{\bf Definition $A_3$\cite{RR} } {\it  Let $u$ be subharmonic on a
region $\Omega$, $u\neq -\infty$, and let $h$ be harmonic on
$\Omega$. We say that $h$ is a harmonic majorant for $u$ if $h\geq
u$ on $\Omega$.  We say that $h$ is a least harmonic majorant for
$u$ if

\noindent {\rm (1)}\
 $h$ is a harmonic majorant for $u$;

 \noindent
 {\rm (2)}\
if $f$ is any harmonic majorant for $u$ in $\Omega$, then $h\leq f$
on $\Omega$. }

\vspace{0.2cm}
 \noindent
{\bf Lemma $A_1$\cite{RR} }  {\it Let $g(x)$ be a nonnegative and
nondecreasing function on $[0,1)$. Let $p(x)$ be any nonnegative
measurable function on $(0,1)$ such that
$$
0<\int_0^a p(t) dt<\infty
$$
for every $a \in (0,1)$ and
$$
\int_0^1 p(t) dt = \infty.
$$
Then
$$
\lim_{x \uparrow 1}g(x) = \sup_{0<r<1} \frac{\int_0^1 g(t)p(\lambda
t)dt}{\int_0^1 p(\lambda t)dt}.
$$}

\vspace{0.2cm}
 \noindent
{\bf Lemma $A_2$\cite{RR} }  {\it Let $V(z)$ be nonnegative and
harmonic on $\Pi$ and have a continuous extension to
$\overline{\Pi}=\{z: \Im z\geq 0\}$. Then
$$
V(z)=cy+\frac{y}{\pi}\int_{-\infty}^{\infty}\frac{V(t)}{(t-x)^2+y^2}dt,\quad
y>0,
$$
where $c$ is given by
$$
c=\lim_{y \rightarrow \infty} \frac{V(iy)}{y}.
$$
}

\vspace{0.2cm}
 \noindent
{\bf Lemma $A_3$\cite{HK} }  {\it The polynomials $a_v(x,\xi)$ are
harmonic in $x$ for fixed $\xi$, and continuous in $x$, $\xi$
jointly for $|\xi|\neq 0$. If $|x|=\rho$, $|\xi|=r>0$, we have the
sharp inequality
$$
|a_v(x,\xi)|\leq \frac{b_{v,m}\rho^v}{r^{m+v-2}},
$$
where $b_{v,m}=1/v$ if $m=2$, $v\geq 1$;
$$
b_{v,m}=(v+m-3)(v+m-4)\cdots (v+1)/(m-1)!, \quad m\geq 3, v\geq 0.
$$
 }

\vspace{0.2cm}
 \noindent
{\bf Lemma $A_4$\cite{HK} }  {\it If $|\xi|=r>0$, then
$K_q(x,\xi)-K(x,\xi)$ is harmonic in ${\bf R}^m$. We set $|x|=\rho$
and have the following estimates
$$
|K_q(x,\xi)|\leq 4^{m+q}\frac{\rho^{q+1}}{r^{m+q-1}} if \rho \leq
\frac{1}{2}r.
$$
If $q=0$, $m=2$, we have
$$
K_0(x,\xi)\leq \log(1+\rho/r),
$$
while in all other cases
$$
K_q(x,\xi,m)\leq
4^{m+q}\frac{\rho^q}{r^{m+q-2}}\inf\{1,\frac{\rho}{r}\}.
$$
}

\vspace{0.2cm}
 \noindent
{\bf Theorem $A_1$\cite{RR} } {\it  For every continuous
complex-valued function $f$ on $\Gamma$ there is a unique continuous
function $h$ on $\overline{D}=D \bigcup \Gamma$ such that the
restriction of $h$ to $\Gamma$ is $f$ and the restriction of $h$ to
$D$ is harmonic. The function $h$ is given on $D$ by
$$
h(z)=\int_{\Gamma}P(z,e^{it})f(e^{it}) d\sigma(e^{it}),\quad z\in D.
$$
}

\vspace{0.2cm}
 \noindent
{\bf Theorem $A_2$\cite{RR} }(Mean value property)  {\it If $h$ is
harmonic on a region $\Omega$ and $\overline{D}(a,R)\subseteq
\Omega$, then
$$
h(a)=\frac{1}{2\pi}\int_0^{2\pi}h(a+Re^{it})dt.
$$
}

\vspace{0.2cm}
 \noindent
{\bf Theorem $A_3$\cite{RR} } (Maximum principle)  {\it A
real-valued harmonic function $h$ on an open connected set $\Omega$
cannot attain either a maximum or a mimimum value in $\Omega$
without reducing to a constant.}

\vspace{0.2cm}
 \noindent
{\bf Theorem $A_4$\cite{RR} } {\it A continuous function $h$ on a
region $\Omega$ has the mean value property if and only if  $h$ is
harmonic on $\Omega$.}

\vspace{0.2cm}
 \noindent
{\bf Theorem $A_5$\cite{RR} } {\it Every nonnegative harmonic
function $h$ on the unit dist $D$ has a respresentation
$$
h(z)=\int_{\Gamma}P(z,e^{it}) d\mu(e^{it}),\quad z\in D,
$$
where $\mu$ is a finite  nonnegative measure on $\Gamma$.}

\vspace{0.2cm}
 \noindent
{\bf Theorem $A_6$\cite{RR} } (Herglotz and Riesz Representation
Theorem) {\it Let $f$ be a analytic and satisfy $\Re f \geq 0$ on
$D$. Then
$$
f(z)=\int_{\Gamma}\frac{e^{i\theta}+z}{e^{i\theta}-z}
d\mu(e^{i\theta})+ic,\quad z\in D,
$$
for some finite nonnegative Borel measure $\mu$ on $\Gamma$ and some
real constant $c$.}

\vspace{0.2cm}
 \noindent
{\bf Theorem $A_7$\cite{RR} } (Stieltjes Inversion Formula)  {\it
Let $\mu$ be a complex Borel measure on $\Gamma$, and let on the
unit dist $D$ has a respresentation
$$
h(z)=\int_{\Gamma}P(z,e^{it}) d\mu(e^{it}),\quad z\in D.
$$
Let $\gamma=\{e^{it}: a<t<b\}$ be an open arc on the unit circle
with endpoints $\alpha=e^{ia}$ and $\beta=e^{ib}$, $0<b-a<2\pi$.
Then
$$
\lim_{r\uparrow 1} \frac{1}{2\pi}\int_a^b h(r e^{i\theta})
d\theta=\mu(\gamma)+\frac{1}{2}\mu(\{\alpha\})+\frac{1}{2}\mu(\{\beta\}).
$$
}

\vspace{0.2cm}
 \noindent
{\bf Theorem $A_8$\cite{RR} } {\it  Let $u$ be usc on a region
$\Omega$ in the complex plane. The following are equivalent:

\noindent
{\rm (1)}\ $u$ is subharmonic on $\Omega$;

 \noindent
 {\rm (2)}\
 for each $a \in \Omega$ and all sufficiently small $R>0$, if $p$
is a polynomial such that  $u\leq \Re p$ on $\partial D(a,R)$, then
$u\leq \Re p$ on $D(a,R)$;

\noindent
 {\rm (3)}\   for each $a \in \Omega$ and all sufficiently small $R>0$,
$$
u(a)\leq \frac{1}{2\pi}\int_0^{2\pi}u(a+Re^{i\theta})d\theta.
$$
In this case, the properties expressed in  {\rm (2)} and  {\rm (3)}
hold for all disks $D(a,R)$ such that $\overline{D}(a,R)\subseteq
\Omega$. }

\vspace{0.2cm}
 \noindent
{\bf Theorem $A_9$\cite{RR} } (Maximum principle)  {\it Assume that
$u$ is subharmonic on a region $\Omega$. If there is a point  $z_0
\in \Omega$ suth that $u(z_0)\geq u(z)$ for all $z \in \Omega$, then
$u\equiv const.$ in $\Omega$.}

\vspace{0.2cm}
 \noindent
{\bf Theorem $A_{10}$\cite{RR} }  {\it Assume that $u$ is
subharmonic in $D(a,R)$ and $u\neq -\infty$. If $0<r_1<r_2<R$, then
$$
-\infty < \frac{1}{2\pi}\int_0^{2\pi}u(a+r_1 e^{i\theta})d\theta
\leq \frac{1}{2\pi}\int_0^{2\pi}u(a+r_2 e^{i\theta})d\theta.
$$
Moreover, whether $u(a)$ is finite or $-\infty$,
$$
\lim_{r \downarrow 0} \frac{1}{2\pi}\int_0^{2\pi}u(a+r
e^{i\theta})d\theta=u(a).
$$
}

\vspace{0.2cm}
 \noindent
{\bf Theorem $A_{11}$\cite{RR} }  {\it Let $u$ be subharmonic in the
unit disk $D$,
 $u\neq -\infty$. There exists a harmonic majorant for $u$ if and
 only if
$$
\sup_{0<r<1} \frac{1}{2\pi}\int_0^{2\pi}u(a+r e^{it})dt<\infty.
$$
In this case there is a least harmonic majorant $h$ for $u$, and $h$
is given by
$$
h(z)=\lim_{r \uparrow 1} \frac{1}{2\pi}\int_0^{2\pi}P(z/r,e^{it})u(r
e^{it})dt
$$
for all $z \in D$.}

 \vspace{0.2cm}
 \noindent
{\bf Theorem $A_{12}$\cite{RR} } (Poisson Representation) {\it Every
nonnegative harmonic function $V(z)$ on $\Pi$ has a representation
$$
V(z)=cy+\frac{y}{\pi}\int_{-\infty}^{\infty}\frac{d\mu(t)}{(t-x)^2+y^2},\quad
y>0,
$$
where $c\geq 0$ and $\mu$ is a nonnegative Borel measure on
$(-\infty,\infty)$ such that
$$
\int_{-\infty}^{\infty}\frac{d\mu(t)}{1+t^2}<\infty.
$$}

 \vspace{0.2cm}
 \noindent
{\bf Theorem $A_{13}$\cite{RR} } ( Nevanlinna Representation) {\it
Every holomorphic function $F(z)$ such that $\Im F(z)\geq 0$ for $z
\in \Pi$ has a representation
$$
F(z)=b+cz+\frac{1}{\pi}\int_{-\infty}^{\infty}\bigg[\frac{1}{t-z}-\frac{t}{1+t^2}\bigg]d\mu(t),\quad
y>0.
$$
where $b=\overline{b}$, $c\geq 0$, and $\mu$ is a nonnegative Borel
measure on $(-\infty,\infty)$ which satisfies
$$
\int_{-\infty}^{\infty}\frac{d\mu(t)}{1+t^2}<\infty.
$$
}

\vspace{0.2cm}
 \noindent
{\bf Theorem $A_{14}$\cite{RR} } (Stieltjes Inversion Formula) {\it
Let $V(z)$ be given by
$$
V(z)=cy+\frac{y}{\pi}\int_{-\infty}^{\infty}\frac{d\mu(t)}{(t-x)^2+y^2},\quad
y>0,
$$
where $c\geq 0$ and $\mu$ is a nonnegative Borel measure satisfying
$$
\int_{-\infty}^{\infty}\frac{d\mu(t)}{1+t^2}<\infty.
$$
If $-\infty<a<b<\infty$, then
$$
\lim_{y \downarrow 0} \int_a^b V(x+iy)
dx=\mu((a,b))+\frac{1}{2}\mu(\{a\})+\frac{1}{2}\mu(\{b\}).
$$
}

\vspace{0.2cm}
 \noindent
{\bf Theorem $A_{15}$\cite{RR} } (Fatou's Theorem)  {\it Let
$$
V(z)=\frac{y}{\pi}\int_{-\infty}^{\infty}\frac{d\mu(t)}{(t-x)^2+y^2},\quad
y>0,
$$
where $\mu$ is a nonnegative Borel measure on $(-\infty,\infty)$
satisfying
$$
\int_{-\infty}^{\infty}\frac{d\mu(t)}{1+t^2}<\infty.
$$
If $d\mu=Fdx+d\mu_s$ is the Lebesgue decomposition of $\mu$, then
$$
\lim_{z\rightarrow x}V(z) =F(x)
$$
nontangentially a.e. on $(-\infty,\infty)$. }

 \vspace{0.2cm}
 \noindent
{\bf Theorem $A_{16}$\cite{RR} } {\it  Let $F$ be holomorphic on
$D_{+}(0,R)$ for some $R>0$, and suppose $F\neq 0$. Then $F \in
N^{+}(D_{+}(0,R))$ if and only if
$$
\log|F(z)|\leq \frac{R^2-|z|^2}{\pi}\int_0^{\pi}\frac{2yR\sin t}{|R
e^{it}-z|^2|R e^{-it}-z|^2}K(R e^{it})dt
$$
$$
+\frac{y}{\pi}\int_{-R}^R\bigg(\frac{1}{|t-z|^2|}-\frac{R^2}{|R^2-tz|^2}\bigg)K(t)dt
$$
for all $z \in D_{+}(0,R)$ and some real-valued Borel function
$K(\zeta)$ on $\Gamma_{+}(0,R)$ such that
$$
\int_0^{\pi}|K(R e^{it})|\sin t dt+\int_{-R}^R
|K(t)|(R^2-t^2)dt<\infty.
$$
}

\vspace{0.2cm}
 \noindent
{\bf Theorem $A_{17}$\cite{ABR} } (Mean value property)  {\it If $u$
is harmonic on $\overline{B}(a,r)$, then $u$ equals the average of
$u$ over $\partial B(a,r)$. More precisely,
$$
u(a)=\int_S u(a+r\zeta)d\sigma (\zeta).
$$
}

\vspace{0.2cm}
 \noindent
{\bf  Theorem $A_{18}$\cite{ABR} } (Solution of the Dirichlet
problem for the ball) {\it Suppose $f$ is continuous on $S$. Define
$u$ on $\overline{B}$ by
$$
 u(x)=\left\{\begin{array}{ll}
 P[f](x)  &   \mbox{if }   x \in B  ,\\
 f(x)&
\mbox{if}\   x \in S.
 \end{array}\right.\eqno{(1.6)}
$$
Then $u$ is continuous on $\overline{B}$ and harmonic on $B$. }

\vspace{0.2cm}
 \noindent
{\bf Theorem $A_{19}$\cite{ABR} } {\it If $u$ is a continuous
function on $\overline{B}$ that is harmonic on $B$, then $u=P[u|_S]$
on $B$.}

\vspace{0.2cm}
 \noindent
{\bf  Theorem $A_{20}$\cite{ABR} } (Solution of the Dirichlet
problem for $H$) {\it Suppose $f$ is continuous and bounded on ${\bf
R}^{n-1}$. Define $u$ on $\overline{H}$ by
$$
 u(z)=\left\{\begin{array}{ll}
 P_H[f](z)  &   \mbox{if }   x \in H  ,\\
 f(z)&
\mbox{if}\   x \in {\bf R}^{n-1}.
 \end{array}\right.\eqno{(1.6)}
$$
Then $u$ is continuous on $\overline{H}$ and harmonic on $H$.
Moreover,
$$
|u|\leq ||f||_{\infty}
$$
on $\overline{H}$. }

\vspace{0.2cm}
 \noindent
{\bf Theorem $A_{21}$\cite{ABR} } {\it Suppose $u$ is a continuous
bounded function on $\overline{H}$ that is harmonic on $H$. Then $u$
is the Poisson integral of its boundary values. More precisely,
$$
u=P_H[u|_{{\bf R}^{n-1}}]
$$
on $H$.}

\vspace{0.2cm}
 \noindent
{\bf Theorem $A_{22}$\cite{R} } (The Schwarz reflection principle)
{\it Suppose $L$ is a segment of the real axis, $\Omega^{+}$ is a
region in $\Pi^{+}$, and every $t \in L$ is the center of an open
disc $D_t$ such that $\Pi^{+} \bigcap D_t$ lies in $\Omega^{+}$. Let
$\Omega^{-}$ be the reflection of $\Omega^{+}$:
$$
\Omega^{-}=\{z:\ \overline{z}\ in \ \Omega^{+}\}.
$$
Suppose $f=u+iv$ is holomorphic in $\Omega^{+}$, and
$$
\lim_{n\rightarrow\infty}v(z_n)=0
$$
for every sequence $\{z_n\}$ in $\Omega^{+}$ which converges to a
point of $L$.

  Then there is a function $F$, holomorphic in $\Omega^{+} \bigcup L \bigcup \Omega^{-}$,
such that $F(z)=f(z)$ in $\Omega^{+}$; this $F$ satisfies the
relation
$$
F(\overline{z})=\overline{F(z)} \quad (z \in \Omega^{+}\cup L\cup
\Omega^{-}).
$$
}

 \vspace{0.2cm}
 \noindent
{\bf Theorem $A_{23}$\cite{ABR} }  {\it  If $u$ is positive and
harmonic on $H$, then there exists a positive Borel measure $\mu$ on
${\bf R}^{n-1}$ and a nonnegative constant $c$ such that
$$
u(x,y)=cy+\int_{{\bf R}^{n-1}} P_H(z,t) d\mu(t)
$$
for all $(x,y) \in H$.}

 \vspace{0.2cm}
 \noindent
{\bf Theorem $A_{24}$\cite{LLST} } (Hayman) {\it  Let
$$
v(z)=\int\int_{{\bf C}_{+}}
\log\bigg|\frac{\zeta-\overline{z}}{\zeta-z}\bigg|d\mu(\zeta)
+\frac{y}{\pi}\int_{-\infty}^{\infty}\frac{d\nu(t)}{(t-x)^2+y^2},
$$
where $d\mu(\zeta)$ and $d\nu(t)$ are nonnegative Borel measures
such that
$$
\int\int_{{\bf C}_{+}} \frac{\Im
\zeta}{1+|\zeta|^2}d\mu(\zeta)<\infty, \quad \int_{{\bf
R}}\frac{d\nu(t)}{1+t^2}<\infty.
$$
The asymptotic relation
$$
v(z)= o(|z|),  \quad    |z|\rightarrow\infty
$$
holds everywhere in ${\bf C}_{+}$ outside some exceptional set of
disks of finite view.}

\vspace{0.2cm}
 \noindent
{\bf Theorem $A_{25}$\cite{HK} } (Green's Theorem) {\it Suppose that
$D$ is an admissible domain with boundary $S$ in ${\bf R}^{m}$ and
that $u \in C^1$ and $v \in C^2$ in $\overline D$. Then
$$
\int_S u(x)\frac{\partial v}{\partial n}d\sigma =-\int_D
\bigg\{\sum_v \frac{\partial u}{\partial x_v}\frac{\partial
v}{\partial x_v}+u\nabla_2 v\bigg\}dx,
$$
where
$$
\nabla_2=\sum_{v=1}^m \frac{\partial^2}{\partial x_v^2}
$$
is Laplace's operator. Hence if $u,v \in C^2$ in $\overline D$ we
have
$$
\int_S \bigg(u\frac{\partial v}{\partial n}-v\frac{\partial
u}{\partial n}\bigg)d\sigma =\int_D (v\nabla^2 u-u\nabla^2 v)dx.
$$
Here $\partial/\partial n$ denotes differentiation along the inward
normal into $D$.}

\vspace{0.2cm}
 \noindent
{\bf Theorem $A_{26}$\cite{HK} } {\it If $D=D(0,R)$ and $\xi$ a
point of $D$, $\xi'=\xi R^2 |\xi|^{-2}$, and if for $m=2$
$$
g(x,\xi,D)=\log \frac{|x-\xi'||\xi|}{|x-\xi|R},\quad \xi \neq 0;
\quad g(x,0,D)=\log \frac{R}{|x|};
$$
while for $m>2$
$$
g(x,\xi,D)=|x-\xi|^{2-m}-\{|\xi||x-\xi'|/R\}^{2-m},\quad \xi \neq 0;
$$
$$
g(x,0,D)=|x|^{2-m}-R^{2-m};
$$
then $g(x,\xi,D)$ is a (classical) Green's function of $D$.}

\vspace{0.2cm}
 \noindent
{\bf Theorem $A_{27}$\cite{HK} } (Poisson's Integral)  {\it If $u$
is harmonic in $D(x_0,R)$ and continuous in $C(x_0,R)$ then for $\xi
\in D(x_0,R)$ we have
$$
u(\xi)=\frac{1}{c_m}\int_{S(x_0,R)}
\frac{R^2-|\xi-x_0|^2}{R|x-\xi|^m}u(x)d\sigma_x,
$$
where $d\sigma_x$ denotes an element of surface area of $S(x_0,R)$
and $c_m=2\pi^{m/2}/\Gamma(m/2)$. }

\vspace{0.2cm}
 \noindent
{\bf Theorem $A_{28}$\cite{HK} }  {\it Suppose that $u(x)$ is s.h.
in $C(x_0,R)$. Then for $\xi \in D(x_0,R)$ we have
$$
u(\xi)\leq \int_{S(x_0,R)} u(x)K(x,\xi)d\sigma_x,
$$
where $K(x,\xi)$ is the Poisson kernel given by
$$
K(x,\xi)=\frac{1}{c_m}\frac{R^2-|\xi-x_0|^2}{R|x-\xi|^m}
$$
and $d\sigma_x$ denotes an element of surface area of $S(x_0,R)$. }

\vspace{0.2cm}
 \noindent
{\bf Theorem $A_{29}$\cite{HK} } (Riesz's Theorem) {\it Suppose that
$u(x)$ is s.h. and not identically $-\infty$, in a domain $D$ in
${\bf R}^m$. Then there exists a unique Borel-measure $\mu$ in $D$
such that for any compact subset $E$ of $D$ we have
$$
u(x)=\int_E u(x)K(x-\xi)d\mu e_{\xi}+h(x),
$$
where $h(x)$ is harmonic in the interior of $E$. }

\vspace{0.2cm}
 \noindent
{\bf Theorem $A_{30}$\cite{HK} }  {\it Suppose that $D$ is a bounded
regular domain in ${\bf R}^m$ whose frontier $F$ has zero
m-dimensional Lebesgue measure, and that $u(x)$ is s.h. and not
identically $-\infty$ on $D \bigcup F$. Then we have for $x \in D$
$$
u(x)=\int_F u(\xi)d\omega (x,e_{\xi})-\int_D g(x,\xi,D)d\mu e_{\xi},
$$
where $\omega (x,e)$ is the harmonic measure of $e$ at $x$,
$g(x,\xi,D)$ is the Green's function of $D$ and $d\mu$ is the Riesz
measure of $u$ in $D$. }

\vspace{0.2cm}
 \noindent
{\bf Theorem $A_{31}$\cite{HK} }  (Weierstrass' Theorem) {\it
Suppose that $\mu$ is a Borel measure in ${\bf R}^m$, let $n(t)$ be
the measure of $D(0,t)$ and let $q(t)$ be a positive integer-valued
increasing function of $t$, continuous on the right, and so chosen
that
$$
\int_1^{\infty} (\frac{t_0}{t})^{q(t)+m-1}dn(t)<\infty
$$
for all positive $t_0$. Then there exists functions $u(x)$, s.h. in
${\bf R}^m$ and with Riesz measure $\mu$, and all such functions
take the form
$$
u(x)=\int_{|\xi|<1} K(x-\xi)d\mu e_{\xi}+\int_{|\xi|\geq1}
K_{q(|\xi|)}(x-\xi)d\mu e_{\xi}+v(x),
$$
where $v(x)$ is harmonic in ${\bf R}^m$. The second integral
converges absolutely near $\infty$ and uniformly for $|x|\leq \rho$
and any fixed positive $\rho$. }

\vspace{0.2cm}
 \noindent
{\bf Proposition $A_1$\cite{ABR} } (Polar coordinates formula)  {\it
The polar coordinates formula for integration on ${\bf R}^n$ states
that for a Borel measurable, integrable function $f$ on  ${\bf
R}^n$,
$$
\frac{1}{nV(B)}\int_{{\bf R}^n} f dV= \int_0^{\infty} r^{n-1} \int_S
f(r\zeta) d\sigma(\zeta) dr,
$$
the constant arises from the normalization of $\sigma$. }

\vspace{0.2cm}
 \noindent
{\bf  Proposition $A_2$\cite{ABR} } {\it Let $\zeta \in S$. Then
$P(\cdot,\zeta)$ is harmonic on ${\bf R}^n -\{\zeta\}$. }

\vspace{0.2cm}
 \noindent
{\bf Proposition $A_3$\cite{ABR} } {\it The Poisson kernel has the
following properties:

\noindent {\rm (a)}\ $P(x,\zeta)>0$ for all $x \in B$ and all $\zeta
\in S$;

 \noindent
 {\rm (b)}\
 $\int_S P(x,\zeta) d\sigma (\zeta)=1$ for all $x \in B$;

\noindent
 {\rm (c)}\
for every $\eta \in S$ and every $\delta >0$,
$$
\int_{|\zeta-\eta|>\delta} P(x,\zeta) d\sigma (\zeta)\rightarrow 0,
\quad {\rm as}  \ x\rightarrow \eta.
$$
 }

\vspace{0.2cm} \noindent {\bf Proposition $A_4$\cite{LLST} } (R.
Nevanlinna's formula for a half-disk) {\it Let $f(z)$ be a
meromorphic function in the half-disk $\overline{D}_R^{+}$, $a_n$ be
its zeros and $b_n$ its poles. Then we obtain

\begin{eqnarray*}
\log|f(z)|
&=& \frac{R^2-|z|^2}{2\pi}\int_0^{\pi}\bigg(\frac{1}{|R
e^{i\theta}-z|^2}-\frac{1}{|R e^{i\theta}-\overline{z}|^2}\bigg)\log
|f(R
e^{i\theta})|d\theta \\
& & +\frac{\Re z}{\pi}\int_{-R}^R
\bigg(\frac{1}{|t-z|^2|}-\frac{R^2}{|R^2-tz|^2}\bigg)\log
|f(t)|dt \\
& & +\sum_{a_n \in D_R^{+}} \log
\bigg|\frac{z-a_n}{z-\overline{a}_n}\cdot\frac{R^2-a_nz}{R^2-\overline{a}_nz}\bigg|-
\sum_{b_n \in D_R^{+}} \log
\bigg|\frac{z-b_n}{z-\overline{b}_n}\cdot\frac{R^2-b_nz}{R^2-\overline{b}_nz}\bigg|.
\end{eqnarray*}
}

\chapter{Growth Estimates for a Class of Subharmonic Functions in the Half Plane}

\section{Introduction and Basic Notations}

  For $z\in{\bf C}\backslash\{0\}$, let \cite{HN}
$$
E(z)=(2\pi)^{-1}\log|z|,
$$
where $|z|$ is the Euclidean norm. We know that $E$ is locally
integrable in ${\bf C} $.

  First, we define the Green function $G(z,\zeta)$ for the upper half plane
 ${\bf C}_{+}$ by \cite{HN}
$$
 G(z,\zeta)=E(z-\zeta)-E(z-\overline{\zeta}),
 \qquad z,\zeta\in\overline{{\bf C}_{+}} ,\  z\neq \zeta, \eqno{(2.1.1)}
$$
then we define the Poisson kernel $P(z,\xi) $ when $z\in {\bf
C}_{+}$ and $\xi\in
\partial {\bf C}_{+} $ by
$$
 P(z,\xi)=-\frac{\partial G(z,\zeta)}{\partial
 \eta}\bigg|_{\eta=0}=\frac{y}{\pi|z-\xi|^2}.\eqno{(2.1.2)}
$$

  The  Dirichlet problem of the upper half plane is to find a function
 $u$ satisfying
$$
 u\in C^2({\bf C}_{+}),\eqno{(2.1.3)}
$$
$$
 \Delta u=0,   z\in {\bf C}_{+}, \eqno{(2.1.4)}
$$
$$
 \lim_{z\rightarrow x}u(z)=f(x)\ {\rm nontangentially  \  a.e.}x\in \partial {\bf C}_{+},
\eqno{(2.1.5)}
$$
where $f$ is a measurable function of ${\bf R}$. The Poisson
integral of the upper half plane is defined by
$$
v(z)=P[f](z)=\int_{{\bf R}}P(z,\xi)f(\xi)d\xi, \eqno{(2.1.6)}
$$
where $P(z,\xi)$ is defined by (2.1.2).

 As we all know, the Poisson integral $P[f]$ exists if
$$
\int_{{\bf R}}\frac{|f(\xi)|}{1+|\xi|^2} d\xi<\infty. \eqno{(2.1.7)}
$$(see \cite{ABR}, \cite{F} and \cite{MS})In this chapter, we replace the condition into
$$
\int_{{\bf R}}\frac{|f(\xi)|^p}{(1+|\xi|)^\gamma}
d\xi<\infty,\eqno{(2.1.8)}
$$
where $1\leq p<\infty$ and $\gamma$ is a real number, then we can
get the asymptotic behaviour of harmonic functions.

  Next, we will generalize these results to subharmonic functions.

\section{Preliminary Lemma}

\vspace{0.3cm}

Let $\mu$ be a positive Borel measure  in ${\bf C},\ \beta\geq0$,
the maximal function $M(d\mu)(z)$ of order $\beta$ is defined by
$$
M(d\mu)(z)=\sup_{ 0<r<\infty}\frac{\mu(B(z,r))}{r^\beta},
$$
then the maximal function $M(d\mu)(z):{\bf C} \rightarrow
[0,\infty)$ is lower semicontinuous, hence measurable. To see this,
for any $\lambda >0 $, let $D(\lambda)=\{z\in{\bf
C}:M(d\mu)(z)>\lambda\}$. Fix $z \in D(\lambda)$, then there exists
 $r>0$ such that $\mu(B(z,r))>tr^\beta$ for some $t>\lambda$, and
there exists $ \delta>0$ satisfying
$(r+\delta)^\beta<\frac{tr^\beta}{\lambda}$. If $|\zeta-z|<\delta$,
then $B(\zeta,r+\delta)\supset B(z,r)$, therefore
$\mu(B(\zeta,r+\delta))\geq
tr^\beta=t(\frac{r}{r+\delta})^\beta(r+\delta)^\beta>\lambda(r+\delta)^\beta$.
Thus $B(z,\delta)\subset D(\lambda)$. This proves that $D(\lambda)$
is open for each $\lambda>0$.

 In order to obtain the results, we
need the lemma below:

\vspace{0.2cm}
 \noindent
{\bf Lemma 2.2.1} {\it Let $\mu$ be a positive Borel measure in
${\bf C},\ \beta\geq0,\ \mu({\bf C})<\infty,$ for any $\lambda \geq
5^{\beta} \mu({\bf C})$, set
$$
E(\lambda)=\{z\in{\bf C}:|z|\geq2,M(d\mu)(z) >
\frac{\lambda}{|z|^{\beta}}\},
$$
then there exists $z_j\in E(\lambda),\ \ \rho_j> 0,\ j=1,2,\cdots$,
such that
$$
E(\lambda) \subset \bigcup_{j=1}^\infty B(z_j,\rho_j) \eqno{(2.2.1)}
$$
and
$$
\sum _{j=1}^{\infty}\frac{\rho_j^{\beta}}{|z_j|^{\beta}}\leq
\frac{3\mu({\bf C})5^{\beta}}{\lambda} .\eqno{(2.2.2)}
$$
}

 Proof: Let $E_k(\lambda)=\{z\in E(\lambda):2^k\leq |z|<2^{k+1}\}$,
then for any $z \in E_k(\lambda),$ there exists $r(z)>0$, such that
$\mu(B(z,r(z))) >\lambda \big(\frac{r(z)}{|z|}\big)^{\beta} $,
therefore $r(z)\leq 2^{k-1}$.
 Since $E_k(\lambda)$ can be covered
by
 the union of a family of balls $\{B(z,r(z)):z\in E_k(\lambda) \}$,
 by the Vitali Lemma \cite{SS}, there exists $\Lambda_k \subset E_k(\lambda)$,
$\Lambda_k$ is at most countable, such that $\{B(z,r(z)):z\in
\Lambda_k \}$ are disjoint and
$$
E_k(\lambda) \subset
 \cup_{z\in \Lambda_k} B(z,5r(z)),
$$
so
$$
E(\lambda)=\cup_{k=1}^\infty E_k(\lambda) \subset \cup_{k=1}^\infty
\cup_{z\in \Lambda_k} B(z,5r(z)).
$$

  On the other hand, note that $ \cup_{z\in \Lambda_k} B(z,r(z)) \subset \{z:2^{k-1}\leq
|z|<2^{k+2}\} $, so that
$$
 \sum_{z \in \Lambda_k}\frac{(5r(z))^{\beta}}{|z|^{\beta}}
\leq 5^\beta\sum_{z\in\Lambda_k}\frac{\mu(B(z,r(z)))}{\lambda} \leq
\frac{5^\beta}{\lambda} \mu\{z:2^{k-1}\leq |z|<2^{k+2}\}.
$$
Hence we obtain
$$
 \sum _{k=1}^{\infty}\sum
_{z \in \Lambda_k}\frac{(5r(z))^{\beta}}{|z|^{\beta}}
 \leq
 \sum _{k=1}^{\infty}\frac{5^\beta}{\lambda} \mu\{z:2^{k-1}\leq |z|<2^{k+2}\}
 \leq
\frac{3\mu({\bf C})5^{\beta}}{\lambda}.
$$
  Rearrange $ \{z:z \in \Lambda_k,k=1,2,\cdots\} $ and $
\{5r(z):z \in \Lambda_k,k=1,2,\cdots\}
 $, we get $\{z_j\}$
and $\{\rho_j\}$
 such that (2.2.1) and
(2.2.2) hold.

\section{$p=1$}

 \subsection*{ 1. Introduction and Main Theorems}

\vspace{0.3cm} In this section, we will consider measurable
functions $f$ in ${\bf R}$  satisfying
$$
\int_{{\bf R}}\frac{|f(\xi)|}{1+|\xi|^{2+m}}
d\xi<\infty,\eqno{(2.3.1)}
$$
where $m$ is a nonnegative integer. This is just (2.1.8) when $p=1$
and $\gamma=2+m$. To obtain a solution of Dirichlet problem for the
boundary date $f$, as in
 \cite{STS}, \cite{STU},  \cite{Y} and \cite{MS}, we use the following
  modified functions defined by
$$
 E_m(z-\zeta)=\left\{\begin{array}{ll}
 E(z-\zeta)  &   \mbox{when }   |\zeta|\leq 1,  \\
 E(z-\zeta) -
 \frac{1}{2\pi}\Re(\log\zeta-\sum_{k=1}^{m-1}\frac{z^k}{k\zeta^{k}})
  & \mbox{when}\   |\zeta|> 1.
 \end{array}\right.
$$
Then we can define the modified Green function $G_m(z,\zeta)$ and
the modified Poisson
 kernel $P_m(z,\xi)$ by (see \cite{MS} and \cite{DI})
$$
 G_m(z,\zeta)=E_{m+1}(z-\zeta)-E_{m+1}(z-\overline{\zeta}),
\qquad z,\zeta\in\overline{{\bf C}_{+}}, \ z\neq
 \zeta;\eqno{(2.3.2)}
$$
$$
 P_m(z,\xi)=\left\{\begin{array}{ll}
 P(z,\xi)  &   \mbox{when }   |\xi|\leq 1  ,\\
 P(z,\xi) - \frac{1}{\pi}\Im\sum_{k=0}^{m}\frac{z^k}{\xi^{1+k}}&
\mbox{when}\   |\xi|> 1,
 \end{array}\right.\eqno{(2.3.3)}
$$
where
 $ z=x+iy, \zeta=\xi+i\eta$.

   Hayman \cite{LLST} has proved the
  following result:

\vspace{0.2cm}
 \noindent
{\bf Theorem B } {\it Let $f$ be a measurable function in ${\bf R}$
satisfying (2.1.7) and $\mu$ be a positive Borel measure satisfying
$$
\int_{{\bf C}_{+}}\frac{\eta}{1+|\zeta|^{2}} d\mu(\zeta)<\infty.
$$
Write the subharmonic function
$$
u(z)= v(z)+h(z), \quad z\in {\bf C}_{+},
$$
where $v(z)$ is the harmonic function defined by (2.1.6), $h(z)$ is
defined by
$$
h(z)= \int_{{\bf C}_{+}} G(z,\zeta)d\mu(\zeta)
$$
and $G(z,\zeta)$ is defined by (2.1.1). Then there exists $z_j\in
{\bf C}_{+},\ \rho_j>0,$ such that
$$
\sum _{j=1}^{\infty}\frac{\rho_j}{|z_j|}<\infty
$$
holds and
$$
u(z)= o(|z|),  \quad  {\rm as}  \ |z|\rightarrow\infty
 $$
holds in ${\bf C}_{+}-G$, where $ G=\bigcup_{j=1}^\infty
B(z_j,\rho_j)$. }

  Our aim in this section is to establish the following theorems.

\vspace{0.2cm}
 \noindent
{\bf Theorem 2.3.1} {\it Let $f$ be a measurable function in ${\bf
R}$ satisfying (2.3.1), and $0< \alpha\leq 2$. Let $v(z)$ be the
harmonic function defined by
$$
v(z)= \int_{{\bf R}}P_m(z,\xi)f(\xi)d\xi, \quad z\in {\bf C}_{+},
 \eqno{(2.3.4)}
$$
where $P_m(z,\xi)$ is defined by (2.3.3). Then there exists $z_j\in
{\bf C}_{+},\ \rho_j>0,$ such that
$$
\sum
_{j=1}^{\infty}\frac{\rho_j^{2-\alpha}}{|z_j|^{2-\alpha}}<\infty
\eqno{(2.3.5)}
$$
holds and
$$
v(z)= o(y^{1-\alpha}|z|^{m+\alpha}),  \quad  {\rm as}  \
|z|\rightarrow\infty   \eqno{(2.3.6)}
 $$
holds in ${\bf C}_{+}-G$, where $ G=\bigcup_{j=1}^\infty
B(z_j,\rho_j)$. }

\vspace{0.2cm}
 \noindent
 {\bf Remark 2.3.1 } {\it If $\alpha=2$, then (2.3.5)
 is a finite sum, the set $G$ is the union of finite disks,
 so (2.3.6) holds in ${\bf C}_{+}$. }

  Next, we will generalize Theorem 2.3.1 to subharmonic functions.

\vspace{0.2cm}
 \noindent
{\bf Theorem 2.3.2 } {\it Let $f$ be a measurable function in ${\bf
R}$ satisfying (2.3.1) and $\mu$ be a positive  Borel  measure
satisfying
$$
\int_{{\bf C}_{+}}\frac{\eta}{1+|\zeta|^{2+m}} d\mu(\zeta)<\infty.
$$
Write the subharmonic function
$$
u(z)= v(z)+h(z), \quad z\in {\bf C}_{+},
$$
where $v(z)$ is the harmonic function defined by (2.3.4), $h(z)$ is
defined by
$$
h(z)= \int_{{\bf C}_{+}} G_m(z,\zeta)d\mu(\zeta)
$$
and $G_m(z,\zeta)$ is defined by (2.3.2). Then there exists $z_j\in
{\bf C}_{+},\ \rho_j>0,$ such that (2.3.5) holds and
$$
u(z)= o(y^{1-\alpha}|z|^{m+\alpha}),  \quad  {\rm as}  \
|z|\rightarrow\infty  \eqno{(2.3.7)}
 $$
holds in ${\bf C}_{+}-G$, where $ G=\bigcup_{j=1}^\infty
B(z_j,\rho_j)$ and $0< \alpha<2$. }

 \vspace{0.2cm}
 \noindent
 {\bf Remark 2.3.2 } {\it If $\alpha=1, m=0$, this is just the result of
Hayman, so our result (2.3.7) is the generalization of Theorem B. }

\vspace{0.4cm}

\subsection*{2.   Main Lemma}

  In order to obtain the results, we need the following lemma:

\vspace{0.2cm}
 \noindent
{\bf Lemma 2.3.1 } {\it The following inequalities hold:\\
{\rm (1)}\ If $|\xi| > 1$, then $|P_m(z,\zeta)-P(z,\zeta)|
           \leq\sum_{k=0}^{m-1}\frac{2^ky|z|^k}{\pi|\xi|^{2+k}} ;$\\
{\rm (2)}\ If $|\xi-z| > 3|z|$, then
            $|P_m(z,\zeta)|
            \leq \frac{2^{m+1}y|z|^m}{\pi|\xi|^{m+2}}; $\\
{\rm (3)}\ If $|\xi| > 1$, then
            $|G_m(z,\zeta)-G(z,\zeta)|\leq
            \frac{1}{\pi}\sum_{k=1}^{m}\frac{ky\eta|z|^{k-1}}{|\zeta|^{1+k}};$\\
{\rm (4)}\ If $|\xi-z| > 3|z|$, then $|G_m(z,\zeta)|\leq
            \frac{1}{\pi}\sum_{k=m+1}^{\infty}\frac{ky\eta|z|^{k-1}}{|\zeta|^{1+k}}.$\\
}

\vspace{0.4cm}

\subsection*{3.   Proof of Theorems }

 \emph{Proof of Theorem 2.3.1}

Define the measure $dm(\xi)$ and the kernel $K(z,\xi)$ by
$$
dm(\xi)=\frac{|f(\xi)|}{1+|\xi|^{2+m}} d\xi ,\ \ K(z,\xi)=
P_m(z,\xi)(1+|\xi|^{2+m}).
$$
  For any $\varepsilon >0$, there exists $R_\varepsilon >2$, such that
$$
\int_{|\xi|\geq
R_\varepsilon}dm(\xi)\leq\frac{\varepsilon}{5^{2-\alpha}}.
$$
For every Lebesgue measurable set $E \subset {\bf R}$ , the measure
$m^{(\varepsilon)}$ defined by $m^{(\varepsilon)}(E)
=m(E\cap\{x\in{\bf R}:|x|\geq R_\varepsilon\}) $ satisfies
$m^{(\varepsilon)}({\bf R})\leq\frac{\varepsilon}{5^{2-\alpha}}$,
write
\begin{eqnarray*}
&v_1(z)& =\int_{|\xi-z| \leq 3|z|} P(z,\xi)(1+|\xi|^{2+m})
dm^{(\varepsilon)}(\xi),\\
&v_2(z)&=\int_{|\xi-z| \leq 3|z|}
(P_m(z,\xi)-P(z,\xi))(1+|\xi|^{2+m})
dm^{(\varepsilon)}(\xi), \\
&v_3(z)&=\int_{|\xi-z| > 3|z|} K(z,\xi)dm^{(\varepsilon)}(\xi), \\
&v_4(z)&=\int_{1<|\xi|<R_\varepsilon}K(z,\xi) dm(\xi), \\
&v_5(z)&=\int_{|\xi|\leq1}K(z,\xi) dm(\xi), \\
\end{eqnarray*}
then
$$
|v(z)| \leq |v_1(z)|+|v_2(z)|+|v_3(z)|+|v_4(z)|+|v_5(z)|.
\eqno{(2.3.8)}
$$
Let $ E_1(\lambda)=\{z\in{\bf C}:|z|\geq2,\exists \ t>0, s.t.
m^{(\varepsilon)}(B(z,t)\cap{\bf R}
)>\lambda(\frac{t}{|z|})^{2-\alpha}\}$, therefore, if $ |z|\geq
2R_\varepsilon$ and $z \notin E_1(\lambda)
 $, then
$$
\forall t>0,\ m^{(\varepsilon)}(B(z,t)\cap{\bf R} )\leq\lambda
\bigg(\frac{t}{|z|}\bigg)^{2-\alpha}.
$$
  So we have
\begin{eqnarray*}
|v_1(z)|
&\leq& \int_{y\leq|\xi-z| \leq
3|z|}\frac{y}{\pi|z-\xi|^2}2|\xi|^{2+m} dm^{(\varepsilon)}(\xi) \\
&\leq& \int_{y\leq|\xi-z| \leq
3|z|}\frac{2y}{\pi|z-\xi|^2}(4|z|)^{2+m}
dm^{(\varepsilon)}(\xi) \\
&=& \frac{2^{2m+5}}{\pi}y|z|^{2+m}\int_{y\leq|\xi-z| \leq
3|z|}\frac{1}{|z-\xi|^2} dm^{(\varepsilon)}(\xi) \\
%
&=& \frac{2^{2m+5}}{\pi}y|z|^{m+2}\int_{y}^{3|z|} \frac{1}{t^2}
dm_z^{(\varepsilon)}(t),
\end{eqnarray*}
where  $m_z^{(\varepsilon)}(t)=\int_{|\xi-z| \leq t}
dm^{(\varepsilon)}(\xi)$, since for $z \notin E_1(\lambda)$,
\begin{eqnarray*}
\int_{y}^{3|z|} \frac{1}{t^2} dm_z^{(\varepsilon)}(t)
&\leq&  \frac{m_z^{(\varepsilon)}(3|z|)}{(3|z|)^2}+2
\int_{y}^{3|z|} \frac{m_z^{(\varepsilon)}(t)}{t^{3}} dt \\
&\leq& \frac{\lambda}{3^\alpha |z|^2}+2 \int_{y}^{3|z|}
\frac{\lambda\frac{t^{2-\alpha}}{|z|^{2-\alpha}}}{t^{3}} dt \\
&\leq& \frac{\lambda}{ |z|^2}\bigg(\frac{1}{3^\alpha}+
\frac{2}{\alpha}\frac{|z|^\alpha}{y^\alpha}\bigg),
\end{eqnarray*}
so that
\begin{eqnarray*}
|v_1(z)|
&\leq& \frac{2^{2m+5}}{\pi}y|z|^{m+2} \frac{\lambda}{|z|^2}
\bigg(\frac{1}{3^\alpha}+\frac{2}{\alpha}\frac{|z|^\alpha}{y^\alpha}\bigg) \\
&\leq& \frac{2^{2m+5}}{\pi}
\bigg(\frac{1}{3^\alpha}+\frac{2}{\alpha}\bigg)\lambda
y^{1-\alpha}|z|^{m+\alpha}. \hspace{52mm}(2.3.9)
\end{eqnarray*}

   By (1) of Lemma 2.3.1, we obtain
\begin{eqnarray*}
|v_2(z)|
&\leq& \int_{y\leq|\xi-z| \leq 3|z|}
\sum_{k=0}^{m-1}\frac{2^ky|z|^k}{\pi}\frac{2|\xi|^{2+m}}{|\xi|^{2+k}} dm^{(\varepsilon)}(\xi) \\
&\leq&\int_{y\leq|\xi-z| \leq 3|z|}
\sum_{k=0}^{m-1}\frac{2^{k+1}y|z|^k}{\pi}(4|z|)^{m-k} dm^{(\varepsilon)}(\xi) \\
&\leq& \frac{2^{2m+1}}{\pi}\sum_{k=0}^{m-1}\frac{1}{2^k}\frac{1}{5^{2-\alpha}}\varepsilon y|z|^m \\
&\leq& \frac{4^{m-1+\alpha}}{\pi}\varepsilon y|z|^m. \hspace{75mm}
(2.3.10)
\end{eqnarray*}

  By (2) of Lemma 2.3.1, we see that \cite{HK}
\begin{eqnarray*}
|v_3(z)|
&\leq& \int_{|\xi-z| > 3|z|}
\frac{2^{m+1}y|z|^m}{\pi|\xi|^{m+2}}2|\xi|^{2+m} dm^{(\varepsilon)}(\xi) \\
&=& \int_{|\xi-z| > 3|z|}
\frac{2^{m+2}y|z|^m}{\pi} dm^{(\varepsilon)}(\xi) \\
&\leq& \frac{2^{m+2}}{\pi}\frac{\varepsilon}{5^{2-\alpha}}
y|z|^m \\
&\leq& \frac{2^{m-2+2\alpha}}{\pi}\varepsilon y|z|^m . \hspace{74mm}
(2.3.11)
\end{eqnarray*}

  Write
\begin{eqnarray*}
v_4(z)
&=& \int_{1<|\xi|<R_\varepsilon}\big[P(z,\xi)+(P_m(z,\zeta)-P(z,\zeta))\big](1+|\xi|^{2+m}) dm(\xi) \\
&=& v_{41}(z)-v_{42}(z),
\end{eqnarray*}
then
\begin{eqnarray*}
|v_{41}(z)|
&\leq& \int_{1<|\xi|<R_\varepsilon}\frac{y}{\pi|z-\xi|^2}
2|\xi|^{2+m} dm(\xi) \\
&\leq& \frac{2R_\varepsilon^{2+m}y}{\pi}
\int_{1<|\xi|<R_\varepsilon}\frac{1}{\big(\frac{|z|}{2}\big)^2}
 dm(\xi) \\
&\leq& \frac{2^3R_\varepsilon^{2+m}m({\bf R})}{\pi}\frac{y}{|z|^2}.
\hspace{84mm} (2.3.12)
\end{eqnarray*}

  Moreover, by (1) of Lemma 2.3.1, we obtain
\begin{eqnarray*}
|v_{42}(z)|
&\leq& \int_{1<|\xi|<R_\varepsilon}
\sum_{k=0}^{m-1}\frac{2^ky|z|^k}{\pi|\xi|^{2+k}}\cdot 2|\xi|^{2+m} dm(\xi) \\
&\leq& \sum_{k=0}^{m-1}\frac{2^{k+1}}{\pi}y|z|^kR_\varepsilon^{m-k}
m({\bf
R}) \\
&\leq& \frac{2^{m+1}R_\varepsilon^m m({\bf R})}{\pi}
y|z|^{m-1}.\hspace{85mm} (2.3.13)
\end{eqnarray*}

  In case $|\xi|\leq 1$, note that
$$
K(z,\xi)=P_m(z,\xi)(1+|\xi|^{2+m}) \leq\frac{2y}{\pi|z-\xi|^2},
$$
so that
$$
|v_5(z)|\leq \int_{|\xi|\leq1}\frac{2y}{\pi|z-\xi|^2} dm(\xi) \leq
\int_{|\xi|\leq1}\frac{2y}{\pi \big(\frac{|z|}{2}\big)^2} dm(\xi)
\leq \frac{2^3m({\bf R})}{\pi}\frac{y}{|z|^2}. \eqno{(2.3.14)}
$$

  Thus, by collecting (2.3.8), (2.3.9), (2.3.10), (2.3.11), (2.3.12), (2.3.13) and
(2.3.14), there exists a positive constant $A$ independent of
$\varepsilon$, such that if $ |z|\geq 2R_\varepsilon$ and $\ z
\notin E_1(\varepsilon)$, we have
$$
|v(z)|\leq A\varepsilon y^{1-\alpha}|z|^{m+\alpha}.
$$

 Let $\mu_\varepsilon$ be a measure in ${\bf C}$ defined by
$ \mu_\varepsilon(E)= m^{(\varepsilon)}(E\cap{\bf R})$ for every
measurable set $E$ in ${\bf C}$. Take
$\varepsilon=\varepsilon_p=\frac{1}{2^{p+2}}, p=1,2,3,\cdots$, then
there exists a sequence $ \{R_p\}$: $1=R_0<R_1<R_2<\cdots$ such that
$$
\mu_{\varepsilon_p}({\bf C})=\int_{|\xi|\geq
R_p}dm(\xi)<\frac{\varepsilon_p}{5^{2-\alpha}}.
$$
Take $\lambda=3\cdot5^{2-\alpha}\cdot2^p\mu_{\varepsilon_p}({\bf
C})$ in Lemma 2.2.1, then $\exists \ z_{j,p}$ and $ \rho_{j,p}$,
where $R_{p-1}\leq |z_{j,p}|<R_p$ such that
$$
\sum
_{j=1}^{\infty}\bigg(\frac{\rho_{j,p}}{|z_{j,p}|}\bigg)^{2-\alpha}
\leq \frac{1}{2^{p}}.
$$
So if $R_{p-1}\leq |z|<R_p$ and $z\notin G_p=\cup_{j=1}^\infty
B(z_{j,p},\rho_{j,p})$, we have
$$
|v(z)|\leq A\varepsilon_py^{1-\alpha}|z|^{m+\alpha},
$$
thereby
$$
\sum _{p=1}^{\infty}
\sum_{j=1}^{\infty}\bigg(\frac{\rho_{j,p}}{|z_{j,p}|}\bigg)^{2-\alpha}
\leq \sum _{p=1}^{\infty}\frac{1}{2^{p}}=1<\infty.
$$
Set $ G=\cup_{p=1}^\infty G_p$, thus Theorem 2.3.1 holds.

 \emph{Proof of Theorem 2.3.2}

  Define the measure $dn(\zeta)$ and the kernel $L(z,\zeta)$ by
$$
dn(\zeta)=\frac{\eta d\mu(\zeta)}{1+|\zeta|^{2+m}},\ \
L(z,\zeta)=G_m(z,\zeta)\frac{1+|\zeta|^{2+m}}{\eta},
$$
then the function $h(z)$ can be written as
$$
h(z)=\int_{{\bf C}_{+}} L(z,\zeta) dn(\zeta).
$$

  For any $\varepsilon >0$, there exists $R_\varepsilon >2$, such that
$$
\int_{|\zeta|\geq
R_\varepsilon}dn(\zeta)<\frac{\varepsilon}{5^{2-\alpha}}.
$$
For every Lebesgue measurable set $E \subset {\bf C}$ , the measure
$n^{(\varepsilon)}$ defined by $n^{(\varepsilon)}(E)
=n(E\cap\{\zeta\in {\bf C}_{+}:|\zeta|\geq R_\varepsilon\}) $
satisfies $n^{(\varepsilon)}({\bf
C}_{+})\leq\frac{\varepsilon}{5^{2-\alpha}}$, write
\begin{eqnarray*}
&h_1(z)&
=\int_{|\zeta-z|\leq\frac{y}{2}}G(z,\zeta)\frac{1+|\zeta|^{2+m}}{\eta}
dn^{(\varepsilon)}(\zeta), \\
&h_2(z)&=\int_{\frac{y}{2}<|\zeta-z|\leq3|z|}G(z,\zeta)\frac{1+|\zeta|^{2+m}}{\eta}
dn^{(\varepsilon)}(\zeta), \\
&h_3(z)&=\int_{|\zeta-z|\leq3|z|}[G_m(z,\zeta)-G(z,\zeta)]\frac{1+|\zeta|^{2+m}}{\eta}
dn^{(\varepsilon)}(\zeta), \\
&h_4(z)&=\int_{|\zeta-z|>3|z|}L(z,\zeta)
dn^{(\varepsilon)}(\zeta), \\
&h_5(z)&=\int_{1<|\zeta|<R_\varepsilon}L(z,\zeta) dn(\zeta), \\
&h_6(z)&=\int_{|\zeta|\leq1}L(z,\zeta) dn(\zeta), \\
\end{eqnarray*}
then
$$
h(z)=h_1(z)+h_2(z)+h_3(z)+h_4(z)+h_5(z)+h_6(z). \eqno{(2.3.15)}
$$
Let $ E_2(\lambda)=\{z\in{\bf C}:|z|\geq2,\exists \ t>0, s.t.
n^{(\varepsilon)}(B(z,t)\cap {\bf C}_{+}
)>\lambda(\frac{t}{|z|})^{2-\alpha}\}, $ therefore, if $ |z|\geq
2R_\varepsilon$ and $z\notin E_2(\lambda)
 $, then
$$\forall t>0, \ n^{(\varepsilon)}(B(z,t)\cap {\bf C}_{+}
)\leq\lambda \bigg(\frac{t}{|z|}\bigg)^{2-\alpha}.$$

  So we have
\begin{eqnarray*}
|h_1(z)|
&\leq&
\int_{|\zeta-z|\leq\frac{y}{2}}\frac{1}{2\pi}\log\bigg|\frac{\zeta-\overline{z}}{\zeta-z}\bigg|
\frac{1+|\zeta|^{2+m}}{\eta}
dn^{(\varepsilon)}(\zeta) \\
&\leq& \int_{|\zeta-z|\leq\frac{y}{2}}\frac{1}{2\pi}
\log\frac{3y}{|\zeta-z|}\frac{2|\zeta|^{2+m}}{\frac{y}{2}}
dn^{(\varepsilon)}(\zeta) \\
&\leq& \frac{2\times (3/2)^{2+m}
}{\pi}\frac{|z|^{2+m}}{y}\int_{|\zeta-z|\leq\frac{y}{2}}\log\frac{3y}{|\zeta-z|}
dn^{(\varepsilon)}(\zeta) \\
&=& \frac{2\times (3/2)^{2+m}
}{\pi}\frac{|z|^{2+m}}{y}\int_0^\frac{y}{2}
\log\frac{3y}{t} dn_z^{(\varepsilon)}(t)\\
&\leq& \frac{2\times (3/2)^{2+m}
}{\pi}\bigg[\frac{\log6}{2^{2-\alpha}}+
\frac{1}{(2-\alpha)2^{2-\alpha}}\bigg]\lambda
y^{1-\alpha}|z|^{m+\alpha}, \hspace{15mm} (2.3.16)
\end{eqnarray*}
where $ n_z^{(\varepsilon)}(t)=\int_{|\zeta-z| \leq t}
dn^{(\varepsilon)}(\zeta)$.\\

  Note that
$$
|G(z,\zeta)|=|E(z-\zeta)-E(z-\overline{\zeta})|\leq
\frac{y\eta}{\pi|z-\zeta|^2}, \eqno{(2.3.17)}
$$
then by (2.3.17), we have
\begin{eqnarray*}
|h_2(z)|
&\leq&
\int_{\frac{y}{2}<|\zeta-z|\leq3|z|}\frac{y\eta}{\pi|z-\zeta|^2}
\frac{2|\zeta|^{2+m}}{\eta}
dn^{(\varepsilon)}(\zeta) \\
&\leq&
\frac{2}{\pi}y(4|z|)^{2+m}\int_{\frac{y}{2}<|\zeta-z|\leq3|z|}
\frac{1}{|z-\zeta|^2} dn^{(\varepsilon)}(\zeta) \\
&=& \frac{2^{2m+5}}{\pi}y|z|^{2+m}\int_\frac{y}{2}^{3|z|}
\frac{1}{t^2} dn_z^{(\varepsilon)}(t)\\
&\leq& \frac{2^{2m+5}}{\pi}y|z|^{2+m}\frac{\lambda}{
|z|^2}\bigg(\frac{1}{3^\alpha}+
\frac{2^{\alpha+1}}{\alpha}\frac{|z|^\alpha}{y^\alpha}\bigg) \\
&\leq& \frac{2^{2m+5}}{\pi}\bigg(\frac{1}{3^\alpha}+
\frac{2^{\alpha+1}}{\alpha}\bigg)\lambda y^{1-\alpha}|z|^{m+\alpha}.
\hspace{40mm} (2.3.18)
\end{eqnarray*}

   By (3) of Lemma 2.3.1, we obtain
\begin{eqnarray*}
|h_3(z)|
&\leq&
\int_{|\zeta-z|\leq3|z|}\frac{1}{\pi}\sum_{k=1}^{m}\frac{ky\eta|z|^{k-1}}{|\zeta|^{1+k}}
\frac{2|\zeta|^{2+m}}{\eta} dn^{(\varepsilon)}(\zeta) \\
&=&
\int_{|\zeta-z|\leq3|z|}\frac{2}{\pi}\sum_{k=1}^{m}ky|z|^{k-1}|\zeta|^{m-k+1}
 dn^{(\varepsilon)}(\zeta) \\
&\leq&
\int_{|\zeta-z|\leq3|z|}\frac{2}{\pi}\sum_{k=1}^{m}ky|z|^{k-1}(4|z|)^{m-k+1}
 dn^{(\varepsilon)}(\zeta) \\
&\leq&
\frac{2}{\pi}\sum_{k=1}^{m}4^{m-k+1}k\frac{1}{5^{2-\alpha}}\varepsilon y|z|^m \\
&\leq& \frac{2^{2m+2\alpha+1}}{9\pi}\varepsilon y|z|^m.
\hspace{73mm} (2.3.19)
\end{eqnarray*}

  By (4) of Lemma 2.3.1, we see that
\begin{eqnarray*}
|h_4(z)|
&\leq&
\int_{|\zeta-z|>3|z|}\frac{1}{\pi}\sum_{k=m+1}^{\infty}\frac{ky\eta|z|^{k-1}}{|\zeta|^{1+k}}|
\frac{2|\zeta|^{2+m}}{\eta} dn^{(\varepsilon)}(\zeta) \\
&=&
\int_{|\zeta-z|>3|z|}\frac{2}{\pi}\sum_{k=m+1}^{\infty}\frac{ky|z|^{k-1}}{|\zeta|^{k-(m+1)}}
 dn^{(\varepsilon)}(\zeta) \\
&\leq& \int_{|\zeta-z|>3|z|}\frac{2}{\pi}\sum_{k=m+1}^{\infty}ky
\frac{|z|^{k-1}}{(2|z|)^{k-m-1}}
dn^{(\varepsilon)}(\zeta) \\
&\leq& \frac{2^{m+2}}{\pi}\sum_{k=m+1}^{\infty}\frac{k}{2^{k}}
 \frac{1}{5^{2-\alpha}}\varepsilon y|z|^m \\
&\leq& \frac{4^{\alpha-1}(m+2)}{\pi}\varepsilon y|z|^m .
\hspace{68mm} (2.3.20)
\end{eqnarray*}

  Write
\begin{eqnarray*}
h_5(z)
&=&
\int_{1<|\zeta|<R_\varepsilon}[G(z,\zeta)+(G_m(z,\zeta)-G(z,\zeta))]\frac{1+|\zeta|^{2+m}}{\eta} dn(\zeta) \\
&=& h_{51}(z)+h_{52}(z),
\end{eqnarray*}
then we obtain by (2.3.17)
\begin{eqnarray*}
|h_{51}(z)|
&\leq& \int_{1<|\zeta|<R_\varepsilon}
\frac{y\eta}{\pi|z-\zeta|^2}\frac{2|\zeta|^{2+m}}{\eta} dn(\zeta) \\
&\leq& \int_{1<|\zeta|<R_\varepsilon}\frac{2}{\pi}
\frac{yR_\varepsilon^{2+m}}{|z-\zeta|^2} dn(\zeta) \\
&\leq&
\frac{2R_\varepsilon^{2+m}}{\pi}y\int_{1<|\zeta|<R_\varepsilon}
\frac{1}{(\frac{|z|}{2})^2} dn(\zeta) \\
&\leq& \frac{2^3R_\varepsilon^{2+m}n({\bf
C}_{+})}{\pi}\frac{y}{|z|^2}. \hspace{85mm} (2.3.21)
\end{eqnarray*}

  Moreover, by (3) of Lemma 2.3.1, we obtain
\begin{eqnarray*}
|h_{52}(z)|
&\leq&
\int_{1<|\zeta|<R_\varepsilon}\frac{1}{\pi}\sum_{k=1}^{m}\frac{ky\eta|z|^{k-1}}{|\zeta|^{1+k}}
\frac{2|\zeta|^{2+m}}{\eta} dn(\zeta) \\
&=&
\int_{1<|\zeta|<R_\varepsilon}\frac{2}{\pi}\sum_{k=1}^{m}ky|z|^{k-1}|\zeta|^{m-k+1}
dn(\zeta) \\
&\leq&
\int_{1<|\zeta|<R_\varepsilon}\frac{2}{\pi}\sum_{k=1}^{m}ky|z|^{k-1}R_\varepsilon^{m-k+1}
 dn(\zeta) \\
&\leq& \frac{m(m+1)R_\varepsilon^{m}n({\bf C}_{+}) }{\pi}
y|z|^{m-1}. \hspace{85mm} (2.3.22)
\end{eqnarray*}

  In case $|\zeta|\leq 1$, by (2.3.17), we have
$$
|L(z,\zeta)|\leq \frac{y\eta}{\pi|z-\zeta|^2}\frac{2}{\eta}
=\frac{2y}{\pi|z-\zeta|^2},
$$
so that
\begin{eqnarray*}
|h_6(z)| \leq \int_{|\zeta|\leq1}\frac{2y}{\pi|z-\zeta|^2} dn(\zeta)
\leq \int_{|\zeta|\leq1}\frac{2y}{\pi(\frac{|z|}{2})^2}
dn(\zeta)\leq \frac{2^3n({\bf C}_{+})}{\pi}\frac{y}{|z|^2}.
\hspace{85mm} (2.3.23)
\end{eqnarray*}

  Thus, by collecting (2.3.15), (2.3.16), (2.3.18), (2.3.19),
(2.3.20), (2.3.21), (2.3.22) and (2.3.23), there exists a positive
constant $A$ independent of $\varepsilon$, such that if $ |z|\geq
2R_\varepsilon$ and $\  z \notin E_2(\varepsilon)$, we have
$$
 |h(z)|\leq A\varepsilon y^{1-\alpha}|z|^{m+\alpha}.
$$

  Similarly, if $z\notin G$, we have
$$
h(z)= o(y^{1-\alpha}|z|^{m+\alpha}), \quad  {\rm as} \
|z|\rightarrow\infty. \eqno{(2.3.24)}
$$

By (2.3.6) and (2.3.24), we obtain that
$$
u(z)=v(z)+h(z)= o(y^{1-\alpha}|z|^{m+\alpha}), \quad  {\rm as} \
|z|\rightarrow\infty
$$
holds in  ${\bf C}_{+}-G$.

\section{$p>1$(General Kernel)}

 \subsection*{ 1. Introduction and Main Theorems}

\vspace{0.3cm}

In this section, we will consider measurable functions $f$ in ${\bf
R}$ satisfying
$$
\int_{{\bf R}}\frac{|f(\xi)|^p}{(1+|\xi|)^\gamma}
d\xi<\infty,\eqno{(2.4.1)}
$$
where $\gamma$ is defined as in Theorem 2.4.1.

  In order to describe the asymptotic behaviour of subharmonic functions
in the upper half plane (see \cite{ML}, \cite{MLJ1}, and
\cite{MLJ2}),
 we establish the following theorems.

\vspace{0.2cm}
 \noindent
{\bf Theorem 2.4.1} {\it Let $1\leq p<\infty,\
\frac{1}{p}+\frac{1}{q}=1$ and
$$
1-p<\gamma <1+p  \quad  {\rm in \  case}  \ p>1;
$$
$$
 0<\gamma \leq 2  \quad  \quad \quad \quad {\rm in \  case}  \ p=1.
$$
If $f$ is a measurable function in ${\bf R}$ satisfying (2.4.1) and
$v(z)$ is the harmonic function defined by (2.1.6), then there
exists $z_j\in {\bf C}_{+},\ \rho_j>0,$ such that
$$
\sum
_{j=1}^{\infty}\frac{\rho_j^{2p-\alpha}}{|z_j|^{2p-\alpha}}<\infty
\eqno{(2.4.2)}
$$
holds and
$$
v(z)=
o(y^{1-\frac{\alpha}{p}}|z|^{\frac{\gamma}{p}+\frac{1}{q}-2+\frac{\alpha}{p}}),
\quad  {\rm as}  \ |z|\rightarrow\infty   \eqno{(2.4.3)}
$$
holds in ${\bf C}_{+}-G$, where $ G=\bigcup_{j=1}^\infty
B(z_j,\rho_j)$ and $0< \alpha\leq 2p$. }

\vspace{0.2cm}
 \noindent
 {\bf Remark 2.4.1 } {\it If $\gamma=1-p$, $p>1$, then
$$
v(z)=
o(y^{1-\frac{\alpha}{p}}(\log|z|)^{\frac{1}{q}}|z|^{\frac{\gamma}{p}+\frac{1}{q}-2+\frac{\alpha}{p}}),
\quad  {\rm as}  \ |z|\rightarrow\infty
$$
holds in ${\bf C}_{+}-G$. }

  Next, we will generalize Theorem 2.4.1 to subharmonic functions.

\vspace{0.2cm}
 \noindent
{\bf Theorem 2.4.2 } {\it Let $p$ and $\gamma$ be as in Theorem
2.4.1. If $f$ is a measurable function in ${\bf R}$ satisfying
(2.4.1) and $\mu$ is a positive Borel measure satisfying
$$
\int_{{\bf C}_{+}}\frac{\eta^p}{(1+|\zeta|)^\gamma}
d\mu(\zeta)<\infty
$$
and
$$
\int_{{\bf C}_{+}}\frac{1}{1+|\zeta|} d\mu(\zeta)<\infty.
$$
Write the subharmonic function
$$
u(z)= v(z)+h(z), \quad z\in {\bf C}_{+},
$$
where $v(z)$ is the harmonic function defined by (2.1.6), $h(z)$ is
defined by
$$
h(z)= \int_{{\bf C}_{+}} G(z,\zeta)d\mu(\zeta)
$$
and $G(z,\zeta)$ is defined by (2.1.1). Then there exists $z_j\in
{\bf C}_{+},\ \rho_j>0,$ such that (2.4.2) holds and
$$
u(z)=
o(y^{1-\frac{\alpha}{p}}|z|^{\frac{\gamma}{p}+\frac{1}{q}-2+\frac{\alpha}{p}}),
\quad  {\rm as}  \ |z|\rightarrow\infty \eqno{(2.4.4)}
$$
holds in ${\bf C}_{+}-G$, where $ G=\bigcup_{j=1}^\infty
B(z_j,\rho_j)$ and $0< \alpha<2p$. }

\vspace{0.2cm}
 \noindent
 {\bf Remark 2.4.2 } {\it If $\gamma=1-p$, $p>1$, then
$$
u(z)=
o(y^{1-\frac{\alpha}{p}}(\log|z|)^{\frac{1}{q}}|z|^{\frac{\gamma}{p}+\frac{1}{q}-2+\frac{\alpha}{p}}),
\quad  {\rm as}  \ |z|\rightarrow\infty
$$
holds in ${\bf C}_{+}-G$. }

\vspace{0.2cm}
 \noindent
 {\bf Remark 2.4.3 } {\it If $\alpha=1$, $p=1$ and $\gamma=2$, then (2.4.2) holds and (2.4.4)
 holds in ${\bf C}_{+}-G$. This is just
 the the result of Hayman, therefore, our result (2.4.4) is the
 generalization of Theorem B.}

\vspace{0.4cm}

\subsection*{2.   Main Lemmas }

 In order to obtain the results, we
need these lemmas below:

\vspace{0.2cm}
 \noindent
{\bf Lemma 2.4.1 } {\it The kernel function $\frac{1}{|z-\zeta|^2}$
has the
following estimates:\\
{\rm (1)}\ If $|\zeta|\leq \frac{|z|}{2}$, then
$\frac{1}{|z-\zeta|^2}\leq
\frac{4}{|z|^2}$;\\
{\rm (2)}\ If $|\zeta|> 2|z|$, then $\frac{1}{|z-\zeta|^2}\leq
\frac{4}{|\zeta|^2}$.\\
}

\vspace{0.2cm}
 \noindent
{\bf Lemma 2.4.2 } {\it The Green function $G(z,\zeta)$ has the
following estimates:\\
{\rm (1)}\ $|G(z,\zeta)|\leq A\log\frac{3y}{|z-\zeta|}$;\\
{\rm (2)}\ $|G(z,\zeta)|\leq \frac{y\eta}{\pi|z-\zeta|^2}$.}

  Proof: (1) is obvious; (2) follows by the Mean Value Theorem for
Derivatives.

\vspace{0.2cm}
 \noindent
{\bf Lemma 2.4.3 }  {\it The following estimate holds:
$$
\int_0^{\frac{y}{2}}
t^{2p-\alpha-1}\bigg(\log\frac{3y}{t}\bigg)^{p-1} dt \leq
\frac{3^{2p-\alpha}}{(2p-\alpha)^p}\Gamma (p) y^{2p-\alpha}.
$$
}

\subsection*{3.   Proof of Theorems }

 \emph{Proof of Theorem 2.4.1}

 We prove only the case $p>1$; the proof of the case $p=1$ is similar.
Define the measure $dm(\xi)$ by
$$
dm(\xi)=\frac{|f(\xi)|^p}{(1+|\xi|)^{\gamma}} d\xi.
$$
  For any $\varepsilon >0$, there exists $R_\varepsilon >2$, such that
$$
\int_{|\xi|\geq
R_\varepsilon}dm(\xi)\leq\frac{\varepsilon^p}{5^{2p-\alpha}}.
$$
For every Lebesgue measurable set $E \subset {\bf R}$  , the measure
$m^{(\varepsilon)}$ defined by $m^{(\varepsilon)}(E)
=m(E\cap\{x\in{\bf R}:|x|\geq R_\varepsilon\}) $ satisfies
$m^{(\varepsilon)}({\bf R})\leq\frac{\varepsilon^p}{5^{2p-\alpha}}$,
write
\begin{eqnarray*}
&v_1(z)& =\int_{G_1} P(z,\xi)f(\xi)
d\xi,\\
&v_2(z)&=\int_{G_2} P(z,\xi)f(\xi)
d\xi, \\
&v_3(z)&=\int_{G_3} P(z,\xi)f(\xi)
d\xi, \\
&v_4(z)&=\int_{G_4} P(z,\xi)f(\xi)
d\xi, \\
\end{eqnarray*}
where
\begin{eqnarray*}
&G_1& =\{\xi\in {\bf R}: R_\varepsilon<|\xi|\leq \frac{|z|}{2}\},\\
&G_2& =\{\xi\in {\bf R}: \frac{|z|}{2}<|\xi| \leq 2|z|\}, \\
&G_3& =\{\xi\in {\bf R}: |\xi|>2|z|\}, \\
&G_4& =\{\xi\in {\bf R}: |\xi|\leq R_\varepsilon\}. \\
\end{eqnarray*}
Then
$$
v(z) =v_1(z)+v_2(z)+v_3(z)+v_4(z). \eqno{(2.4.5)}
$$

  First, if $\gamma >1-p$, then $\frac{\gamma q}{p}
  +1>0$, so that we obtain by (1) of Lemma 2.4.1 and H\"{o}lder's
inequality
\begin{eqnarray*}
|v_1(z)|
&\leq& \int_{G_1}\frac{ y}{\pi}\frac{4}{|z|^2}|f(\xi)| d\xi \\
&\leq& \frac{4}{\pi}\frac{y}{|z|^2}
\bigg(\int_{G_1}\frac{|f(\xi)|^p}{|\xi|^\gamma}
d\xi\bigg)^{1/p}\bigg(\int_{G_1}|\xi|^{\frac{\gamma q}{p}}
d\xi\bigg)^{1/q}, \\
\end{eqnarray*}
since
$$
\int_{G_1}|\xi|^{\frac{\gamma q}{p}} d\xi\leq \frac{2}{\frac{\gamma
q}{p}+1}\bigg(\frac{|z|}{2}\bigg)^{\frac{\gamma q}{p}+1} ,
$$
so that
$$
|v_1(z)|\leq A \varepsilon y|z|^{\frac{\gamma}{p}+\frac{1}{q}-2}.
 \eqno{(2.4.6)}
$$

  Let $ E_1(\lambda)=\{z\in{\bf C}:|z|\geq2,\exists \ t>0,
s.t.m^{(\varepsilon)}(B(z,t)\cap{\bf R}
)>\lambda^p(\frac{t}{|z|})^{2p-\alpha}\}$, therefore, if $ |z|\geq
2R_\varepsilon$ and $z \notin E_1(\lambda)
 $, then we have
$$
\forall t>0,\ m^{(\varepsilon)}(B(z,t)\cap{\bf R} )\leq\lambda^p
\bigg(\frac{t}{|z|}\bigg)^{2p-\alpha}.
$$

  If $\gamma >1-p$, then $\frac{\gamma q}{p}
  +1>0$, so that we obtain by H\"{o}lder's inequality
\begin{eqnarray*}
|v_2(z)|
&\leq& \frac{y}{\pi}
\bigg(\int_{G_2}\frac{|f(\xi)|^p}{|z-(\xi,0)|^{2p}|\xi|^\gamma}
d\xi\bigg)^{1/p}\bigg(\int_{G_2}|\xi|^{\frac{\gamma q}{p}}
d\xi\bigg)^{1/q} \\
&\leq& Ay|z|^{\frac{\gamma}{p}+\frac{1}{q}}
\bigg(\int_{G_2}\frac{|f(\xi)|^p}{|z-(\xi,0)|^{2p}|\xi|^\gamma}
d\xi\bigg)^{1/p},\\
\end{eqnarray*}
since
\begin{eqnarray*}
\int_{G_2}\frac{|f(\xi)|^p}{|z-(\xi,0)|^{2p}|\xi|^\gamma} d\xi
&\leq& \int_y^{3|z|}
\frac{2^\gamma+1}{t^{2p}} dm_z^{(\varepsilon)}(t) \\
&\leq& \frac{\lambda^p}{
|z|^{2p}}(2^\gamma+1)\bigg(\frac{1}{3^\alpha}+
\frac{2p}{\alpha}\bigg)\frac{|z|^\alpha}{y^\alpha}, \\
\end{eqnarray*}
where  $m_z^{(\varepsilon)}(t)=\int_{|z-(\xi,0)| \leq t}
dm^{(\varepsilon)}(\xi)$.\\
Hence we have
$$
|v_2(z)|\leq A \lambda
y^{1-\frac{\alpha}{p}}|z|^{\frac{\gamma}{p}+\frac{1}{q}-2+\frac{\alpha}{p}}.
\eqno{(2.4.7)}
$$

  If $\gamma <1+p$, then $(\frac{\gamma}{p}-2)q+1<0$,
so that we obtain by (2) of Lemma 2.4.1 and H\"{o}lder's inequality

\begin{eqnarray*}
|v_3(z)|
&\leq& \int_{G_3}\frac{y}{\pi}\frac{4}{|\xi|^2}|f(\xi)| d\xi \\
&\leq& \frac{4}{\pi}y
\bigg(\int_{G_3}\frac{|f(\xi)|^p}{|\xi|^\gamma}
d\xi\bigg)^{1/p}\bigg(\int_{G_3}|\xi|^{(\frac{\gamma }{p}-2)q}
d\xi\bigg)^{1/q} \\
&\leq& A \varepsilon y|z|^{\frac{\gamma}{p}+\frac{1}{q}-2}.
\hspace{80mm} (2.4.8)
\end{eqnarray*}

  Finally, by (1) of Lemma 2.4.1, we obtain
$$
|v_4(z)|\leq \frac{4}{\pi}\frac{y}{|z|^2} \int_{G_4}{|f(\xi)|}d\xi,
$$
which implies by $\gamma >1-p$ that
$$
|v_4(z)|\leq A \varepsilon y|z|^{\frac{\gamma}{p}+\frac{1}{q}-2}.
\eqno{(2.4.9)}
$$

  Thus, by collecting (2.4.5), (2.4.6), (2.4.7), (2.4.8) and
(2.4.9), there exists a positive constant $A$ independent of
$\varepsilon$, such that if $ |z|\geq 2R_\varepsilon$ and $\  z
\notin E_1(\varepsilon)$, we have
$$
|v(z)|\leq A\varepsilon
y^{1-\frac{\alpha}{p}}|z|^{\frac{\gamma}{p}+\frac{1}{q}-2+\frac{\alpha}{p}}.
$$

 Let $\mu_\varepsilon$ be a measure in ${\bf C}$ defined by
$ \mu_\varepsilon(E)= m^{(\varepsilon)}(E\cap{\bf R})$ for every
measurable set $E$ in ${\bf C}$. Take
$\varepsilon=\varepsilon_p=\frac{1}{2^{p+2}}, p=1,2,3,\cdots$, then
there exists a sequence $ \{R_p\}$: $1=R_0<R_1<R_2<\cdots$ such that
$$
\mu_{\varepsilon_p}({\bf C})=\int_{|\xi|\geq
R_p}dm(\xi)<\frac{\varepsilon_p^p}{5^{2p-\alpha}}.
$$
Take $\lambda=3\cdot5^{2p-\alpha}\cdot2^p\mu_{\varepsilon_p}({\bf
C})$ in Lemma 2.2.1, then there exists $z_{j,p}$ and $ \rho_{j,p}$,
where $R_{p-1}\leq |z_{j,p}|<R_p$, such that
$$
\sum
_{j=1}^{\infty}\bigg(\frac{\rho_{j,p}}{|z_{j,p}|}\bigg)^{2p-\alpha}
\leq \frac{1}{2^{p}}.
$$
If $R_{p-1}\leq |z|<R_p$ and $z\notin G_p=\cup_{j=1}^\infty
B(z_{j,p},\rho_{j,p})$, we have
$$
|v(z)|\leq
A\varepsilon_py^{1-\frac{\alpha}{p}}|z|^{\frac{\gamma}{p}+\frac{1}{q}-2+\frac{\alpha}{p}}.
$$
Thereby
$$
\sum _{p=1}^{\infty}
\sum_{j=1}^{\infty}\bigg(\frac{\rho_{j,p}}{|z_{j,p}|}\bigg)^{2p-\alpha}
\leq \sum _{p=1}^{\infty}\frac{1}{2^{p}}=1<\infty.
$$
Set $ G=\cup_{p=1}^\infty G_p$, thus Theorem 2.4.1 holds.

 \emph{Proof of Theorem 2.4.2}

 We prove only the case $p>1$; the remaining case $p=1$ can be proved similarly.
Define the measure $dn(\zeta)$ by
$$
dn(\zeta)=\frac{\eta^p}{(1+|\zeta|)^{\gamma}} d\mu(\zeta).
$$

  For any $\varepsilon >0$, there exists $R_\varepsilon >2$, such that
$$
\int_{|\zeta|\geq
R_\varepsilon}dn(\zeta)<\frac{\varepsilon^p}{5^{2p-\alpha}}.
$$
For every Lebesgue measurable set $E \subset {\bf C}$, the measure
$n^{(\varepsilon)}$ defined by $n^{(\varepsilon)}(E)
=n(E\cap\{\zeta\in {\bf C}_{+}:|\zeta|\geq R_\varepsilon\}) $
satisfies $n^{(\varepsilon)}({\bf
C}_{+})\leq\frac{\varepsilon^p}{5^{2p-\alpha}}$, write
\begin{eqnarray*}
&h_1(z)& = \int_{F_1} G(z,\zeta)d\mu(\zeta),\\
&h_2(z)&=\int_{F_2} G(z,\zeta)d\mu(\zeta),\\
&h_3(z)&=\int_{F_3} G(z,\zeta)d\mu(\zeta), \\
&h_4(z)&=\int_{F_4} G(z,\zeta)d\mu(\zeta), \\
\end{eqnarray*}
where
\begin{eqnarray*}
&F_1& =\{\zeta\in {\bf C}_{+}: R_\varepsilon<|\zeta|\leq \frac{|z|}{2}\},\\
&F_2& =\{\zeta\in {\bf C}_{+}: \frac{|z|}{2}<|\zeta| \leq 2|z|\}, \\
&F_3& =\{\zeta\in {\bf C}_{+}: |\zeta|>2|z|\}, \\
&F_4& =\{\zeta\in {\bf C}_{+}: |\zeta|\leq R_\varepsilon\}. \\
\end{eqnarray*}
Then
$$
h(z) =h_1(z)+h_2(z)+h_3(z)+h_4(z). \eqno{(2.4.10)}
$$

  First, if $\gamma >1-p$, then $\frac{\gamma q}{p}
  +1>0$, so that we obtain by (1) of Lemma 2.4.1, (2) of Lemma
2.4.2 and H\"{o}lder's inequality
\begin{eqnarray*}
|h_1(z)|
&\leq& \int_{F_1}\frac{y\eta}{\pi|z-\zeta|^2} d\mu(\zeta) \\
&\leq& \int_{F_1}\frac{y\eta}{\pi}\frac{4}{|z|^2} d\mu(\zeta) \\
&\leq& \frac{4}{\pi}\frac{y}{|z|^2}
\bigg(\int_{F_1}\frac{\eta^p}{|\zeta|^\gamma}
d\mu(\zeta)\bigg)^{1/p}\bigg(\int_{F_1}|\zeta|^{\frac{\gamma q}{p}}
d\mu(\zeta)\bigg)^{1/q}, \\
\end{eqnarray*}
since
$$
\int_{F_1}|\zeta|^{\frac{\gamma q}{p}} d\mu(\zeta)\leq
2\bigg(\frac{|z|}{2}\bigg)^{\frac{\gamma
q}{p}+1}\int_{H}\frac{1}{(1+|\zeta|)} d\mu(\zeta),
$$
so that
$$
|h_1(z)|\leq A \varepsilon y|z|^{\frac{\gamma}{p}+\frac{1}{q}-2}.
\eqno{(2.4.11)}
$$

  Let $ E_2(\lambda)=\{z\in{\bf C}:|z|\geq2,\exists \ t>0, s.t.
n^{(\varepsilon)}(B(z,t)\cap {\bf C}_{+}
)>\lambda^p\big(\frac{t}{|z|}\big)^{2p-\alpha}\},$ therefore, if $
|z|\geq 2R_\varepsilon$ and $z\notin E_2(\lambda)
 $, then we have
$$
\forall t>0, \ n^{(\varepsilon)}(B(z,t)\cap H
)\leq\lambda^p(\frac{t}{|z|})^{2p-\alpha}.
$$

  If $\gamma >1-p$, then $\frac{\gamma q}{p}
  +1>0$, so that we obtain by H\"{o}lder's inequality
\begin{eqnarray*}
|h_2(z)|
&\leq& \bigg(\int_{F_2}\frac{|G(z,\zeta)|^p}{|\zeta|^\gamma}
d\mu(\zeta)\bigg)^{1/p}\bigg(\int_{F_2}|\zeta|^{\frac{\gamma q}{p}}
d\mu(\zeta)\bigg)^{1/q} \\
&\leq& \bigg((2^\gamma +1)\int_{F_2}\frac{|G(z,\zeta)|^p}{\eta^p}
dn(\zeta)\bigg)^{1/p}\bigg(\int_{F_2}|\zeta|^{\frac{\gamma q}{p}}
d\mu(\zeta)\bigg)^{1/q}\\
&\leq& A|z|^{\frac{\gamma}{p}+\frac{1}{q}}\bigg(\int_{F_2}\frac{|G(z,\zeta)|^p}{\eta^p}dn(\zeta)\bigg)^{1/p},\\
\end{eqnarray*}
since
\begin{eqnarray*}
\int_{F_2}\frac{|G(z,\zeta)|^p}{\eta^p} dn(\zeta)
&\leq&  \int_{|z-\zeta|\leq 3|z|}\frac{|G(z,\zeta)|^p}{\eta^p} dn^{(\varepsilon)}(\zeta) \\
&=&  \int_{|z-\zeta|\leq \frac{y}{2}}\frac{|G(z,\zeta)|^p}{\eta^p}
dn^{(\varepsilon)}(\zeta)+
 \int_{\frac{y}{2}<|z-\zeta|\leq 3|z|}\frac{|G(z,\zeta)|^p}{\eta^p} dn^{(\varepsilon)}(\zeta)\\
&=& h_{21}(z)+h_{22}(z),
\end{eqnarray*}
so that we have by (1) of Lemma 2.4.2 and Lemma 2.4.3
\begin{eqnarray*}
h_{21}(z)
&\leq& \int_{|z-\zeta|\leq
\frac{y}{2}}\bigg(\frac{A}{y}\log\frac{3y}{|z-\zeta|}\bigg)^p
dn^{(\varepsilon)}(\zeta) \\
&=& \frac{A}{y^p}\int_0^{\frac{y}{2}}
\bigg(\log\frac{3y}{t}\bigg)^p dn_z^{(\varepsilon)}(t) \\
&\leq& A\lambda^p\frac{y^{p-\alpha}}{|z|^{2p-\alpha}}
+A\lambda^p\frac{1}{y^p|z|^{2p-\alpha}}\int_0^{\frac{y}{2}}
t^{2p-\alpha-1}\bigg(\log\frac{3y}{t}\bigg)^{p-1} dt \\
&\leq& A\lambda^p\frac{y^{p-\alpha}}{|z|^{2p-\alpha}}.\\
\end{eqnarray*}

  Moreover, we have by (2) of Lemma 2.4.2
\begin{eqnarray*}
h_{22}(z)
&\leq& \int_{\frac{y}{2}<|z-\zeta|\leq 3|z|
}\bigg(\frac{y}{\pi|z-\zeta|^2}\bigg)^p
dn^{(\varepsilon)}(\zeta) \\
&=& \bigg(\frac{y}{\pi}\bigg)^p\int_{\frac{y}{2}}^{3|z|}
\frac{1}{t^{2p}} dn_z^{(\varepsilon)}(t) \\
&\leq&
\bigg(\frac{1}{\pi}\bigg)^p\bigg(\frac{1}{3^\alpha}+\frac{2^{\alpha+1}p}{\alpha}\bigg)
\lambda^p\frac{y^{p-\alpha}}{|z|^{2p-\alpha}},\\
\end{eqnarray*}
where  $n_z^{(\varepsilon)}(t)=\int_{|z-\zeta| \leq t}
dn^{(\varepsilon)}(\zeta)$.\\
Hence we have
$$
|h_2(z)|\leq A \lambda
y^{1-\frac{\alpha}{p}}|z|^{\frac{\gamma}{p}+\frac{1}{q}-2+\frac{\alpha}{p}}.
\eqno{(2.4.12)}
$$

  If $\gamma <1+p$, then $(\frac{\gamma}{p}-2)q
  +1<0$, so that we obtain by (2) of Lemma 2.4.1, (2) of Lemma 2.4.2 and H\"{o}lder's inequality

\begin{eqnarray*}
|h_3(z)|
&\leq& \int_{F_3}\frac{y\eta}{\pi|z-\zeta|^2} d\mu(\zeta) \\
&\leq& \int_{F_3}\frac{y\eta}{\pi}\frac{4}{|\zeta|^2} d\mu(\zeta) \\
&\leq& \frac{4}{\pi}y \bigg(\int_{F_3}\frac{\eta^p}{|\zeta|^\gamma}
d\mu(\zeta)\bigg)^{1/p}\bigg(\int_{F_3}|\zeta|^{(\frac{\gamma
}{p}-2)q}
d\mu(\zeta)\bigg)^{1/q} \\
&\leq& A\varepsilon y|z|^{\frac{\gamma}{p}+\frac{1}{q}-2}.
\hspace{80mm} (2.4.13)
\end{eqnarray*}

  Finally, by (1) of Lemma 2.4.1 and (2) of Lemma 2.4.2, we obtain
$$
|h_4(z)|\leq \int_{F_4}\frac{y\eta}{\pi|z-\zeta|^2} d\mu(\zeta) \leq
\frac{4}{\pi}\frac{y}{|z|^2} \int_{F_4}\eta d\mu(\zeta),
$$
which implies by $\gamma >1-p$ that
$$
|h_4(z)|\leq A\varepsilon y|z|^{\frac{\gamma}{p}+\frac{1}{q}-2}.
\eqno{(2.4.14)}
$$

  Thus, by collecting (2.4.10), (2.4.11),
(2.4.12), (2.4.13) and (2.4.14), there exists a positive constant
$A$ independent of $\varepsilon$, such that if $ |z|\geq
2R_\varepsilon$ and $\  z \notin E_2(\varepsilon)$, we have
$$
 |h(z)|\leq A\varepsilon y^{1-\frac{\alpha}{p}}|z|^{\frac{\gamma}{p}+\frac{1}{q}-2+\frac{\alpha}{p}}.
$$

  Similarly, if $z\notin G$, we have
$$
h(z)=
o(y^{1-\frac{\alpha}{p}}|z|^{\frac{\gamma}{p}+\frac{1}{q}-2+\frac{\alpha}{p}}),\quad
{\rm as} \ |z|\rightarrow\infty. \eqno{(2.4.15)}
$$

  By (2.4.3) and  (2.4.15), we obtain that
$$
u(z)=v(z)+h(z)=
o(y^{1-\frac{\alpha}{p}}|z|^{\frac{\gamma}{p}+\frac{1}{q}-2+\frac{\alpha}{p}}),\quad
{\rm as} \ |z|\rightarrow\infty
$$
holds in ${\bf C}_{+}-G$, thus we complete the proof of Theorem
2.4.2.

\section{$p>1$(Modified Kernel)}

\subsection*{ 1. Introduction and Main Theorems}

\vspace{0.3cm}

  In this section, we will consider
measurable functions $f$ in ${\bf R}$  satisfying
$$
\int_{{\bf R}}\frac{|f(\xi)|^p}{(1+|\xi|)^\gamma}
d\xi<\infty,\eqno{(2.5.1)}
$$
where $\gamma$ is defined as in Theorem 2.5.1.

  In order to describe the asymptotic behaviour of subharmonic
functions represented by the modified kernel in the upper half plane
(see \cite{LDG}, \cite{ZYT}, \cite{ZA}, \cite{ML}, \cite{MLJ1}, and
\cite{MLJ2}), we establish the following theorems.

\vspace{0.2cm}
 \noindent
{\bf Theorem 2.5.1} {\it Let $1\leq p<\infty,\
\frac{1}{p}+\frac{1}{q}=1$ and
$$
1+mp<\gamma <1+(m+1)p  \qquad  {\rm in \  case}  \  p>1;
$$
$$
\quad \ \    m+1<\gamma \leq m+2  \qquad  \qquad {\rm in \ case}  \
p=1.
$$
If $f$ is a measurable function in ${\bf R}$ satisfying (2.5.1) and
$v(z)$ is the harmonic function defined by
$$
v(z)=\int_{{\bf R}}P_m(z,\xi)f(\xi)d\xi,  \eqno{(2.5.2)}
$$
then there exists $z_j\in {\bf C}_{+},\ \rho_j>0,$ such that
$$
\sum
_{j=1}^{\infty}\frac{\rho_j^{2p-\alpha}}{|z_j|^{2p-\alpha}}<\infty
\eqno{(2.5.3)}
$$
holds and
$$
v(z)=
o(y^{1-\frac{\alpha}{p}}|z|^{\frac{\gamma}{p}+\frac{1}{q}-2+\frac{\alpha}{p}}),
\quad  {\rm as}  \ |z|\rightarrow\infty   \eqno{(2.5.4)}
$$
holds in ${{\bf C}_{+}}-G$, where $ G=\bigcup_{j=1}^\infty
B(z_j,\rho_j)$ and $0< \alpha\leq 2p$. }

 \vspace{0.2cm}
 \noindent
 {\bf Remark 2.5.1 } {\it If $\gamma=1+mp$, $p>1$, then
$$
v(z)=
o(y^{1-\frac{\alpha}{p}}(\log|z|)^{\frac{1}{q}}
|z|^{\frac{\gamma}{p}+\frac{1}{q}-2+\frac{\alpha}{p}}),
\quad  {\rm as}  \ |z|\rightarrow\infty
$$
holds in ${{\bf C}_{+}}-G$. }

  Next, we will generalize Theorem 2.5.1 to subharmonic functions.

\vspace{0.2cm}
 \noindent
{\bf Theorem 2.5.2 } {\it Let $p$ and $\gamma$ be as in Theorem
2.5.1. If $f$ is a measurable function in ${\bf R}$ satisfying
(2.5.1) and $\mu$ is a positive Borel measure satisfying
$$
\int_{{\bf C}_{+}}\frac{\eta^p}{(1+|\zeta|)^\gamma}
d\mu(\zeta)<\infty
$$
and
$$
\int_{{\bf C}_{+}}\frac{1}{1+|\zeta|} d\mu(\zeta)<\infty.
$$
Write the subharmonic function
$$
u(z)= v(z)+h(z), \quad z\in {{\bf C}_{+}},
$$
where $v(z)$ is the harmonic function defined by (2.5.2), $h(z)$ is
defined by
$$
h(z)= \int_{{\bf C}_{+}} G_m(z,\zeta)d\mu(\zeta)
$$
and $G_m(z,\zeta)$ is defined by (2.3.2). Then there exists $z_j\in
{{\bf C}_{+}},\ \rho_j>0,$ such that (2.5.3) holds and
$$
u(z)=
o(y^{1-\frac{\alpha}{p}}|z|^{\frac{\gamma}{p}+\frac{1}{q}-2+\frac{\alpha}{p}}),
\quad  {\rm as}  \ |z|\rightarrow\infty
 \eqno{(2.5.5)}
$$
holds in ${{\bf C}_{+}}-G$, where $G=\bigcup_{j=1}^\infty
B(z_j,\rho_j)$ and $0< \alpha<2p$. }

\vspace{0.2cm}
 \noindent
 {\bf Remark 2.5.2 } {\it If $\gamma=1+mp$, $p>1$, then
$$
u(z)=
o(y^{1-\frac{\alpha}{p}}(\log|z|)^{\frac{1}{q}}|z|^{\frac{\gamma}{p}+\frac{1}{q}-2+\frac{\alpha}{p}}),
\quad  {\rm as}  \ |z|\rightarrow\infty
$$
holds in ${{\bf C}_{+}}-G$. }

\vspace{0.2cm}
 \noindent
 {\bf Remark 2.5.3} {\it If $\alpha=1$, $p=1$, $m=0$ and $\gamma=2$,
 then (2.5.3) holds and (2.5.5) holds in ${\bf C}_{+}-G$.
 This is just the result of Hayman, therefore,
 our result (2.5.5) is the generalization of Theorem B.
 }

\subsection*{ 2.  Main Lemmas}

 In order to obtain the results, we
need these lemmas below:

\vspace{0.2cm}
 \noindent
{\bf Lemma 2.5.1 } {\it The modified Poisson
 kernel $P_m(z,\xi)$ has the
following estimates:\\
{\rm (1)}\ If $1<|\xi|\leq \frac{|z|}{2}$, then $|P_m(z,\xi)|\leq
\frac{Ay|z|^{m-1}}{|\xi|^{m+1}}$;\\
{\rm (2)}\ If $\frac{|z|}{2}<|\xi| \leq 2|z|$, then
$|P_m(z,\xi)|\leq
\frac{Ay}{|z-(\xi,0)|^2}$;\\
{\rm (3)}\ If $|\xi|>2|z|$, then $|P_m(z,\xi)|\leq
\frac{Ay|z|^{m}}{|\xi|^{m+2}}$;\\
{\rm (4)}\ If $|\xi|\leq 1$, then $|P_m(z,\xi)|\leq
\frac{Ay}{|z|^{2}}$.\\
}

\vspace{0.2cm}
 \noindent
{\bf Lemma 2.5.2 }  {\it The modified Green function $G_m(z,\zeta)$
has the
following estimates:\\
{\rm (1)}\ If $1<|\zeta|\leq \frac{|z|}{2}$, then
$|G_m(z,\zeta)|\leq
\frac{Ay\eta|z|^{m-1}}{|\zeta|^{m+1}}$;\\
{\rm (2)}\ If $\frac{|z|}{2}<|\zeta| \leq 2|z|$, then
$|G_m(z,\zeta)|\leq
\frac{Ay\eta}{|z-\zeta|^2}$;\\
{\rm (3)}\ If $|\zeta|>2|z|$, then $|G_m(z,\zeta)|\leq
\frac{Ay\eta|z|^{m}}{|\zeta|^{m+2}}$;\\
{\rm (4)}\ If $|\zeta|\leq 1$, then $|G_m(z,\zeta)|\leq
\frac{y\eta}{\pi|z-\zeta|^2}\leq
\frac{Ay\eta}{|z|^{2}}$;\\
{\rm (5)}\ If $|\zeta-z| \leq \frac{y}{2}$, then $|G_m(z,\zeta)|\leq
A\log\frac{3y}{|z-\zeta|}$.\\
}

\subsection*{3.   Proof of Theorems}

 \emph{Proof of Theorem 2.5.1}

 We prove only the case $p>1$; the proof of the case $p=1$ is similar.
Define the measure $dm(\xi)$ by
$$
dm(\xi)=\frac{|f(\xi)|^p}{(1+|\xi|)^{\gamma}} d\xi.
$$
  For any $\varepsilon >0$, there exists $R_\varepsilon >2$, such that
$$
\int_{|\xi|\geq
R_\varepsilon}dm(\xi)\leq\frac{\varepsilon^p}{5^{2p-\alpha}}.
$$
For every Lebesgue measurable set $E \subset {\bf R}$, the measure
$m^{(\varepsilon)}$ defined by $m^{(\varepsilon)}(E)
=m(E\cap\{x\in{\bf R}:|x|\geq R_\varepsilon\}) $ satisfies
$m^{(\varepsilon)}({\bf R})\leq\frac{\varepsilon^p}{5^{2p-\alpha}}$,
write
\begin{eqnarray*}
&v_1(z)& =\int_{G_1} P_m(z,\xi)f(\xi)
d\xi,\\
&v_2(z)&=\int_{G_2} P_m(z,\xi)f(\xi)
d\xi, \\
&v_3(z)&=\int_{G_3} P_m(z,\xi)f(\xi)
d\xi, \\
&v_4(z)&=\int_{G_4} P_m(z,\xi)f(\xi)
d\xi, \\
\end{eqnarray*}
where
\begin{eqnarray*}
&G_1& =\{\xi\in {\bf R}: 1<|\xi|\leq \frac{|z|}{2}\},\\
&G_2& =\{\xi\in {\bf R}: \frac{|z|}{2}<|\xi| \leq 2|z|\}, \\
&G_3& =\{\xi\in {\bf R}: |\xi|>2|z|\}, \\
&G_4& =\{\xi\in {\bf R}: |\xi|\leq 1\}. \\
\end{eqnarray*}
Then
$$
v(z) =v_1(z)+v_2(z)+v_3(z)+v_4(z). \eqno{(2.5.6)}
$$

  First, if $\gamma >1+mp$, then $(\frac{\gamma}{p}-m-1)q
  +1>0$. For $R_\varepsilon >2$, we have
$$
v_1(z)=\int_{1<|\xi|\leq R_\varepsilon} P_m(z,\xi)f(\xi) d\xi
+\int_{R_\varepsilon<|\xi|\leq \frac{|z|}{2} } P_m(z,\xi)f(\xi)
d\xi=v_{11}(z)+v_{12}(z),
$$
if $|z|>2R_\varepsilon$, then we obtain by (1) of Lemma 2.5.1 and
H\"{o}lder's inequality
\begin{eqnarray*}
|v_{11}(z)|
&\leq& \int_{1<|\xi|\leq R_\varepsilon}\frac{Ay|z|^{m-1}}{|\xi|^{m+1}}|f(\xi)| d\xi \\
&\leq& Ay|z|^{m-1} \bigg(\int_{1<|\xi|\leq
R_\varepsilon}\frac{|f(\xi)|^p}{|\xi|^\gamma}
d\xi\bigg)^{1/p}\bigg(\int_{1<|\xi|\leq
R_\varepsilon}|\xi|^{(\frac{\gamma}{p}-m-1)q}
d\xi\bigg)^{1/q}, \\
\end{eqnarray*}
since
$$
\int_{1<|\xi|\leq R_\varepsilon}|\xi|^{(\frac{\gamma}{p}-m-1)q}
d\xi\leq AR_\varepsilon^{(\frac{\gamma}{p}-m-1)q+1},
$$
so that
$$
|v_{11}(z)|\leq
Ay|z|^{m-1}R_\varepsilon^{(\frac{\gamma}{p}-m-1)+\frac{1}{q}}.
\eqno{(2.5.7)}
$$

  Moreover, we have similarly

\begin{eqnarray*}
|v_{12}(z)|
&\leq& Ay|z|^{m-1} \bigg(\int_{R_\varepsilon<|\xi|\leq
\frac{|z|}{2}}\frac{|f(\xi)|^p}{|\xi|^\gamma}
d\xi\bigg)^{1/p}\bigg(\int_{R_\varepsilon<|\xi|\leq
\frac{|z|}{2}}|\xi|^{(\frac{\gamma}{p}-m-1)q}
d\xi\bigg)^{1/q} \\
&\leq& Ay|z|^{\frac{\gamma}{p}+\frac{1}{q}-2}
\bigg(\int_{R_\varepsilon<|\xi|\leq
\frac{|z|}{2}}\frac{|f(\xi)|^p}{|\xi|^\gamma}
d\xi\bigg)^{1/p}, \\
\end{eqnarray*}
which implies by arbitrariness of $R_\varepsilon$ that
$$
|v_{12}(z)|\leq A\varepsilon y|z|^{\frac{\gamma}{p}+\frac{1}{q}-2}.
\eqno{(2.5.8)}
$$

  Let $ E_1(\lambda)=\{z\in{\bf C}:|z|\geq2,\exists \
t>0, {\rm s.t.}\  m^{(\varepsilon)}(B(z,t)\cap{\bf R}
)>\lambda^p(\frac{t}{|z|})^{2p-\alpha}\}$, therefore, if $ |z|\geq
2R_\varepsilon$ and $z \notin E_1(\lambda)
 $, then we have
$$
\forall t>0,\ m^{(\varepsilon)}(B(z,t)\cap{\bf R}
)\leq\lambda^p\bigg(\frac{t}{|z|}\bigg)^{2p-\alpha}.
$$
If $\gamma >1+mp$, then $(\frac{\gamma}{p}-m-1)q
  +1>0$, so that we obtain by (2) of Lemma 2.5.1 and
H\"{o}lder's inequality
\begin{eqnarray*}
|v_2(z)|
&\leq& \int_{G_2}\frac{Ay}{|z-(\xi,0)|^2}|f(\xi)| d\xi \\
&\leq& Ay
\bigg(\int_{G_2}\frac{|f(\xi)|^p}{|z-(\xi,0)|^{2p}|\xi|^\gamma}
d\xi\bigg)^{1/p}\bigg(\int_{G_2}|\xi|^{\frac{\gamma q}{p}}
d\xi\bigg)^{1/q} \\
&\leq& Ay|z|^{\frac{\gamma}{p}+\frac{1}{q}}
\bigg(\int_{G_2}\frac{|f(\xi)|^p}{|z-(\xi,0)|^{2p}|\xi|^\gamma}
d\xi\bigg)^{1/p},\\
\end{eqnarray*}
since
\begin{eqnarray*}
\int_{G_2}\frac{|f(\xi)|^p}{|z-(\xi,0)|^{2p}|\xi|^\gamma} d\xi
&\leq& \int_y^{3|z|}
\frac{2^\gamma+1}{t^{2p}} dm_z^{(\varepsilon)}(t) \\
&\leq& \frac{\lambda^p}{
|z|^{2p}}(2^\gamma+1)\bigg(\frac{1}{3^\alpha}+
\frac{2p}{\alpha}\bigg)\frac{|z|^\alpha}{y^\alpha}, \\
\end{eqnarray*}
where  $m_z^{(\varepsilon)}(t)=\int_{|z-(\xi,0)| \leq t}
dm^{(\varepsilon)}(\xi)$.\\
Hence we have
$$
|v_2(z)|\leq A\lambda
y^{1-\frac{\alpha}{p}}|z|^{\frac{\gamma}{p}+\frac{1}{q}-2+\frac{\alpha}{p}}.\eqno{(2.5.9)}
$$

  If $\gamma <1+(m+1)p$, then $(\frac{\gamma}{p}-m-2)q+1<0$,
so that we obtain by (3) of Lemma 2.5.1 and H\"{o}lder's inequality
\begin{eqnarray*}
|v_3(z)|
&\leq& \int_{G_3}\frac{Ay|z|^{m}}{|\xi|^{m+2}}|f(\xi)| d\xi \\
&\leq& Ay|z|^{m} \bigg(\int_{G_3}\frac{|f(\xi)|^p}{|\xi|^\gamma}
d\xi\bigg)^{1/p}\bigg(\int_{G_3}|\xi|^{(\frac{\gamma }{p}-m-2)q}
d\xi\bigg)^{1/q} \\
&\leq& A\varepsilon y|z|^{\frac{\gamma}{p}+\frac{1}{q}-2}.
\hspace{80mm} (2.5.10)
\end{eqnarray*}

  Finally, by (4) of Lemma 2.5.1, we obtain
$$
|v_4(z)|\leq \frac{Ay}{|z|^{2}}
\int_{G_4}{|f(\xi)|}d\xi.\eqno{(2.5.11)}
$$

  Thus, by collecting (2.5.6), (2.5.7), (2.5.8), (2.5.9), (2.5.10) and
(2.5.11), there exists a positive constant $A$ independent of
$\varepsilon$, such that if $ |z|\geq 2R_\varepsilon$ and $\  z
\notin E_1(\varepsilon)$, we have
$$
|v(z)|\leq A\varepsilon
y^{1-\frac{\alpha}{p}}|z|^{\frac{\gamma}{p}+\frac{1}{q}-2+\frac{\alpha}{p}}.
$$

 Let $\mu_\varepsilon$ be a measure in ${\bf C}$ defined by
$ \mu_\varepsilon(E)= m^{(\varepsilon)}(E\cap{\bf R})$ for every
measurable set $E$ in ${\bf C}$. Take
$\varepsilon=\varepsilon_p=\frac{1}{2^{p+2}}, p=1,2,3,\cdots$, then
there exists a sequence $ \{R_p\}$: $1=R_0<R_1<R_2<\cdots$ such that
$$
\mu_{\varepsilon_p}({\bf C})=\int_{|\xi|\geq
R_p}dm(\xi)<\frac{\varepsilon_p^p}{5^{2p-\alpha}}.
$$
Take $\lambda=3\cdot5^{2p-\alpha}\cdot2^p\mu_{\varepsilon_p}({\bf
C})$ in Lemma 2.2.1, then there exists $z_{j,p}$ and $ \rho_{j,p}$,
where $R_{p-1}\leq |z_{j,p}|<R_p$, such that
$$
\sum
_{j=1}^{\infty}\bigg(\frac{\rho_{j,p}}{|z_{j,p}|}\bigg)^{2p-\alpha}
\leq \frac{1}{2^{p}}.
$$
If $R_{p-1}\leq |z|<R_p$ and $z\notin G_p=\cup_{j=1}^\infty
B(z_{j,p},\rho_{j,p})$, we have
$$
|v(z)|\leq
A\varepsilon_py^{1-\frac{\alpha}{p}}|z|^{\frac{\gamma}{p}+\frac{1}{q}-2+\frac{\alpha}{p}}.
$$
Thereby
$$
\sum _{p=1}^{\infty}
\sum_{j=1}^{\infty}\bigg(\frac{\rho_{j,p}}{|z_{j,p}|}\bigg)^{2p-\alpha}
\leq \sum _{p=1}^{\infty}\frac{1}{2^{p}}=1<\infty.
$$
Set $ G=\cup_{p=1}^\infty G_p$, thus Theorem 2.5.1 holds.

 \emph{Proof of Theorem 2.5.2}

 We prove only the case $p>1$; the remaining case $p=1$
can be proved similarly. Define the measure $dn(\zeta)$ by
$$
dn(\zeta)=\frac{\eta^p}{(1+|\zeta|)^{\gamma}} d\mu(\zeta).
$$

  For any $\varepsilon >0$, there exists $R_\varepsilon >2$, such that
$$
\int_{|\zeta|\geq
R_\varepsilon}dn(\zeta)<\frac{\varepsilon^p}{5^{2p-\alpha}}.
$$
For every Lebesgue measurable set $E \subset {\bf C}$, the measure
$n^{(\varepsilon)}$ defined by $n^{(\varepsilon)}(E)
=n(E\cap\{\zeta\in {{\bf C}_{+}}:|\zeta|\geq R_\varepsilon\}) $
satisfies $n^{(\varepsilon)}({{\bf
C}_{+}})\leq\frac{\varepsilon^p}{5^{2p-\alpha}}$, write
\begin{eqnarray*}
&h_1(z)& = \int_{F_1} G_m(z,\zeta)d\mu(\zeta),\\
&h_2(z)&=\int_{F_2} G_m(z,\zeta)d\mu(\zeta),\\
&h_3(z)&=\int_{F_3} G_m(z,\zeta)d\mu(\zeta), \\
&h_4(z)&=\int_{F_4} G_m(z,\zeta)d\mu(\zeta), \\
\end{eqnarray*}
where
\begin{eqnarray*}
&F_1& =\{\zeta\in {{\bf C}_{+}}: 1<|\zeta|\leq \frac{|z|}{2}\},\\
&F_2& =\{\zeta\in {{\bf C}_{+}}: \frac{|z|}{2}<|\zeta| \leq 2|z|\}, \\
&F_3& =\{\zeta\in {{\bf C}_{+}}: |\zeta|>2|z|\}, \\
&F_4& =\{\zeta\in {{\bf C}_{+}}: |\zeta|\leq 1\}. \\
\end{eqnarray*}
Then
$$
h(z) =h_1(z)+h_2(z)+h_3(z)+h_4(z). \eqno{(2.5.12)}
$$

  First, if $\gamma >1+mp$, then $(\frac{\gamma}{p}-m-1)q
  +1>0$. For $R_\varepsilon >2$, we have
$$
h_1(z)=\int_{1<|\zeta|\leq R_\varepsilon} G_m(z,\zeta) d\mu(\zeta)
+\int_{R_\varepsilon<|\zeta|\leq \frac{|z|}{2} } G_m(z,\zeta)
d\mu(\zeta)=h_{11}(z)+h_{12}(z),
$$
if $|z|>2R_\varepsilon$, then we obtain by (1) of Lemma 2.5.2 and
H\"{o}lder's inequality
\begin{eqnarray*}
|h_{11}(z)|
&\leq& \int_{1<|\zeta|\leq R_\varepsilon}\frac{Ay\eta|z|^{m-1}}{|\zeta|^{m+1}} d\mu(\zeta) \\
&\leq& Ay|z|^{m-1} \bigg(\int_{1<|\zeta|\leq
R_\varepsilon}\frac{\eta^p}{|\zeta|^\gamma}
d\mu(\zeta)\bigg)^{1/p}\bigg(\int_{1<|\zeta|\leq
R_\varepsilon}|\zeta|^{(\frac{\gamma}{p}-m-1)q}
d\mu(\zeta)\bigg)^{1/q}, \\
\end{eqnarray*}
since
$$
\int_{1<|\zeta|\leq R_\varepsilon}|\zeta|^{(\frac{\gamma}{p}-m-1)q}
d\mu(\zeta)\leq AR_\varepsilon^{(\frac{\gamma}{p}-m-1)q+1},
$$
so that
$$
|h_{11}(z)|\leq
Ay|z|^{m-1}R_\varepsilon^{(\frac{\gamma}{p}-m-1)+\frac{1}{q}}.
\eqno{(2.5.13)}
$$

  Moreover, we have similarly

\begin{eqnarray*}
|h_{12}(z)|
&\leq& Ay|z|^{m-1} \bigg(\int_{R_\varepsilon<|\zeta|\leq
\frac{|z|}{2}}\frac{\eta^p}{|\zeta|^\gamma}
d\mu(\zeta)\bigg)^{1/p}\bigg(\int_{R_\varepsilon<|\zeta|\leq
\frac{|z|}{2}}|\zeta|^{(\frac{\gamma}{p}-m-1)q}
d\mu(\zeta)\bigg)^{1/q} \\
&\leq& Ay|z|^{\frac{\gamma}{p}+\frac{1}{q}-2}
\bigg(\int_{R_\varepsilon<|\zeta|\leq
\frac{|z|}{2}}\frac{\eta^p}{|\zeta|^\gamma}
d\mu(\zeta)\bigg)^{1/p}, \\
\end{eqnarray*}
which implies by arbitrariness of $R_\varepsilon$ that
$$
|h_{12}(z)|\leq A\varepsilon y|z|^{\frac{\gamma}{p}+\frac{1}{q}-2}.
\eqno{(2.5.14)}
$$

  Let $ E_2(\lambda)=\{z\in{\bf C}:|z|\geq2,\exists \
t>0,{\rm s.t.}\ n^{(\varepsilon)}(B(z,t)\cap {{\bf C}_{+}}
)>\lambda^p(\frac{t}{|z|})^{2p-\alpha}\}, $ therefore, if $ |z|\geq
2R_\varepsilon$ and $z\notin E_2(\lambda)
 $, then we have
$$
\forall t>0, \ n^{(\varepsilon)}(B(z,t)\cap {{\bf C}_{+}}
)\leq\lambda^p\bigg(\frac{t}{|z|}\bigg)^{2p-\alpha}.
$$
If $\gamma >1+mp$, then $(\frac{\gamma}{p}-m-1)q
  +1>0$, so that we obtain by H\"{o}lder's inequality
\begin{eqnarray*}
|h_2(z)|
&\leq& \bigg(\int_{F_2}\frac{|G_m(z,\zeta)|^p}{|\zeta|^\gamma}
d\mu(\zeta)\bigg)^{1/p}\bigg(\int_{F_2}|\zeta|^{\frac{\gamma q}{p}}
d\mu(\zeta)\bigg)^{1/q} \\
&\leq& \bigg((2^\gamma +1)\int_{F_2}\frac{|G_m(z,\zeta)|^p}{\eta^p}
dn(\zeta)\bigg)^{1/p}\bigg(\int_{F_2}|\zeta|^{\frac{\gamma q}{p}}
d\mu(\zeta)\bigg)^{1/q}\\
&\leq& A|z|^{\frac{\gamma}{p}+\frac{1}{q}}\bigg(\int_{F_2}\frac{|G_m(z,\zeta)|^p}{\eta^p}dn(\zeta)\bigg)^{1/p},\\
\end{eqnarray*}
since
\begin{eqnarray*}
\int_{F_2}\frac{|G_m(z,\zeta)|^p}{\eta^p} dn(\zeta)
&\leq&  \int_{|z-\zeta|\leq 3|z|}\frac{|G_m(z,\zeta)|^p}{\eta^p} dn^{(\varepsilon)}(\zeta) \\
&=&  \int_{|z-\zeta|\leq \frac{y}{2}}\frac{|G_m(z,\zeta)|^p}{\eta^p}
dn^{(\varepsilon)}(\zeta)+
 \int_{\frac{y}{2}<|z-\zeta|\leq 3|z|}\frac{|G_m(z,\zeta)|^p}{\eta^p} dn^{(\varepsilon)}(\zeta)\\
&=& h_{21}(z)+h_{22}(z),
\end{eqnarray*}
so that we have by (5) of Lemma 2.5.2 and Lemma 2.4.3
\begin{eqnarray*}
h_{21}(z)
&\leq& \int_{|z-\zeta|\leq
\frac{y}{2}}\bigg(\frac{A}{y}\log\frac{3y}{|z-\zeta|}\bigg)^p
dn^{(\varepsilon)}(\zeta) \\
&=& \frac{A}{y^p}\int_0^{\frac{y}{2}}
\bigg(\log\frac{3y}{t}\bigg)^p dn_z^{(\varepsilon)}(t) \\
&\leq& A\lambda^p\frac{y^{p-\alpha}}{|z|^{2p-\alpha}}
+A\lambda^p\frac{1}{y^p|z|^{2p-\alpha}}\int_0^{\frac{y}{2}}
t^{2p-\alpha-1}\bigg(\log\frac{3y}{t}\bigg)^{p-1} dt \\
&\leq& A\lambda^p\frac{y^{p-\alpha}}{|z|^{2p-\alpha}}.\\
\end{eqnarray*}

  Moreover, we have by (2) of Lemma 2.5.2
\begin{eqnarray*}
h_{22}(z)
&\leq& \int_{\frac{y}{2}<|z-\zeta|\leq 3|z|
}\bigg(\frac{Ay}{|z-\zeta|^{2}}\bigg)^p
dn^{(\varepsilon)}(\zeta) \\
&=& (Ay)^p\int_{\frac{y}{2}}^{3|z|}
\frac{1}{t^{2p}} dn_z^{(\varepsilon)}(t) \\
&\leq& A\bigg(\frac{1}{3^\alpha}+\frac{2p2^\alpha}{\alpha}\bigg)
\lambda^p\frac{y^{p-\alpha}}{|z|^{2p-\alpha}},\\
\end{eqnarray*}
where  $n_z^{(\varepsilon)}(t)=\int_{|z-\zeta| \leq t}
dn^{(\varepsilon)}(\zeta)$.\\
Hence we have
$$
|h_2(z)|\leq A\lambda
y^{1-\frac{\alpha}{p}}|z|^{\frac{\gamma}{p}+\frac{1}{q}-2+\frac{\alpha}{p}}.
\eqno{(2.5.15)}
$$

 If $\gamma <1+(m+1)p$, then $(\frac{\gamma}{p}-m-2)q+1<0$,
so that we obtain by (3) of Lemma 2.5.2 and H\"{o}lder's inequality
\begin{eqnarray*}
|h_3(z)|
&\leq& \int_{F_3}\frac{Ay\eta|z|^{m}}{|\zeta|^{m+2}} d\mu(\zeta) \\
&\leq& Ay|z|^{m} \bigg(\int_{F_3}\frac{\eta^p}{|\zeta|^\gamma}
d\mu(\zeta)\bigg)^{1/p}\bigg(\int_{F_3}|\zeta|^{(\frac{\gamma
}{p}-m-2)q}
d\mu(\zeta)\bigg)^{1/q} \\
&\leq& A\varepsilon y|z|^{\frac{\gamma}{p}+\frac{1}{q}-2}.
\hspace{78mm} (2.5.16)
\end{eqnarray*}

  Finally, by (4) of Lemma 2.5.2, we obtain
$$
|h_4(z)|\leq \frac{Ay}{|z|^{2}}
\int_{F_4}{\eta}d\mu(\zeta).\eqno{(2.5.17)}
$$

  Thus, by collecting (2.5.12), (2.5.13), (2.5.14),
(2.5.15), (2.5.16) and (2.5.17), there exists a positive constant
$A$ independent of $\varepsilon$, such that if $ |z|\geq
2R_\varepsilon$ and $\  z \notin E_2(\varepsilon)$, we have
$$
 |h(z)|\leq A\varepsilon y^{1-\frac{\alpha}{p}}|z|^{\frac{\gamma}{p}+\frac{1}{q}-2+\frac{\alpha}{p}}.
$$

  Similarly, if $z\notin G$, we have
$$
h(z)=
o(y^{1-\frac{\alpha}{p}}|z|^{\frac{\gamma}{p}+\frac{1}{q}-2+\frac{\alpha}{p}}),
\quad {\rm as} \ |z|\rightarrow\infty. \eqno{(2.5.18)}
$$

  By (2.5.4) and (2.5.18), we obtain that
$$
u(z)=v(z)+h(z)=
o(y^{1-\frac{\alpha}{p}}|z|^{\frac{\gamma}{p}+\frac{1}{q}-2+\frac{\alpha}{p}}),
\quad {\rm as} \ |z|\rightarrow\infty
$$
holds in ${{\bf C}_{+}}-G$, thus we complete the proof of Theorem
2.5.2.

\newpage
\ \thispagestyle{empty}
\newpage

\chapter{Growth Estimates for a Class of Subharmonic Functions in the Half Space}

\section{Introduction and Basic Notations}

  For $x\in{\bf R}^{n}\backslash\{0\}$, let \cite{HN}
$$
 E(x)=-r_n|x|^{2-n},
$$
where $|x|$ is the Euclidean norm, $r_n=\frac{1}{(n-2)\omega_{n}}$
and $\omega_{n}=\frac{2\pi^{\frac{n}{2}}}{\Gamma(\frac{n}{2})}$ is
the surface area of the unit sphere in ${\bf R}^{n} $. We know that
$E$ is locally integrable in ${\bf R}^{n} $.

 The Green function $G(x,y)$ for the upper half space
 $H$ is given by \cite{HN}
$$
 G(x,y)=E(x-y)-E(x-y^{\ast}), \qquad x,y\in\overline{H} ,\  x\neq y,
\eqno{(3.1.1)}
$$
where $^{\ast}$ denotes the reflection in the boundary plane
$\partial H$ just as $y^{\ast}=(y_1,y_2,\cdots,y_{n-1},-y_n)$, then
we define the Poisson kernel $P(x,y')$ when $x\in H$ and $y'\in
\partial H $ by
$$
 P(x,y')=-\frac{\partial G(x,y)}{\partial
 y_n}\bigg|_{y_n=0}=\frac{2x_n}{\omega_n|x-(y',0)|^n}.
\eqno{(3.1.2)}
$$

  The  Dirichlet problem of the upper half space is to find a function
 $u$ satisfying
$$
 u\in C^2(H), \eqno{(3.1.3)}
$$
$$
 \Delta u=0,   x\in H, \eqno{(3.1.4)}
$$
$$
 \lim_{x\rightarrow x'}u(x)=f(x')
 \ {\rm nontangentially  \  a.e.}x'\in \partial H, \eqno{(3.1.5)}
$$
where $f$ is a measurable function of ${\bf R}^{n-1} $. The Poisson
integral of the upper half space is defined by
$$
u(x)=P[f](x)=\int_{{\bf R}^{n-1}}P(x,y')f(y')dy', \eqno{(3.1.6)}
$$
where $P(x,y')$ is defined by (3.1.2).

  As we all know, the Poisson integral $P[f]$ exists if
$$
\int_{{\bf R}^{n-1}}\frac{|f(y')|}{1+|y'|^n} dy'<\infty.
$$(see \cite{ABR}, \cite{F} and \cite{MS})In this chapter,
we replace the condition into
$$
\int_{{\bf R}^{n-1}}\frac{|f(y')|^p}{(1+|y'|)^{\gamma}} dy'<\infty,
\eqno{(3.1.7)}
$$
where $1\leq p<\infty$ and $\gamma$ is a real number, then we can
get the asymptotic behaviour of harmonic functions.

  Next, we will generalize these results to subharmonic functions.

\section{Preliminary Lemma}

Let $\mu$ be a positive Borel measure  in ${\bf R}^n,\ \beta\geq0$,
the maximal function $M(d\mu)(x)$ of order $\beta$ is defined by
$$
M(d\mu)(x)=\sup_{ 0<r<\infty}\frac{\mu(B(x,r))}{r^\beta},
$$
then the maximal function $M(d\mu)(x):{\bf R}^n \rightarrow
[0,\infty)$ is lower semicontinuous, hence measurable. To see this,
for any $ \lambda >0 $, let $D(\lambda)=\{x\in{\bf
R}^{n}:M(d\mu)(x)>\lambda\}$. Fix $x \in D(\lambda)$, then there
exists
 $r>0$ such that $\mu(B(x,r))>tr^\beta$ for some $t>\lambda$, and
there exists $ \delta>0$ satisfying
$(r+\delta)^\beta<\frac{tr^\beta}{\lambda}$. If $|y-x|<\delta$, then
$B(y,r+\delta)\supset B(x,r)$, therefore $\mu(B(y,r+\delta))\geq
tr^\beta >\lambda(r+\delta)^\beta$. Thus $B(x,\delta)\subset
D(\lambda)$. This proves that $D(\lambda)$ is open for each
$\lambda>0$.

 In order to obtain the results, we
need the lemma below:

\vspace{0.2cm}
 \noindent
{\bf Lemma 3.2.1 } {\it Let $\mu$ be a positive Borel measure  in
${\bf R}^n,\ \beta\geq0,\ \mu({\bf R}^n)<\infty,$ for any $ \lambda
\geq 5^{\beta} \mu({\bf R}^n)$, set
$$
E(\lambda)=\{x\in{\bf R}^{n}:|x|\geq2,M(d\mu)(x) >
\frac{\lambda}{|x|^{\beta}}\},
$$
then there exists $ x_j\in E(\lambda),\  \ \rho_j> 0,\
j=1,2,\cdots$, such that
$$
E(\lambda) \subset \bigcup_{j=1}^\infty B(x_j,\rho_j) \eqno{(3.2.1)}
$$
and
$$
\sum _{j=1}^{\infty}\frac{\rho_j^{\beta}}{|x_j|^{\beta}}\leq
\frac{3\mu({\bf R}^n)5^{\beta}}{\lambda} .\eqno{(3.2.2)}
$$
}

 Proof: Let $E_k(\lambda)=\{x\in E(\lambda):2^k\leq |x|<2^{k+1}\}$,
then  for any $ x \in E_k(\lambda),$ there exists $ r(x)>0$, such
that $\mu(B(x,r(x))) >\lambda \big(\frac{r(x)}{|x|}\big)^{\beta} $,
therefore $r(x)\leq 2^{k-1}$.
 Since $E_k(\lambda)$ can be covered
by
 the union of a family of balls $\{B(x,r(x)):x\in E_k(\lambda) \}$,
 by the Vitali Lemma \cite{SS}, there exists $  \Lambda_k\subset E_k(\lambda)$,
$\Lambda_k$ is at most countable, such that $\{B(x,r(x)):x\in
\Lambda_k \}$ are disjoint and
$$
E_k(\lambda) \subset
 \cup_{x\in \Lambda_k} B(x,5r(x)),
$$
so
$$
E(\lambda)=\cup_{k=1}^\infty E_k(\lambda) \subset \cup_{k=1}^\infty
\cup_{x\in \Lambda_k} B(x,5r(x)).
$$

  On the other hand, note that $ \cup_{x\in \Lambda_k} B(x,r(x)) \subset \{x:2^{k-1}\leq
|x|<2^{k+2}\} $, so that
$$
 \sum_{x \in \Lambda_k}\frac{(5r(x))^{\beta}}{|x|^{\beta}}
\leq 5^\beta\sum_{x\in\Lambda_k}\frac{\mu(B(x,r(x)))}{\lambda} \leq
\frac{5^\beta}{\lambda} \mu\{x:2^{k-1}\leq |x|<2^{k+2}\}.
$$
Hence we obtain
$$
 \sum _{k=1}^{\infty}\sum
_{x \in \Lambda_k}\frac{(5r(x))^{\beta}}{|x|^{\beta}}
 \leq
 \sum _{k=1}^{\infty}\frac{5^\beta}{\lambda} \mu\{x:2^{k-1}\leq |x|<2^{k+2}\}
 \leq
\frac{3\mu({\bf R}^n)5^{\beta}}{\lambda}.
$$
  Rearrange $ \{x:x \in \Lambda_k,k=1,2,\cdots\} $ and $
\{5r(x):x \in \Lambda_k,k=1,2,\cdots\}
 $, we get $\{x_j\}$
and $\{\rho_j\}$ such that (3.2.1) and (3.2.2) hold.

\section{$p=1$}

 \section*{ 1. Introduction and Main Theorems}

\vspace{0.3cm}

 In this section, we will
consider measurable functions $f$ in ${\bf R}^{n-1}$ satisfying (see
\cite{ABR}, \cite{F} and \cite{MS})
$$
\int_{{\bf R}^{n-1}}\frac{|f(y')|}{1+|y'|^{n+m}} dy'<\infty,
\eqno{(3.3.1)}
$$
where $m$ is a nonnegative integer. This is just (3.1.7) when $p=1$
and $\gamma=n+m$. It is well known that the Poisson kernel $P(x,y')$
has a series expansion in terms of the ultraspherical ( or
Gegenbauer ) polynomials $ C^{\lambda}_{k}(t)\ (\lambda
=\frac{n}{2})$(see \cite{SW} and \cite{S}). The latter can be
defined by a generating function
$$
(1-2tr+r^2)^{-\lambda } = \sum _{k=0}^{\infty}
C_{k}^{\lambda}(t)r^k, \eqno{(3.3.2)}
$$
where $|r|<1$, $ |t|\leq 1$ and $ \lambda > 0$. The coefficients $
C^{\lambda}_{k}(t) $ is called the ultraspherical ( or Gegenbauer )
polynomial of degree $ k $  associated with $ \lambda $, the
function $ C^{\lambda}_{k}(t) $ is a  polynomial of degree $ k $ in
$ t $.
 To obtain a solution of  Dirichlet problem for the boundary date $f$, as in
  \cite{STS}, \cite{STU},  \cite{Y} and \cite{MS}, we use the following modified functions defined by
$$
 E_m(x-y)=\left\{\begin{array}{ll}
 E(x-y)  &   \mbox{when }   |y|\leq 1,  \\
 E(x-y) +\sum_{k=0}^{m-1}\frac{r_n|x|^k}{|y|^{n-2+k}}C^{\frac{n-2}{2}}_{k}
 \left( \frac{x\cdot y}{|x||y|}\right ) &
\mbox{when}\   |y|> 1.
 \end{array}\right.
$$
Then we can define the modified Green function $G_m(x,y)$ and the
modified Poisson
 kernel $P_m(x,y')$ by (see  \cite{HRH}, \cite{HRS},  \cite{DI}, \cite{J} and \cite{MS})
$$
 G_m(x,y)=E_{m+1}(x-y)-E_{m+1}(x-y^{\ast}), \qquad x,y\in\overline{H}, \ x\neq
 y; \eqno{(3.3.3)}
$$
$$
 P_m(x,y')=\left\{\begin{array}{ll}
 P(x,y')  &   \mbox{when }   |y'|\leq 1  ,\\
 P(x,y') - \sum_{k=0}^{m-1}\frac{2x_n|x|^k}{\omega_n|y'|^{n+k}}C^{n/2}_{k}
 \left( \frac{x\cdot (y',0)}{|x||y'|}\right )&
\mbox{when}\   |y'|> 1.
 \end{array}\right.\eqno{(3.3.4)}
$$

   Siegel-Talvila \cite{STS} have proved the
  following result:

\vspace{0.2cm}
 \noindent
{\bf Theorem C }  {\it Let $f$ be a measurable function in ${\bf
R}^{n-1}$ satisfying (3.3.1), then the harmonic function
$$
v(x)= \int_{{\bf R}^{n-1}}P_m(x,y')f(y')dy', \quad x\in H
\eqno{(3.3.5)}
$$
satisfies (3.1.3), (3.1.4), (3.1.5) and
$$
v(x)= o(x_n^{1-n}|x|^{m+n}),  \quad  {\rm as}  \
|x|\rightarrow\infty,  \eqno{(3.3.6)}
$$
where $P_m(x,y')$ is defined by (3.3.4). }

  In order to describe the asymptotic behaviour
of subharmonic functions in the half space (see \cite{ZDG},
\cite{ZG}, \cite{ML}, \cite{MLJ1} and \cite{MLJ2}), we establish the
following theorems.

\vspace{0.2cm}
 \noindent
{\bf Theorem 3.3.1 }  {\it Let $f$ be a measurable function in ${\bf
R}^{n-1}$ satisfying (3.3.1), and $0< \alpha\leq n$. Let $v(x)$ be
the harmonic function defined by
 (3.3.5). Then there exists $x_j\in H,\ \rho_j>0,$ such that
$$
\sum
_{j=1}^{\infty}\frac{\rho_j^{n-\alpha}}{|x_j|^{n-\alpha}}<\infty
\eqno{(3.3.7)}
$$
holds and
$$
v(x)= o(x_n^{1-\alpha}|x|^{m+\alpha}),  \quad  {\rm as}  \
|x|\rightarrow\infty   \eqno{(3.3.8)}
$$
holds in $H-G$, where $ G=\bigcup_{j=1}^\infty B(x_j,\rho_j)$. }

\vspace{0.2cm}
 \noindent
 {\bf Remark 3.3.1 }  {\it If $\alpha=n$, then (3.3.7) is a finite sum,
 the set $G$ is the union of finite balls, so (3.3.6) holds in $H$.
 This is just the result of Siegel-Talvila, therefore,
 our result (3.3.8) is the generalization of Theorem C.
}

  Next, we will generalize Theorem 3.3.1 to subharmonic functions.

\vspace{0.2cm}
 \noindent
{\bf Theorem 3.3.2 }  {\it Let $f$ be a measurable function in ${\bf
R}^{n-1}$ satisfying (3.3.1) and $\mu$ be a positive  Borel measure
satisfying
$$
\int_H\frac{y_n}{1+|y|^{n+m}} d\mu(y)<\infty. \eqno{(3.3.9)}
$$
Write the subharmonic function
$$
u(x)= v(x)+h(x), \quad x\in H,
$$
where $v(x)$ is the harmonic function defined by (3.3.5), $h(x)$ is
defined by
$$
h(x)= \int_H G_m(x,y)d\mu(y)
$$
and $G_m(x,y)$ is defined by (3.3.3). Then there exists $x_j\in H,\
\rho_j>0,$ such that (3.3.7) holds and
$$
u(x)= o(x_n^{1-\alpha}|x|^{m+\alpha}),  \quad  {\rm as}  \
|x|\rightarrow\infty
 $$
holds in $H-G$, where $ G=\bigcup_{j=1}^\infty B(x_j,\rho_j)$ and
$0< \alpha<2$. }

  Next we are concerned with minimal thinness \cite{AE} at
infinity for $v(x)$ and $h(x)$, for a set $E \subset H$ and an open
set $F \subset {\bf R}^{n-1}$, we consider the capacity
$$
C(E;F)= \inf\int_{{\bf R}^{n-1}}g(y')dy',
$$
where the infimum is taken over all nonnegative measurable functions
$g$ such that $g=0$ outside $F$ and
$$
\int_{{\bf R}^{n-1}}\frac{g(y')}{|x-(y',0)|^n}dy'\geq 1, \qquad
{\rm for \ all}\ x\in E.
$$

  We say that $E \subset H$ is minimally thin at infinity if
$$
\sum _{i=1}^{\infty}2^{-in}C(E_i;F_i)<\infty,
$$
where $E_i=\{x\in E:\; 2^i\leq |x|<2^{i+1}\}$ and $F_i=\{x\in {\bf
R}^{n-1}:\; 2^i< |x|<2^{i+3}\}$.

\vspace{0.2cm}
 \noindent
{\bf Theorem 3.3.3 }  {\it Let $f$ be a measurable function in ${\bf
R}^{n-1}$ satisfying (3.3.1), then there exists a set $E \subset H$
such that $E$ is minimally thin at infinity and
$$
 \lim_{|x|\rightarrow\infty,x\in H-E}\frac{v(x)}{x_n|x|^m}=0.
$$
}

  Similarly, for $h(x)$, we can also conclude the following:

\vspace{0.2cm}
 \noindent
{\bf Corollary 3.3.1 }  {\it Let $\mu$ be a positive  Borel  measure
satisfying (3.3.9), then there exists a set $E \subset H$ such that
$E$ is minimally thin at infinity and
$$
 \lim_{|x|\rightarrow\infty,x\in H-E}\frac{h(x)}{x_n|x|^m}=0.
$$
}

Finally we are concerned with rarefiedness \cite{AE} at infinity for
$v(x)$ and $h(x)$, for a set $E \subset H$ and an open set $F
\subset H$, we consider the capacity
$$
C(E;F)= \inf\int_{H}g(y)d\mu(y),
$$
where the infimum is taken over all nonnegative measurable functions
$g$ such that $g=0$ outside $F$ and
$$
\int_{H}\frac{g(y)}{|x-y|^{n-1}}d\mu(y)\geq 1, \qquad  {\rm for \
all}\  x\in E.
$$

  We say that $E \subset H$ is rarefied at infinity if
$$
\sum _{i=1}^{\infty}2^{-i(n-1)}C(E_i;F_i)<\infty,
$$
where $E_i$ is as in Theorem 3.3.3 and $F_i=\{x\in H:\; 2^i<
|x|<2^{i+3}\}$.

\vspace{0.2cm}
 \noindent
{\bf Theorem 3.3.4 }  {\it Let $\mu$ be a positive  Borel  measure
satisfying (3.3.9), then there exists a set $E \subset H$ such that
$E$ is rarefied at infinity and
$$
 \lim_{|x|\rightarrow\infty,x\in H-E}\frac{h(x)}{|x|^{m+1}}=0.
$$
}

  Similarly, for $v(x)$, we can also conclude the following:

\vspace{0.2cm}
 \noindent
{\bf Corollary 3.3.2 }  {\it Let $f$ be a measurable function in
${\bf R}^{n-1}$ satisfying (3.3.1), then there exists a set $E
\subset H$ such that $E$ is rarefied at infinity and
$$
 \lim_{|x|\rightarrow\infty,x\in H-E}\frac{v(x)}{|x|^{m+1}}=0.
$$
}

\section*{ 2.  Main Lemmas}

  In order to obtain the results, we need the following lemmas:

\vspace{0.2cm}
 \noindent
{\bf Lemma 3.3.1 }  {\it Gegenbauer polynomials have the following
properties:\\
{\rm (1)}\ $|C_k^{\lambda }(t)|\leq C_k^{\lambda }(1)=\frac{\Gamma
           (2\lambda +k)}{\Gamma (2\lambda )\Gamma(k+1)}, \ \ |t|\leq 1 ;$\\
{\rm (2)}\ $\frac{d}{dt}C_k^\lambda (t)=2\lambda C_{k-1}^{\lambda+1}(t),\ \ k \geq 1;$\\
{\rm (3)}\ $ \sum _{k=0}^{\infty} C_k^\lambda (1)r^k=(1-r)^{-2\lambda};$\\
{\rm (4)}\ $ |C^{\frac{n-2}{2}}_{k}\left(t\right)-C^{\frac{n-2}{2}}_{k}
            \left( t^{\ast}\right )| \leq(n-2)C^{n/2}_{k-1} \left( 1\right)
           |t-t^{\ast}|,\ \ |t|\leq 1, \ \ |t^{\ast}|\leq 1$.\\

Proof: (1) and (2) can be derived from \cite{SW} and \cite{G}; (3)
follows by taking $t=1$ in (3.3.2); (4) follows
 by (1), (2) and the Mean Value Theorem for Derivatives.

\vspace{0.2cm}
 \noindent
{\bf Lemma 3.3.2 }  {\it The Green function $G(x,y)$ has the
following
estimates:\\
{\rm (1)}\ $|G(x,y)|\leq \frac{r_n}{|x-y|^{n-2}};$\\
{\rm (2)}\ $|G(x,y)|\leq \frac{2x_ny_n}{\omega_n|x-y|^n};$\\
{\rm (3)}\ $|G(x,y)|\leq \frac{Ax_ny_n}{|x-y|^{n-2}|x-y^{\ast}|^2}.$\\

Proof: (1) is obvious; (2) follows by the Mean Value Theorem for
Derivatives; (3) can be derived from \cite{AE}.

\section*{3.   Proof of Theorems }

\vspace{0.3cm}

 \emph{Proof of Theorem 3.3.1}

Define the measure $dm(y')$ and the kernel $K(x,y')$ by
$$
dm(y')=\frac{|f(y')|}{1+|y'|^{n+m}} dy' ,\ \ K(x,y')=
P_m(x,y')(1+|y'|^{n+m}).
$$
  For any $\varepsilon >0$, there exists $R_\varepsilon >2$, such that
$$
\int_{|y'|\geq
R_\varepsilon}dm(y')\leq\frac{\varepsilon}{5^{n-\alpha}}.
$$
For every Lebesgue measurable set $E \subset {\bf R}^{n-1}$, the
measure $m^{(\varepsilon)}$ defined by $m^{(\varepsilon)}(E)
=m(E\cap\{x'\in{\bf R}^{n-1}:|x'|\geq R_\varepsilon\}) $ satisfies
$m^{(\varepsilon)}({\bf
R}^{n-1})\leq\frac{\varepsilon}{5^{n-\alpha}}$, write
\begin{eqnarray*}
& &v_1(x)=\int_{|x-(y',0)| \leq 3|x|} P(x,y')(1+|y'|^{n+m})
dm^{(\varepsilon)}(y'),\\
& &v_2(x)=\int_{|x-(y',0)| \leq 3|x|}
(P_m(x,y')-P(x,y'))(1+|y'|^{n+m})
dm^{(\varepsilon)}(y'), \\
& &v_3(x)=\int_{|x-(y',0)| > 3|x|} K(x,y')dm^{(\varepsilon)}(y'), \\
& &v_4(x)=\int_{1<|y'|<R_\varepsilon}K(x,y') dm(y'), \\
& &v_5(x)=\int_{|y'|\leq1}K(x,y') dm(y'), \\
\end{eqnarray*}
then
$$
|v(x)| \leq |v_1(x)|+|v_2(x)|+|v_3(x)|+|v_4(x)|+|v_5(x)|.
\eqno{(3.3.10)}
$$
Let $ E_1(\lambda)=\{x\in{\bf R}^{n}:|x|\geq2,\exists \ t>0, s.t.
m^{(\varepsilon)}(B(x,t)\cap{\bf R}^{n-1}
)>\lambda(\frac{t}{|x|})^{n-\alpha}\}$, therefore, if $ |x|\geq
2R_\varepsilon$ and $ x \notin E_1(\lambda)$, then we have

\begin{eqnarray*}
|v_1(x)|
&\leq& \int_{x_n\leq|x-(y',0)| \leq
3|x|}\frac{2x_n}{\omega_n|x-(y',0)|^n}2|y'|^{n+m} dm^{(\varepsilon)}(y') \\
&\leq& \frac{4^{n+m+1}}{\omega_n}x_n|x|^{m+n}\int_{x_n}^{3|x|}
\frac{1}{t^n} dm_x^{(\varepsilon)}(t) \\
&\leq& \frac{4^{n+m+1}}{\omega_n}
\bigg(\frac{1}{3^\alpha}+\frac{n}{\alpha}\bigg)\lambda
x_n^{1-\alpha}|x|^{m+\alpha}, \hspace{45mm} (3.3.11)
\end{eqnarray*}
where  $m_x^{(\varepsilon)}(t)=\int_{|x-(y',0)| \leq t}
dm^{(\varepsilon)}(y')$.

   By (1) and (3) of Lemma 3.3.1, we obtain
\begin{eqnarray*}
|v_2(x)|
&\leq& \int_{x_n\leq|x-(y',0)| \leq 3|x|}
\sum_{k=0}^{m-1}\frac{2x_n|x|^k}{\omega_n}C^{n/2}_{k}
 \left( 1\right )\frac{2|y'|^{n+m}}{|y'|^{n+k}} dm^{(\varepsilon)}(y') \\
&\leq&
\frac{4^{m+1}}{\omega_n}\sum_{k=0}^{m-1}\frac{1}{4^k}C^{n/2}_{k}
 \left( 1\right )\frac{1}{5^{n-\alpha}}\varepsilon x_n|x|^m \\
&\leq& \frac{4^{m+1+\alpha}}{\omega_n\cdot 3^n}\varepsilon x_n|x|^m.
\hspace{67mm} (3.3.12)
\end{eqnarray*}

  By (1) and (3) of Lemma 3.3.1, we see that \cite{HK}
\begin{eqnarray*}
|v_3(x)|
&\leq& \int_{|x-(y',0)| > 3|x|}
\sum_{k=m}^{\infty}\frac{4x_n|x|^k}{\omega_n(2|x|)^{k-m}}C^{n/2}_{k}
\left( 1\right ) dm^{(\varepsilon)}(y') \\
&\leq& \frac{2^{m+2}}{\omega_n}\frac{\varepsilon}{5^{n-\alpha}}
\sum_{k=m}^{\infty}\frac{1}{2^k}C^{n/2}_{k}
\left( 1\right )x_n|x|^m \\
&\leq& \frac{2^{m-n+2\alpha+2}}{\omega_n}\varepsilon x_n|x|^m .
\hspace{62mm} (3.3.13)
\end{eqnarray*}

  Write
\begin{eqnarray*}
v_4(x)
&=& \int_{1<|y'|<R_\varepsilon}[P(x,y') +(P_m(x,y')-P(x,y'))](1+|y'|^{n+m}) dm(y') \\
&=& v_{41}(x)+v_{42}(x),
\end{eqnarray*}
then
\begin{eqnarray*}
|v_{41}(x)|
&\leq& \int_{1<|y'|<R_\varepsilon}\frac{2x_n}{\omega_n|x-(y',0)|^n}
2|y'|^{n+m} dm(y') \\
&\leq& \frac{4R_\varepsilon^{n+m}x_n}{\omega_n}
\int_{1<|y'|<R_\varepsilon}\frac{1}{(\frac{|x|}{2})^n}
 dm(y') \\
&\leq& \frac{2^{n+2}R_\varepsilon^{n+m} m({\bf R}^{n-1})}{\omega_n}
\frac{x_n}{|x|^n}. \hspace{50mm} (3.3.14)
\end{eqnarray*}

  Moreover, by (1) and (3) of Lemma 3.3.1, we obtain
\begin{eqnarray*}
|v_{42}(x)|
&\leq& \int_{1<|y'|<R_\varepsilon}
\sum_{k=0}^{m-1}\frac{2x_n|x|^k}{\omega_n|y'|^{n+k}}C^{n/2}_{k}
 \left( 1\right )\cdot 2|y'|^{n+m} dm(y') \\
&\leq& \sum_{k=0}^{m-1}\frac{4}{\omega_n}C^{n/2}_{k} \left( 1\right
)x_n|x|^kR_\varepsilon^{m-k} m({\bf
R}^{n-1}) \\
&\leq& \frac{2^{n+m+1}R_\varepsilon^m m({\bf R}^{n-1})}{\omega_n}
x_n|x|^{m-1}.\hspace{41mm} (3.3.15)
\end{eqnarray*}

  In case $|y'|\leq 1$, note that
$$
K(x,y')=P_m(x,y')(1+|y'|^{n+m})
\leq\frac{4x_n}{\omega_n|x-(y',0)|^n},
$$
so that
$$
|v_5(x)|\leq \int_{|y'|\leq1}\frac{4x_n}{\omega_n
\big(\frac{|x|}{2}\big)^n} dm(y') \leq\frac{2^{n+2} m({\bf
R}^{n-1})}{\omega_n} \frac{x_n}{|x|^n}.\eqno{(3.3.16)}
$$

  Thus, by collecting (3.3.10), (3.3.11), (3.3.12), (3.3.13), (3.3.14), (3.3.15) and
(3.3.16), there exists a positive constant $A$ independent of
$\varepsilon$, such that if $ |x|\geq 2R_\varepsilon$ and $\  x
\notin E_1(\varepsilon)$, we have
$$
|v(x)|\leq A\varepsilon x_n^{1-\alpha}|x|^{m+\alpha}.
$$

 Let $\mu_\varepsilon$ be a measure in ${\bf R}^n$ defined by
$ \mu_\varepsilon(E)= m^{(\varepsilon)}(E\cap{\bf R}^{n-1})$ for
every measurable set $E$ in ${\bf R}^n$. Take
$\varepsilon=\varepsilon_p=\frac{1}{2^{p+2}}, p=1,2,3,\cdots$, then
there exists a sequence $ \{R_p\}$: $1=R_0<R_1<R_2<\cdots$ such that
$$
\mu_{\varepsilon_p}({\bf R}^n)=\int_{|y'|\geq
R_p}dm(y')<\frac{\varepsilon_p}{5^{n-\alpha}}.
$$
Take $\lambda=3\cdot5^{n-\alpha}\cdot2^p\mu_{\varepsilon_p}({\bf
R}^n)$ in Lemma 3.2.1, then there exists $ x_{j,p}$ and $
\rho_{j,p}$, where $R_{p-1}\leq |x_{j,p}|<R_p,$ such that
$$
\sum
_{j=1}^{\infty}\bigg(\frac{\rho_{j,p}}{|x_{j,p}|}\bigg)^{n-\alpha}
\leq \frac{1}{2^{p}}.
$$
So if $R_{p-1}\leq |x|<R_p$ and $ x\notin G_p=\cup_{j=1}^\infty
B(x_{j,p},\rho_{j,p})$, we have
$$
|v(x)|\leq A\varepsilon_px_n^{1-\alpha}|x|^{m+\alpha},
$$
thereby
$$
\sum _{p=1}^{\infty}
\sum_{j=1}^{\infty}\bigg(\frac{\rho_{j,p}}{|x_{j,p}|}\bigg)^{n-\alpha}
\leq \sum _{p=1}^{\infty}\frac{1}{2^{p}}=1<\infty.
$$
Set $ G=\cup_{p=1}^\infty G_p$, thus Theorem 3.3.1 holds.

 \emph{Proof of Theorem 3.3.2}

  Define the measure $dn(y)$ and the kernel $L(x,y)$ by
$$
dn(y)=\frac{y_n d\mu(y)}{1+|y|^{n+m}},\ \
L(x,y)=G_m(x,y)\frac{1+|y|^{n+m}}{y_n},
$$
then the function $h(x)$ can be written as
$$
h(x)=\int_H L(x,y) dn(y).
$$

  For any $\varepsilon >0$, there exists $R_\varepsilon >2$, such that
$$
\int_{|y|\geq R_\varepsilon}dn(y)<\frac{\varepsilon}{5^{n-\alpha}}.
$$
For every Lebesgue measurable set $E \subset {\bf R}^{n}$, the
measure $n^{(\varepsilon)}$ defined by $n^{(\varepsilon)}(E)
=n(E\cap\{y\in H:|y|\geq R_\varepsilon\}) $ satisfies
$n^{(\varepsilon)}(H)\leq\frac{\varepsilon}{5^{n-\alpha}}$, write
\begin{eqnarray*}
&h_1(x)& =\int_{|x-y|\leq\frac{x_n}{2}}G(x,y)\frac{1+|y|^{n+m}}{y_n}
dn^{(\varepsilon)}(y), \\
&h_2(x)&=\int_{\frac{x_n}{2}<|x-y|\leq3|x|}G(x,y)\frac{1+|y|^{n+m}}{y_n}
dn^{(\varepsilon)}(y), \\
&h_3(x)&=\int_{|x-y|\leq3|x|}(G_m(x,y)-G(x,y))\frac{1+|y|^{n+m}}{y_n}
dn^{(\varepsilon)}(y), \\
&h_4(x)&=\int_{|x-y|>3|x|}L(x,y)
dn^{(\varepsilon)}(y), \\
&h_5(x)&=\int_{1<|y|<R_\varepsilon}L(x,y) dn(y), \\
&h_6(x)&=\int_{|y|\leq1}L(x,y) dn(y),\\
\end{eqnarray*}
then
$$
h(x)=h_1(x)+h_2(x)+h_3(x)+h_4(x)+h_5(x)+h_6(x). \eqno{(3.3.17)}
$$
Let $ E_2(\lambda)=\{x\in{\bf R}^{n}:|x|\geq2,\exists \ t>0, s.t.
n^{(\varepsilon)}(B(x,t)\cap H
)>\lambda(\frac{t}{|x|})^{n-\alpha}\}, $ therefore, if $ |x|\geq
2R_\varepsilon$ and $ x \notin E_2(\lambda)$, then we have by (1) of
Lemma 3.3.2
\begin{eqnarray*}
|h_1(x)|
&\leq& \int_{|x-y|\leq\frac{x_n}{2}}
\frac{r_n}{|x-y|^{n-2}}\frac{2|y|^{n+m}}{\frac{x_n}{2}}
dn^{(\varepsilon)}(y) \\
&\leq& 4\times
(3/2)^{n+m}r_n\frac{|x|^{n+m}}{x_n}\int_0^\frac{x_n}{2}
\frac{1}{t^{n-2}} dn_x^{(\varepsilon)}(t)\\
&\leq& 4\times (3/2)^{n+m}r_n\bigg[\frac{1}{2^{2-\alpha}}+
\frac{n-2}{(2-\alpha)2^{2-\alpha}}\bigg]\lambda
x_n^{1-\alpha}|x|^{m+\alpha}, \hspace{12mm} (3.3.18)
\end{eqnarray*}
where $ n_x^{(\varepsilon)}(t)=\int_{|x-y| \leq t}
dn^{(\varepsilon)}(y)$.

 By (2) of Lemma 3.3.2, we have
\begin{eqnarray*}
|h_2(x)|
&\leq&
\int_{\frac{x_n}{2}<|x-y|\leq3|x|}\frac{2x_ny_n}{\omega_n|x-y|^n}
\frac{2|y|^{n+m}}{y_n}
dn^{(\varepsilon)}(y) \\
&\leq&
\frac{4^{n+m+1}}{\omega_n}x_n|x|^{n+m}\int_\frac{x_n}{2}^{3|x|}
\frac{1}{t^n} dn_x^{(\varepsilon)}(t)\\
&\leq& \frac{4^{n+m+1}}{\omega_n}\bigg(\frac{1}{3^\alpha}+
\frac{n2^\alpha}{\alpha}\bigg)\lambda x_n^{1-\alpha}|x|^{m+\alpha}.
\hspace{40mm} (3.3.19)
\end{eqnarray*}

  First note $C^{\lambda}_{0}
\left( t\right )\equiv 1$ \cite{SW}, then we obtain by (1), (3) and
(4) of Lemma 3.3.1 and taking $t=\frac{x\cdot y}{|x||y|},\
t^{\ast}=\frac{x\cdot y^{\ast}}{|x||y^{\ast}|}$ in (4) of Lemma
3.3.1

\begin{eqnarray*}
|h_3(x)|
&\leq&
\int_{|x-y|\leq3|x|}\sum_{k=1}^{m}\frac{r_n|x|^k}{|y|^{n-2+k}}
2(n-2)C^{n/2}_{k-1} \left( 1\right )\frac{x_ny_n}{|x||y|}
\frac{2|y|^{n+m}}{y_n} dn^{(\varepsilon)}(y) \\
&\leq&
\frac{4^{m+1}}{\omega_n}\sum_{k=1}^{m}\frac{1}{4^{k-1}}C^{n/2}_{k-1}
 \left( 1\right )\frac{1}{5^{n-\alpha}}\varepsilon x_n|x|^m \\
&\leq& \frac{4^{m+1+\alpha}}{\omega_n\cdot 3^n}\varepsilon x_n|x|^m.
\hspace{73mm} (3.3.20)
\end{eqnarray*}

  By (1), (3) and (4) of Lemma 3.3.1, we see that
\begin{eqnarray*}
|h_4(x)|
&\leq&
\int_{|x-y|>3|x|}\sum_{k=m+1}^{\infty}\frac{r_n|x|^k}{|y|^{n-2+k}}
2(n-2)C^{n/2}_{k-1} \left( 1\right )\frac{x_ny_n}{|x||y|}
\frac{2|y|^{n+m}}{y_n} dn^{(\varepsilon)}(y) \\
&\leq&
\frac{2^{m+2}}{\omega_n}\sum_{k=m+1}^{\infty}\frac{1}{2^{k-1}}C^{n/2}_{k-1}
 \left( 1\right )\frac{1}{5^{n-\alpha}}\varepsilon x_n|x|^m \\
&\leq& \frac{2^{m-n+2\alpha+2}}{\omega_n}\varepsilon x_n|x|^m .
\hspace{67mm} (3.3.21)
\end{eqnarray*}

  Write
\begin{eqnarray*}
h_5(x)
&=&
\int_{1<|y|<R_\varepsilon}[G(x,y)+(G_m(x,y)-G(x,y))]\frac{1+|y|^{n+m}}{y_n} dn(y) \\
&=& h_{51}(x)+h_{52}(x),
\end{eqnarray*}
then we obtain by (2) of Lemma 3.3.2
\begin{eqnarray*}
|h_{51}(x)|
&\leq& \int_{1<|y|<R_\varepsilon}
\frac{2x_ny_n}{\omega_n|x-y|^n}\frac{2|y|^{n+m}}{y_n} dn(y) \\
&\leq&
\frac{4R_\varepsilon^{n+m}}{\omega_n}x_n\int_{1<|y|<R_\varepsilon}
\frac{1}{(\frac{|x|}{2})^n} dn(y) \\
&\leq& \frac{2^{n+2}R_\varepsilon^{n+m}n(H) }{\omega_n}
\frac{x_n}{|x|^n} .\hspace{65mm} (3.3.22)
\end{eqnarray*}

  Moreover, by (1), (3) and (4) of Lemma 3.3.1, we obtain
\begin{eqnarray*}
|h_{52}(x)|
&\leq&
\int_{1<|y|<R_\varepsilon}\sum_{k=1}^{m}\frac{r_n|x|^k}{|y|^{n-2+k}}
2(n-2)C^{n/2}_{k-1} \left( 1\right )\frac{x_ny_n}{|x||y|}
\frac{2|y|^{n+m}}{y_n} dn(y) \\
&\leq& \sum_{k=1}^{m}\frac{4}{\omega_n}C^{n/2}_{k-1}
 \left( 1\right )x_n|x|^{k-1}R_\varepsilon^{m-k+1}n(H) \\
&\leq& \frac{2^{n+m+1}R_\varepsilon^{m}n(H) }{\omega_n} x_n|x|^{m-1}
.\hspace{56mm} (3.3.23)
\end{eqnarray*}

  In case $|y|\leq 1$, by (2) of Lemma 3.3.2, we have
$$
|L(x,y)|\leq \frac{2x_ny_n}{\omega_n|x-y|^n}\frac{2}{y_n}
=\frac{4x_n}{\omega_n|x-y|^n},
$$
so that
$$
|h_6(x)|\leq \int_{|y|\leq1}\frac{4x_n}{\omega_n(\frac{|x|}{2})^n}
dn(y)\leq \frac{2^{n+2}n(H) }{\omega_n}
\frac{x_n}{|x|^n}.\eqno{(3.3.24)}
$$

  Thus, by collecting (3.3.17), (3.3.18), (3.3.19), (3.3.20),
(3.3.21), (3.3.22), (3.3.23) and (3.3.24), there exists a positive
constant $A$ independent of $\varepsilon$, such that if $ |x|\geq
2R_\varepsilon$ and $\  x \notin E_2(\varepsilon)$, we have
$$
 |h(x)|\leq A\varepsilon x_n^{1-\alpha}|x|^{m+\alpha}.
$$

  Similarly, if $x\notin G$, we have
$$
h(x)= o(x_n^{1-\alpha}|x|^{m+\alpha}), \quad  {\rm as} \
|x|\rightarrow\infty. \eqno{(3.3.25)}
$$

By (3.3.8) and  (3.3.25), we obtain that
$$
u(x)=v(x)+h(x)= o(x_n^{1-\alpha}|x|^{m+\alpha}),  \quad  {\rm as} \
|x|\rightarrow\infty
$$
holds in $H-G$, thus we complete the proof of Theorem 3.3.2.

 \emph{Proof of Theorem 3.3.3 and 3.3.4}

We prove only Theorem 3.3.4, the proof of Theorem 3.3.3 is similar.
By (3.3.20), (3.3.21), (3.3.22), (3.3.23) and (3.3.24) we have
$$
 \lim_{|x|\rightarrow\infty,x\in H}\frac{h_3(x)+h_4(x)+h_5(x)+h_6(x)}{|x|^{m+1}}=0.
 \eqno{(3.3.26)}
$$
In view of (3.3.9), we can find a sequence $\{a_i\}$ of positive
numbers such that $\lim_{i\rightarrow\infty}a_i=\infty$ and
$$
\sum
_{i=1}^{\infty}a_i\int_{F_i}\frac{y_n}{|y|^{n+m}}d\mu(y)<\infty.
$$
Consider the sets
$$
E_i=\{x\in H:\; 2^i\leq |x|<2^{i+1},|h_1(x)+h_2(x)|\geq
a_i^{-1}2^{im}|x|\}
$$
for $i=1,2,\cdots$. If $x\in E_i$, then we obtain by (3) of Lemma
3.3.2
$$
a_i^{-1}\leq 2^{-im}|x|^{-1}|h_1(x)+h_2(x)|\leq
A2^{-i(m+1)}\int_{F_i}\frac{y_n}{|x-y|^{n-1}}d\mu(y),
$$
so that it follows from the definition of $C(E_i;F_i)$ that
$$
C(E_i;F_i)\leq Aa_i2^{-i(m+1)}\int_{F_i}y_nd\mu(y)\leq
Aa_i2^{i(n-1)}\int_{F_i}\frac{y_n}{|y|^{n+m}}d\mu(y).
$$
Define $E=\bigcup_{i=1}^\infty E_i$, then
$$
\sum _{i=1}^{\infty}2^{-i(n-1)}C(E_i;F_i)<\infty.
$$
Clearly,
$$
 \lim_{|x|\rightarrow\infty,x\in H-E}\frac{h_1(x)+h_2(x)}{|x|^{m+1}}=0.
 \eqno{(3.3.27)}
$$
Thus, by collecting (3.3.26) and (3.3.27), the proof of Theorem
3.3.4 is completed.

\section{$p>1$(General Kernel)}

 \section*{ 1. Introduction and Main Theorems}

\vspace{0.3cm}

  In this section, we will consider measurable functions $f$ in ${\bf
R}^{n-1}$ satisfying
$$
\int_{{\bf R}^{n-1}}\frac{|f(y')|^p}{(1+|y'|)^{\gamma}} dy'<\infty,
\eqno{(3.4.1)}
$$
where $\gamma$ is defined as in Theorem 3.4.1.

  In order to describe the asymptotic behaviour of subharmonic functions
in the upper half space (see \cite{ML}, \cite{MLJ1}, and
\cite{MLJ2}),
 we establish the following theorems.

\vspace{0.2cm}
 \noindent
{\bf Theorem 3.4.1 } {\it Let $1\leq p<\infty,\
\frac{1}{p}+\frac{1}{q}=1$ and
$$
-(n-1)(p-1)<\gamma <(n-1)+p  \quad  {\rm in \  case}  \ p>1;
$$
$$
 0<\gamma \leq n  \quad  {\rm in \  case}  \ p=1.
$$
If $f$ is a measurable function in ${\bf R}^{n-1}$ satisfying
(3.4.1) and $v(x)$ is the harmonic function defined by (3.1.6), then
there exists $x_j\in H,\ \rho_j>0,$ such that
$$
\sum
_{j=1}^{\infty}\frac{\rho_j^{pn-\alpha}}{|x_j|^{pn-\alpha}}<\infty
\eqno{(3.4.2)}
$$
holds and
$$
v(x)=
o(x_n^{1-\frac{\alpha}{p}}|x|^{\frac{\gamma}{p}+\frac{n-1}{q}-n+\frac{\alpha}{p}}),
\quad  {\rm as}  \ |x|\rightarrow\infty   \eqno{(3.4.3)}
$$
holds in $H-G$, where $ G=\bigcup_{j=1}^\infty B(x_j,\rho_j)$ and
$0< \alpha\leq np$. }

\vspace{0.2cm}
 \noindent
 {\bf Remark 3.4.1 } {\it If $\alpha=n$, $p=1$ and $\gamma=n$, then (3.4.2) is a finite sum,
 the set $G$ is the union of finite balls, so (3.4.3) holds in $H$. This is just
 the case $m=0$ of the result of Siegel-Talvila.
}

\vspace{0.2cm}
 \noindent
 {\bf Remark 3.4.2 } {\it If $\gamma=-(n-1)(p-1)$, $p>1$, then
$$
v(x)=
o(x_n^{1-\frac{\alpha}{p}}(\log|x|)^{\frac{1}{q}}
|x|^{\frac{\gamma}{p}+\frac{n-1}{q}-n+\frac{\alpha}{p}}),
\quad  {\rm as}  \ |x|\rightarrow\infty
$$
holds in $H-G$. }

  Next, we will generalize Theorem 3.4.1 to subharmonic functions.

\vspace{0.2cm}
 \noindent
{\bf Theorem 3.4.2 } {\it Let $p$ and $\gamma$ be as in Theorem
3.4.1. If $f$ is a measurable function in ${\bf R}^{n-1}$ satisfying
(3.4.1) and $\mu$ is a positive Borel measure satisfying
$$
\int_H\frac{y_n^p}{(1+|y|)^\gamma} d\mu(y)<\infty
$$
and
$$
\int_H\frac{1}{(1+|y|)^{n-1}} d\mu(y)<\infty.
$$
Write the subharmonic function
$$
u(x)= v(x)+h(x), \quad x\in H,
$$
where $v(x)$ is the harmonic function defined by (3.1.6), $h(x)$ is
defined by
$$
h(x)= \int_H G(x,y)d\mu(y)
$$
and $G(x,y)$ is defined by (3.1.1). Then there exists $x_j\in H,\
\rho_j>0,$ such that (3.4.2) holds and
$$
u(x)=
o(x_n^{1-\frac{\alpha}{p}}|x|^{\frac{\gamma}{p}+\frac{n-1}{q}-n+\frac{\alpha}{p}}),
\quad  {\rm as}  \ |x|\rightarrow\infty
$$
holds in $H-G$, where $ G=\bigcup_{j=1}^\infty B(x_j,\rho_j)$ and
$0< \alpha<2p$. }

\vspace{0.2cm}
 \noindent
 {\bf Remark 3.4.3 } {\it If $\gamma=-(n-1)(p-1)$, $p>1$, then
$$
u(x)= o(x_n^{1-\frac{\alpha}{p}}(\log|x|)^{\frac{1}{q}}
|x|^{\frac{\gamma}{p}+\frac{n-1}{q}-n+\frac{\alpha}{p}}),
 \quad {\rm as}  \ |x|\rightarrow\infty
$$
holds in $H-G$. }

\section*{ 2. Main Lemmas}

\vspace{0.3cm}

 In order to obtain the results, we
need these lemmas below:

\vspace{0.2cm}
 \noindent
{\bf Lemma 3.4.1 } {\it The kernel function $\frac{1}{|x-y|^n}$ has
the
following estimates:\\
{\rm (1)}\ If $|y|\leq \frac{|x|}{2}$, then $\frac{1}{|x-y|^n}\leq
\frac{2^n}{|x|^n}$;\\
{\rm (2)}\ If $|y|> 2|x|$, then $\frac{1}{|x-y|^n}\leq
\frac{2^n}{|y|^n}$.\\
}

\vspace{0.2cm}
 \noindent
{\bf Lemma 3.4.2 } {\it The Green function $G(x,y)$ has the
following estimates:\\
{\rm (1)}\ $|G(x,y)|\leq \frac{r_n}{|x-y|^{n-2}}$;\\
{\rm (2)}\ $|G(x,y)|\leq \frac{2x_ny_n}{\omega_n|x-y|^n}$.\\
}

Proof: (1) is obvious; (2) follows by the Mean Value Theorem for
Derivatives.

\vspace{0.4cm}

\section*{3.   Proof of Theorems}

 \emph{Proof of Theorem 3.4.1}

 We prove only the case $p>1$; the proof of the case $p=1$ is similar.
Define the measure $dm(y')$ by
$$
dm(y')=\frac{|f(y')|^p}{(1+|y'|)^{\gamma}} dy'.
$$
  For any $\varepsilon >0$, there exists $R_\varepsilon >2$, such that
$$
\int_{|y'|\geq
R_\varepsilon}dm(y')\leq\frac{\varepsilon^p}{5^{pn-\alpha}}.
$$
For every Lebesgue measurable set $E \subset {\bf R}^{n-1}$, the
measure $m^{(\varepsilon)}$ defined by $m^{(\varepsilon)}(E)
=m(E\cap\{x'\in{\bf R}^{n-1}:|x'|\geq R_\varepsilon\}) $ satisfies
$m^{(\varepsilon)}({\bf
R}^{n-1})\leq\frac{\varepsilon^p}{5^{pn-\alpha}}$, write
\begin{eqnarray*}
&v_1(x)& =\int_{G_1} P(x,y')f(y')
dy',\\
&v_2(x)&=\int_{G_2} P(x,y')f(y')
dy', \\
&v_3(x)&=\int_{G_3} P(x,y')f(y')
dy', \\
&v_4(x)&=\int_{G_4} P(x,y')f(y')
dy', \\
\end{eqnarray*}
where
\begin{eqnarray*}
&G_1& =\{y'\in {\bf R}^{n-1}: R_\varepsilon<|y'|\leq \frac{|x|}{2}\},\\
&G_2& =\{y'\in {\bf R}^{n-1}: \frac{|x|}{2}<|y'| \leq 2|x|\}, \\
&G_3& =\{y'\in {\bf R}^{n-1}: |y'|>2|x|\}, \\
&G_4& =\{y'\in {\bf R}^{n-1}: |y'|\leq R_\varepsilon\}. \\
\end{eqnarray*}
Then
$$
v(x) =v_1(x)+v_2(x)+v_3(x)+v_4(x). \eqno{(3.4.4)}
$$

  First, if $\gamma >-(n-1)(p-1)$, then $\frac{\gamma q}{p}
  +(n-1)>0$, so that we obtain by (1) of Lemma 3.4.1 and H\"{o}lder's
inequality
\begin{eqnarray*}
|v_1(x)|
&\leq& \int_{G_1}\frac{ 2x_n}{\omega_n}\frac{2^n}{|x|^n}|f(y')| dy' \\
&\leq& \frac{ 2^{n+1}}{\omega_n}\frac{x_n}{|x|^n}
\bigg(\int_{G_1}\frac{|f(y')|^p}{|y'|^\gamma}
dy'\bigg)^{1/p}\bigg(\int_{G_1}|y'|^{\frac{\gamma q}{p}}
dy'\bigg)^{1/q}, \\
\end{eqnarray*}
since
$$
\int_{G_1}|y'|^{\frac{\gamma q}{p}} dy'\leq
\omega_{n-1}\frac{1}{\frac{\gamma
q}{p}+n-1}\bigg(\frac{|x|}{2}\bigg)^{\frac{\gamma q}{p}+n-1} ,
$$
so that
$$
|v_1(x)|\leq A \varepsilon
x_n|x|^{\frac{\gamma}{p}+\frac{n-1}{q}-n}. \eqno{(3.4.5)}
$$

  Let $ E_1(\lambda)=\{x\in{\bf R}^{n}:|x|\geq2,\exists \ t>0, s.t.
m^{(\varepsilon)}(B(x,t)\cap{\bf R}^{n-1}
)>\lambda^p(\frac{t}{|x|})^{pn-\alpha}\}$, therefore, if $ |x|\geq
2R_\varepsilon$ and $x \notin E_1(\lambda)
 $, then we have
$$
\forall t>0,\ m^{(\varepsilon)}(B(x,t)\cap{\bf R}^{n-1}
)\leq\lambda^p \bigg(\frac{t}{|x|}\bigg)^{pn-\alpha}.
$$

  If $\gamma >-(n-1)(p-1)$, then $\frac{\gamma q}{p}
  +(n-1)>0$, so that we obtain by H\"{o}lder's inequality
\begin{eqnarray*}
|v_2(x)|
&\leq& \frac{2x_n}{\omega_n}
\bigg(\int_{G_2}\frac{|f(y')|^p}{|x-(y',0)|^{pn}|y'|^\gamma}
dy'\bigg)^{1/p}\bigg(\int_{G_2}|y'|^{\frac{\gamma q}{p}}
dy'\bigg)^{1/q} \\
&\leq& Ax_n|x|^{\frac{\gamma}{p}+\frac{n-1}{q}}
\bigg(\int_{G_2}\frac{|f(y')|^p}{|x-(y',0)|^{pn}|y'|^\gamma}
dy'\bigg)^{1/p},\\
\end{eqnarray*}
since
\begin{eqnarray*}
\int_{G_2}\frac{|f(y')|^p}{|x-(y',0)|^{pn}|y'|^\gamma} dy'
&\leq& \int_{x_n}^{3|x|}
\frac{2^\gamma+1}{t^{pn}} dm_x^{(\varepsilon)}(t) \\
&\leq& \frac{\lambda^p}{
|x|^{pn}}(2^\gamma+1)\bigg(\frac{1}{3^\alpha}+
\frac{pn}{\alpha}\bigg)\frac{|x|^\alpha}{x_n^\alpha}, \\
\end{eqnarray*}
where  $m_x^{(\varepsilon)}(t)=\int_{|x-(y',0)| \leq t}
dm^{(\varepsilon)}(y')$.\\
Hence we have
$$
|v_2(x)|\leq A \lambda
x_n^{1-\frac{\alpha}{p}}|x|^{\frac{\gamma}{p}+\frac{n-1}{q}-n+\frac{\alpha}{p}}.
\eqno{(3.4.6)}
$$

  If $\gamma <(n-1)+p$, then $(\frac{\gamma}{p}-n)q+(n-1)<0$,
so that we obtain by (2) of Lemma 3.4.1 and H\"{o}lder's inequality

\begin{eqnarray*}
|v_3(x)|
&\leq& \int_{G_3}\frac{ 2x_n}{\omega_n}\frac{2^n}{|y'|^n}|f(y')| dy' \\
&\leq& \frac{ 2^{n+1}}{\omega_n}x_n
\bigg(\int_{G_3}\frac{|f(y')|^p}{|y'|^\gamma}
dy'\bigg)^{1/p}\bigg(\int_{G_3}|y'|^{(\frac{\gamma }{p}-n)q}
dy'\bigg)^{1/q} \\
&\leq& A\varepsilon x_n|x|^{\frac{\gamma}{p}+\frac{n-1}{q}-n}.
 \hspace{96mm} (3.4.7)
\end{eqnarray*}

  Finally, by (1) of Lemma 3.4.1, we obtain
$$
|v_4(x)|\leq \frac{ 2^{n+1}}{\omega_n}\frac{x_n}{|x|^n}
\int_{G_4}{|f(y')|}dy',
$$
which implies by $\gamma >-(n-1)(p-1)$ that
$$
|v_4(x)|\leq A \varepsilon
x_n|x|^{\frac{\gamma}{p}+\frac{n-1}{q}-n}. \eqno{(3.4.8)}
$$

  Thus, by collecting (3.4.4), (3.4.5), (3.4.6), (3.4.7) and
(3.4.8), there exists a positive constant $A$ independent of
$\varepsilon$, such that if $ |x|\geq 2R_\varepsilon$ and $\  x
\notin E_1(\varepsilon)$, we have
$$
|v(x)|\leq A\varepsilon
x_n^{1-\frac{\alpha}{p}}|x|^{\frac{\gamma}{p}+\frac{n-1}{q}-n+\frac{\alpha}{p}}.
$$

 Let $\mu_\varepsilon$ be a measure in ${\bf R}^n$ defined by
$ \mu_\varepsilon(E)= m^{(\varepsilon)}(E\cap{\bf R}^{n-1})$ for
every measurable set $E$ in ${\bf R}^n$. Take
$\varepsilon=\varepsilon_p=\frac{1}{2^{p+2}}, p=1,2,3,\cdots$, then
there exists a sequence $ \{R_p\}$: $1=R_0<R_1<R_2<\cdots$ such that
$$
\mu_{\varepsilon_p}({\bf R}^n)=\int_{|y'|\geq
R_p}dm(y')<\frac{\varepsilon_p^p}{5^{pn-\alpha}}.
$$
Take $\lambda=3\cdot5^{pn-\alpha}\cdot2^p\mu_{\varepsilon_p}({\bf
R}^n)$ in Lemma 3.2.1, then there exists $x_{j,p}$ and $
\rho_{j,p}$, where $R_{p-1}\leq |x_{j,p}|<R_p$, such that
$$
\sum
_{j=1}^{\infty}\bigg(\frac{\rho_{j,p}}{|x_{j,p}|}\bigg)^{pn-\alpha}
\leq \frac{1}{2^{p}}.
$$
If $R_{p-1}\leq |x|<R_p$ and $x\notin G_p=\cup_{j=1}^\infty
B(x_{j,p},\rho_{j,p})$, we have
$$
|v(x)|\leq
A\varepsilon_px_n^{1-\frac{\alpha}{p}}|x|^{\frac{\gamma}{p}+\frac{n-1}{q}-n+\frac{\alpha}{p}}.
$$
Thereby
$$
\sum _{p=1}^{\infty}
\sum_{j=1}^{\infty}\bigg(\frac{\rho_{j,p}}{|x_{j,p}|}\bigg)^{pn-\alpha}
\leq \sum _{p=1}^{\infty}\frac{1}{2^{p}}=1<\infty.
$$
Set $ G=\cup_{p=1}^\infty G_p$, thus Theorem 3.4.1 holds.

 \emph{Proof of Theorem 3.4.2}

 We prove only the case $p>1$; the remaining case $p=1$ can be proved similarly.
Define the measure $dn(y)$ by
$$
dn(y)=\frac{y_n^p}{(1+|y|)^{\gamma}} d\mu(y).
$$

  For any $\varepsilon >0$, there exists $R_\varepsilon >2$, such that
$$
\int_{|y|\geq
R_\varepsilon}dn(y)<\frac{\varepsilon^p}{5^{pn-\alpha}}.
$$
For every Lebesgue measurable set $E \subset {\bf R}^{n}$, the
measure $n^{(\varepsilon)}$ defined by $n^{(\varepsilon)}(E)
=n(E\cap\{y\in H:|y|\geq R_\varepsilon\}) $ satisfies
$n^{(\varepsilon)}(H)\leq\frac{\varepsilon^p}{5^{pn-\alpha}}$, write
\begin{eqnarray*}
&h_1(x)& = \int_{F_1} G(x,y)d\mu(y),\\
&h_2(x)&=\int_{F_2} G(x,y)d\mu(y),\\
&h_3(x)&=\int_{F_3} G(x,y)d\mu(y), \\
&h_4(x)&=\int_{F_4} G(x,y)d\mu(y), \\
\end{eqnarray*}
where
\begin{eqnarray*}
&F_1& =\{y\in H: R_\varepsilon<|y|\leq \frac{|x|}{2}\},\\
&F_2& =\{y\in H: \frac{|x|}{2}<|y| \leq 2|x|\}, \\
&F_3& =\{y\in H: |y|>2|x|\}, \\
&F_4& =\{y\in H: |y|\leq R_\varepsilon\}. \\
\end{eqnarray*}
Then
$$
h(x) =h_1(x)+h_2(x)+h_3(x)+h_4(x). \eqno{(3.4.9)}
$$

  First, if $\gamma >-(n-1)(p-1)$, then $\frac{\gamma q}{p}
  +(n-1)>0$, so that we obtain by (1) of Lemma 3.4.1, (2) of Lemma
3.4.2 and H\"{o}lder's inequality
\begin{eqnarray*}
|h_1(x)|
&\leq& \int_{F_1}\frac{2x_ny_n}{\omega_n|x-y|^n} d\mu(y) \\
&\leq& \int_{F_1}\frac{2x_ny_n}{\omega_n}\frac{2^n}{|x|^n} d\mu(y) \\
&\leq& \frac{ 2^{n+1}}{\omega_n}\frac{x_n}{|x|^n}
\bigg(\int_{F_1}\frac{y_n^p}{|y|^\gamma}
d\mu(y)\bigg)^{1/p}\bigg(\int_{F_1}|y|^{\frac{\gamma q}{p}}
d\mu(y)\bigg)^{1/q}, \\
\end{eqnarray*}
since
$$
\int_{F_1}|y|^{\frac{\gamma q}{p}} d\mu(y)\leq
 2^{n-1}\bigg(\frac{|x|}{2}\bigg)^{\frac{\gamma
q}{p}+n-1}\int_{H}\frac{1}{(1+|y|)^{n-1}} d\mu(y),
$$
so that
$$
|h_1(x)|\leq A\varepsilon x_n|x|^{\frac{\gamma}{p}+\frac{n-1}{q}-n}.
\eqno{(3.4.10)}
$$

  Let $ E_2(\lambda)=\{x\in{\bf R}^{n}:|x|\geq2,\exists \ t>0, s.t.
n^{(\varepsilon)}(B(x,t)\cap H
)>\lambda^p(\frac{t}{|x|})^{pn-\alpha}\}, $ therefore, if $ |x|\geq
2R_\varepsilon$ and $x\notin E_2(\lambda)
 $, then we have
$$
\forall t>0, \ n^{(\varepsilon)}(B(x,t)\cap H
)\leq\lambda^p\bigg(\frac{t}{|x|}\bigg)^{pn-\alpha}.
$$

  If $\gamma >-(n-1)(p-1)$, then $\frac{\gamma q}{p}
  +(n-1)>0$, so that we obtain by H\"{o}lder's inequality
\begin{eqnarray*}
|h_2(x)|
&\leq& \bigg(\int_{F_2}\frac{|G(x,y)|^p}{|y|^\gamma}
d\mu(y)\bigg)^{1/p}\bigg(\int_{F_2}|y|^{\frac{\gamma q}{p}}
d\mu(y)\bigg)^{1/q} \\
&\leq& \bigg((2^\gamma +1)\int_{F_2}\frac{|G(x,y)|^p}{y_n^p}
dn(y)\bigg)^{1/p}\bigg(\int_{F_2}|y|^{\frac{\gamma q}{p}}
d\mu(y)\bigg)^{1/q}\\
&\leq& A|x|^{\frac{\gamma}{p}+\frac{n-1}{q}}\bigg(\int_{F_2}\frac{|G(x,y)|^p}{y_n^p}dn(y)\bigg)^{1/p},\\
\end{eqnarray*}
since
\begin{eqnarray*}
\int_{F_2}\frac{|G(x,y)|^p}{y_n^p} dn(y)
&\leq&  \int_{|x-y|\leq 3|x|}\frac{|G(x,y)|^p}{y_n^p} dn^{(\varepsilon)}(y) \\
&=&  \int_{|x-y|\leq \frac{x_n}{2}}\frac{|G(x,y)|^p}{y_n^p}
dn^{(\varepsilon)}(y)+
 \int_{\frac{x_n}{2}<|x-y|\leq 3|x|}\frac{|G(x,y)|^p}{y_n^p} dn^{(\varepsilon)}(y)\\
&=& h_{21}(x)+h_{22}(x),
\end{eqnarray*}
so that we have by (1) of Lemma 3.4.2
\begin{eqnarray*}
h_{21}(x)
&\leq& \int_{|x-y|\leq
\frac{x_n}{2}}\bigg(\frac{2r_n}{x_n|x-y|^{n-2}}\bigg)^p
dn^{(\varepsilon)}(y) \\
&=& \bigg(\frac{2r_n}{x_n}\bigg)^p\int_0^{\frac{x_n}{2}}
\frac{1}{t^{p(n-2)}} dn_x^{(\varepsilon)}(t) \\
&\leq& (2r_n)^p\frac{np-\alpha}{(2p-\alpha)2^{2p-\alpha}}
\lambda^p\frac{x_n^{p-\alpha}}{|x|^{np-\alpha}}.\\
\end{eqnarray*}

  Moreover, we have by (2) of Lemma 3.4.2
\begin{eqnarray*}
h_{22}(x)
&\leq& \int_{\frac{x_n}{2}<|x-y|\leq 3|x|
}\bigg(\frac{2x_n}{\omega_n|x-y|^{n}}\bigg)^p
dn^{(\varepsilon)}(y) \\
&=& \bigg(\frac{2x_n}{\omega_n}\bigg)^p\int_{\frac{x_n}{2}}^{3|x|}
\frac{1}{t^{pn}} dn_x^{(\varepsilon)}(t) \\
&\leq&
\bigg(\frac{2}{\omega_n}\bigg)^p\bigg(\frac{1}{3^\alpha}+\frac{np2^\alpha}{\alpha}\bigg)
\lambda^p\frac{x_n^{p-\alpha}}{|x|^{np-\alpha}},\\
\end{eqnarray*}
where  $n_x^{(\varepsilon)}(t)=\int_{|x-y| \leq t}
dn^{(\varepsilon)}(y)$.\\
Hence we have
$$
|h_2(x)|\leq A \lambda
x_n^{1-\frac{\alpha}{p}}|x|^{\frac{\gamma}{p}+\frac{n-1}{q}-n+\frac{\alpha}{p}}.
\eqno{(3.4.11)}
$$

  If $\gamma <(n-1)+p$, then $(\frac{\gamma}{p}-n)q
  +(n-1)<0$, so that we obtain by (2) of Lemma 3.4.1, (2) of Lemma 3.4.2 and H\"{o}lder's inequality

\begin{eqnarray*}
|h_3(x)|
&\leq& \int_{F_3}\frac{2x_ny_n}{\omega_n|x-y|^n} d\mu(y) \\
&\leq& \int_{F_3}\frac{ 2x_ny_n}{\omega_n}\frac{2^n}{|y|^n} d\mu(y) \\
&\leq& \frac{ 2^{n+1}}{\omega_n}x_n
\bigg(\int_{F_3}\frac{y_n^p}{|y|^\gamma}
d\mu(y)\bigg)^{1/p}\bigg(\int_{F_3}|y|^{(\frac{\gamma }{p}-n)q}
d\mu(y)\bigg)^{1/q} \\
&\leq& A\varepsilon x_n|x|^{\frac{\gamma}{p}+\frac{n-1}{q}-n}.
 \hspace{96mm} (3.4.12)
\end{eqnarray*}

  Finally, by (1) of Lemma 3.4.1 and (2) of Lemma 3.4.2, we obtain
$$
|h_4(x)|\leq \int_{F_4}\frac{2x_ny_n}{\omega_n|x-y|^n} d\mu(y) \leq
\frac{ 2^{n+1}}{\omega_n}\frac{x_n}{|x|^n} \int_{F_4}y_n d\mu(y),
$$
which implies by $\gamma >-(n-1)(p-1)$ that
$$
|h_4(x)|\leq A \varepsilon
x_n|x|^{\frac{\gamma}{p}+\frac{n-1}{q}-n}.  \eqno{(3.4.13)}
$$

  Thus, by collecting (3.4.9), (3.4.10),
(3.4.11), (3.4.12) and (3.4.13), there exists a positive constant
$A$ independent of $\varepsilon$, such that if $ |x|\geq
2R_\varepsilon$ and $\  x \notin E_2(\varepsilon)$, we have
$$
 |h(x)|\leq A\varepsilon x_n^{1-\frac{\alpha}{p}}|x|^{\frac{\gamma}{p}+\frac{n-1}{q}-n+\frac{\alpha}{p}}.
$$

  Similarly, if $x\notin G$, we have
$$
h(x)=
o(x_n^{1-\frac{\alpha}{p}}|x|^{\frac{\gamma}{p}+\frac{n-1}{q}-n+\frac{\alpha}{p}}),
\quad {\rm as} \ |x|\rightarrow\infty. \eqno{(3.4.14)}
$$

  By (3.4.3) and  (3.4.14), we obtain that
$$
u(x)=v(x)+h(x)=
o(x_n^{1-\frac{\alpha}{p}}|x|^{\frac{\gamma}{p}+\frac{n-1}{q}-n+\frac{\alpha}{p}}),
\quad {\rm as} \ |x|\rightarrow\infty
$$
holds in $H-G$, thus we complete the proof of Theorem 3.4.2.

\section{the Estimates for the Modified Poisson Kernel
 and Green Function}

 \section*{ 1. Introduction and Main Theorems}

\vspace{0.3cm}

 Recall that the modified Poisson kernel
$P_m(x,y')$ and the modified Green function  $G_m(x,y)$ (see
\cite{HRH}, \cite{HRS},  \cite{DI}, \cite{J} and \cite{MS}) are
defined respectively by
$$
 P_m(x,y')=\left\{\begin{array}{ll}
 P(x,y')  &   \mbox{when }   |y'|\leq 1  ,\\
 P(x,y') - \sum_{k=0}^{m-1}\frac{2x_n|x|^k}{\omega_n|y'|^{n+k}}C^{n/2}_{k}
 \left( \frac{x\cdot (y',0)}{|x||y'|}\right )&
\mbox{when}\   |y'|> 1
 \end{array}\right.
$$
and
$$
 G_m(x,y)=E_{m+1}(x-y)-E_{m+1}(x-y^{\ast}), \qquad x,y\in\overline{H}, \ x\neq
 y,
$$
where
$$
 E_m(x-y)=\left\{\begin{array}{ll}
 E(x-y)  &   \mbox{when }   |y|\leq 1,  \\
 E(x-y) +\sum_{k=0}^{m-1}\frac{r_n|x|^k}{|y|^{n-2+k}}C^{\frac{n-2}{2}}_{k}
 \left( \frac{x\cdot y}{|x||y|}\right ) &
\mbox{when}\   |y|> 1.
 \end{array}\right.
$$

  In our discussions, the estimates for the modified Poisson kernel
$P_m(x,y')$ and the modified Green function  $G_m(x,y)$ are
fundamental, therefore, we establish the following theorems.

 \vspace{0.2cm}
 \noindent
{\bf Theorem 3.5.1 } {\it Suppose $|y'|>1$, then we have the
estimates:
$$
 |P_m(x,y')|\leq\left\{\begin{array}{ll}
 A\frac{x_n}{|x-y'|^n}s'^{m+n-1} , &   \mbox{when }   s'>1  ,\\
   A\frac{x_n}{|x-y'|^n}s'^m,& \mbox{when}\   s'\leq1,
 \end{array}\right.
$$
where $s'=\frac{|x|}{|y'|}$. }

\vspace{0.2cm}
 \noindent
{\bf Theorem 3.5.2 }  {\it Suppose $|y|>1$, then we have the
estimates:
$$
|G_m(x,y)| \leq\left\{\begin{array}{ll}
 A\frac{x_ny_n}{|x-y|^{n-2}}\frac{|x|^m}{|y|^{m+1}}
 \big(\frac{|x|}{|y|^2}
 +\frac{1}{|y|}+\frac{|x|}{|x-y|^2} \big),
  &   \mbox{when }   s\leq1  ,\\
 A\frac{x_ny_n}{|x-y|^{n-2}}\frac{|x|^{m+n-4}}{|y|^{m+n-2}}
  \big(1 +\frac{|x|}{|y|}+\frac{|x|^2}{|x-y|^2} \big),& \mbox{when}\   s>1,
 \end{array}\right.
$$
where $s=\frac{|x|}{|y|}$. }

 \section*{ 2. Main Lemma}

\vspace{0.3cm}

   In order to obtain the results, we need the lemma below:

\vspace{0.2cm}
 \noindent
{\bf Lemma 3.5.1 } {\it Suppose $|y'|>1$, set $s'=\frac{|x|}{|y'|}$
and $t'= \frac{x\cdot y'}{|x||y'|} $, then
$$
P_m(x,y')=P(x,y')[mC_m^{n/2}(t')I_{m-1}^{(n)}(s',t')-(n+m-1)C_{m-1}^{n/2}(t')I_m^{(n)}(s',t')],
$$
where
$$
I_m^{(n)}(s',t')=\int_0^{s'}(1-2t'\xi+\xi^2)^{n/2-1}\xi^m d\xi,
\quad s'>0,\ |t'|<1.
$$
}

\vspace{0.4cm}

\section*{3.   Proof of Theorems }

 \emph{Proof of Theorem 3.5.1}

Suppose $|y'|>1$, since
\begin{eqnarray*}
I_m^{(n)}(s',t')
&=& \int_0^{s'}(1-2t'\xi+\xi^2)^{n/2-1}\xi^m d\xi \\
&\leq& \left\{\begin{array}{ll}
 As'^{m+1} , &   \mbox{when }   s'\leq1  ,\\
 As'^{m+n-1},& \mbox{when}\   s'>1,
 \end{array}\right.
\end{eqnarray*}
we can obtain by Lemma 3.5.1
\begin{eqnarray*}
|P_m(x,y')|
&\leq& P(x,y')[m|C_m^{n/2}(t')||I_{m-1}^{(n)}(s',t')|
+(n+m-1)|C_{m-1}^{n/2}(t')||I_m^{(n)}(s',t')|] \\
&\leq& P(x,y')[mA|I_{m-1}^{(n)}(s',t')|+(n+m-1)A|I_m^{(n)}(s',t')|].
\end{eqnarray*}
When $s'>1$,
$$
|P_m(x,y')|\leq A\frac{x_n}{|x-y'|^n}s'^{m+n-1};
$$
when $s'\leq1$,
$$
|P_m(x,y')|\leq A\frac{x_n}{|x-y'|^n}s'^m.
$$
Thus
$$
 |P_m(x,y')|\leq\left\{\begin{array}{ll}
 A\frac{x_n}{|x-y'|^n}s'^{m+n-1} , &   \mbox{when }   s'>1  ,\\
   A\frac{x_n}{|x-y'|^n}s'^m,& \mbox{when}\   s'\leq1.
 \end{array}\right.
$$

\emph{Proof of Theorem 3.5.2}

Suppose $|y|>1$, by Lemma 3.5.1, we obtain
$$
P_m(x,y)=P(x,y)[mC_m^{n/2}(t)I_{m-1}^{(n)}(s,t)
-(n+m-1)C_{m-1}^{n/2}(t)I_m^{(n)}(s,t)],
$$
where $s=\frac{|x|}{|y|}$, $t= \frac{x\cdot y}{|x||y|}.$\\
Thus
$$
E_m(x-y)=E(x-y)[mC_m^{\frac{n-2}{2}}(t)I_{m-1}^{(n-2)}(s,t)
-(n+m-3)C_{m-1}^{\frac{n-2}{2}}(t)I_m^{(n-2)}(s,t)].
$$
Similarly, we can obtain
$$
E_m(x-y^{\ast})=E(x-y^{\ast})[mC_m^{\frac{n-2}{2}}(t^{\ast})I_{m-1}^{(n-2)}(s^{\ast},t^{\ast})
-(n+m-3)C_{m-1}^{\frac{n-2}{2}}(t^{\ast})I_m^{(n-2)}(s^{\ast},t^{\ast})],
$$
so that
\begin{eqnarray*}
G_m(x,y)
&=& (m+1)[E(x-y)C_{m+1}^{\frac{n-2}{2}}(t)I_m^{(n-2)}(s,t)\\
& & -E(x-y^{\ast})C_{m+1}^{\frac{n-2}{2}}(t^{\ast})I_m^{(n-2)}(s^{\ast},t^{\ast})] \\
& & -(n+m-2)[E(x-y)C_m^{\frac{n-2}{2}}(t)I_{m+1}^{(n-2)}(s,t)\\
& & -E(x-y^{\ast})C_m^{\frac{n-2}{2}}(t^{\ast})I_{m+1}^{(n-2)}(s^{\ast},t^{\ast})] \\
&=& (m+1)[f-f^{\ast}]-(n+m-2)[g-g^{\ast}], \hspace{33mm} (3.5.1)
\end{eqnarray*}
where
\begin{eqnarray*}
&  & f-f^{\ast}\\
& =& E(x-y)C_{m+1}^{\frac{n-2}{2}}(t)I_m^{(n-2)}(s,t)
-E(x-y)C_{m+1}^{\frac{n-2}{2}}(t)I_m^{(n-2)}(s,t^{\ast}) \\
& & +E(x-y)C_{m+1}^{\frac{n-2}{2}}(t)I_m^{(n-2)}(s^{\ast},t^{\ast})
-E(x-y)C_{m+1}^{\frac{n-2}{2}}(t^{\ast})I_m^{(n-2)}(s^{\ast},t^{\ast}) \\
&
&+E(x-y)C_{m+1}^{\frac{n-2}{2}}(t^{\ast})I_m^{(n-2)}(s^{\ast},t^{\ast})
-E(x-y^{\ast})C_{m+1}^{\frac{n-2}{2}}(t^{\ast})I_m^{(n-2)}(s^{\ast},t^{\ast})\\
&=& I_1+I_2+I_3.  \hspace{96mm} (3.5.2)
\end{eqnarray*}

For the first term, we have
\begin{eqnarray*}
I_1
&=& E(x-y)C_{m+1}^{\frac{n-2}{2}}(t)[I_m^{(n-2)}(s,t)
-I_m^{(n-2)}(s,t^{\ast})] \\
&=& E(x-y)C_{m+1}^{\frac{n-2}{2}}(t)
\bigg[\int_0^s(1-2t\xi+\xi^2)^{\frac{n-2}{2}-1}\xi^m d\xi
-\int_0^s(1-2t^{\ast}\xi+\xi^2)^{\frac{n-2}{2}-1}\xi^m d\xi \bigg] \\
&=&\frac{-2(n-4)x_ny_n}{|x||y|}
E(x-y)C_{m+1}^{\frac{n-2}{2}}(t)I_{m+1}^{(n-4)}(s,t_0),
\end{eqnarray*}
where $t^{\ast}<t_0<t$, thus
$$
|I_1|\leq\left\{\begin{array}{ll}
 A\frac{x_ny_n}{|x-y|^{n-2}}\frac{|x|^{m+1}}{|y|^{m+3}} , &   \mbox{when }   s\leq1  ,\\
  A\frac{x_ny_n}{|x-y|^{n-2}}\frac{|x|^{m+n-5}}{|y|^{m+n-3}},& \mbox{when}\
  s>1;
 \end{array}\right. \eqno{(3.5.3)}
$$
for the second term, we have
$$
I_2=E(x-y)I_m^{(n-2)}(s^{\ast},t^{\ast})[C_{m+1}^{\frac{n-2}{2}}(t)-C_{m+1}^{\frac{n-2}{2}}(t^{\ast})],
$$
by (4) of Lemma 3.3.1, we have
\begin{eqnarray*}
|I_2|
&\leq& \frac{1}{(n-2)\omega_n}\frac{1}{|x-y|^{n-2}}
|I_m^{(n-2)}(s^{\ast},t^{\ast})|(n-2)C_{m}^{\frac{n}{2}}(1)|t-t^{\ast}|\\
&\leq & \left\{\begin{array}{ll}
 A\frac{x_ny_n}{|x-y|^{n-2}}\frac{|x|^m}{|y|^{m+2}} , &   \mbox{when }   s\leq1  ,\\
  A\frac{x_ny_n}{|x-y|^{n-2}}\frac{|x|^{m+n-4}}{|y|^{m+n-2}},& \mbox{when}\
  s>1;
 \end{array}\right.
\hspace{44mm} (3.5.4)
\end{eqnarray*}
for the third term, we have
$$
I_3=C_{m+1}^{\frac{n-2}{2}}(t^{\ast})I_m^{(n-2)}(s^{\ast},t^{\ast})[E(x-y)-E(x-y^{\ast})],
$$
thus
\begin{eqnarray*}
|I_3|
&\leq & A\frac{1}{(n-2)\omega_n}\frac{2(n-2)x_ny_n}{|x-y|^n}|I_m^{(n-2)}(s,t^{\ast})| \\
&\leq & \left\{\begin{array}{ll}
 A\frac{x_ny_n}{|x-y|^n}\frac{|x|^{m+1}}{|y|^{m+1}} , &   \mbox{when }   s\leq1  ,\\
  A\frac{x_ny_n}{|x-y|^n}\frac{|x|^{m+n-3}}{|y|^{m+n-3}},& \mbox{when}\   s>1.
 \end{array}\right.
\hspace{49mm} (3.5.5)
\end{eqnarray*}
So we have by (3.5.2), (3.5.3), (3.5.4) and (3.5.5)
\begin{eqnarray*}
&  & |f-f^{\ast}|\\
&\leq& |I_1|+|I_2|+|I_3| \\
&\leq& \left\{\begin{array}{ll}
 A\frac{x_ny_n}{|x-y|^{n-2}}\frac{|x|^m}{|y|^{m+1}}
 \big(\frac{|x|}{|y|^2} +\frac{1}{|y|}+\frac{|x|}{|x-y|^2}\big), &   \mbox{when }   s\leq1  ,\\
   A\frac{x_ny_n}{|x-y|^{n-2}}\frac{|x|^{m+n-5}}{|y|^{m+n-3}}
 \big(1 +\frac{|x|}{|y|}+\frac{|x|^2}{|x-y|^2}\big),& \mbox{when}\   s>1.
 \end{array}\right.
\hspace{25mm} (3.5.6)
\end{eqnarray*}

Similarly,
\begin{eqnarray*}
&  & g-g^{\ast}\\
& =& E(x-y)C_m^{\frac{n-2}{2}}(t)I_{m+1}^{(n-2)}(s,t)
-E(x-y)C_m^{\frac{n-2}{2}}(t)I_{m+1}^{(n-2)}(s,t^{\ast}) \\
& & +E(x-y)C_m^{\frac{n-2}{2}}(t)I_{m+1}^{(n-2)}(s^{\ast},t^{\ast})
-E(x-y)C_m^{\frac{n-2}{2}}(t^{\ast})I_{m+1}^{(n-2)}(s^{\ast},t^{\ast}) \\
&
&+E(x-y)C_m^{\frac{n-2}{2}}(t^{\ast})I_{m+1}^{(n-2)}(s^{\ast},t^{\ast})
-E(x-y^{\ast})C_m^{\frac{n-2}{2}}(t^{\ast})I_{m+1}^{(n-2)}(s^{\ast},t^{\ast})\\
&=& J_1+J_2+J_3, \hspace{93mm} (3.5.7)
\end{eqnarray*}
and we have the similar estimates:
$$
|J_1|\leq\left\{\begin{array}{ll}
 A\frac{x_ny_n}{|x-y|^{n-2}}\frac{|x|^{m+2}}{|y|^{m+4}} , &   \mbox{when }   s\leq1  ,\\
  A\frac{x_ny_n}{|x-y|^{n-2}}\frac{|x|^{m+n-4}}{|y|^{m+n-2}},& \mbox{when}\ s>1;
 \end{array}\right.\eqno{(3.5.8)}
$$
$$
|J_2|\leq\left\{\begin{array}{ll}
 A\frac{x_ny_n}{|x-y|^{n-2}}\frac{|x|^{m+1}}{|y|^{m+3}} , &   \mbox{when }   s\leq1  ,\\
 A\frac{x_ny_n}{|x-y|^{n-2}}\frac{|x|^{m+n-3}}{|y|^{m+n-1}},& \mbox{when}\ s>1;
 \end{array}\right.\eqno{(3.5.9)}
$$
$$
|J_3|\leq\left\{\begin{array}{ll}
 A\frac{x_ny_n}{|x-y|^n}\frac{|x|^{m+2}}{|y|^{m+2}} , &   \mbox{when }   s\leq1  ,\\
 A\frac{x_ny_n}{|x-y|^n}\frac{|x|^{m+n-2}}{|y|^{m+n-2}},& \mbox{when}\   s>1.
 \end{array}\right. \eqno{(3.5.10)}
$$
So we have by (3.5.7), (3.5.8), (3.5.9) and (3.5.10)
\begin{eqnarray*}
&  & |g-g^{\ast}|\\
&\leq&  |J_1|+|J_2|+|J_3| \\
&\leq& \left\{\begin{array}{ll}
 A\frac{x_ny_n}{|x-y|^{n-2}}\frac{|x|^{m+1}}{|y|^{m+2}}
 \big(\frac{|x|}{|y|^2} +\frac{1}{|y|}+\frac{|x|}{|x-y|^2}\big), &   \mbox{when }   s\leq1  ,\\
   A\frac{x_ny_n}{|x-y|^{n-2}}\frac{|x|^{m+n-4}}{|y|^{m+n-2}}
 \big(1 +\frac{|x|}{|y|}+\frac{|x|^2}{|x-y|^2}\big),& \mbox{when}\   s>1.
 \end{array}\right.
\hspace{24mm} (3.5.11)
\end{eqnarray*}

Hence we finally obtain by (3.5.1), (3.5.6) and (3.5.11)
\begin{eqnarray*}
&  & |G_m(x,y)|\\
&\leq& (m+1)|f-f^{\ast}|+(n+m-2)|g-g^{\ast}| \\
&\leq& \left\{\begin{array}{ll}
 A\frac{x_ny_n}{|x-y|^{n-2}}\frac{|x|^m}{|y|^{m+1}}
  \big(\frac{|x|}{|y|^2}+\frac{1}{|y|}
  +\frac{|x|}{|x-y|^2}\big), &   \mbox{when }   s\leq1  ,\\
 A\frac{x_ny_n}{|x-y|^{n-2}}\frac{|x|^{m+n-4}}{|y|^{m+n-2}}
  \big(1 +\frac{|x|}{|y|}+\frac{|x|^2}{|x-y|^2}
  \big),& \mbox{when}\   s>1.
 \end{array}\right.
\end{eqnarray*}

\section{$p>1$(Modified Kernel)}

 \section*{ 1. Introduction and Main Theorems}

\vspace{0.3cm}

  In this section, we will
consider measurable functions $f$ in ${\bf R}^{n-1}$ satisfying
$$
\int_{{\bf R}^{n-1}}\frac{|f(y')|^p}{(1+|y'|)^{\gamma}} dy'<\infty,
\eqno{(3.6.1)}
$$
where $\gamma$ is defined as in Theorem 3.6.1.

  In order to describe the asymptotic behaviour of subharmonic functions
represented by the modified kernel in
  the upper half space (see \cite{ML}, \cite{MLJ1}
and \cite{MLJ2}), we
   establish the following theorems.

\vspace{0.2cm}
 \noindent
{\bf Theorem 3.6.1} {\it Let $1\leq p<\infty,\
\frac{1}{p}+\frac{1}{q}=1$ and
$$
(n-1)+mp<\gamma <(n-1)+(m+1)p  \qquad  {\rm in \  case}  \  p>1;
$$
$$
\quad \ \    m+n-1<\gamma \leq m+n  \qquad  \qquad \quad \qquad {\rm
in \ case}  \ p=1.
$$
If $f$ is a measurable function in ${\bf R}^{n-1}$ satisfying
(3.6.1) and  $v(x)$ is the harmonic function defined by
$$
v(x)= \int_{{\bf R}^{n-1}}P_m(x,y')f(y')dy', \quad x\in H,
\eqno{(3.6.2)}
$$
then there exists $x_j\in H,\ \rho_j>0,$ such that
$$
\sum
_{j=1}^{\infty}\frac{\rho_j^{pn-\alpha}}{|x_j|^{pn-\alpha}}<\infty
\eqno{(3.6.3)}
$$
holds and
$$
v(x)=
o(x_n^{1-\frac{\alpha}{p}}|x|^{\frac{\gamma}{p}+\frac{n-1}{q}-n+\frac{\alpha}{p}}),
\quad  {\rm as}  \ |x|\rightarrow\infty   \eqno{(3.6.4)}
$$
holds in $H-G$, where $ G=\bigcup_{j=1}^\infty B(x_j,\rho_j)$ and
$0< \alpha\leq np$. }

\vspace{0.2cm}
 \noindent
 {\bf Remark 3.6.1 } {\it If $\alpha=n$, $p=1$ and $\gamma=n+m$, then (3.6.3) is a finite sum,
 the set $G$ is the union of finite balls, so (3.6.4) holds in $H$.
 This is just the result of Siegel-Talvila, therefore,
 our result (3.6.4) is the generalization of Theorem C.
}

 \vspace{0.2cm}
 \noindent
 {\bf Remark 3.6.2 } {\it If $\gamma=(n-1)+mp$, $p>1$, then
$$
v(x)=
o(x_n^{1-\frac{\alpha}{p}}(\log|x|)^{\frac{1}{q}}|x|^{\frac{\gamma}{p}+\frac{n-1}{q}-n+\frac{\alpha}{p}}),
\quad  {\rm as}  \ |x|\rightarrow\infty
$$
holds in $H-G$. }

  Next, we will generalize Theorem 3.6.1 to subharmonic functions.

\vspace{0.2cm}
 \noindent
{\bf Theorem 3.6.2 } {\it Let $p$ and $\gamma$ be as in Theorem
3.6.1. If $f$ is a measurable function in ${\bf R}^{n-1}$ satisfying
(3.6.1) and $\mu$ is a positive Borel measure satisfying
$$
\int_H\frac{y_n^p}{(1+|y|)^\gamma} d\mu(y)<\infty
$$
and
$$
\int_H\frac{1}{(1+|y|)^{n-1}} d\mu(y)<\infty.
$$
Write the subharmonic function
$$
u(x)= v(x)+h(x), \quad x\in H,
$$
where $v(x)$ is the harmonic function defined by (3.6.2), $h(x)$ is
defined by
$$
h(x)= \int_H G_m(x,y)d\mu(y)
$$
and $G_m(x,y)$ is defined by (3.3.3). Then there exists $x_j\in H,\
\rho_j>0,$ such that (3.6.3) holds and
$$
u(x)=
o(x_n^{1-\frac{\alpha}{p}}|x|^{\frac{\gamma}{p}+\frac{n-1}{q}-n+\frac{\alpha}{p}}),
\quad  {\rm as}  \ |x|\rightarrow\infty
$$
holds in $H-G$, where $ G=\bigcup_{j=1}^\infty B(x_j,\rho_j)$ and
$0< \alpha<2p$. }

\vspace{0.2cm}
 \noindent
 {\bf Remark 3.6.3 } {\it If $\gamma=(n-1)+mp$, $p>1$, then
$$
u(x)=
o(x_n^{1-\frac{\alpha}{p}}(\log|x|)^{\frac{1}{q}}|x|^{\frac{\gamma}{p}+\frac{n-1}{q}-n+\frac{\alpha}{p}}),
\quad  {\rm as}  \ |x|\rightarrow\infty
$$
holds in $H-G$. }

 \section*{ 2. Main Lemmas}

 In order to obtain the results, we
need these lemmas below:

\vspace{0.2cm}
 \noindent
{\bf Lemma 3.6.1 }  {\it The modified Poisson
 kernel $P_m(x,y')$ has the
following estimates:\\
{\rm (1)}\ If $1<|y'|\leq \frac{|x|}{2}$, then $|P_m(x,y')|\leq
\frac{Ax_n|x|^{m-1}}{|y'|^{m+n-1}}$;\\
{\rm (2)}\ If $\frac{|x|}{2}<|y'| \leq 2|x|$, then $|P_m(x,y')|\leq
\frac{Ax_n}{|x-(y',0)|^n}$;\\
{\rm (3)}\ If $|y'|>2|x|$, then $|P_m(x,y')|\leq
\frac{Ax_n|x|^{m}}{|y'|^{m+n}}$;\\
{\rm (4)}\ If $|y'|\leq 1$, then $|P_m(x,y')|\leq
\frac{Ax_n}{|x|^{n}}$.\\
}

\vspace{0.2cm}
 \noindent
{\bf Lemma 3.6.2 } {\it The modified Green function $G_m(x,y)$ has
the
following estimates:\\
{\rm (1)}\ If $1<|y|\leq \frac{|x|}{2}$, then $|G_m(x,y)|\leq
\frac{Ax_ny_n|x|^{m-1}}{|y|^{m+n-1}}$;\\
{\rm (2)}\ If $\frac{|x|}{2}<|y| \leq 2|x|$, then $|G_m(x,y)|\leq
\frac{Ax_ny_n}{|x-y|^n}$;\\
{\rm (3)}\ If $|y|>2|x|$, then $|G_m(x,y)|\leq
\frac{Ax_ny_n|x|^{m}}{|y|^{m+n}}$;\\
{\rm (4)}\ If $|y|\leq 1$, then $|G_m(x,y)|\leq
\frac{2x_ny_n}{\omega_n|x-y|^n}\leq
\frac{Ax_ny_n}{|x|^{n}}$;\\
{\rm (5)}\ If $|x-y| \leq \frac{x_n}{2}$, then $|G_m(x,y)|\leq
\frac{A}{|x-y|^{n-2}}$.\\
}

\vspace{0.4cm}

\section*{3.   Proof of Theorems }

 \emph{Proof of Theorem 3.6.1}

 We prove only the case $p>1$; the proof of the case $p=1$ is similar.
Define the measure $dm(y')$ by
$$
dm(y')=\frac{|f(y')|^p}{(1+|y'|)^{\gamma}} dy'.
$$
  For any $\varepsilon >0$, there exists $R_\varepsilon >2$, such that
$$
\int_{|y'|\geq
R_\varepsilon}dm(y')\leq\frac{\varepsilon^p}{5^{pn-\alpha}}.
$$
For every Lebesgue measurable set $E \subset {\bf R}^{n-1}$, the
measure $m^{(\varepsilon)}$ defined by $m^{(\varepsilon)}(E)
=m(E\cap\{x'\in{\bf R}^{n-1}:|x'|\geq R_\varepsilon\}) $ satisfies
$m^{(\varepsilon)}({\bf
R}^{n-1})\leq\frac{\varepsilon^p}{5^{pn-\alpha}}$, write
\begin{eqnarray*}
&v_1(x)& =\int_{G_1} P_m(x,y')f(y')
dy',\\
&v_2(x)&=\int_{G_2} P_m(x,y')f(y')
dy', \\
&v_3(x)&=\int_{G_3} P_m(x,y')f(y')
dy', \\
&v_4(x)&=\int_{G_4} P_m(x,y')f(y')
dy', \\
\end{eqnarray*}
where
\begin{eqnarray*}
&G_1& =\{y'\in {\bf R}^{n-1}: 1<|y'|\leq \frac{|x|}{2}\},\\
&G_2& =\{y'\in {\bf R}^{n-1}: \frac{|x|}{2}<|y'| \leq 2|x|\}, \\
&G_3& =\{y'\in {\bf R}^{n-1}: |y'|>2|x|\}, \\
&G_4& =\{y'\in {\bf R}^{n-1}: |y'|\leq 1\}. \\
\end{eqnarray*}
Then
$$
v(x) =v_1(x)+v_2(x)+v_3(x)+v_4(x). \eqno{(3.6.5)}
$$

  First, if $\gamma >(n-1)+mp$, then $(\frac{\gamma}{p}-m-n+1)q
  +(n-1)>0$. For $R_\varepsilon >2$, we have
$$
v_1(x)=\int_{1<|y'|\leq R_\varepsilon} P_m(x,y')f(y') dy'
+\int_{R_\varepsilon<|y'|\leq \frac{|x|}{2} } P_m(x,y')f(y')
dy'=v_{11}(x)+v_{12}(x),
$$
if $|x|>2R_\varepsilon$, then we obtain by (1) of Lemma 3.6.1 and
H\"{o}lder's inequality
\begin{eqnarray*}
|v_{11}(x)|
&\leq& \int_{1<|y'|\leq R_\varepsilon}\frac{Ax_n|x|^{m-1}}{|y'|^{m+n-1}}|f(y')| dy' \\
&\leq& Ax_n|x|^{m-1} \bigg(\int_{1<|y'|\leq
R_\varepsilon}\frac{|f(y')|^p}{|y'|^\gamma}
dy'\bigg)^{1/p}\bigg(\int_{1<|y'|\leq
R_\varepsilon}|y'|^{(\frac{\gamma}{p}-m-n+1)q}
dy'\bigg)^{1/q}, \\
\end{eqnarray*}
since
$$
\int_{1<|y'|\leq R_\varepsilon}|y'|^{(\frac{\gamma}{p}-m-n+1)q}
dy'\leq AR_\varepsilon^{(\frac{\gamma}{p}-m-n+1)q+(n-1)},
$$
so that
$$
|v_{11}(x)|\leq
Ax_n|x|^{m-1}R_\varepsilon^{(\frac{\gamma}{p}-m-n+1)+\frac{n-1}{q}}.
\eqno{(3.6.6)}
$$

  Moreover, we have similarly

\begin{eqnarray*}
|v_{12}(x)|
&\leq& Ax_n|x|^{m-1} \bigg(\int_{R_\varepsilon<|y'|\leq
\frac{|x|}{2}}\frac{|f(y')|^p}{|y'|^\gamma}
dy'\bigg)^{1/p}\bigg(\int_{R_\varepsilon<|y'|\leq
\frac{|x|}{2}}|y'|^{(\frac{\gamma}{p}-m-n+1)q}
dy'\bigg)^{1/q} \\
&\leq& Ax_n|x|^{\frac{\gamma}{p}+\frac{n-1}{q}-n}
\bigg(\int_{R_\varepsilon<|y'|\leq
\frac{|x|}{2}}\frac{|f(y')|^p}{|y'|^\gamma}
dy'\bigg)^{1/p}, \\
\end{eqnarray*}
which implies by arbitrariness of $R_\varepsilon$ that
$$
|v_{12}(x)|\leq A\varepsilon
x_n|x|^{\frac{\gamma}{p}+\frac{n-1}{q}-n}. \eqno{(3.6.7)}
$$

  Let $ E_1(\lambda)=\{x\in{\bf R}^{n}:|x|\geq2,\exists \
t>0, {\rm s.t.}\  m^{(\varepsilon)}(B(x,t)\cap{\bf R}^{n-1}
)>\lambda^p(\frac{t}{|x|})^{pn-\alpha}\}$, therefore, if $ |x|\geq
2R_\varepsilon$ and $x \notin E_1(\lambda)
 $, then we have
$$
\forall t>0,\ m^{(\varepsilon)}(B(x,t)\cap{\bf R}^{n-1}
)\leq\lambda^p\bigg(\frac{t}{|x|}\bigg)^{pn-\alpha}.
$$
If $\gamma >(n-1)+mp$, then $(\frac{\gamma}{p}-m-n+1)q
  +(n-1)>0$, so that we obtain by (2) of Lemma 3.6.1 and
H\"{o}lder's inequality
\begin{eqnarray*}
|v_2(x)|
&\leq& \int_{G_2}\frac{Ax_n}{|x-(y',0)|^n}|f(y')| dy' \\
&\leq& Ax_n
\bigg(\int_{G_2}\frac{|f(y')|^p}{|x-(y',0)|^{pn}|y'|^\gamma}
dy'\bigg)^{1/p}\bigg(\int_{G_2}|y'|^{\frac{\gamma q}{p}}
dy'\bigg)^{1/q} \\
&\leq& Ax_n|x|^{\frac{\gamma}{p}+\frac{n-1}{q}}
\bigg(\int_{G_2}\frac{|f(y')|^p}{|x-(y',0)|^{pn}|y'|^\gamma}
dy'\bigg)^{1/p},\\
\end{eqnarray*}
since
\begin{eqnarray*}
\int_{G_2}\frac{|f(y')|^p}{|x-(y',0)|^{pn}|y'|^\gamma} dy'
&\leq& \int_{x_n}^{3|x|}
\frac{2^\gamma+1}{t^{pn}} dm_x^{(\varepsilon)}(t) \\
&\leq& \frac{\lambda^p}{
|x|^{pn}}(2^\gamma+1)\bigg(\frac{1}{3^\alpha}+
\frac{pn}{\alpha}\bigg)\frac{|x|^\alpha}{x_n^\alpha}, \\
\end{eqnarray*}
where  $m_x^{(\varepsilon)}(t)=\int_{|x-(y',0)| \leq t}
dm^{(\varepsilon)}(y')$.\\
Hence we have
$$
|v_2(x)|\leq A\lambda
x_n^{1-\frac{\alpha}{p}}|x|^{\frac{\gamma}{p}+\frac{n-1}{q}-n+\frac{\alpha}{p}}.\eqno{(3.6.8)}
$$

  If $\gamma <(n-1)+(m+1)p$, then $(\frac{\gamma}{p}-m-n)q+(n-1)<0$,
so that we obtain by (3) of Lemma 3.6.1 and H\"{o}lder's inequality
\begin{eqnarray*}
|v_3(x)|
&\leq& \int_{G_3}\frac{Ax_n|x|^{m}}{|y'|^{m+n}}|f(y')| dy' \\
&\leq& Ax_n|x|^{m} \bigg(\int_{G_3}\frac{|f(y')|^p}{|y'|^\gamma}
dy'\bigg)^{1/p}\bigg(\int_{G_3}|y'|^{(\frac{\gamma }{p}-m-n)q}
dy'\bigg)^{1/q} \\
&\leq& A\varepsilon x_n|x|^{\frac{\gamma}{p}+\frac{n-1}{q}-n}.
\hspace{75mm} (3.6.9)
\end{eqnarray*}

  Finally, by (4) of Lemma 3.6.1, we obtain
$$
|v_4(x)|\leq \frac{Ax_n}{|x|^{n}}
\int_{G_4}{|f(y')|}dy'.\eqno{(3.6.10)}
$$

  Thus, by collecting (3.6.5), (3.6.6), (3.6.7), (3.6.8), (3.6.9) and
(3.6.10), there exists a positive constant $A$ independent of
$\varepsilon$, such that if $ |x|\geq 2R_\varepsilon$ and $\  x
\notin E_1(\varepsilon)$, we have
$$
|v(x)|\leq A\varepsilon
x_n^{1-\frac{\alpha}{p}}|x|^{\frac{\gamma}{p}+\frac{n-1}{q}-n+\frac{\alpha}{p}}.
$$

 Let $\mu_\varepsilon$ be a measure in ${\bf R}^n$ defined by
$ \mu_\varepsilon(E)= m^{(\varepsilon)}(E\cap{\bf R}^{n-1})$ for
every measurable set $E$ in ${\bf R}^n$. Take
$\varepsilon=\varepsilon_p=\frac{1}{2^{p+2}}, p=1,2,3,\cdots$, then
there exists a sequence $ \{R_p\}$: $1=R_0<R_1<R_2<\cdots$ such that
$$
\mu_{\varepsilon_p}({\bf R}^n)=\int_{|y'|\geq
R_p}dm(y')<\frac{\varepsilon_p^p}{5^{pn-\alpha}}.
$$
Take $\lambda=3\cdot5^{pn-\alpha}\cdot2^p\mu_{\varepsilon_p}({\bf
R}^n)$ in Lemma 3.2.1, then there exists $x_{j,p}$ and $
\rho_{j,p}$, where $R_{p-1}\leq |x_{j,p}|<R_p$, such that
$$
\sum
_{j=1}^{\infty}\bigg(\frac{\rho_{j,p}}{|x_{j,p}|}\bigg)^{pn-\alpha}
\leq \frac{1}{2^{p}}.
$$
If $R_{p-1}\leq |x|<R_p$ and $x\notin G_p=\cup_{j=1}^\infty
B(x_{j,p},\rho_{j,p})$, we have
$$
|v(x)|\leq
A\varepsilon_px_n^{1-\frac{\alpha}{p}}|x|^{\frac{\gamma}{p}+\frac{n-1}{q}-n+\frac{\alpha}{p}}.
$$
Thereby
$$
\sum _{p=1}^{\infty}
\sum_{j=1}^{\infty}\bigg(\frac{\rho_{j,p}}{|x_{j,p}|}\bigg)^{pn-\alpha}
\leq \sum _{p=1}^{\infty}\frac{1}{2^{p}}=1<\infty.
$$
Set $ G=\cup_{p=1}^\infty G_p$, thus Theorem 3.6.1 holds.

 \emph{Proof of Theorem 3.6.2}

 We prove only the case $p>1$; the remaining case $p=1$
can be proved similarly. Define the measure $dn(y)$ by
$$
dn(y)=\frac{y_n^p}{(1+|y|)^{\gamma}} d\mu(y).
$$

  For any $\varepsilon >0$, there exists $R_\varepsilon >2$, such that
$$
\int_{|y|\geq
R_\varepsilon}dn(y)<\frac{\varepsilon^p}{5^{pn-\alpha}}.
$$
For every Lebesgue measurable set $E \subset {\bf R}^{n}$, the
measure $n^{(\varepsilon)}$ defined by $n^{(\varepsilon)}(E)
=n(E\cap\{y\in H:|y|\geq R_\varepsilon\}) $ satisfies
$n^{(\varepsilon)}(H)\leq\frac{\varepsilon^p}{5^{pn-\alpha}}$, write
\begin{eqnarray*}
&h_1(x)& = \int_{F_1} G_m(x,y)d\mu(y),\\
&h_2(x)&=\int_{F_2} G_m(x,y)d\mu(y),\\
&h_3(x)&=\int_{F_3} G_m(x,y)d\mu(y), \\
&h_4(x)&=\int_{F_4} G_m(x,y)d\mu(y), \\
\end{eqnarray*}
where
\begin{eqnarray*}
&F_1& =\{y\in H: 1<|y|\leq \frac{|x|}{2}\},\\
&F_2& =\{y\in H: \frac{|x|}{2}<|y| \leq 2|x|\}, \\
&F_3& =\{y\in H: |y|>2|x|\}, \\
&F_4& =\{y\in H: |y|\leq 1\}. \\
\end{eqnarray*}
Then
$$
h(x) =h_1(x)+h_2(x)+h_3(x)+h_4(x). \eqno{(3.6.11)}
$$

  First, if $\gamma >(n-1)+mp$, then $(\frac{\gamma}{p}-m-n+1)q
  +(n-1)>0$. For $R_\varepsilon >2$, we have
$$
h_1(x)=\int_{1<|y|\leq R_\varepsilon} G_m(x,y) d\mu(y)
+\int_{R_\varepsilon<|y|\leq \frac{|x|}{2} } G_m(x,y)
d\mu(y)=h_{11}(x)+h_{12}(x),
$$
if $|x|>2R_\varepsilon$, then we obtain by (1) of Lemma 3.6.2 and
H\"{o}lder's inequality
\begin{eqnarray*}
|h_{11}(x)|
&\leq& \int_{1<|y|\leq R_\varepsilon}\frac{Ax_ny_n|x|^{m-1}}{|y|^{m+n-1}} d\mu(y) \\
&\leq& Ax_n|x|^{m-1} \bigg(\int_{1<|y|\leq
R_\varepsilon}\frac{y_n^p}{|y|^\gamma}
d\mu(y)\bigg)^{1/p}\bigg(\int_{1<|y|\leq
R_\varepsilon}|y|^{(\frac{\gamma}{p}-m-n+1)q}
d\mu(y)\bigg)^{1/q}, \\
\end{eqnarray*}
since
$$
\int_{1<|y|\leq R_\varepsilon}|y|^{(\frac{\gamma}{p}-m-n+1)q}
d\mu(y)\leq AR_\varepsilon^{(\frac{\gamma}{p}-m-n+1)q+(n-1)},
$$
so that
$$
|h_{11}(x)|\leq
Ax_n|x|^{m-1}R_\varepsilon^{(\frac{\gamma}{p}-m-n+1)+\frac{n-1}{q}}.
\eqno{(3.6.12)}
$$

  Moreover, we have similarly

\begin{eqnarray*}
|h_{12}(x)|
&\leq& Ax_n|x|^{m-1} \bigg(\int_{R_\varepsilon<|y|\leq
\frac{|x|}{2}}\frac{y_n^p}{|y|^\gamma}
d\mu(y)\bigg)^{1/p}\bigg(\int_{R_\varepsilon<|y|\leq
\frac{|x|}{2}}|y|^{(\frac{\gamma}{p}-m-n+1)q}
d\mu(y)\bigg)^{1/q} \\
&\leq& Ax_n|x|^{\frac{\gamma}{p}+\frac{n-1}{q}-n}
\bigg(\int_{R_\varepsilon<|y|\leq
\frac{|x|}{2}}\frac{y_n^p}{|y|^\gamma}
d\mu(y)\bigg)^{1/p}, \\
\end{eqnarray*}
which implies by arbitrariness of $R_\varepsilon$ that
$$
|h_{12}(x)|\leq A\varepsilon
x_n|x|^{\frac{\gamma}{p}+\frac{n-1}{q}-n}. \eqno{(3.6.13)}
$$

  Let $ E_2(\lambda)=\{x\in{\bf R}^{n}:|x|\geq2,\exists \
t>0,{\rm s.t.}\ n^{(\varepsilon)}(B(x,t)\cap H
)>\lambda^p(\frac{t}{|x|})^{pn-\alpha}\}, $ therefore, if $ |x|\geq
2R_\varepsilon$ and $x\notin E_2(\lambda)
 $, then we have
$$
\forall t>0, \ n^{(\varepsilon)}(B(x,t)\cap H
)\leq\lambda^p\bigg(\frac{t}{|x|}\bigg)^{pn-\alpha}.
$$
If $\gamma >(n-1)+mp$, then $(\frac{\gamma}{p}-m-n+1)q
  +(n-1)>0$, so that we obtain by H\"{o}lder's inequality
\begin{eqnarray*}
|h_2(x)|
&\leq& \bigg(\int_{F_2}\frac{|G_m(x,y)|^p}{|y|^\gamma}
d\mu(y)\bigg)^{1/p}\bigg(\int_{F_2}|y|^{\frac{\gamma q}{p}}
d\mu(y)\bigg)^{1/q} \\
&\leq& \bigg((2^\gamma +1)\int_{F_2}\frac{|G_m(x,y)|^p}{y_n^p}
dn(y)\bigg)^{1/p}\bigg(\int_{F_2}|y|^{\frac{\gamma q}{p}}
d\mu(y)\bigg)^{1/q}\\
&\leq& A|x|^{\frac{\gamma}{p}+\frac{n-1}{q}}\bigg(\int_{F_2}\frac{|G_m(x,y)|^p}{y_n^p}dn(y)\bigg)^{1/p},\\
\end{eqnarray*}
since
\begin{eqnarray*}
\int_{F_2}\frac{|G_m(x,y)|^p}{y_n^p} dn(y)
&\leq&  \int_{|x-y|\leq 3|x|}\frac{|G_m(x,y)|^p}{y_n^p} dn^{(\varepsilon)}(y) \\
&=&  \int_{|x-y|\leq \frac{x_n}{2}}\frac{|G_m(x,y)|^p}{y_n^p}
dn^{(\varepsilon)}(y)\\
& &+ \int_{\frac{x_n}{2}<|x-y|\leq 3|x|}\frac{|G_m(x,y)|^p}{y_n^p} dn^{(\varepsilon)}(y)\\
&=& h_{21}(x)+h_{22}(x),
\end{eqnarray*}
so that we have by (5) of Lemma 3.6.2
\begin{eqnarray*}
h_{21}(x)
&\leq& \int_{|x-y|\leq
\frac{x_n}{2}}\bigg(\frac{A}{x_n|x-y|^{n-2}}\bigg)^p
dn^{(\varepsilon)}(y) \\
&=& \bigg(\frac{A}{x_n}\bigg)^p\int_0^{\frac{x_n}{2}}
\frac{1}{t^{p(n-2)}} dn_x^{(\varepsilon)}(t) \\
&\leq& A\frac{np-\alpha}{(2p-\alpha)2^{2p-\alpha}}
\lambda^p\frac{x_n^{p-\alpha}}{|x|^{np-\alpha}}.\\
\end{eqnarray*}

  Moreover, we have by (2) of Lemma 3.6.2
\begin{eqnarray*}
h_{22}(x)
&\leq& \int_{\frac{x_n}{2}<|x-y|\leq 3|x|
}\bigg(\frac{Ax_n}{|x-y|^{n}}\bigg)^p
dn^{(\varepsilon)}(y) \\
&=& (Ax_n)^p\int_{\frac{x_n}{2}}^{3|x|}
\frac{1}{t^{pn}} dn_x^{(\varepsilon)}(t) \\
&\leq& A\bigg(\frac{1}{3^\alpha}+\frac{np2^\alpha}{\alpha}\bigg)
\lambda^p\frac{x_n^{p-\alpha}}{|x|^{np-\alpha}},\\
\end{eqnarray*}
where  $n_x^{(\varepsilon)}(t)=\int_{|x-y| \leq t}
dn^{(\varepsilon)}(y)$.\\
Hence we have
$$
|h_2(x)|\leq A\lambda
x_n^{1-\frac{\alpha}{p}}|x|^{\frac{\gamma}{p}+\frac{n-1}{q}-n+\frac{\alpha}{p}}.
\eqno{(3.6.14)}
$$

 If $\gamma <(n-1)+(m+1)p$, then $(\frac{\gamma}{p}-m-n)q+(n-1)<0$,
so that we obtain by (3) of Lemma 3.6.2 and H\"{o}lder's inequality
\begin{eqnarray*}
|h_3(x)|
&\leq& \int_{F_3}\frac{Ax_ny_n|x|^{m}}{|y|^{m+n}} d\mu(y) \\
&\leq& Ax_n|x|^{m} \bigg(\int_{F_3}\frac{y_n^p}{|y|^\gamma}
d\mu(y)\bigg)^{1/p}\bigg(\int_{F_3}|y|^{(\frac{\gamma }{p}-m-n)q}
d\mu(y)\bigg)^{1/q} \\
&\leq& A\varepsilon x_n|x|^{\frac{\gamma}{p}+\frac{n-1}{q}-n}.
\hspace{73mm} (3.6.15)
\end{eqnarray*}

  Finally, by (4) of Lemma 3.6.2, we obtain
$$
|h_4(x)|\leq \frac{Ax_n}{|x|^{n}}
\int_{F_4}{y_n}d\mu(y).\eqno{(3.6.16)}
$$

  Thus, by collecting (3.6.11), (3.6.12), (3.6.13),
(3.6.14), (3.6.15) and (3.6.16), there exists a positive constant
$A$ independent of $\varepsilon$, such that if $ |x|\geq
2R_\varepsilon$ and $\  x \notin E_2(\varepsilon)$, we have
$$
 |h(x)|\leq A\varepsilon x_n^{1-\frac{\alpha}{p}}|x|^{\frac{\gamma}{p}+\frac{n-1}{q}-n+\frac{\alpha}{p}}.
$$

  Similarly, if $x\notin G$, we have
$$
h(x)=
o(x_n^{1-\frac{\alpha}{p}}|x|^{\frac{\gamma}{p}+\frac{n-1}{q}-n+\frac{\alpha}{p}}),\quad
{\rm as} \ |x|\rightarrow\infty. \eqno{(3.6.17)}
$$

  By (3.6.4) and (3.6.17), we obtain that
$$
u(x)=v(x)+h(x)=
o(x_n^{1-\frac{\alpha}{p}}|x|^{\frac{\gamma}{p}+\frac{n-1}{q}-n+\frac{\alpha}{p}}),\quad
{\rm as} \ |x|\rightarrow\infty
$$
holds in $H-G$, thus we complete the proof of Theorem 3.6.2.

\chapter{a Generalization of Harmonic Majorants}

\section{a Generalization of Harmonic Majorants in the Upper Half Plane}

\section*{ 1. Introduction and Main Theorem}

\vspace{0.2cm}

  The Poisson kernel for the
half plane ${\bf C}_+=\{z=x+iy \in {\bf C}:\; y>0\}$ is the
 function
$$
 P(z,t)=
\frac{y}{\pi|z-t|^2},
$$
where $z \in {\bf C}_+$ and $t\in {\bf R}$.

If $p\geq 0$ is an integer, we define a modified
 Cauchy kernel of order $ p$ for $z\in {\bf C}_+-\{t\}$ by
$$
 C_p(z,t)=\left\{\begin{array}{ll}
 \frac{1}{\pi}\frac{1}{t-z} , &   \mbox{when }   |t|\leq 1  ,\\
 \frac{1}{\pi}\frac{1}{t-z} - \frac{1}{\pi}\sum_{k=0}^{p}\frac{z^k}{t^{k+1}},&
\mbox{when}\   |t|> 1,
 \end{array}\right.
$$
then we define a modified
 Poisson kernel of order $p$ for the upper half plane  by
$$
P_p(z,t)=\Im C_p(z,t).
$$

 Flett and Kuran \cite{RR} proved the following theorem:

 \vspace{0.2cm}
 \noindent
{\bf Theorem D } {\it  Let $G(z)$ be nonnegative and subharmonic in
${\bf C}_+$. Then $G(z)$ has a harmonic majorant in ${\bf C}_+$ if
and only if
$$
\sup_{y>0}\int_{-\infty}^{\infty}\frac{G(x+iy)}{x^2+(y+1)^2}dx<\infty.
$$
}

\vspace{0.2cm}
 \noindent
 {\bf Remark 4.1.1 } {\it If $G(z)$ has a harmonic majorant in ${\bf C}_+$,
 then there exists a harmonic function
$$
H(z)=cy+\frac{y}{\pi}\int_{-\infty}^{\infty}\frac{d\mu(t)}{(t-x)^2+y^2},\quad
y>0,
$$
where $c\geq 0$ and $\mu$ is a nonnegative Borel measure on
$(-\infty,\infty)$ such that
$$
\int_{-\infty}^{\infty}\frac{d\mu(t)}{1+t^2}<\infty,
$$
and
$$
G(z)\leq H(z).
$$
}

  In this section, We will generalize Theorem D partly to the modified kernel.

 \vspace{0.2cm}
 \noindent
{\bf Theorem 4.1.1 } {\it  Let
$$
H(z)=\Im[Q_p(z)+\frac{1}{\pi}\int_{-\infty}^{\infty}C_p(z,t)d\mu(t)],\quad
z=x+iy,\  y>0,
$$
where
$$
Q_p(z)=\sum _{k=0}^{p}a_k z^k, \quad a_k\in {\bf R},\
k=0,1,2,\cdots,p
$$
and $\mu$ is a nonnegative Borel measure on $(-\infty,\infty)$ such
that
$$
\int_{-\infty}^{\infty}\frac{1}{1+|t|^{p+1}}d\mu(t)<\infty.
$$
If  $G(z)$ is subharmonic in ${\bf C}_+$ and
$$
G(z)\leq H(z),
$$
then
$$
\sup_{y>0}\int_{-\infty}^{\infty}
\frac{G(x+iy)}{[x^2+(y+1)^2]^{\frac{p+1}{2}}}dx<\infty.
$$
}

\vspace{0.2cm}
 \noindent
 {\bf Remark 4.1.2 } {\it If $p=1$,  this is just the result of
 Flett and Kuran, therefore, our result is partly the generalization of
Theorem D. }

\vspace{0.4cm}

\section*{2.   Main lemmas }

\vspace{0.3cm}

In order to obtain the result, we need these lemmas below:

\vspace{0.2cm}
 \noindent
{\bf Lemma 4.1.1 } {\it For any $|t|>1$,  the following equality
$$
\Im C_p(z,t)=\Im \frac{tz^{p+1}-|z|^2z^p} {|t-z|^2 t^{p+1}}
\eqno{(4.1.1)}
$$
holds . }

\vspace{0.2cm} Proof: For $|t|>1$, since
$$
 C_p(z,t)=\frac{1}{t-z} - \sum_{k=0}^{p}\frac{z^k}{t^{k+1}}
 =\frac{z^{p+1}}{(t-z)t^{p+1}},
$$
then
$$
\Im C_p(z,t)=\Im \frac{z^{p+1}}{(t-z)t^{p+1}} =\Im
\bigg[\frac{z^{p+1}(t-\overline{z})}{|t-z|^2 t^{p+1}}\bigg] =\Im
\frac{tz^{p+1}-|z|^2z^p} {|t-z|^2 t^{p+1}}.
$$
This proves the equality (4.1.1).

\vspace{0.2cm}
 \noindent
{\bf Lemma 4.1.2 } {\it There exists $A>0$, such that the inequality
$$
\Im(tz^{p+1}-|z|^2z^p)\leq Ay(t^2+y^2)(x^2+y^2)^{\frac{p-1}{2}}
$$
holds in the following conditions:\\
{\rm (1)}\ $ p=2m-1,m=1,2,\cdots;$\\
{\rm (2)}\ $ p=2m,m=1,2,\cdots,x\geq0;$\\
{\rm (3)}\ $ p=2m,m=1,2,\cdots,x<0,|t|\geq |x|.$\\
 }

\vspace{0.4cm}

\section*{3.   Proof of Theorem }

\vspace{0.3cm}

  First from $G(z)\leq H(z)$, we obtain
\begin{eqnarray*}
& & \int_{-\infty}^{\infty}
\frac{G(x+iy)}{[x^2+(y+1)^2]^{\frac{p+1}{2}}}dx \\
&\leq & \int_{-\infty}^{\infty} \frac{\Im
Q_p(z)}{[x^2+(y+1)^2]^{\frac{p+1}{2}}}dx\\
& &+
\frac{1}{\pi}\int_{-\infty}^{\infty}\frac{1}{[x^2+(y+1)^2]^{\frac{p+1}{2}}}dx
\int_{-\infty}^{\infty}\Im C_p(z,t)d\mu(t) \\
&= & I_1+I_2. \hspace{103mm} (4.1.2)
\end{eqnarray*}
For the first term, we have
\begin{eqnarray*}
I_1
&=& \int_{-\infty}^{\infty} \frac{\sum _{k=0}^{p}a_k \Im
(x+iy)^k}{[x^2+(y+1)^2]^{\frac{p+1}{2}}}dx \\
&=& y \int_{-\infty}^{\infty} \frac{\sum _{k=0}^{p}a_k
[C_k^1x^{k-1}-C_k^3x^{k-3}y^2+\cdots]}{[x^2+(y+1)^2]^{\frac{p+1}{2}}}dx \\
&\leq & y \int_{-\infty}^{\infty} \frac{\sum _{k=0}^{p}|a_k|\sum
_{i=1}^{[\frac{k+1}{2}]}
C_k^{2i-1}|x|^{(k-1)-(2i-2)}y^{2i-2}}{[x^2+(y+1)^2]^{\frac{p+1}{2}}}dx \\
&\leq& A y \int_{-\infty}^{\infty}
\frac{1}{[x^2+(y+1)^2]^{\frac{p-k}{2}+1}}dx \\
&\leq& A y \int_{-\infty}^{\infty} \frac{1}{x^2+(y+1)^2}dx \\
&\leq& A\pi \frac{y}{y+1}\leq A\pi; \hspace{86mm} (4.1.3)
\end{eqnarray*}
for the second term, we will discuss in the following conditions:

\vspace{0.2cm}
\noindent
{\rm (1)}\ $ p=2m-1,m=1,2,\cdots; $

\begin{eqnarray*}
I_2
& =& \frac{1}{\pi}\int_{-\infty}^{\infty}\int_{-\infty}^{\infty}\Im
C_p(z,t)d\mu(t)\frac{1}{[x^2+(y+1)^2]^{\frac{p+1}{2}}}dx \\
&= & \frac{1}{\pi} \int_{-\infty}^{\infty} \int_{|t|\leq
1}^{\infty}\frac{ y}{|t-z|^2}d\mu(t)\cdot\frac{1}{[x^2+(y+1)^2]^{\frac{p+1}{2}}}dx\\
& & + \frac{1}{\pi} \int_{-\infty}^{\infty} \int_{|t|> 1}\frac{
\Im(tz^{p+1}-|z|^2z^p)}{|t-z|^2t^{p+1}}d\mu(t) \cdot \frac{1}{[x^2+(y+1)^2]^{\frac{p+1}{2}}}dx \\
&= & J_{11}+J_{12}.\hspace{95mm} (4.1.4)
\end{eqnarray*}
Note that
$$
\frac{y}{\pi} \int_{-\infty}^{\infty}\frac{
1}{(t-x)^2+y^2}\frac{y+1}{x^2+(y+1)^2}dx
=\frac{2y+1}{t^2+(2y+1)^2},\eqno{(4.1.5)}
$$
we have
\begin{eqnarray*}
J_{11}
& =& \frac{1}{\pi}\int_{|t|\leq 1} \int_{-\infty}^{\infty}\frac{
y}{(t-x)^2+y^2}\frac{1}{[x^2+(y+1)^2]^{\frac{p+1}{2}}}dxd\mu(t) \\
&\leq &  \int_{|t|\leq 1}\frac{1}{y+1}[\frac{y}{\pi}
\int_{-\infty}^{\infty}\frac{
1}{(t-x)^2+y^2}\frac{y+1}{x^2+(y+1)^2}dx]d\mu(t) \\
&= &  \int_{|t|\leq 1}\frac{1}{y+1}\frac{2y+1}{t^2+(2y+1)^2}d\mu(t) \\
&\leq &  2\int_{|t|\leq 1}\frac{1}{t^2+1}d\mu(t) \\
&\leq & 4\int_{-\infty}^{\infty}\frac{1}{1+|t|^{p+1}}d\mu(t)<\infty.
\hspace{65mm}(4.1.6)
\end{eqnarray*}
Moreover,
\begin{eqnarray*}
J_{12}
& \leq&  \frac{1}{\pi} \int_{-\infty}^{\infty}\int_{|t|> 1}\frac{
My(t^2+y^2)(x^2+y^2)^{\frac{p-1}{2}}}{[(t-x)^2+y^2]|t|^{p+1}}\frac{1}{[x^2+(y+1)^2]^{\frac{p+1}{2}}}d\mu(t)dx \\
&= &  \frac{My}{\pi} \int_{-\infty}^{\infty}\int_{|t|> 1}\frac{
(t^2+y^2)}{[(t-x)^2+y^2][x^2+(y+1)^2]|t|^{p+1}}d\mu(t)dx \\
&= &  M \int_{|t|> 1}\frac{1}{y+1}\frac{y}{\pi}
\int_{-\infty}^{\infty}\frac{
1}{(t-x)^2+y^2}\frac{y+1}{x^2+(y+1)^2}dx\frac{t^2+y^2}{|t|^{p+1}}d\mu(t),
\end{eqnarray*}
again by (4.1.5), we have
\begin{eqnarray*}
J_{12}
& \leq & M \int_{|t|>
1}\frac{1}{y+1}\frac{2y+1}{t^2+(2y+1)^2}\frac{t^2+y^2}{|t|^{p+1}}d\mu(t) \\
&\leq & 4M \int_{|t|> 1}\frac{1}{1+|t|^{p+1}}d\mu(t)
<\infty;\hspace{57mm}(4.1.7)
\end{eqnarray*}

\vspace{0.2cm}
 \noindent
 {\rm (2)}\
$ p=2m,m=1,2,\cdots. $

\begin{eqnarray*}
I_2
& = & \frac{1}{\pi}\int_{-\infty}^{\infty}\int_{-\infty}^{\infty}\Im
C_p(z,t)d\mu(t)\frac{1}{[x^2+(y+1)^2]^{\frac{p+1}{2}}}dx \\
& = & \int\int_{\{(x,t):\; x\geq 0\}}+\int\int_{\{(x,t):\;x<
0,|t|\geq
-x\}} +\int\int_{\{(x,t):\;x< 0,|t|< -x\}} \\
&= & J_{21}+J_{22}+J_{23}. \hspace{86mm}(4.1.8)
\end{eqnarray*}
Similarly, we can obtain in the same method as {\rm (1)} that
$J_{21}<\infty $ and $ J_{22}<\infty.$

 Write
\begin{eqnarray*}
J_{23}
& = & \frac{1}{\pi}\int\int_{\{(x,t):\;x< 0,|t|< -x\}}\Im
C_p(z,t)d\mu(t)\frac{1}{[x^2+(y+1)^2]^{\frac{p+1}{2}}}dx \\
& = & \int\int_{\{(x,t):\;x< 0,|t|< -x\}\bigcap\{(x,t):\;|t|\leq
1\}}
+\int\int_{\{(x,t):\;x< 0,|t|< -x\}\bigcap\{(x,t):\;|t|>1\}} \\
&= & K_1+K_2, \hspace{95mm} (4.1.9)
\end{eqnarray*}
again, we can obtain in the same method as {\rm (1)} that $
K_1<\infty$.

  In the following, we will show that
$$
K_2 <\infty. \eqno{(4.1.10)}
$$
Write $D=\{(x,t):\;x< 0,|t|< -x\}\bigcap\{(x,t):\;|t|>1\}$, then
\begin{eqnarray*}
K_2
& = & \frac{1}{\pi}\int\int_{D}\Im
C_p(z,t)d\mu(t)\frac{1}{[x^2+(y+1)^2]^{\frac{p+1}{2}}}dx\\
& = & \frac{1}{\pi}\int\int_{D} \bigg \{\frac{t\sum
_{i=1}^{\frac{p}{2}+1}(-1)^{i+1}C_{p+1}^{2i-1}
x^{(p+1)-(2i-1)}y^{2i-1}} {[(t-x)^2+y^2]|t|^{p+1}} \\
&   & -\frac{(x^2+y^2)\sum
_{i=1}^{\frac{p}{2}}(-1)^{i+1}C_{p}^{2i-1} x^{p-(2i-1)}y^{2i-1}}
{[(t-x)^2+y^2]|t|^{p+1}} \bigg \}
 \cdot
\frac{1}{[x^2+(y+1)^2]^{\frac{p+1}{2}}}d\mu(t)dx \\
&\leq&  \frac{1}{\pi}\int\int_{D} \frac{\sum
_{i=1}^{\frac{p}{2}+1}C_{p+1}^{2i-1} |x|^{p-2i+3}y^{2i-2}+
(x^2+y^2)\sum _{i=1}^{\frac{p}{2}}C_{p}^{2i-1} |x|^{p-2i+1}y^{2i-2}}
{[x^2+(y+1)^2]^{\frac{p+1}{2}}} \\
&  & \times \frac{y}{[(t-x)^2+y^2]|t|^{p+1}}d\mu(t)dx.
\end{eqnarray*}
Note that
$$
\frac{\sum _{i=1}^{\frac{p}{2}+1}C_{p+1}^{2i-1}
|x|^{p-2i+3}y^{2i-2}+ (x^2+y^2)\sum _{i=1}^{\frac{p}{2}}C_{p}^{2i-1}
|x|^{p-2i+1}y^{2i-2}} {[x^2+(y+1)^2]^{\frac{p+1}{2}}} \leq A,
$$
then we have
\begin{eqnarray*}
K_2
&\leq& \frac{Ay}{\pi}\int\int_{D}\frac{1}{[(t-x)^2+y^2]|t|^{p+1}}d\mu(t)dx \\
&\leq& A\int_{-\infty}^{\infty}\frac{y}{\pi}\int_{-\infty}^{\infty}
\frac{1}{(t-x)^2+y^2}dx \frac{2}{1+|t|^{p+1}}d\mu(t).
\end{eqnarray*}

Note that
$$
\frac{y}{\pi}\int_{-\infty}^{\infty}\frac{1}{(t-x)^2+y^2}dx=1,
$$
then we have
$$
K_2\leq
2A\int_{-\infty}^{\infty}\frac{1}{1+|t|^{p+1}}d\mu(t)<\infty.
$$
So the result follows by collecting (4.1.2), (4.1.3), (4.1.4),
(4.1.6), (4.1.7), (4.1.8), (4.1.9) and (4.1.10).

\section{a Generalization of Harmonic Majorants in the Upper Half Space}

 \section*{ 1. Introduction and Main Theorem}

\vspace{0.3cm}

\vspace{0.2cm}
  The Poisson kernel for the
half space $H $ is the
 function
$$
 P(x,y')=
\frac{2x_n}{\omega_n|x-y'|^n},
$$
where $x \in H$, $y\in \partial{H}$ and $ \omega_n = \frac{2\pi
^{\frac{n}{2}}}{\Gamma (\frac{n}{2})}$ is the area of the unit
sphere in ${\bf R}^n$.

  In this section, We will generalize Theorem D to the
upper half space.

 \vspace{0.2cm}
 \noindent
{\bf Theorem 4.2.1 } {\it  Let
$$
H(x)=cx_n+\int_{{\bf R}^{n-1}}P(x,y')d\mu(y'),
$$
where $c\geq 0$ and $\mu$ is a nonnegative Borel measure on ${\bf
R}^{n-1}$ such that
$$
\int_{{\bf R}^{n-1}}\frac{1}{(1+|y'|^2)^{n/2}}d\mu(y')<\infty.
$$
If $G(x)$ is nonnegative and subharmonic in $H$. Then
$$
G(x)\leq H(x)
$$
if and only if
$$
\sup_{x_n>0}\int_{{\bf
R}^{n-1}}\frac{G(x)}{[|x'|^2+(x_n+1)^2]^{\frac{n}{2}}}dx'<\infty.
\eqno{(4.2.1)}
$$
}

\vspace{0.2cm}
 \noindent
 {\bf Remark 4.2.1 } {\it If $n=2$,  this is just the result of
 Flett and Kuran, therefore, our result is the generalization of
Theorem D. }

\vspace{0.4cm}

\section*{2.   Main Lemmas }

\vspace{0.3cm}

   In order to obtain the result, we need these lemmas below:

\vspace{0.2cm}
 \noindent
{\bf Lemma 4.2.1 } {\it Let $H(x)$ be nonnegative and harmonic in
$H$ and have a continuous extension to $\overline{H}$. Then
$$
H(x)=cx_n+\frac{2x_n}{\omega_n}\int_{{\bf
R}^{n-1}}\frac{H(y')}{[|y'-x'|^2+x_n^2]^{\frac{n}{2}}}dy',
$$
where $c$ is given by
$$
c=\lim_{x_n\rightarrow\infty} \frac{H(0,x_n)}{x_n}.
\eqno{(4.2.2)}
$$
}

\vspace{0.2cm}
 \noindent
{\bf Lemma 4.2.2 } {\it Let
\begin{displaymath}
\mathbf{x^T} = \left( \begin{array}{ccc}
x_1  \\
x_2  \\
\vdots  \\
x_n
\end{array} \right),
\end{displaymath}
then we have
$$
D_n=|xx^TE-2x^Tx|=-|x|^{2n}.
$$
}

Proof:

\qquad \qquad  \quad $D_n =||x|^2E-2x^Tx|$
\begin{displaymath}
  =\left| \begin{array}{cccc}
|x|^2-2x_1^2 & -2x_1x_2 & \ldots & -2x_1x_n \\
\vdots & \vdots & \vdots & \vdots\\
-2x_nx_1 & -2x_nx_2 & \ldots &|x|^2-2x_n^2
\end{array} \right|
\end{displaymath}
\qquad \qquad \qquad \qquad $=|x|^2D_{n-1}-2x_n^2|x|^{2(n-1)},$\\
therefore,
\begin{eqnarray*}
D_n
& = & |x|^2D_{n-1}-2x_n^2|x|^{2(n-1)} \\
& = & |x|^2[|x|^2D_{n-2}-2x_{n-1}^2|x|^{2(n-2)}]-2x_n^2|x|^{2(n-1)}\\
& = & |x|^{2\times 2}D_{n-2}-2(x_{n-1}^2+x_n^2)|x|^{2(n-1)} \\
& = & \cdots\cdots \\
& = & |x|^{2\times (n-1)}D_1-2(x_2^2+x_3^2+\cdots +x_n^2)|x|^{2(n-1)}\\
& = & |x|^{2\times (n-1)}(|x|^2-2x_1^2)-2(x_2^2+x_3^2+\cdots
+x_n^2)|x|^{2(n-1)} \\
& = & |x|^{2n}-2(x_1^2+x_2^2+\cdots +x_n^2)|x|^{2(n-1)}\\
& = & |x|^{2n}-2|x|^{2n}=-|x|^{2n}.
\end{eqnarray*}

\vspace{0.2cm}
 \noindent
{\bf Lemma 4.2.3 } {\it Let
$x=(x_1,x_2,\cdots,x_{n-1},x_{n})=(x',x_n)$, where $x' \in {\bf
R}^{n-1}$ and $x_{n} \in {\bf R}$. $S=(0,0,\cdots,0,-1)=(0,-1)$,
where $0 \in {\bf R}^{n-1}$. Suppose
$$
u=\Phi (x)=2(x-S)^*+S,
$$
where $x^*=\frac{x}{|x|^2}$,\ then we have
$$
J_\Phi (x)=-\frac{2^n}{[|x'|^2+(x_n+1)^2]^n},
$$
where $J_\Phi (x)$ is given by
$$
J_\Phi (x)=\frac{\partial
(u_1,u_2,\cdots,u_n)}{\partial(x_1,x_2,\cdots,x_n)}.
$$
}

Proof: Since
\begin{eqnarray*}
u
& = & \Phi (x)=2(x-S)^*+S \\
& = & \frac{2(x-S)}{|x-S|^2}+S\\
& = & \frac{2(x',x_n+1)}{|x'|^2+(x_n+1)^2}+(0,-1) \\
& = & \frac{(2x',2(x_n+1))+(0,-|x'|^2-(x_n+1)^2)}{|x'|^2+(x_n+1)^2} \\
& = & \frac{(2x',1-|x'|^2-x_n^2)}{|x'|^2+(x_n+1)^2}, \hspace{75mm}
(4.2.3)
\end{eqnarray*}
then we have
\begin{eqnarray*}
|\Phi (x)|^2
& = & \bigg|\frac{(2x',1-|x'|^2-x_n^2)}{|x'|^2+(x_n+1)^2}\bigg|^2 \\
& = & \frac{4x'^2+(1-|x'|^2-x_n^2)^2}{[|x'|^2+(x_n+1)^2]^2}\\
& = & \frac{x'^4+2x'^2(1+x_n^2)+(1+x_n)^2(1-x_n)^2}{[|x'|^2+(x_n+1)^2]^2} \\
& = & \frac{|x'|^2+(1-x_n)^2}{|x'|^2+(1+x_n)^2}.
\end{eqnarray*}
Let
$$
u=(u_1,u_2,\cdots,u_n),
$$
by(4.2.3), we obtain
$$
u_i=\left\{\begin{array}{ll}
 \frac{2x_i}{|x'|^2+(x_n+1)^2} , &   \mbox{when }   i=1,2,\cdots,n-1 ,\\
 \frac{2(x_n+1)}{|x'|^2+(x_n+1)^2}-1,& \mbox{when}\   i=n,
 \end{array}\right.
$$
for $i=1,2,\cdots,n-1$, we have
\begin{eqnarray*}
&  & \frac{\partial u_i}{\partial
x_i}=2\frac{|x'|^2+(x_n+1)^2-2x_i^2}{[|x'|^2+(x_n+1)^2]^2}, \\
&  & \frac{\partial u_i}{\partial
x_j}=\frac{-4x_ix_j}{[|x'|^2+(x_n+1)^2]^2},i\neq j,j=1,2,\cdots,n-1\\
&  & \frac{\partial u_i}{\partial
x_n}=\frac{-4x_i(x_n+1)}{[|x'|^2+(x_n+1)^2]^2}, \\
&  & \frac{\partial u_n}{\partial
x_i}=\frac{-4(x_n+1)x_i}{[|x'|^2+(x_n+1)^2]^2},\\
&  & \frac{\partial u_n}{\partial
x_n}=2\frac{|x'|^2+(x_n+1)^2-2(x_n+1)^2}{[|x'|^2+(x_n+1)^2]^2}.\\
\end{eqnarray*}

So we get
\begin{displaymath}
{J_{\phi}(x)} = \left| \begin{array}{cccc}
 \frac{\partial u_1}{\partial x_1} &
  \ldots &  \frac{\partial u_1}{\partial x_{n-1}} &
   \frac{\partial u_1}{\partial x_n}  \\
\vdots & \vdots & \vdots & \vdots\\
 \frac{\partial u_n}{\partial x_1} &
  \ldots &  \frac{\partial u_n}{\partial x_{n-1}} &
   \frac{\partial u_n}{\partial x_n}
\end{array} \right|\qquad \ \ \qquad \ \ \qquad \ \ \qquad \ \ \qquad \ \ \qquad \ \
\quad
\end{displaymath}

\begin{displaymath}
\mathbf{ }=\frac{2^n}{[|x'|^2+(x_n+1)^2]^{2n}} \left|
\begin{array}{cccc}
|x|^2-2x_1^2 & -2x_1x_2 & \ldots & -2x_1x_n \\
\vdots & \vdots & \vdots & \vdots\\
-2x_nx_1 & -2x_nx_2 & \ldots &|x|^2-2x_n^2
\end{array} \right|
\end{displaymath}

\qquad \ \ $=-\frac{2^n}{[|x'|^2+(x_n+1)^2]^n}.$

\vspace{0.4cm}

\section*{3.   Proof of Theorem }

 \vspace{0.2cm}
 We first prove necessity.

First applying Lemma 4.2.1  with $ H(x)\equiv1$, by (4.2.1), we have
$c=0$,
 so we obtain
$$
1=\frac{2x_n}{\omega_n}\int_{{\bf
R}^{n-1}}\frac{1}{[|y'-x'|^2+x_n^2]^{\frac{n}{2}}}dy'.
\eqno{(4.2.4)}
$$
For $a>0,$\ consider the function
$$
H(x)= \frac{x_n+a}{[|x'|^2+(x_n+a)^2]^{\frac{n}{2}}},
$$
it is clear that $H(x)$ is nonnegative and harmonic in $H$, then
applying Lemma 4.2.1 with
$$
H(x)= \frac{x_n+a}{[|x'|^2+(x_n+a)^2]^{\frac{n}{2}}},
$$
by (4.2.2) and
$$
H((0,x_n))= \frac{x_n+a}{(x_n+a)^n} = \frac{1}{(x_n+a)^{n-1}},
$$
we have $c=0$, so we obtain by
$$
H((y',0))= \frac{a}{[|y'|^2+a^2]^{\frac{n}{2}}}
$$
that
$$
\frac{x_n+a}{[|x'|^2+(x_n+a)^2]^{\frac{n}{2}}}
=\frac{2x_n}{\omega_n}\int_{{\bf
R}^{n-1}}\frac{1}{[|y'-x'|^2+x_n^2]^{\frac{n}{2}}}
\frac{a}{[|y'|^2+a^2]^{\frac{n}{2}}}dy'. \eqno{(4.2.5)}
$$
In these two formulas (4.2.4) and (4.2.5), we interchange the roles
of $x'$ and $y'$ and choose $a=x_n+1$, then we get
$$
1=\frac{2x_n}{\omega_n}\int_{{\bf
R}^{n-1}}\frac{1}{[|x'-y'|^2+x_n^2]^{\frac{n}{2}}}dx',
 \eqno{(4.2.6)}
$$
and
$$
\frac{2x_n+1}{[|y'|^2+(2x_n+1)^2]^{\frac{n}{2}}}
=\frac{2x_n}{\omega_n}\int_{{\bf
R}^{n-1}}\frac{1}{[|x'-y'|^2+x_n^2]^{\frac{n}{2}}}
\frac{x_n+1}{[|x'|^2+(x_n+1)^2]^{\frac{n}{2}}}dx',
\eqno{(4.2.7)}
$$
by (4.2.6), we can also get
\begin{eqnarray*}
1
& = & \frac{2(x_n+1)}{\omega_n}\int_{{\bf
R}^{n-1}}\frac{1}{[|x'-y'|^2+(x_n+1)^2]^{\frac{n}{2}}}dx' \\
& = & \frac{2(x_n+1)}{\omega_n}\int_{{\bf
R}^{n-1}}\frac{1}{[|x'|^2+(x_n+1)^2]^{\frac{n}{2}}}dx'.
\hspace{46mm} (4.2.8)
\end{eqnarray*}
Thus from
$$
G(x)\leq cx_n+\int_{{\bf R}^{n-1}}P(x,y')d\mu(y'),
$$
we obtain by (4.2.7) and (4.2.8)
\begin{eqnarray*}
&  & \int_{{\bf R}^{n-1}}
\frac{G(x)}{[|x'|^2+(x_n+1)^2]^{\frac{n}{2}}}dx' \\
& \leq & \int_{{\bf R}^{n-1}}
\frac{cx_n}{[|x'|^2+(x_n+1)^2]^{\frac{n}{2}}}dx' \\
& & + \frac{2}{\omega_n}\int_{{\bf
R}^{n-1}}\frac{x_n}{[|x'|^2+(x_n+1)^2]^{\frac{n}{2}}}dx'\int_{{\bf
R}^{n-1}}\frac{1}{[|y'-x'|^2+x_n^2]^{\frac{n}{2}}}d\mu(y') \\
& = & \frac{c\omega_nx_n}{2(x_n+1)} \\
& & + \frac{1}{x_n+1}\int_{{\bf
R}^{n-1}}\frac{2x_n}{\omega_n}\int_{{\bf
R}^{n-1}}\frac{1}{[|x'-y'|^2+x_n^2]^{\frac{n}{2}}}
\frac{x_n+1}{[|x'|^2+(x_n+1)^2]^{\frac{n}{2}}}dx'd\mu(y')\\
& = & \frac{c\omega_nx_n}{2(x_n+1)}+ \frac{2x_n+1}{x_n+1}\int_{{\bf
R}^{n-1}}\frac{1}{[|y'|^2+(2x_n+1)^2]^{\frac{n}{2}}}d\mu(y')\\
& \leq &\frac{c\omega_n}{2}+ 2\int_{{\bf
R}^{n-1}}\frac{1}{[1+|y'|^2]^{\frac{n}{2}}}d\mu(y').
\end{eqnarray*}
Hence (4.2.1) holds.

 In the other direction, assume that (4.2.1) holds. We show that
$G(x)$ has a harmonic majorant in $H$, or what is the same thing,
$G(\Phi^{-1} (u))$ has a harmonic majorant in $B_n$. It is
sufficient to show that
$$
g(r)=r^{n-1}\int_S G(\Phi^{-1} (r\xi))d\sigma (\xi)
$$
remains bounded as $r \uparrow 1$. By Lemma $A_1$, it is the same
thing to show that
$$
\int_0^1\frac{g(r)}{1-\lambda^nr^n}dr \leq
A\int_0^1\frac{1}{1-\lambda^nt^n}dt
$$
for all $\lambda \in (0,1)$ and some positive constant. Calculate as
follows:

\begin{eqnarray*}
&  & \int_0^1\frac{g(r)}{1-\lambda^nr^n}dr \\
& = & \int_0^1\int_S\frac{G(\Phi^{-1}
(r\xi))}{1-\lambda^nr^n}d\sigma
(\xi) r^{n-1}dr\\
&= & \frac{1}{nV(B)}\int_{B_n}\frac{G(\Phi^{-1}
(u))}{1-\lambda^n|u|^n}dV \\
&=&  \frac{1}{nV(B)}\int_{{\bf
R}_n^{+}}\frac{G(x)}{1-\lambda^n|\Phi(x)|^n}|J_\Phi (x)|dx \\
&=&\frac{1}{nV(B)}\int_{{\bf R}_n^{+}}\frac{G(x)}{1-\lambda^n
\bigg[\frac{|x'|^2+(1-x_n)^2}{|x'|^2+(x_n+1)^2}\bigg]^\frac{n}{2}}
\bigg|-\frac{2^n}{[|x'|^2+(x_n+1)^2]^n}\bigg|dx\\
&=& \frac{2^n}{nV(B)}\int_{{\bf R}_n^{+}}
\frac{G(x)}{[|x'|^2+(x_n+1)^2]^\frac{n}{2}}\frac{1}{[|x'|^2+(x_n+1)^2]^\frac{n}{2}-
\{\lambda^2[|x'|^2+(x_n-1)^2]\}^\frac{n}{2}}dx\\
& \le & \frac{2^n}{nV(B)}\int_0^{\infty}\bigg[\int_{{\bf
R}_{n-1}}\frac{G(x)}{[|x'|^2+(x_n+1)^2]^\frac{n}{2}}dx'\bigg] \cdot
\frac{1}{(x_n+1)^n- \lambda^n|x_n-1|^n}dx_n,
\end{eqnarray*}

by (4.2.1), we have
\begin{eqnarray*}
&  & \int_0^1\frac{g(r)}{1-\lambda^nr^n}dr \\
&\leq & A\int_0^{\infty} \frac{1}{(x_n+1)^n- \lambda^n|x_n-1|^n}dx_n \\
&\leq&  A\int_0^{\infty} \frac{(x_n+1)^{n-2}}{(x_n+1)^n-
\lambda^n|x_n-1|^n}dx_n \\
&=& A\int_0^{\infty} \frac{1}{1- \lambda^n|\frac{x_n-1}{x_n+1}|^n}
\frac{1}{(x_n+1)^2}dx_n\\
& = & \frac{1}{2}A\int_{-1}^{1} \frac{1}{1- \lambda^n|t|^n} dt\\
& = & A\int_0^{1} \frac{1}{1- \lambda^nt^n} dt
\end{eqnarray*}
The change of variables is made with the substitution
$t=\frac{x_n-1}{x_n+1}$. So the result follows.

\chapter{Properties of Limit for Poisson Integral}

\section{Properties of Limit for Poisson Integral in the Upper Half Plane}

 \section*{ 1. Introduction and Main Theorem}

\vspace{0.2cm}

  The Poisson kernel for the upper
half plane ${\bf C}_+=\{z=x+iy \in {\bf C}:\; y>0\}$ is the
 function
$$
 P(z,t)=
\frac{y}{\pi|z-t|^2},
$$
where $z \in {\bf C}_+$ and $t\in {\bf R}$.

If $p\geq 0$ is an integer, we define a modified
 Cauchy kernel of order $ p$ for $z\in {\bf C}_+-\{t\}$ by
$$
 C_p(z,t)=\left\{\begin{array}{ll}
 \frac{1}{\pi}\frac{1}{t-z} , &   \mbox{when }   |t|\leq 1  ,\\
 \frac{1}{\pi}\frac{1}{t-z} -\frac{1}{\pi}\sum_{k=0}^{p}\frac{z^k}{t^{k+1}},&
\mbox{when}\   |t|> 1,
 \end{array}\right.
$$
then we define a modified
 Poisson kernel of order $p$ for the upper half plane  by
$$
P_p(z,t)=\Im C_p(z,t).
$$

  For any $|t|>1$,  the following equality
$$
\Im C_p(z,t)=\frac{tR^{p+1}\sin (p+1)\theta-R^{p+2}\sin p\theta }
{|t-z|^2 t^{p+1}} \eqno{(5.1.1)}
$$
holds, where $z=Re^{i\theta}$.

 Marvin Rosenblum and James Rovnyak \cite{RR}
 proved the following theorem:

 \vspace{0.2cm}
 \noindent
{\bf Theorem E } {\it  If
$$
H(z)=cy+\frac{y}{\pi}\int_{-\infty}^{\infty}\frac{d\mu(t)}{(t-x)^2+y^2},\quad
y>0,
$$
where $c$ is a real number and $\mu$ is a nonnegative Borel measure
on $(-\infty,\infty)$ such that
$$
\int_{-\infty}^{\infty}\frac{d\mu(t)}{1+t^2}<\infty.
$$
Then for every $\theta  \in (0,\pi)$:

\noindent {\rm (1)}\
$$
\lim_{R\rightarrow\infty} \frac{1}{R}H(R e^{i\theta})=c\sin\theta;
$$

 \noindent
 {\rm (2)}\
$$
\lim_{R\rightarrow\infty} \frac{2}{\pi R}\int_0^\pi H(R
e^{i\theta})\sin\theta d\theta=c.
$$
}

  In this section, We will generalize Theorem E to the modified kernel.

 \vspace{0.2cm}
 \noindent
{\bf Theorem 5.1.1 } {\it  If
$$
H(z)=\Im
\bigg[Q_p(z)+\frac{1}{\pi}\int_{-\infty}^{\infty}C_p(z,t)d\mu(t)\bigg],\quad
z=x+iy, \ y>0
$$
where
$$
Q_p(z)=\sum _{k=0}^{p}a_k z^k, \quad a_k\in {\bf R},\
k=0,1,2,\cdots,p
$$
and $\mu$ is a nonnegative Borel measure on $(-\infty,\infty)$ such
that
$$
\int_{-\infty}^{\infty}\frac{1}{1+|t|^{p+1}}d\mu(t)<\infty.
$$
Then for every $\theta  \in (0,\pi):$

\noindent
 {\rm (1)}\
$$
\lim_{R\rightarrow\infty} \frac{1}{R^p}H(R e^{i\theta})
 =\bigg[a_p-\frac{1}{\pi}
 \int_{|t|>1}\frac{1}{t^{p+1}}d\mu(t)\bigg]\sin p \theta;
$$

\noindent
 {\rm (2)}\
$$
\lim_{R\rightarrow\infty} \frac{2}{\pi R^p}\int_0^\pi H(R
e^{i\theta})\sin p\theta d\theta=a_p-\frac{1}{\pi}
 \int_{|t|>1}\frac{1}{t^{p+1}}d\mu(t).
$$
}

\vspace{0.2cm}
 \noindent
{\bf Remark 5.1.1 } {\it If $p=1$,  this is just the result of
 Marvin Rosenblum and James Rovnyak, therefore, our result is the generalization of
Theorem E. }

\vspace{0.4cm}

\section*{2.   Proof of Theorem }

\vspace{0.3cm}

 We first prove the equality (5.1.1). Since
$$
 C_p(z,t)=\frac{1}{t-z} - \sum_{k=0}^{p}\frac{z^k}{t^{k+1}}
 =\frac{z^{p+1}}{(t-z)t^{p+1}},
$$
then
\begin{eqnarray*}
\Im C_p(z,t)
& =& \Im \frac{z^{p+1}}{(t-z)t^{p+1}}
=\Im  \bigg[\frac{z^{p+1}(t-\overline{z})}{|t-z|^2 t^{p+1}} \bigg] \\
&=& \Im  \bigg[\frac{t(R e^{i\theta})^{p+1}-|z|^2(R e^{i\theta})^{p}
}
{|t-z|^2 t^{p+1}} \bigg]\\
&=& \frac{tR^{p+1}\sin (p+1)\theta-R^{p+2}\sin p\theta } {|t-z|^2
t^{p+1}}. \\
\end{eqnarray*}
This proves the equality (5.1.1).

Since
\begin{eqnarray*}
H(z)
& =& \Im  \bigg(\sum _{k=0}^{p}a_k
z^k  \bigg)+\frac{1}{\pi}\int_{-\infty}^{\infty}\Im C_p(z,t)d\mu(t) \\
& =& \Im  \bigg(\sum _{k=0}^{p}a_k R^k e^{ik\theta}
\bigg)+\frac{1}{\pi}  \bigg[\int_{|t|\leq 1}\Im
C_p(z,t)d\mu(t)+\int_{|t|> 1}\Im C_p(z,t)d\mu(t)  \bigg] \\
& =& \sum _{k=0}^{p}a_k R^k \sin k\theta+
\frac{1}{\pi}\int_{|t|\leq1}\Im\frac{1}{t-z}d\mu(t)+
\frac{1}{\pi}\int_{|t|>1}\Im\frac{z^{p+1}}{(t-z)t^{p+1}}d\mu(t) \\
&=& \sum _{k=0}^{p}a_k R^k \sin k\theta+
\frac{1}{\pi}\int_{|t|\leq1}\frac{y}{|t-z|^2}d\mu(t)\\
& & +\frac{1}{\pi}\int_{|t|>1}\frac{tR^{p+1}\sin
(p+1)\theta-R^{p+2}\sin
p\theta } {|t-z|^2 t^{p+1}}d\mu(t), \\
\end{eqnarray*}
for every $\theta \in (0,\pi),$
\begin{eqnarray*}
\frac{H(R e^{i\theta})}{R^p}
& =& \sum _{k=0}^{p}a_k R^{k-p} \sin
k\theta+ \frac{1}{\pi}\int_{|t|\leq1}\frac{y}{|t-z|^2 R^p}d\mu(t) \\
&  & +\frac{1}{\pi}\int_{|t|>1}\frac{R[t\sin (p+1)\theta-R\sin
p\theta]
} {|t-z|^2 t^{p+1}}d\mu(t) \\
&=& I_1+I_2+I_3, \hspace{80mm} (5.1.2)
\end{eqnarray*}
then
$$
\lim_{R\rightarrow\infty} I_1=\lim_{R\rightarrow\infty}
\sum_{k=0}^{p}a_k R^{k-p} \sin k\theta
 =a_p\sin p \theta. \eqno{(5.1.3)}
$$
Moreover,
\begin{eqnarray*}
I_2
& =& \frac{1}{\pi}\int_{|t|\leq1}\frac{y}{|t-z|^2 R^p}d\mu(t) \\
&=& \frac{1}{\pi}\int_{|t|\leq1}\frac{\sin \theta
(1+|t|^{p+1})}{R^{p-1}|t-R e^{i\theta}|^2}
\frac{d\mu(t)}{1+|t|^{p+1}}.
\end{eqnarray*}
Since
$$
\frac{\sin \theta (1+|t|^{p+1})}{R^{p-1}|t-R e^{i\theta}|^2}
<\frac{2}{2^{p-1}(R-|t|)^2} <\frac{2}{2^{p-1}(R-1)^2} <2^{2-p},
$$
by the dominated convergence theorem, we have
$$
\lim_{R\rightarrow\infty} I_2=0. \eqno{(5.1.4)}
$$
Write
\begin{eqnarray*}
I_3
& =& \frac{1}{\pi}\int_{|t|>1}\frac{R[t\sin (p+1)\theta-R\sin
p\theta] } {|t-z|^2 }\frac{d\mu(t)}{t^{p+1}} \\
&=& \frac{1}{\pi}\int_{|t|>1}J\frac{d\mu(t)}{t^{p+1}}. \\
\end{eqnarray*}

  Multiplying (5.1.2) by $2\pi^{-1}\sin p\theta$ and integrating
with respect to $\theta$, we obtain
\begin{eqnarray*}
& & \frac{2}{\pi R^p}\int_0^\pi H(R e^{i\theta})\sin p\theta d\theta \\
& =& \frac{2}{\pi}\int_0^\pi \sum_{k=0}^{p}a_k R^{k-p} \sin
k\theta\sin p\theta d\theta \\
& & + \frac{2}{\pi^2}\int_0^\pi \int_{|t|\leq1}\frac{y\sin
p\theta}{|t-z|^2 R^p}d\mu(t) d\theta \\
& &+ \frac{2}{\pi^2}\int_0^\pi \int_{|t|>1}\frac{R\sin p\theta[t\sin
(p+1)\theta-R\sin p\theta] } {|t-z|^2 }\frac{d\mu(t)}{t^{p+1}}
d\theta \\
&=& I_1'+I_2'+I_3'. \hspace{94mm} (5.1.5)
\end{eqnarray*}
For the first term, we have
\begin{eqnarray*}
I_1'
& =& \frac{2}{\pi}\sum_{k=0}^{p}a_k R^{k-p} \int_0^\pi \sin
k\theta\sin p\theta d\theta \\
& =& \frac{1}{\pi}\sum_{k=0}^{p}a_k R^{k-p} \int_0^\pi [\cos
(k-p)\theta-\cos (k+p)\theta] d\theta \\
& =& \frac{1}{\pi} \bigg[\sum_{k=0}^{p-1}+\sum_{k=p} \bigg] a_k
R^{k-p} \int_0^\pi [\cos
(k-p)\theta-\cos (k+p)\theta] d\theta\\
&= &\frac{1}{\pi}a_p \int_0^\pi (1-\cos 2p\theta) d\theta \\
&=& a_p; \hspace{105mm} (5.1.6)
\end{eqnarray*}
for the second term, we have
$$
I_2'=\frac{2}{\pi^2}\int_0^\pi \int_{|t|\leq1}\frac{\sin \theta\sin
p\theta(1+|t|^{p+1})}{R^{p-1}|t-R
e^{i\theta}|^2}\frac{d\mu(t)}{1+|t|^{p+1}} d\theta
$$
for all $\theta  \in (0,\pi)$ and $R>2$.

If $|t|\leq 1$, then
\begin{eqnarray*}
& &  \bigg|\frac{\sin \theta\sin p\theta(1+|t|^{p+1})}{R^{p-1}|t-R
e^{i\theta}|^2} \bigg| \\
&<& \frac{2}{2^{p-1}(R-|t|)^2}<\frac{2}{2^{p-1}(R-1)^2} \\
&<&\frac{2}{2^{p-1}}=2^{2-p}.
\end{eqnarray*}
Since
$$
\frac{2}{\pi^2}\int_0^\pi
\int_{|t|\leq1}2^{2-p}\frac{d\mu(t)}{1+|t|^{p+1}} d\theta<\infty,
$$
by the dominated convergence theorem, we have
$$
\lim_{R\rightarrow\infty} I_2'=0; \eqno{(5.1.7)}
$$
for the third term, we have
\begin{eqnarray*}
I_3'
&=&  \frac{2}{\pi^2}\int_0^\pi \int_{|t|>1}\frac{R\sin p\theta[t\sin
(p+1)\theta-R\sin p\theta] } {|t-z|^2 }\frac{d\mu(t)}{t^{p+1}}
d\theta \\
&=& \frac{2}{\pi^2}\int_0^\pi \int_{|t|>1}J'\frac{d\mu(t)}{t^{p+1}}
d\theta.
\end{eqnarray*}

In the following, we will show that
$$
J' \leq2p(p+1) \eqno{(5.1.8)}
$$
for all $\theta  \in (0,\pi)$ and $R>2$.

  Since
$$
|t-z|^2=t^2-2Rt\cos\theta+R^2=(t-R\cos\theta)^2+R^2\sin^2\theta \geq
R^2\sin^2\theta
$$
and
$$
|t-z|^2=\left\{\begin{array}{ll}
 (t-R)^2+2Rt(1-\cos\theta)=(t-R)^2+4Rt\sin^2\frac{\theta}{2}
\geq 4R|t|\sin^2\frac{\theta}{2} , &   \mbox{when }   t>0  ,\\
(t+R)^2-2Rt(1+\cos\theta)=(t+R)^2-4Rt\cos^2\frac{\theta}{2} \geq
4R|t|\cos^2\frac{\theta}{2},& \mbox{when}\   t<0,
 \end{array}\right.
$$
then
\begin{eqnarray*}
|J'|
&\leq  & \frac{R|\sin p\theta|(|t\sin (p+1)\theta|-R|\sin p\theta|)
} {|t-z|^2 } \\
&\leq & \frac{|t|R p(p+1)\sin^2 \theta+R^2p^2\sin^2\theta } {|t-z|^2 } \\
&\leq & p(p+1)\frac{|t|R\sin^2 \theta+R^2\sin^2\theta } {|t-z|^2 } \\
&= & p(p+1) \bigg(\frac{|t|R\sin^2 \theta } {|t-z|^2
}+\frac{R^2\sin^2\theta} {|t-z|^2 } \bigg) \\
&\leq &  2p(p+1).
\end{eqnarray*}
This proves (5.1.8).

Since
$$ \frac{2}{\pi^2}\int_0^\pi
\int_{|t|>1}2p(p+1)\frac{d\mu(t)}{t^{p+1}} d\theta
=\frac{4p(p+1)}{\pi} \int_{|t|>1}\frac{d\mu(t)}{t^{p+1}} <\infty,
$$
by the dominated convergence theorem, we have
$$
\lim_{R\rightarrow\infty}
I_3'=\lim_{R\rightarrow\infty}\frac{2}{\pi^2}\int_0^\pi
\int_{|t|>1}J'\frac{d\mu(t)}{t^{p+1}} d\theta.
$$
Note that
\begin{eqnarray*}
J'
&= & J'+ \sin^2 p\theta- \sin^2 p\theta \\
&= & \frac{Rt\sin (p+1)\theta\sin p\theta-R^2\sin^2 p\theta }
{|t-z|^2}
+\frac{|t-z|^2 \sin^2 p\theta} {|t-z|^2}- \sin^2 p\theta \\
&= & \frac{Rt\sin (p+1)\theta\sin p\theta-R^2\sin^2
p\theta+(t^2-2Rt\cos\theta+R^2)\sin^2 p\theta } {|t-z|^2} - \sin^2
p\theta\\
&= & \frac{t^2\sin^2 p\theta +Rt\sin p\theta[\sin
(p+1)\theta-2\cos\theta\sin p\theta]} {|t-z|^2} - \sin^2 p\theta \\
&= &  \frac{t^2\sin^2 p\theta -Rt\sin p\theta\sin (p-1)\theta}
{|t-z|^2} - \sin^2 p\theta,
\end{eqnarray*}
then we obtain
$$
\lim_{R\rightarrow\infty} J'=- \sin^2 p\theta.
$$
Therefore
\begin{eqnarray*}
\lim_{R\rightarrow\infty} I_3'
&= & \frac{2}{\pi^2}\int_0^\pi
\int_{|t|>1}\lim_{R\rightarrow\infty}J'\frac{d\mu(t)}{t^{p+1}}
d\theta \\
&= & \frac{2}{\pi^2}\int_0^\pi \int_{|t|>1}- \sin^2
p\theta\frac{d\mu(t)}{t^{p+1}} d\theta \\
&= & -\frac{2}{\pi^2}\int_0^\pi \frac{1-\cos 2p\theta}{2}d\theta
\int_{|t|>1}\frac{d\mu(t)}{t^{p+1}} \\
&= &  -\frac{1}{\pi} \int_{|t|>1}\frac{d\mu(t)}{t^{p+1}}.
\hspace{70mm} (5.1.9)
\end{eqnarray*}
Thus, (2) holds by collecting (5.1.5), (5.1.6), (5.1.7) and (5.1.9).

  Similarly, we have
$$
|p\sin \theta J|\leq 2p(p+1),
$$
so
$$
|J|\leq \frac{2(p+1)}{\sin \theta}.
$$
Since
$$
\frac{1}{\pi} \int_{|t|>1} \frac{2(p+1)}{\sin
\theta}\frac{d\mu(t)}{t^{p+1}} =\frac{2(p+1)}{\pi \sin \theta}
\int_{|t|>1} \frac{d\mu(t)}{t^{p+1}}<\infty
$$
and
$$
\lim_{R\rightarrow\infty} J=- \sin p\theta,
$$
by the dominated convergence theorem, we have
$$
\lim_{R\rightarrow\infty} I_3=\frac{1}{\pi}
\int_{|t|>1}\lim_{R\rightarrow\infty}J\frac{d\mu(t)}{t^{p+1}}
=-\frac{\sin p\theta}{\pi} \int_{|t|>1}\frac{d\mu(t)}{t^{p+1}}.
\eqno{(5.1.10)}
$$
So (1) follows by (5.1.2), (5.1.3), (5.1.4) and (5.1.10).

\section{Properties of Limit for Poisson Integral in the Upper Half Space}

 \section*{ 1. Introduction and Main Theorem}

\vspace{0.2cm}

  The Poisson kernel for the upper
half space $H $ is the
 function
$$
 P(x,y')=
\frac{2x_n}{\omega_n|x-y'|^n},
$$
where $x \in H$, $y'\in \partial{H}$ and $ \omega_n = \frac{2\pi
^{\frac{n}{2}}}{\Gamma (\frac{n}{2})}$ is the area of the unit
sphere in ${\bf R}^n$.

   In this section, We will generalize Theorem E to the
upper half space.

 \vspace{0.2cm}
 \noindent
{\bf Theorem 5.2.1} {\it  If
$$
H(x)=cx_n+\int_{{\bf R}^{n-1}}P(x,y')d\mu(y'),
$$
where $c$ is a real number and $\mu$ is a nonnegative Borel measure
on ${\bf R}^{n-1}$ such that
$$
\int_{{\bf
R}^{n-1}}\frac{1}{(1+|y'|^2)^{\frac{n}{2}}}d\mu(y')<\infty.
$$
Then for every $\theta_j  \in (0,\pi),j=1,2,\cdots,n-1,$

\noindent {\rm (1)}\
$$
\lim_{R\rightarrow\infty}
\frac{1}{R}H(x)=c\sin\theta_1\sin\theta_2\cdots\sin\theta_{n-1};
$$

\noindent
 {\rm (2)}\
$$
\lim_{R\rightarrow\infty}
\frac{1}{R}\int_0^\pi\int_0^\pi\cdots\int_0^\pi
H(x)(\sin\theta_1\sin\theta_2\cdots\sin\theta_{n-1})^{n-1}
d\theta_1d\theta_2\cdots d\theta_{n-1}=2^{n-1}I_n^{n-1}c,
$$
where $R=|x|$ and
$$
 I_n=\left\{\begin{array}{ll}
 \frac{(2k-1)!!}{(2k)!!}\frac{\pi}{2} , &   \mbox{when }   n=2k  ,\\
 \frac{(2k)!!}{(2k+1)!!},&
\mbox{when}\   n=2k+1 .
 \end{array}\right.
$$
}

\vspace{0.2cm}
 \noindent
{\bf Remark 5.2.1 } {\it If $n=2$,  this is just the result of
 Marvin Rosenblum and James Rovnyak, therefore, our result is the generalization of
Theorem E. }

\vspace{0.4cm}

\section*{2.   Proof of Theorem }

\vspace{0.3cm}

Write $x'=R\xi,\ x_n=R\eta$, by the formula of polar coordinates, we
have
$$
\eta=\sin\theta_1\sin\theta_2\cdots\sin\theta_{n-1}.  \eqno{(5.2.1)}
$$
 For every $\theta_j  \in (0,\pi),j=1,2,\cdots,n-1,$
\begin{eqnarray*}
\frac{1}{R}H(x)
&= & c\eta+\frac{2}{\omega_n}\int_{{\bf
R}^{n-1}}\frac{2R\eta}{\omega_n|x-y'|^n}d\mu(y') \\
&= & c\eta+\frac{2}{\omega_n}\int_{{\bf
R}^{n-1}}\eta\frac{(1+|y'|^2)^{n/2}}{|x-y'|^n}\frac{1}{(1+|y'|^2)^{n/2}}d\mu(y') \\
&=& c\eta+\frac{2}{\omega_n}\int_{{\bf
R}^{n-1}}J_1\frac{1}{(1+|y'|^2)^{n/2}}d\mu(y').\hspace{40mm} (5.2.2)
\end{eqnarray*}

  Multiplying this by $\eta^{n-1}$ and integrating
with respect to $\theta_j$, we obtain
\begin{eqnarray*}
& & \frac{1}{R}\int_0^\pi\int_0^\pi\cdots\int_0^\pi
H(x)\eta^{n-1} d\theta_1d\theta_2\cdots d\theta_{n-1} \\
&= & \int_0^\pi\int_0^\pi\cdots\int_0^\pi c\eta^{n}
d\theta_1d\theta_2\cdots d\theta_{n-1} \\
& & +\frac{2}{\omega_n}\int_0^\pi\int_0^\pi\cdots\int_0^\pi
\int_{{\bf
R}^{n-1}}\eta^n\frac{(1+|y'|^2)^{n/2}}{|x-y'|^n}\frac{1}{(1+|y'|^2)^{n/2}}d\mu(y')
d\theta_1d\theta_2\cdots d\theta_{n-1} \\
&= & K_1+K_2. \hspace{100mm} (5.2.3)
\end{eqnarray*}
By (5.2.1), we have
\begin{eqnarray*}
K_1
&= & c\int_0^\pi\int_0^\pi\cdots\int_0^\pi
(\sin\theta_1\sin\theta_2\cdots\sin\theta_{n-1})^{n}
d\theta_1d\theta_2\cdots d\theta_{n-1} \\
&=& c \bigg(\int_0^\pi \sin^n \theta d\theta \bigg)^{n-1}.
\end{eqnarray*}
Since
$$
\int_0^\pi \sin^n \theta d\theta=2\int_0^{\pi/2} \sin^n \theta
d\theta=2I_n,
$$
we obtain
$$
K_1=2^{n-1}I_n^{n-1}c.\eqno{(5.2.4)}
$$
Moreover,
$$
K_2=\frac{2}{\omega_n}\int_0^\pi\int_0^\pi\cdots\int_0^\pi
\int_{{\bf
R}^{n-1}}\eta^n\frac{(1+|y'|^2)^{n/2}}{|x-y'|^n}\frac{1}{(1+|y'|^2)^{n/2}}d\mu(y')
d\theta_1d\theta_2\cdots d\theta_{n-1}
$$
$$
=\frac{2}{\omega_n}\int_0^\pi\int_0^\pi\cdots\int_0^\pi \int_{{\bf
R}^{n-1}}J_2\frac{1}{(1+|y'|^2)^{n/2}}d\mu(y')
d\theta_1d\theta_2\cdots d\theta_{n-1}.
$$

  In the following, we will show that
$$
J_2 \leq2^{n/2} \eqno{(5.2.5)}
$$
for all $y',\theta_j  \in (0,\pi),j=1,2,\cdots,n-1,$ and $R>2$.

If $|y'|<1$, then since $R>2$,
$$
J_2\leq \frac{2^{n/2}}{(R-|y'|)^n}
 \leq \frac{2^{n/2}}{(R-1)^n}
 \leq2^{n/2};
$$

if $|y'|\geq1$, since
\begin{eqnarray*}
& & |x-y'|^2=|(x',x_n)-(y',0)|^2=|x'-y'|^2+x_n^2 \\
&= & |y'|^2-2x'\cdot y'+R^2
=|y'|^2(|\xi|^2+\eta^2)-2y'\cdot R\xi+R^2 \\
&=& |y'|^2\eta^2+(|y'|^2|\xi|^2-|y'\cdot \xi|^2)+(y'\cdot \xi-R)^2
\geq|y'|^2\eta^2,
\end{eqnarray*}

then
\begin{eqnarray*}
J_2
& =&\eta^n\frac{(1+|y'|^2)^{n/2}}{|x-y'|^n}
=\bigg[\eta^2\frac{(1+|y'|^2)}{|x-y'|^2} \bigg]^{n/2} \\
&\leq &  \bigg[\eta^2\frac{(1+|y'|^2)}{|y'|^2\eta^2}\bigg]^{n/2}
=\bigg(1+\frac{1}{|y'|^2}\bigg)^{n/2}\leq 2^{n/2}.
\end{eqnarray*}
This proves (5.2.5).

Since
\begin{eqnarray*}
& & \frac{2}{\omega_n}\int_0^\pi\int_0^\pi\cdots\int_0^\pi
\int_{{\bf R}^{n-1}}2^{n/2}\frac{1}{(1+|y'|^2)^{n/2}}d\mu(y')
d\theta_1d\theta_2\cdots d\theta_{n-1} \\
&=& \frac{2}{\omega_n}2^{n/2}\pi^{n-1} \int_{{\bf
R}^{n-1}}\frac{1}{(1+|y'|^2)^{n/2}}d\mu(y')< \infty,
\end{eqnarray*}
by the dominated convergence theorem, we have
\begin{eqnarray*}
\lim_{R\rightarrow\infty} K_2
&= & \lim_{R\rightarrow\infty}
\frac{2}{\omega_n}\int_0^\pi\int_0^\pi\cdots\int_0^\pi \int_{{\bf
R}^{n-1}}J_2\frac{1}{(1+|y'|^2)^{n/2}}d\mu(y')
d\theta_1d\theta_2\cdots d\theta_{n-1} \\
&=&\lim_{R\rightarrow\infty}
\frac{2}{\omega_n}\int_0^\pi\int_0^\pi\cdots\int_0^\pi \int_{{\bf
R}^{n-1}}\lim_{R\rightarrow\infty}J_2\frac{1}{(1+|y'|^2)^{n/2}}d\mu(y')
d\theta_1d\theta_2\cdots d\theta_{n-1} \\
&=& 0.  \hspace{100mm} (5.2.6)
\end{eqnarray*}
Thus, (2) holds by collecting (5.2.3), (5.2.4) and (5.2.6).

 Moreover,
\begin{eqnarray*}
J_1
& =& \eta\frac{(1+|y'|^2)^{n/2}}{|x-y'|^n}=\frac{J_2}{\eta^{n-1}} \\
&\leq &
\frac{2^{n/2}}{(\sin\theta_1\sin\theta_2\cdots\sin\theta_{n-1})^{n-1}}.
\end{eqnarray*}

Since
\begin{eqnarray*}
& & \frac{2}{\omega_n}\int_{{\bf R}^{n-1}}\frac{2^{n/2}}
{(\sin\theta_1\sin\theta_2\cdots\sin\theta_{n-1})^{n-1}}
\frac{1}{(1+|y'|^2)^{n/2}}d\mu(y') \\
&=& \frac{2}{\omega_n}\frac{2^{n/2}}
{(\sin\theta_1\sin\theta_2\cdots\sin\theta_{n-1})^{n-1}}\int_{{\bf
R}^{n-1}} \frac{1}{(1+|y'|^2)^{n/2}}d\mu(y')<\infty,
\end{eqnarray*}
by the dominated convergence theorem, we have
\begin{eqnarray*}
& & \lim_{R\rightarrow\infty} \frac{2}{\omega_n}\int_{{\bf
R}^{n-1}}J_1\frac{1}{(1+|y'|^2)^{n/2}}d\mu(y') \\
&=& \frac{2}{\omega_n}\int_{{\bf
R}^{n-1}}\lim_{R\rightarrow\infty}J_1\frac{1}{(1+|y'|^2)^{n/2}}d\mu(y')=0.
\hspace{45mm} (5.2.7)
\end{eqnarray*}
So (1) follows by (5.2.2) and (5.2.7).

\newpage
\ \thispagestyle{empty}
\newpage
\chapter{a Lower Bound for a Class of  Harmonic Functions in the Half Space}

\section{Introduction and Main Theorem}

   B.Ya.Levin \cite{LLST} has proved the
  following result:

\vspace{0.2cm}
 \noindent
{\bf Theorem F } {\it Let $u(z)$ be a harmonic function in the upper
half plane ${\bf C}_{+}=\{z=x+iy=Re^{i\theta}, y>0\}$ with
continuous boundary values on the real axis. Suppose that
$$
u(z)\leq KR^\rho, \quad z\in {\bf C}_{+},R=|z|>1,\rho>1,
$$
and
$$
|u(z)|\leq K, \quad z\in \overline{{\bf C}_{+}}, R=|z|\leq1,\Im
z\geq0.
$$
Then
$$
u(z)\geq -cK\frac{1+R^\rho}{\sin\theta}, \quad z\in {\bf C}_{+},
$$
where $c$ does not depend on $K,R,\theta$ and the function $u(z)$. }

  Our aim in this chapter is to establish the following main theorem.

\vspace{0.2cm}
 \noindent
{\bf Theorem 6.1.1} {\it Let $u(x)$ be a harmonic function in the
upper half space $H$ with continuous boundary values on the boundary
$\partial{H}$, write $|x'|=|x|\cos\theta, x_n=|x|\sin\theta$
$(0<\theta\leq \pi/2)$. Suppose that
$$
u(x)\leq KR^{\rho(R)}, \quad x\in H,R=|x|>1,\rho(R)>1,
 \eqno{(6.1.1)}
$$
and
$$
u(x)\geq -K, \quad x\in \overline{H}, R=|x|\leq1,x_n\geq0.
\eqno{(6.1.2)}
$$
Then
$$
u(x)\geq -cK\frac{1+(2R)^{\rho(R)}}{\sin^{n-1}\theta}, \quad x\in
H,\eqno{(6.1.3)}
$$
where $c$ does not depend on $K,R,\theta$ and the function $u(x)$,
$\rho(R)$ is nondecreasing in $[1,+\infty)$. }

\vspace{0.2cm}
 \noindent
 {\bf Remark 6.1.1} {\it If $n=2, \rho(R)\equiv \rho$,  this is just the result of
B.Ya.Levin, therefore, our result (6.1.3) is the generalization of
Theorem F. }

\section{Main Lemmas}

 In order to obtain the result, we
need these lemmas below:

\vspace{0.2cm}
 \noindent
{\bf Lemma 6.2.1 } {\it Let $u(x)$ be a harmonic function in the
upper half space $H=\{x=(x',x_n)\in {\bf R}^{n}:\; x_n>0\}$ with
continuous boundary values on the boundary $\partial{H}$, $R>1$.
Then we have

\begin{eqnarray*}
& & \int_{\{x\in {\bf
R}^{n}:\;|x|=R,x_n>0\}}u(x)\frac{nx_n}{R^{n+1}}d\sigma(x) \\
&+ &  \int_{\{x\in {\bf R}^{n}:\;1<|x'|<R,x_n=0\}}u(x')
\bigg(\frac{1}{|x'|^{n}}-\frac{1}{R^{n}} \bigg)dx'
=c_1+\frac{c_2}{R^{n}},
\end{eqnarray*}

where
\begin{eqnarray*}
& & c_1=\int_{\{x\in {\bf R}^{n}:\;|x|=1,x_n>0\}}
\bigg[(n-1)x_nu(x)+x_n\frac{\partial
u(x)}{\partial n} \bigg]d\sigma(x), \\
& &c_2=\int_{\{x\in {\bf R}^{n}:\;|x|=1,x_n>0\}}
\bigg[x_nu(x)-x_n\frac{\partial u(x)}{\partial n} \bigg]d\sigma(x).
\end{eqnarray*}
}

\vspace{0.2cm}
 \noindent
{\bf Lemma 6.2.2 } {\it Let $u(x)$ be a harmonic function in the
upper half space $H=\{x=(x',x_n)\in {\bf R}^{n}:\; x_n>0\}$ with
continuous boundary values on the boundary $\partial{H}$, $R>1$.
Then on the closed half ball
$\overline{B}_{R}^{+}=\overline{B}_{R}\cap H=\{x\in \overline{H}:\;
|x|\leq R\}$, we have
\begin{eqnarray*}
u(x)
& =& \int_{\{y\in H:\;|y|=R, y_n>0\}}\frac{R^2-|x|^2}{\omega_nR}
\bigg(\frac{1}{|y-x|^n}-\frac{1}{|y-x^{\ast}|^n}
\bigg)u(y)d\sigma(y) \\
& & +\frac{2x_n}{\omega_n}\int_{\{y\in
\overline{H}:\;|y'|<R,y_n=0\}}
 \bigg(\frac{1}{|y'-x|^n}-\frac{R^n}{|x|^n}
 \frac{1}{|y'-\widetilde{x}|^n} \bigg)u(y')dy',
\end{eqnarray*}
where $\widetilde{x}=R^2x/|x|^2,\  x^{\ast}=(x',-x_n)$, and
$\omega_{n}=\frac{\pi^{\frac{n}{2}}}{\Gamma(1+\frac{n}{2})}$ is the
volume of the unit $n$-ball in ${\bf R}^{n} $. }

\vspace{0.2cm}
 \noindent
 {\bf Remark 6.2.1 } {\it Lemma 6.2.1 is the generalization of the Carleman formula
 for harmonic functions in the upper half plane to the upper half
space; Lemma 6.2.2 is the generalization of the Nevanlinna formula
 for harmonic functions in the upper half disk to the upper half
ball. }

\section{ Proof of Theorem}

  We use Lemma 6.2.1 to the harmonic function $u(x)$,
\begin{eqnarray*}
& & \int_{\{x\in {\bf
R}^{n}:\;|x|=R,x_n>0\}}u^{-}(x)\frac{nx_n}{R^{n+1}}d\sigma(x) \\
& & +\int_{\{x\in {\bf R}^{n}:\;1<|x'|<R,x_n=0\}}u^{-}(x')
\bigg(\frac{1}{|x'|^{n}}-\frac{1}{R^{n}} \bigg)dx' \\
& =& \int_{\{x\in {\bf
R}^{n}:\;|x|=R,x_n>0\}}u^{+}(x)\frac{nx_n}{R^{n+1}}d\sigma(x)\\
& &+\int_{\{x\in {\bf R}^{n}:\;1<|x'|<R,x_n=0\}}u^{+}(x')
\bigg(\frac{1}{|x'|^{n}}-\frac{1}{R^{n}} \bigg)dx'+
c_1+\frac{c_2}{R^{n}},\hspace{16mm} (6.3.1)
\end{eqnarray*}
where $u^{+}(x)=\max\{u(x), 0\}, u^{-}(x)=(-u(x))^{+}$ and
$u(x)=u^{+}(x)-u^{-}(x)$.

  The terms on the right-hand of (6.3.1) can be estimated by using (6.1.1):
\begin{eqnarray*}
& & \int_{\{x\in {\bf R}^{n}:\;|x|=R,x_n>0\}}
u^{+}(x)\frac{nx_n}{R^{n+1}}d\sigma(x) \leq AKR^{\rho(R)-1},
\hspace{34mm} (6.3.2)\\
& &\int_{\{x\in {\bf R}^{n}:\;1<|x'|<R,x_n=0\}} u^{+}(x')
\bigg(\frac{1}{|x'|^{n}}-\frac{1}{R^{n}} \bigg)dx' \leq
AKR^{\rho(R)-1}. \hspace{19mm} (6.3.3)
\end{eqnarray*}

Thus, for $R>1$, we can obtain by (6.3.1), (6.3.2) and (6.3.3)
\begin{eqnarray*}
& & \int_{\{x\in {\bf R}^{n}:\;|x|=R,x_n>0\}}
u^{-}(x)\frac{nx_n}{R^{n+1}}d\sigma(x) \leq AKR^{\rho(R)-1},
\hspace{34mm} (6.3.4)\\
& &\int_{\{x\in {\bf R}^{n}:\;1<|x'|<R,x_n=0\}} u^{-}(x')
\bigg(\frac{1}{|x'|^{n}}-\frac{1}{R^{n}} \bigg)dx' \leq
AKR^{\rho(R)-1}.  \hspace{19mm} (6.3.5)
\end{eqnarray*}

Note that
\begin{eqnarray*}
& & \int_{\{x\in {\bf R}^{n}:\;1<|x'|<R,x_n=0\}}
\frac{u^{-}(x')}{|x'|^{n}}dx' \\
&\leq& \frac{2^n}{2^n-1}\int_{\{x\in {\bf R}^{n}:\;1<|x'|<R,x_n=0\}}
u^{-}(x') \bigg(\frac{1}{|x'|^{n}}-\frac{1}{(2R)^{n}} \bigg)dx' \\
&\leq& AK(2R)^{\rho(R)-1}.\hspace{88mm} (6.3.6)
\end{eqnarray*}

  We use Lemma 6.2.2 to the harmonic function $-u(x)$, and note that $-u(x)\leq
  u^{-}(x)$, we have
\begin{eqnarray*}
-u(x)
&=& \int_{\{y\in H:\;|y|=R, y_n>0\}}\frac{R^2-|x|^2}{\omega_nR}
\bigg(\frac{1}{|y-x|^n}-\frac{1}{|y-x^{\ast}|^n} \bigg)(-u(y))d\sigma(y) \\
& & +\frac{2x_n}{\omega_n}\int_{\{y\in
\overline{H}:\;|y'|<R,y_n=0\}}
 \bigg(\frac{1}{|y'-x|^n}
-\frac{R^n}{|x|^n}\frac{1}{|y'-\widetilde{x}|^n} \bigg)(-u(y'))dy' \\
&\leq &  \int_{\{y\in H:\;|y|=R, y_n>0\}}\frac{R^2-|x|^2}{\omega_nR}
\bigg(\frac{1}{|y-x|^n}-\frac{1}{|y-x^{\ast}|^n} \bigg)u^{-}(y)d\sigma(y) \\
& & +\frac{2x_n}{\omega_n}\int_{\{y\in
\overline{H}:\;|y'|<R,y_n=0\}}
 \bigg(\frac{1}{|y'-x|^n}-\frac{R^n}{|x|^n}\frac{1}{|y'-\widetilde{x}|^n} \bigg)u^{-}(y')dy'\\
&=&I_1+I_2 .\hspace{91mm} (6.3.7)
\end{eqnarray*}

  Note that the following estimates:
$$
\frac{1}{|y-x|^n}-\frac{1}{|y-x^{\ast}|^n}\leq
\frac{2nx_ny_n}{\omega_n|y-x|^{n+2}}, \eqno{(6.3.8)}
$$
$$
|y-x|^n\leq x_n^n=|x|^n\sin^n\theta, \quad x\in H,
y_n=0.\eqno{(6.3.9)}
$$
Put $|x|=r>1/2,\ R=2r$ in (6.3.7), then by (6.3.4), (6.3.8)and
(6.3.9), we have
\begin{eqnarray*}
I_1
&\leq &\int_{\{y\in H:\;|y|=R,
y_n>0\}}\frac{R^2-r^2}{\omega_nR}\frac{2nx_ny_n}{\omega_n|y-x|^{n+2}}u^{-}(y)d\sigma(y) \\
&\leq& AKR^{\rho(R)}  \hspace{94mm} (6.3.10)
\end{eqnarray*}
and
\begin{eqnarray*}
I_2
&\leq & \frac{2x_n}{\omega_n}\int_{\{y\in
\overline{H}:\;|y'|<R,y_n=0\}} \frac{1}{x_n^n}u^{-}(y')dy' \\
&= & \frac{2}{\omega_nx_n^{n-1}}\int_{\{y\in
\overline{H}:\;|y'|<R,y_n=0\}} u^{-}(y')dy' \\
&= & \frac{2}{\omega_nx_n^{n-1}}\int_{\{y\in
\overline{H}:\;1<|y'|<R,y_n=0\}} u^{-}(y')dy'+
\frac{2}{\omega_nx_n^{n-1}}\int_{\{y\in \overline{H}:\;|y'|\leq
1,y_n=0\}} u^{-}(y')dy' \\
&=& I_{21}+I_{22}.\hspace{95mm} (6.3.11)
\end{eqnarray*}
For the first integral we have by (6.3.6)
\begin{eqnarray*}
I_{21}
&\leq &\frac{2R^n}{\omega_nx_n^{n-1}}\int_{\{y\in
\overline{H}:\;1<|y'|<R,y_n=0\}} \frac{u^{-}(y')}{|y'|^n}dy' \\
&\leq& AK\frac{(2R)^{\rho(R)}}{\sin^{n-1}\theta},\hspace{86mm}
(6.3.12)
\end{eqnarray*}
for the second integral we have by (6.1.2)
\begin{eqnarray*}
I_{22}
&\leq &\frac{2K}{\omega_nx_n^{n-1}}\int_{\{y\in
\overline{H}:\;1<|y'|<R,y_n=0\}} dy' \\
&\leq& AK\frac{1}{\sin^{n-1}\theta}.\hspace{87mm} (6.3.13)
\end{eqnarray*}
By collecting (6.3.7), (6.3.10), (6.3.11), (6.3.12) and (6.3.13), we
have for $|x|>1/2$,
$$
-u(x)\leq AK\frac{1+(2R)^{\rho(R)}}{\sin^{n-1}\theta},
\eqno{(6.3.14)}
$$
for $|x\leq1/2$, we can get by (6.1.2)
$$
-u(x)\leq K\leq K\frac{1+(2R)^{\rho(R)}}{\sin^{n-1}\theta},
\eqno{(6.3.15)}
$$
so we obtain by (6.3.14) and (6.3.15)
$$
u(x)\geq -cK\frac{1+(2R)^{\rho(R)}}{\sin^{n-1}\theta}, \quad x\in H.
$$

\vspace{0.2cm} \noindent
 {\bf Remark 6.3.1 } {\it By modifying (6.3.6):
\begin{eqnarray*}
& & \int_{\{x\in {\bf R}^{n}:\;1<|x'|<R,x_n=0\}}
\frac{u^{-}(x')}{|x'|^{n}}dx' \\
&\leq & \frac{(N+1)^n}{(N+1)^n-N^n}\int_{\{x\in {\bf
R}^{n}:\;1<|x'|<R,x_n=0\}} u^{-}(x')
\bigg(\frac{1}{|x'|^{n}}-\frac{1}{(\frac{N+1}{N}R)^{n}} \bigg)dx'\\
&\leq & AK \bigg(\frac{N+1}{N}R  \bigg)^{\rho(R)-1},
\end{eqnarray*}

we can get:
$$
u(x)\geq -cK\frac{1+ \bigg(\frac{N+1}{N}R
 \bigg)^{\rho(R)}}{\sin^{n-1}\theta}, \quad x\in H.
$$}
{\bf Remark 6.3.2 } {\it A example: suppose $u(z)=\Re
e^{-iz}=e^y\cos x$ is a harmonic function in the upper half plane
${\bf C}_{+}$ with continuous boundary values on the real axis,
write $|x|=R\cos\theta, y=R\sin\theta (0<\theta\leq \pi/2)$. Let
$K=1, \rho(R)=\frac{R}{\log R}$ , then $u(z)$ satisfies
$$
u(z)\leq e^R \leq KR^{\rho(R)}.
$$
Thus
$$
u(z)\geq -e^R \geq-cK\frac{1+(2R)^{\rho(R)}}{\sin^{n-1}\theta},
\quad x\in H.
$$
}

\chapter{the Carleman Formula of Subharmonic Functions in the Half Space}

\section{Introduction and Main Theorem}

  B.Ya.Levin \cite{LLST} has proved the
following result:

\vspace{0.2cm}
 \noindent
{\bf Theorem G } (Carleman's formula) {\it Let $f(z)$ be a
meromorphic function in a closed sector $\overline{S}=\{z:\; \rho
\leq |z| \leq R,\  \Im z\geq 0\}$ whose zeros and poles do not lie
on the boundary $\partial{S}$. Then we obtain
\begin{eqnarray*}
& & \sum_{\rho <|a_n|< R}
\bigg(\frac{1}{|a_n|}-\frac{|a_n|}{R^2}\bigg)\sin
\alpha_n-\sum_{\rho <|b_n|< R}
\bigg(\frac{1}{|b_n|}-\frac{|b_n|}{R^2}\bigg)\sin \beta_n \\
&=&
\frac{1}{2\pi}\int_{\rho}^{R}\bigg(\frac{1}{t^2}-\frac{1}{R^2}\bigg)\log
|f(t)f(-t)|dt+\frac{1}{\pi R}\int_0^{\pi}\log |f(R
e^{i\varphi})|\sin \varphi d\varphi-A_f(\rho,R),
\end{eqnarray*}
where $a_n=|a_n|e^{i\alpha_n}$, $b_n=|b_n|e^{i\beta_n}$ are zeros
and poles of the function $f(z)$, and the remainder term
$A_f(\rho,R)$ is expressed by
$$
A_f(\rho,R)=-\frac{1}{2\pi}\int_0^{\pi}\bigg[\bigg(\frac{1}{\rho^2}+\frac{1}{R^2}\bigg)\log
|f(\rho
e^{i\varphi})|-\bigg(\frac{1}{\rho}-\frac{\rho}{R^2}\bigg)\frac{\partial
}{\partial \rho}\log |f(\rho e^{i\varphi})|\bigg]\rho \sin \varphi
d\varphi.
$$ }

 The object of this chapter is to generalize the Carleman formula
 for meromorphic functions in the upper half plane to subharmonic functions
 in the upper half space. we derive the following main theorem.

\vspace{0.2cm}
 \noindent
{\bf Theorem 7.1.1 } {\it Let $u(x)$ be a subharmonic function in
the upper half space $H$ with continuous boundary values on the
boundary $\partial{H}$, for $R>r>0$, we have
\begin{eqnarray*}
& &\int_{\{x\in {\bf
R}^{n}:\;|x|=R,x_n>0\}}u(x)\frac{nx_n}{R^{n+1}}d\sigma(x) \\
&+& \int_{{\{x\in {\bf R}^{n}:\ r<|x'|<R,x_n=0\}}}u(x')
\bigg(\frac{1}{|x'|^n}-\frac{1}{R^n} \bigg)dx' \geq A_u(r,R),
\end{eqnarray*}
where
$$
A_u(r,R)=c_1(r)+\frac{c_2(r)}{R^n}
$$
is a function depending on $r$ and $R$ and $ c_1(r)$, $ c_2(r)$ are
functions depending only on $r$, they are denoted by
$$
c_1(r)=\int_{\{x\in {\bf R}^{n}:\;|x|=r,x_n>0\}}
\bigg[\frac{(n-1)x_n}{r^{n+1}}u(x)+\frac{x_n}{r^n}\frac{\partial
u(x)}{\partial n} \bigg]d\sigma(x),
$$
and
$$
c_2(r)=\int_{\{x\in {\bf R}^{n}:\;|x|=r,x_n>0\}}
\bigg[\frac{x_n}{r}u(x)-x_n\frac{\partial u(x)}{\partial n}
\bigg]d\sigma(x).
$$
}

\section{Main Lemma}

 In order to obtain the result, we
need the lemma below:

\vspace{0.2cm}
 \noindent
{\bf Lemma 7.2.1 } {\it Suppose that $D$ is an admissible domain
with boundary $S$ in ${\bf R}^{n}$. If $u,v \in C^2$ in $\overline
D$, then we have
$$
\int_S  \bigg[u(x)\frac{\partial v(x)}{\partial
n}-v(x)\frac{\partial u(x)}{\partial n} \bigg]d\sigma(x)
= \int_D
[v(x)\triangle u(x)-u(x)\triangle v(x)]dx.
$$
Here $\partial/\partial n$ denotes differentiation along the inward
normal into $D$. }

\vspace{0.2cm}
 \noindent
 {\bf Remark 7.2.1} {\it Lemma 7.2.1 is just called the second Green's formula.
 }

\section{Proof of Theorem}

  Apply the second Green's formula  to the subharmonic function $u(x)$
and $v(x)=\frac{x_n}{|x|^n}-\frac{x_n}{R^n}$ in the resulting sphere
$$
B_{r,R}^{+}={\{x\in {\bf R}^{n}:\\
;r<|x|<R,x_n>0\}},
$$ we obtain
$$
\int_{\partial B_{r,R}^{+}} \bigg(u(x)\frac{\partial v(x)}{\partial
n}-v(x)\frac{\partial u(x)}{\partial n} \bigg)d\sigma(x) =\int_{
B_{r,R}^{+}} v(x)\triangle u(x)dx\geq 0. \eqno{(7.3.1)}
$$
The function  $v(x)$ is harmonic in $H$, the equations
$$
v(x)=0, \quad \frac{\partial v(x)}{\partial n}=\frac{nx_n}{R^{n+1}}
\eqno{(7.3.2)}
$$
hold on the half sphere $\{x\in {\bf R}^{n}:\;|x|=R,x_n>0\}$.

While the equation
$$
\frac{\partial v(x)}{\partial n}=-\frac{x_n}{r}
\bigg(\frac{n-1}{r^n}+\frac{1}{R^n} \bigg) \eqno{(7.3.3)}
$$
holds on the half sphere $\{x\in {\bf R}^{n}:\;|x|=r,x_n>0\}$.

Moreover, the equations
$$
v(x)=0,\quad \frac{\partial v(x)}{\partial
n}=\frac{1}{|x|^n}-\frac{1}{R^n} \eqno{(7.3.4)}
$$
hold on ${\{x\in {\bf R}^{n}:\;r<|x'|<R,x_n=0\}}$.

Thus
\begin{eqnarray*}
0
&\leq& \int_{\partial B_{r,R}^{+}} \bigg[u(x)\frac{\partial
v(x)}{\partial
n}-v(x)\frac{\partial u(x)}{\partial n} \bigg]d\sigma(x) \\
&=& \int_{\{x\in {\bf R}^{n}:\;|x|=R,x_n>0\}}
\bigg[u(x)\frac{\partial
v(x)}{\partial n}-v(x)\frac{\partial u(x)}{\partial n} \bigg]d\sigma(x) \\
& & +\int_{\{x\in {\bf R}^{n}:\;|x|=r,x_n>0\}}
\bigg[u(x)\frac{\partial
v(x)}{\partial n}-v(x)\frac{\partial u(x)}{\partial n} \bigg]d\sigma(x) \\
& & +\int_{{\{x\in{\bf R}^{n}:\;r<|x'|<R,x_n=0\}}}
\bigg[u(x)\frac{\partial v(x)}
{\partial n}-v(x)\frac{\partial u(x)}{\partial n} \bigg]d\sigma(x) \\
&=&I_1+I_2+I_3 .\hspace{92mm}(7.3.5)
\end{eqnarray*}
For the first term we have by (7.3.2)
$$
I_1=\int_{\{x\in {\bf
R}^{n}:\;|x|=R,x_n>0\}}u(x)\frac{nx_n}{R^{n+1}}d\sigma(x);\eqno{(7.3.6)}
$$
for the second term we have by (7.3.3)
\begin{eqnarray*}
I_2
&= & \int_{\{x\in {\bf R}^{n}:\;|x|=r,x_n>0\}}
\bigg[-u(x)\frac{x_n}{r} \bigg(\frac{n-1}{r^n}+\frac{1}{R^n} \bigg)
- \bigg(\frac{x_n}{|x|^n} -\frac{x_n}{R^n} \bigg)\frac{\partial
u(x)}{\partial n} \bigg]d\sigma(x) \\
&=&-c_1(r)-\frac{c_2(r)}{R^n}; \hspace{84mm} (7.3.7)
\end{eqnarray*}
for the third term we have by (7.3.4)
$$
I_3=\int_{{\{x\in {\bf R}^{n}:\ r<|x'|<R,x_n=0\}}}u(x')
\bigg(\frac{1}{|x'|^n}-\frac{1}{R^n} \bigg)dx'.\eqno{(7.3.8)}
$$
By collecting (7.3.5), (7.3.6), (7.3.7) and (7.3.8), we have
\begin{eqnarray*}
& &\int_{\{x\in {\bf
R}^{n}:\;|x|=R,x_n>0\}}u(x)\frac{nx_n}{R^{n+1}}d\sigma(x) \\
&+& \int_{{\{x\in {\bf R}^{n}:\ r<|x'|<R,x_n=0\}}}u(x')
\bigg(\frac{1}{|x'|^n}-\frac{1}{R^n} \bigg)dx' \geq A_u(r,R).
\end{eqnarray*}
This completes the proof of Theorem.

\chapter{a Generalization of the Nevanlinna Formula for Analytic Functions in the Right Half Plane}

\section{Introduction and Main Theorem}

 { Recall that ${\bf C}$  denote the complex plane
with points $z=x+iy$, where $x, y \in {\bf R}$.  The boundary and
closure of an open  $\Omega$ of ${\bf C}$  are denoted by
$\partial{\Omega}$
 and $\overline{\Omega}$ respectively.
 The right half plane is the set
 ${\bf C}_+=\{z=x+iy \in {\bf C}:\; x>0\}$, whose boundary is
 $\partial{{\bf C}_+}$.
 We identify ${\bf C}$ with ${\bf R} \times {\bf R}$ and
${\bf R} $ with $ {\bf R} \times \{0\}$,
  with this convention we then have $ \partial {{\bf C}_+}={\bf R}$.

\vspace{0.2cm}
    Suppose $R>1$, We  write $B_{+}(0,R)=\{z:\; |z|<R,\Re z>0\}$
    for the open right half
  disk of radius $R$ in ${\bf C}$ centered at the origin, whose boundary is
$\partial B_{+}(0,R)=\{z:\;z=it,|t|\leq R\}\bigcup \{z:\;z=R
e^{i\theta},|\theta|\leq \frac{\pi}{2}\}$.

\vspace{0.3cm}
  Let $\rho>1$, if the function $f(x)$ is analytic in the open right
half plane ${\bf C}_+$, continuous in the closed right half plane
$\overline{{\bf C}_+}$, and satisfies the following conditions:
$$
 \int_{-\infty}^{\infty}\frac{ \log^{+}|F(it)|}{1+|t|^{\rho+1}}dt <\infty
\eqno{(8.1.1)}
$$
and
$$
\int\int_{{\bf C}_+}\frac{x\log^{+}|F(z)|}{1+|z|^{\rho +3}}dm(z)
<\infty , \eqno{(8.1.2)}
$$
then a number of results have been achieved in \cite{DFF},
\cite{DFO}, \cite{DO}, \cite{GD}, \cite{DZ}, in this chapter, we
replace the first condition (8.1.1) into
$$
\lim_{\overline{\varepsilon\rightarrow 0}}
\int_{-\infty}^{\infty}\frac{
\log^{+}|F(it+\varepsilon)|}{1+|t|^{\rho+1}}dt <\infty,
\eqno{(8.1.3)}
$$
and that the function $f(x)$ is continuous in the boundary $\partial
{\bf C}_+$ is not needed, we can get the silimar results as
\cite{YDR}.

 \vspace{0.2cm}
 \noindent
{\bf Theorem 8.1.1 } {\it  Suppose $R'>R>1$, $F \in
N^{+}(B_{+}(0,R))$, let $\Lambda_R$ is the set of zeros of $F$ in
$B_{+}(0,R)$ and $\Lambda$ is the set of zeros of $F$ in ${\bf C}_+$
(including repetitions for multiplicities). If the conditions
(8.1.2) and (8.1.3) are satisfied, then

\noindent {\rm (1)}\
$$
\lim_{\overline{\varepsilon\rightarrow 0}}
\int_{-\infty}^{\infty}\frac{
|\log|F(it+\varepsilon)||}{1+|t|^{\rho+1}}dt <\infty;
$$

\noindent {\rm (2)}\
$$
\lim_{\overline{R\rightarrow \infty}}
\frac{1}{R^{\rho}}\int_{-\pi/2}^{\pi/2}|\log|F(R
e^{i\theta})||\cos\theta d\theta =0;
$$

\noindent {\rm (3)}\
$$
\sum_{\lambda_n \in \bigwedge}\frac{\Re \lambda_n}{1+
|\lambda_n|^{\rho+1}}<\infty.
$$
}

\section{Proof of Theorem}

$\forall z \in B_{+}(0,R)$ and $z \notin \Lambda_R$, write
$F_{\varepsilon}(z)=F(z+\varepsilon)$, then (see \cite{LDT},
\cite{PD} and \cite{RR})
\begin{eqnarray*}
\log|F_\varepsilon (z)|
&\leq& \frac{1}{2\pi}\int_{-\pi/2}^{\pi/2} \frac{4Rx \cos \theta
(R^2-|z|^2)}{|R e^{i\theta}-z|^2 |R e^{-i\theta}+z|^2}
\log|F_{\varepsilon}(R e^{i\theta})|d\theta \\
& & +\frac{x}{\pi}\int_{-R}^{R}
\bigg(\frac{1}{|it-z|^2}-\frac{R^2}{|R^2+itz|^2}\bigg)
\log|F_{\varepsilon}(it)|dt. \hspace{16mm} (8.2.1)
\end{eqnarray*}

  Without loss of generality we may assume that $F(1)\neq 0$, then
there exists $ \varepsilon_0 >0$, such that for any
$0<\varepsilon<\varepsilon_0$, we have $F(1+\varepsilon)\neq 0$ and
$|F(1+\varepsilon)|> \frac{|F(1)|}{2}$.

  Suppose $R>2$, $z=1$, by (8.2.1), we have
\begin{eqnarray*}
& & \log\frac{|F(1)|}{2} +
\frac{2R(R^2-1)}{\pi}\int_{-\pi/2}^{\pi/2} \frac{ \cos \theta }{|R
e^{i\theta}-1|^2 |R e^{-i\theta}+1|^2} \log^{-}|F_{\varepsilon}(R
e^{i\theta})|d\theta \\
& & +\frac{1}{\pi}\int_{|t|\leq R/2}
\bigg(\frac{1}{t^2+1}-\frac{R^2}{|R^4+t^2}\bigg)
\log^{-}|F_{\varepsilon}(it)|dt \\
&\leq& \frac{2R(R^2-1)}{\pi}\int_{-\pi/2}^{\pi/2} \frac{ \cos \theta
}{|R e^{i\theta}-1|^2 |R e^{-i\theta}+1|^2}
\log^{+}|F_{\varepsilon}(R e^{i\theta})|d\theta \\
& & +\frac{1}{\pi}\int_{|t|\leq R}
\bigg(\frac{1}{t^2+1}-\frac{R^2}{R^4+t^2}\bigg)
\log^{+}|F_{\varepsilon}(it)|dt .
\end{eqnarray*}

Set
$$
m_{+}^{(\varepsilon)}(R)=\frac{1}{ R}\int_{-\pi/2}^{\pi/2}
\log^{+}|F_{\varepsilon}(R e^{i\theta})|\cos\theta d\theta,
$$
$$
m_{-}^{(\varepsilon)}(R)=\frac{1}{ R}\int_{-\pi/2}^{\pi/2}
\log^{-}|F_{\varepsilon}(R e^{i\theta})|\cos\theta d\theta,
$$
$$
g_{+}^{(\varepsilon)}(t)=\log^{+}|F_{\varepsilon}(it)|+\log^{+}|F_{\varepsilon}(-it)|,
$$
$$
g_{-}^{(\varepsilon)}(t)=\log^{-}|F_{\varepsilon}(it)|+\log^{-}|F_{\varepsilon}(-it)|.
$$

  Note that when $|t|\leq R/2$,
$$
\frac{1}{t^2+1}-\frac{R^2}{R^4+t^2} \geq
\frac{9}{32}\frac{1}{t^2+1};
$$
when $|t|\leq R$,
$$
\frac{1}{t^2+1}-\frac{R^2}{R^4+t^2} \leq \frac{1}{t^2+1}.
$$
So we obtain
\begin{eqnarray*}
& &\frac{8}{27\pi}m_{-}^{(\varepsilon)}(R)
   +\frac{9}{32\pi}\bigg[\int_1^{R/2}
   \frac{1}{2t^2}g_{-}^{(\varepsilon)}(t)dt+\frac{1}{2}\int_{|t|<1}
   \log^{-}|F_{\varepsilon}(it)|dt\bigg] \\
&\leq& \frac{32}{\pi}m_{+}^{(\varepsilon)}(R)
+\frac{1}{\pi}\bigg[\int_1^{R}
 \frac{1}{t^2}g_{+}^{(\varepsilon)}(t)dt+\int_{|t|<1}
\frac{1}{t^2+1}\log^{+}|F_{\varepsilon}(it)|dt\bigg]-
\log\frac{|F(1)|}{2}.\hspace{2mm} (8.2.2)
\end{eqnarray*}

  Multiplying (8.2.2) by $\frac{1}{R^{\rho}}$ and integrating with
respect to $R$, we obtain
\begin{eqnarray*}
& &\frac{8}{27\pi}\int_2^{\infty}
\frac{m_{-}^{(\varepsilon)}(R)}{R^{\rho}}dR
+\frac{9}{32\pi}\int_2^{\infty}\frac{1}{R^{\rho}}\int_1^{R/2}
\frac{1}{2t^2}g_{-}^{(\varepsilon)}(t)dtdR \\
& & +\frac{9}{64\pi}\int_{|t|<1}
\log^{-}|F_{\varepsilon}(it)|dt\cdot\int_2^{\infty}\frac{1}{R^{\rho}}dR \\
&\leq&
\frac{32}{\pi}\int_2^{\infty}\frac{m_{+}^{(\varepsilon)}(R)}{R^{\rho}}dR
+\frac{1}{\pi}\int_2^{\infty}\frac{1}{R^{\rho}}\int_1^{R}
 \frac{1}{t^2}g_{+}^{(\varepsilon)}(t)dtdR \\
& & +\bigg(\frac{1}{\pi}\int_{|t|<1}
\frac{1}{t^2+1}\log^{+}|F_{\varepsilon}(it)|dt-
\log\frac{|F(1)|}{2}\bigg)\cdot\int_2^{\infty}\frac{1}{R^{\rho}}dR.
\end{eqnarray*}
After some elementary calculations, we get
\begin{eqnarray*}
& &\frac{8}{27\pi}\int\int_D
\frac{x\log^{-}|F_{\varepsilon}(z)|}{|z|^{\rho +3}}dm(z)
+\frac{9}{64\pi}\frac{1}{2^{\rho-1}(\rho-1)} \\
& & \times \bigg[\int_1^{\infty}
 \frac{g_{-}^{(\varepsilon)}(t)}{t^{\rho+1}}dt+\int_{|t|<1}
\log^{-}|F_{\varepsilon}(it)|dt\bigg] \\
&\leq& \frac{32}{\pi}\int\int_D
\frac{x\log^{+}|F_{\varepsilon}(z)|}{|z|^{\rho +3}}dm(z)
+\frac{1}{\pi}\frac{1}{\rho-1}\int_1^{\infty}
 \frac{g_{+}^{(\varepsilon)}(t)}{t^{\rho+1}}dt \\
& & +\frac{1}{\pi}\frac{1}{\rho-1}\int_{|t|<1}
\frac{1}{t^2+1}\log^{+}|F_{\varepsilon}(it)|dt-
\frac{1}{2^{\rho-1}(\rho-1)}\log\frac{|F(1)|}{2},
\end{eqnarray*}
where $D=\{(x,y): \; x^2+y^2\geq 4,x\geq0\}$.

 By (8.1.3), there exists a sequence $\{\varepsilon_n\}$,
$\varepsilon_n \rightarrow 0$ and $M>0$, such that
$$
\int_{-\infty}^{\infty}\frac{
\log^{+}|F(it+\varepsilon_n)|}{1+|t|^{\rho+1}}dt \leq M<\infty,
$$
so
\begin{eqnarray*}
& &\int_1^{\infty}
 \frac{g_{+}^{(\varepsilon_n)}(t)}{t^{\rho+1}}dt+\int_{|t|<1}
\frac{1}{t^2+1}\log^{+}|F_{\varepsilon_n}(it)|dt \\
&\leq&  2\int_{-\infty}^{\infty}\frac{
\log^{+}|F(it+\varepsilon_n)|}{1+|t|^{\rho+1}}dt \leq 2M<\infty.
\end{eqnarray*}
When $\varepsilon<\frac{3}{4}|\omega|$, limit
$\varepsilon<\frac{6}{7}$, then we have
$$
\int\int_D \frac{x\log^{+}|F_{\varepsilon}(z)|}{|z|^{\rho +3}}dm(z)
\leq 2\times 4^{\rho +3}\int\int_{{\bf C}_+} \frac{\Re
z\log^{+}|F(z)|}{1+|z|^{\rho +3}}dm(z)\leq M,
$$
where $D'=D+\varepsilon$, $\omega \in D'$ and
$\omega=z+\varepsilon$, so we obtain
$$
\int_1^{\infty}
 \frac{g_{+}^{(\varepsilon_n)}(t)+g_{-}^{(\varepsilon_n)}(t)}{1+t^{\rho+1}}dt\leq
 M,
$$
$$
\int\int_D \frac{x|\log^{+}|F_{\varepsilon}(z)||}{|z|^{\rho
+3}}dm(z)\leq M,
$$
$$
\int_{|t|<1} |\log|F_{\varepsilon_n}(it)||dt \leq M,
$$
and
$$
\int_{-\infty}^{\infty}\frac{|\log|F_{\varepsilon_n}
(it)||}{1+|t|^{\rho+1}}dt \leq M.
$$
Therefore, there exists $ M>0$, such that
$$
\sup_n\int_{-\infty}^{\infty}\frac{|\log|F
(it+\varepsilon_n)||}{1+|t|^{\rho+1}}dt \leq M <\infty.
\eqno{(8.2.3)}
$$
$\forall s(t) \in C_0 (-\infty, +\infty)$, set
$$
T_n (s)=\int_{-\infty}^{\infty}s(t) \frac{\log|F
(it+\varepsilon_n)|}{1+|t|^{\rho+1}}dt,
$$
by (8.2.3), we obtain that $T_n$ is a bounded linear functional in
$C_0 (-\infty, +\infty)$ and
$$
\sup_n\|T_n\|=\sup_n\int_{-\infty}^{\infty} \frac{|\log|F
(it+\varepsilon_n)||}{1+|t|^{\rho+1}}dt\leq M.
$$
Hence, there exists a subsequence $\{T_{n_k}\}$ of $\{T_n\}$, such
that $T_{n_k}$ ${\rm weakly^*}$ converges to $T$, that is to say,
$$
T(s)=\lim_{n_k \rightarrow \infty}T_{n_k} (s),\quad \forall s(t) \in
C_0 (-\infty, +\infty).
$$
By Riesz Representation Theorem \cite{R}, there exists a Radon
measure $\nu$ such that
$$
T(s)=\int_{-\infty}^{\infty}s(t)d\nu.
$$
Set $d\nu=\frac{1}{1+|t|^{\rho+1}}d\mu$, then
$$
T(s)=\int_{-\infty}^{\infty}\frac{s(t)}{1+|t|^{\rho+1}}d\mu,\quad
\forall s(t) \in C_0 (-\infty, +\infty).
$$
Define
$$
T_{\varepsilon} (s)=\int_{-\infty}^{\infty}s(t) \frac{\log|F
(it+\varepsilon)|}{1+|t|^{\rho+1}}dt,
$$
then
$$
\|T_{\varepsilon}\|=\int_{-\infty}^{\infty} \frac{|\log|F
(it+\varepsilon)||}{1+|t|^{\rho+1}}dt,
$$
and
$$
\lim_{\overline{\varepsilon\rightarrow 0}} \|T_{\varepsilon}\| \leq
M,
$$
so
$$
\lim_{\overline{\varepsilon\rightarrow 0}}\int_{-\infty}^{\infty}
\frac{|\log|F (it+\varepsilon)||}{1+|t|^{\rho+1}}dt\leq M <\infty.
$$
Hence (1) holds;

Since
\begin{eqnarray*}
& & \int_2^{\infty}\frac{1}{R^{\rho+1}} \int_{-\pi/2}^{\pi/2}
|\log|F_{\varepsilon_n}(R e^{i\theta})||\cos\theta d\theta dR \\
&=& \int\int_D \frac{x\log^{+}|F_{\varepsilon_n}(z)|}{|z|^{\rho
+3}}dm(z)\leq M,
\end{eqnarray*}
we obtain
\begin{eqnarray*}
& & \int_2^{\infty}\frac{1}{R^{\rho+1}} \int_{-\pi/2}^{\pi/2}
|\log|F(R e^{i\theta})||\cos\theta d\theta dR \\
&\leq&  \lim_{\overline{n\rightarrow
\infty}}\int_2^{\infty}\frac{1}{R^{\rho+1}} \int_{-\pi/2}^{\pi/2}
|\log|F_{\varepsilon_n}(R e^{i\theta})||\cos\theta d\theta dR \leq
M,
\end{eqnarray*}
so
$$
\lim_{\overline{R\rightarrow \infty}}\frac{1}{R^{\rho}}
\int_{-\pi/2}^{\pi/2} |\log|F(R e^{i\theta})||\cos\theta d\theta =0.
$$
Thus (2) holds;

 Write $\lambda_{\varepsilon}=\lambda_n -\varepsilon$,
$\forall z \in B_{+}(0,R)$ and $z \notin \Lambda_R$, then
$$
F_{\varepsilon}(\lambda_{\varepsilon})
=F_{\varepsilon}(\lambda_n-\varepsilon) =F(\lambda_n) =0,
$$
where $\lambda_n$, $\lambda_{\varepsilon}$ are denoted the zeros of
$F$, $F_{\varepsilon}$ respectively. So we have
\begin{eqnarray*}
\log|F_\varepsilon (z)|
&=& \frac{1}{2\pi}\int_{-\pi/2}^{\pi/2} \bigg(\frac{R^2-|z|^2}{|R
e^{i\theta}-z|^2}-\frac{R^2-|z|^2}{|R e^{-i\theta}+z|^2}\bigg)
\log|F_{\varepsilon}(R e^{i\theta})|d\theta \\
& &  +\frac{1}{\pi}\int_{-R}^{R} \bigg(\frac{\Re
z}{|it-z|^2}-\frac{R^2 \Re z}{|R^2+itz|^2}\bigg)
\log|F_{\varepsilon}(it)|dt \\
& & + \sum_{\lambda_{\varepsilon} \in \bigwedge_R}
\log\bigg|\frac{z-\lambda_{\varepsilon}}{R^2-\overline{\lambda_{\varepsilon}}z}
\frac{R^2+\lambda_{\varepsilon}z}{z+\overline{\lambda_{\varepsilon}}}\bigg|.
\hspace{44mm} (8.2.4)
\end{eqnarray*}

 Without loss of generality we may assume that $F(1)\neq 0$, then
there exists $ \varepsilon_0 >0$, such that for any
$0<\varepsilon<\varepsilon_0$, we have $F(1+\varepsilon)\neq 0$ and
$|F(1+\varepsilon)|> \frac{|F(1)|}{2}$.

  Suppose $R>2$, $z=1$, by (8.2.4), we have
\begin{eqnarray*}
\log\frac{|F(1)|}{2}
&\leq & \frac{2R(R^2-1)}{\pi (R-1)^4}\int_{-\pi/2}^{\pi/2}
 \log^{+}|F_{\varepsilon}(R e^{i\theta})|\cos \theta d\theta \\
& & +\frac{1}{\pi}\int_{1}^{R} \frac{1}{t^2}[
\log^{+}|F_{\varepsilon}(it)|+\log^{+}|F_{\varepsilon}(it)|]dt \\
& & +\frac{1}{\pi}\int_{|t|<1} \frac{1}{1+t^2}
\log^{+}|F_{\varepsilon}(it)|dt\\
& & + \sum_{\lambda_{\varepsilon} \in \bigwedge_R}
\log\bigg|\frac{1-\lambda_{\varepsilon}}{R^2-\overline{\lambda_{\varepsilon}}}
\frac{R^2+\lambda_{\varepsilon}}{1+\overline{\lambda_{\varepsilon}}}\bigg|.
\end{eqnarray*}

Note that
$$
\log x <\frac{1}{2}(x^2-1), \quad \forall x \in (0, 1),
$$
then we have
\begin{eqnarray*}
& &\log\bigg|\frac{1-\lambda_{\varepsilon}}
{R^2-\overline{\lambda_{\varepsilon}}}
\frac{R^2+\lambda_{\varepsilon}}{1+\overline{\lambda_{\varepsilon}}}\bigg|\\
&\leq
&\frac{|(R^2-\lambda_{\varepsilon}\overline{\lambda_{\varepsilon}})
+(\overline{\lambda_{\varepsilon}}-R^2
\lambda_{\varepsilon})|^2}{2|1+\lambda_{\varepsilon}|^2
|R^2-\overline{\lambda_{\varepsilon}}|^2} \\
& &-\frac{|(R^2-\lambda_{\varepsilon}
\overline{\lambda_{\varepsilon}})
-(\overline{\lambda_{\varepsilon}}-R^2
\lambda_{\varepsilon})|^2}{2|1+\lambda_{\varepsilon}|^2
|R^2-\overline{\lambda_{\varepsilon}}|^2} \\
& =&-\frac{2(R^2-|\lambda_{\varepsilon}|^2)(R^2-1)\Re
\lambda_{\varepsilon}}{|1+\lambda_{\varepsilon}|^2
|R^2-\overline{\lambda_{\varepsilon}}|^2}.
\end{eqnarray*}
Since
$$
\sum_{1\leq |\lambda_{\varepsilon}|<R/2 } \frac{\Re
\lambda_{\varepsilon}}{|\lambda_{\varepsilon}|^2} =\int_1^{R/2}
\frac{1}{t} dN_0^{(\varepsilon)}(t),
$$
where
$$
N_0^{(\varepsilon)}(t)=\sum_{1\leq |\lambda_{\varepsilon}|<t
}\cos\theta_\varepsilon,
$$
then we obtain
\begin{eqnarray*}
& & \sum_{|\lambda_{\varepsilon}|<1 } \Re \lambda_{\varepsilon}
+\int_1^{R/2} \frac{1}{t} dN_0^{(\varepsilon)}(t)\\
&\leq & \frac{256}{\pi}m_{+}^{(\varepsilon)}(R)
+\frac{8}{\pi}\int_1^{R}
 \frac{1}{t^2}g_{+}^{(\varepsilon)}(t)dt \\
& & +\frac{8}{\pi}\int_{|t|<1}
\frac{1}{t^2+1}\log^{+}|F_{\varepsilon}(it)|dt-
8\log\frac{|F(1)|}{2}.
\end{eqnarray*}

Multiplying this by $\frac{1}{R^{\rho}}$ and integrating with
respect to $R$, we obtain
\begin{eqnarray*}
& & \sum_{|\lambda_{\varepsilon}|<1 } \Re \lambda_{\varepsilon}
\cdot\int_2^{\infty}\frac{1}{R^{\rho}}dR
+\int_2^{\infty}\frac{1}{R^{\rho}}\int_1^{R/2} \frac{1}{t}
dN_0^{(\varepsilon)}(t)dR \\
&\leq
&\frac{256}{\pi}\int_2^{\infty}\frac{m_{+}^{(\varepsilon)}(R)}{R^{\rho}}dR
+\frac{8}{\pi}\int_1^{\infty}\frac{1}{R^{\rho}}\int_1^{R}
 \frac{1}{t^2}g_{+}^{(\varepsilon)}(t)dtdR\\
& &+\int_2^{\infty}\frac{1}{R^{\rho}}
\cdot8\bigg[\frac{1}{\pi}\int_{|t|<1}
\frac{1}{t^2+1}\log^{+}|F_{\varepsilon}(it)|dt-
\log\frac{|F(1)|}{2}\bigg].
\end{eqnarray*}
By some elementary calculations, we get
\begin{eqnarray*}
&  & \sum_{\lambda_{\varepsilon} \in \bigwedge}\frac{\Re
\lambda_{\varepsilon}}{1+ |\lambda_{\varepsilon}|^{\rho+1}}\\
&\leq& \frac{2^{\rho+7}(\rho-1)}{\pi}
\int_2^{\infty}\frac{m_{+}^{(\varepsilon)}(R)}{R^{\rho}}dR \\
& & +\frac{2^{\rho+3}}{\pi}\int_{-\infty}^{\infty} \frac{|\log|F
(it+\varepsilon)||}{1+|t|^{\rho+1}}dt\\
& &- 8\log\frac{|F(1)|}{2}.
\end{eqnarray*}
So we have
\begin{eqnarray*}
&  & \lim_{\overline{\varepsilon\rightarrow
0}}\sum_{\lambda_{\varepsilon} \in \bigwedge}\frac{\Re
\lambda_{\varepsilon}}{1+ |\lambda_{\varepsilon}|^{\rho+1}} \\
&\leq& \frac{2^{\rho+7}(\rho-1)}{\pi}
\overline{\lim_{\varepsilon\rightarrow
0}}\int_2^{\infty}\frac{m_{+}^{(\varepsilon)}(R)}{R^{\rho}}dR \\
& & +\frac{2^{\rho+3}}{\pi}\lim_{\overline{\varepsilon\rightarrow
0}}\int_{-\infty}^{\infty} \frac{|\log|F
(it+\varepsilon)||}{1+|t|^{\rho+1}}dt \\
& & - 8\log\frac{|F(1)|}{2}.
\end{eqnarray*}
Hence
$$
\lim_{\overline{\varepsilon\rightarrow
0}}\sum_{\lambda_{\varepsilon} \in \bigwedge}\frac{\Re
\lambda_{\varepsilon}}{1+ |\lambda_{\varepsilon}|^{\rho+1}} \leq
M<\infty.
$$
Since
$$
\sum_{\lambda_n \in \bigwedge}\frac{\Re \lambda_n}{1+
|\lambda_n|^{\rho+1}} \leq \lim_{\overline{\varepsilon\rightarrow
0}}\sum_{\lambda_{\varepsilon} \in \bigwedge}\frac{\Re
\lambda_{\varepsilon}}{1+ |\lambda_{\varepsilon}|^{\rho+1}} \leq
M<\infty,
$$
thus (3) holds. This completes the proof of Theorem.

\newpage
\ \thispagestyle{empty}
\newpage

\chapter{Integral Representations of Harmonic Functions in the Half Plane}

\section{Introduction and Main Theorem}

\vspace{0.3cm}
  Let $\rho(R)\geq1$ is nondecreasing in $[0, +\infty)$ satisfying
$$
\varepsilon_0=\limsup_{R\rightarrow\infty}\frac{\rho'(R)R \log
R}{\rho(R)}<1. \eqno{(9.1.1)}
$$
For any real number  $ \alpha > 0 $, we denote by $(LU)_{\alpha }$
the space of all measurable  functions $f(x+iy)$ in the  upper half
plane ${\bf C}_+$ which satisfy the following inequality:
$$
\int\int_{{\bf
C}_+}\frac{y|f(x+iy)|}{1+(x^2+y^2)^\frac{{\rho(|z|)+\alpha
+3}}{2}}dxdy <\infty ;\eqno{(9.1.2)}
$$
and $(LV)_{\alpha} $ the set of all measurable functions $g(x) $ in
${\bf R} $ which satisfy the following inequality:
$$
 \int_{-\infty}^{\infty}\frac{ |g(x)|}{1+|x|^{\rho(|x|)+\alpha +1}}dx <\infty
.\eqno{(9.1.3)}
$$
We also denote by $(CH)_{\alpha }$ the set  of all continuous
functions $u(x+iy)$ in the closed upper half plane $\overline{{\bf
C}_+}$, harmonic in the open upper half plane ${\bf C}_+$ with the
positive part $u^{+}(x+iy)=\max\{u(x+iy),0\}\in (LU)_{\alpha}$ and $
u^+(x) \in (LV)_{\alpha }$.

 The Poisson kernel for the upper
half plane ${\bf C}_+$ is the
 function
$$
 P(z,t)=
\frac{y}{\pi|z-t|^2},
$$
where $z \in {\bf C}_+$, $t\in {\bf R}$.

If $u(z)\leq 0  $ is harmonic in the open upper half plane ${\bf
C}_+$, continuous in the closed upper half plane $\overline{{\bf
C}_+}$, then (see \cite{HN}, \cite{SH} and \cite{ABR}) $ u\in
(CH)_{\alpha}$ for each $\alpha > 0 $ and there exists a  constant
$c\leq 0 $ such that
$$
u(z)=cy+  \int_{-\infty}^{\infty}P(z,t)u(t)dt \eqno{(9.1.4)}
 $$
for all $ z \in  {\bf C}_+ $, the integral in (9.1.4) is absolutely
convergent. Motivated by this result, we will prove that if $u\in
(CH)_{\alpha }$, then $ u^{+}(x+iy) \in (LU)_{\alpha }, u(x)\in
(LV)_{\alpha} $ and  a similar representation to (9.1.4)  for the
function $ u \in (CH)_{\alpha }$ holds by modifying the Poisson
kernel $
 P_m(z,t)$.  It is well known (see
\cite{G}, \cite{HN} and \cite{SW}) that the Poisson kernel $
 P(z,t)$
is harmonic in $ z\in {\bf C}-\{t\}$ and has a series expansion:
$$
 P(z,t)=\frac{1}{\pi}\Im \sum_{k=0}^{\infty}\frac{z^k}{ t^{k+1}},
 $$
this series converges for $ |z|<|t|$. So if $m\geq 0$ is an integer,
we define a modified
 Cauchy kernel of order $m$ for $z\in {\bf C}_+$ by
$$
 C_m(z,t)=\left\{\begin{array}{ll}
 \frac{1}{\pi}\frac{1}{t-z} , &   \mbox{when }   |t|\leq 1  ,\\
 \frac{1}{\pi}\frac{1}{t-z} - \frac{1}{\pi}\sum_{k=0}^{m}\frac{z^k}{ t^{k+1}},&
\mbox{when}\   |t|> 1,
 \end{array}\right.\eqno{(9.1.5)}
$$

then we define a modified
 Poisson kernel of order $m$ for the upper half plane  by
$$
P_m(z,t)=\Im C_m(z,t).
$$

 That is to say,
 $$
 P_m(z,t)=\left\{\begin{array}{ll}
 P(z,t) , &   \mbox{when }   |t|\leq 1  ,\\
 P(z,t) - \frac{1}{\pi}\Im\sum_{k=0}^{m}\frac{z^k}{t^{k+1}},&
\mbox{when}\   |t|> 1.
 \end{array}\right.
  $$
The modified
 Poisson kernel $
 P_m(z,t)$ is harmonic in $ z\in {\bf C}_+ $.

  Up to now, a number of results about integral representations
have been achieved in \cite{YDI}, \cite{HD}, \cite{YDR}, \cite{DIH},
\cite{DIAW}, \cite{DZI}, \cite{DIAH}, \cite{ZDI}, \cite{DIAO},
\cite{ZYI}, in this chapter, we will establish the following
theorem.

 \vspace{0.2cm}
 \noindent
{\bf Theorem 9.1.1 } {\it  If  $ u\in (CH)_{\alpha}\ (\alpha > 0) $
, then the following properties hold:

\noindent {\rm (1)}\
$$
 \int_{-\infty}^{\infty}\frac{ |u(x)|}{1+|x|^{\rho(|x|)+\alpha +1}}dx <\infty
;
$$

 \noindent
 {\rm (2)}\   the integral
$$
\int_{-\infty}^{\infty}P_{[\rho(|t|)+\alpha]}(z,t)u(t)dt
 $$
is absolutely convergent, it represents a harmonic function $u_{{\bf
C}_+}(z)$ in ${\bf C}_+$ and can be continuously extended to
$\overline{{\bf C}_+} $ such that $u_{{\bf C}_+}(t)=u(t)$ for $ t\in
{\bf R}$;

\noindent
 {\rm (3)}\   There exists an entire function $Q_{{\bf C}_+}(z)$  which
  satisfies on the boundary ${\bf R}$ that $\Im Q_{{\bf C}_+}(x)=0$
  such that $u(z)=\Im Q_{{\bf C}_+}(z)+u_{{\bf
C}_+}(z)$ for all $ z \in \overline{{\bf C}_+} $. }

\section{Main Lemma}

In order to obtain the result, we need the lemma below:

\vspace{0.2cm}
 \noindent
{\bf Lemma 9.2.1 }  {\it For any $ t\in {\bf R}$ and $|z|>1, y>0 $,
the
 following inequalities
$$
 |C_m(z,t)|\leq \left\{\begin{array}{ll}
 \frac{|z|^{m+1}}
 {\pi y|t|^{m+1}} , &   \mbox{when }
 1<|t|\leq 2|z| ,  \\
\frac{2|z|^{m+1}}{\pi|t|^{m+2}}, & \mbox{when }
 |t|> \max \{1,2|z|\},\\
\frac{1}{\pi y}, & \mbox{when }
 |t|\leq 1
 \end{array}\right.
$$
hold. }

\vspace{0.2cm}
 Proof:  When $ t\in
{\bf R}, |t|\leq 1 $, we have $ |t-z| \geq y$ and so
$$
|C_m(z,t)|\leq\frac{1}{\pi y};
$$
when $ t\in {\bf R}, 1<|t|\leq 2|z|$, we also have $ |t-z| \geq y $
and so by (9.1.5)
\begin{eqnarray*}
|C_m(z,t)|
&=& \frac{1}{\pi} \bigg|\frac{1}{t-z}-
\frac{1-(\frac{z}{t})^{m+1}}{t-z} \bigg| \\
&=&\frac{1}{\pi}\frac{|\frac{z}{t}|^{m+1}}{|t-z|} \leq
\frac{|z|^{m+1}}
 {\pi y|t|^{m+1}} ;
\end{eqnarray*}
when $|t|> \max \{1,2|z|\}$, we have by (9.1.5)
$$
|C_m(z,t)|=\frac{1}{\pi} \bigg|\sum_{k=m+1}^{\infty}\frac{z^k}{
t^{k+1}} \bigg|
 \leq \frac{1}{\pi}\sum_{k=m+1}^{\infty}\frac{|z|^k}{
|t|^{k+1}} \leq\frac{2|z|^{m+1}}{\pi|t|^{m+2}}.
$$
This proves the inequalities.

\section{Proof of Theorem}

   If  $ u\in (CH)_{\alpha}\ (\alpha > 0) $ , suppose $ R>1$, then by the
Carleman formula for harmonic functions in the upper half plane,
\begin{eqnarray*}
& & \frac{1}{\pi R}\int_0^{\pi}u(R e^{i\theta})\sin \theta d\theta \\
&+&\frac{1}{2\pi}\int_{1<|x|<R}u(x)
\bigg(\frac{1}{x^2}-\frac{1}{R^2} \bigg)dx =c_1+\frac{c_2}{R^2},
\end{eqnarray*}
where
$$
c_1=\frac{1}{2\pi}\int_0^{\pi} \bigg[u(R e^{i\theta})+\frac{\partial
u(R e^{i\theta})}{\partial n} \bigg]\sin \theta d\theta,
$$
$$
c_2=\frac{1}{2\pi}\int_0^{\pi} \bigg[u(R e^{i\theta})-\frac{\partial
u(R e^{i\theta})}{\partial n} \bigg]\sin \theta d\theta.
$$
Set
$$
m_{+}(R)=\frac{1}{\pi R}\int_0^{\pi} u^{+}(R e^{i\theta})\sin\theta
d\theta,
$$
$$
m_{-}(R)=\frac{1}{\pi R}\int_0^{\pi} u^{-}(r e^{i\theta})\sin\theta
d\theta,
$$
$$
g_{+}(x)=u^{+}(x)+u^{+}(-x),
$$
$$
g_{-}(x)=u^{-}(x)+u^{-}(-x),
$$
then
\begin{eqnarray*}
& & m_{-}(R)+\frac{1}{2\pi}\int_1^R
 \bigg(\frac{1}{x^2}-\frac{1}{R^2} \bigg)g_{-}(x)dx \\
&= &m_{+}(R)+\frac{1}{2\pi}\int_1^R
 \bigg(\frac{1}{x^2}-\frac{1}{R^2} \bigg)g_{+}(x)dx-
c_1-\frac{c_2}{R^2},\hspace{31mm} (9.3.1)
\end{eqnarray*}
where $u^{+}(z)=\max\{u(z), 0\}, u^{-}(z)=(-u(z))^{+}$ and
$u(z)=u^{+}(z)-u^{-}(z)$.

Since $u\in (CH)_{\alpha}$, by (9.1.2), we obtain
$$
\int_1^\infty\frac{m_{+}(R)}{R^{\rho(R)+\alpha }}dR =
 \frac{1}{\pi}\int\int_{D}\frac{yu^{+}(x+iy)}{(x^2+y^2)^\frac{{\rho(|z|)+\alpha +3}}{2}}dxdy
 <\infty, \eqno{(9.3.2)}
$$
where $D=\{z\in {\bf C}_+ :\;|z|>1\}$.

By (9.1.3), we can also obtain
\begin{eqnarray*}
& &\int_1^\infty\frac{1}{R^{\rho(R)+\alpha }}\int_1^R
 \bigg(\frac{1}{x^2}-\frac{1}{R^2} \bigg)g_{+}(x)dxdR \\
& =& \int_1^\infty g_{+}(x)\int_x^\infty\frac{1}{R^{\rho(R)+\alpha}}
\bigg(\frac{1}{x^2}-\frac{1}{R^2} \bigg)dRdx \\
&\leq&  \frac{2}{3}\int_1^{\infty}\frac{
 g_{+}(x)}{x^{\rho(x)+\alpha +1}}dx < \infty.
\hspace{73mm} (9.3.3)
\end{eqnarray*}
Similarly, we have
\begin{eqnarray*}
& & \int_1^\infty\frac{1}{R^{\rho(R)+\alpha/2 }}\int_1^R
 \bigg(\frac{1}{x^2}-\frac{1}{R^2} \bigg)g_{-}(x)dxdR \\
& =& \int_1^\infty g_{-}(x)\int_x^\infty\frac{1}{R^{\rho(R)+\alpha/2
}} \bigg(\frac{1}{x^2}-\frac{1}{R^2} \bigg)dRdx.\hspace{10mm}(9.3.4)
\end{eqnarray*}
So we have by (9.3.1), (9.3.2) and (9.3.3)
\begin{eqnarray*}
& & \int_1^\infty\frac{1}{R^{\rho(R)+\alpha/2 }}\int_1^R
 \bigg(\frac{1}{x^2}-\frac{1}{R^2} \bigg)g_{-}(x)dxdR \\
&\leq& \int_1^\infty\frac{1}{R^{\rho(R)+\alpha/2 }}  \bigg[2\pi
m_{+}(R)+\int_1^R  \bigg(\frac{1}{x^2}-\frac{1}{R^2}
\bigg)g_{+}(x)dx-
2\pi  \bigg(c_1+\frac{c_2}{R^2} \bigg) \bigg] dR \\
&<& \infty. \hspace{73mm} (9.3.5)
\end{eqnarray*}

  $\forall \alpha >0$, set
$$
I(\alpha)=\lim_{x\rightarrow\infty}
\frac{\int_x^\infty\frac{1}{R^{\rho(R)+\alpha/2 }}
\big(\frac{1}{x^2}-\frac{1}{R^2} \big)dR }{x^{-\big[\rho(x)+\alpha
+1\big]}},
$$
by the L'hospital's rule and (9.1.1), we have
$$
I(\alpha)=+\infty.
$$
Therefore, there exists $ \varepsilon_1>0$, such that for any $
x\geq 1$,
$$
\int_x^\infty\frac{1}{R^{\rho(R)+\alpha/2
}}\bigg(\frac{1}{x^2}-\frac{1}{R^2}\bigg)dR \geq  \frac{
\varepsilon_1}{x^{\rho(x)+\alpha +1}}.
$$
Multiplying this by $g_{-}(x)$ and integrating with respect to $x$,
we can obtain by (9.3.4) and (9.3.5)
$$
\varepsilon_1\int_1^\infty \frac{ g_{-}(x)}{x^{\rho(x)+\alpha +1}}dx
\leq \int_1^\infty
   g_{-}(x)\int_x^\infty\frac{1}{R^{\rho(R)+\alpha/2}}
   \bigg(\frac{1}{x^2}-\frac{1}{R^2}\bigg)dRdx
 <\infty.
$$
Thus
$$
\int_1^\infty \frac{ g_{-}(x)}{x^{\rho(x)+\alpha +1}}dx<\infty.
$$
by (9.3.3), we have
$$
\int_1^\infty \frac{ g_{+}(x)}{x^{\rho(x)+\alpha +1}}dx<\infty.
$$
Hence (1) holds.

 $\forall \alpha >0, R>1$, $\exists M(R)>0$, such that
for any $k>k_R=[2R]+1$, we have
$$
\frac{R^{\rho(k+1)+\alpha+1 } }{k^{\alpha/2}}\leq M(R),
$$
so $\forall \alpha >0, R>1$, if $|z|\leq R, k>k_R=[2R]+1$, then
$|t|\geq 2|z|$ and
\begin{eqnarray*}
& &\sum_{k=k_R}^\infty\int_{k\leq|t|<k+1}
\frac{|z|^{[\rho(|t|)+\alpha]+1} }
{|t|^{[\rho(|t|)+\alpha]+2 }}|u(t)|dt \\
&\leq& \sum_{k=k_R}^\infty\frac{R^{\rho(k+1)+\alpha+1 }
}{k^{\alpha/2}}\int_{k\leq|t|<k+1}\frac{2|u(t)|}{1+|t|^{\rho(|t|)+\alpha/2+1
}}dt \\
&\leq& 2M(R)\int_{|t|\geq
k_R}\frac{|u(t)|}{1+|t|^{\rho(|t|)+\alpha/2+1 }}dt.
\end{eqnarray*}
So the integral is absolutely convergent.

   To verify the boundary behavior of $ u_{{\bf C}_+}(z)$, choose a
large $T>2$, and write
\begin{eqnarray*}
u_{{\bf C}_+}(z)
&= &\int_{|t|\leq 2T}P(z,t)u(t)dt \\
& & -\Im\sum_{k=0}^{[\rho(|t|+\alpha)]}\int_{1< |t|\leq
2T}\frac{z^k}{\pi t^{k+1}}u(t)dt\\
& & +\int_{|t|> 2T}P_{[\rho(|t|+\alpha)]}(z,t)u(t)dt\\
& =& X(z)-Y(z)+Z(z).
\end{eqnarray*}

  Consider $ z\rightarrow x_0$, the first term $X(z)$ approaches $
u(x_0)$ because it is the Poisson integral of $
u(t)\chi_{[-2T,2T]}(t)$, where $\chi_{[-2T,2T]}$ is the
characteristic function of the interval $[-2T,2T]$; the second term
$Y(z)$ is a polynomial times $y$ and tends to $0$; and the third
term $Z(z)$ is $O(y)$ and therefore also to $0$. So the function
$u_{{\bf C}_+}(z) $ can be continuously extended to $\overline{{\bf
C}_+} $ such that $u_{{\bf C}_+}(t)=u(t)$; consequently, $
u(z)-u_{{\bf C}_+}(z) $ is harmonic in ${\bf C}_+$ and can be
continuously extended to
 $\overline{{\bf C}_+} $ with  $ 0$ in the boundary ${\bf R}$
 of ${\bf C}_+ $.
  The Schwarz reflection  principle (\cite{ABR}, p.68 and \cite{GT}, p.28)
  applied to $u(z)-u_{{\bf C}_+}(z)$
shows that there exists an entire function $Q_{{\bf C}_+}(z)$ which
satisfies $Q_{{\bf C}_+}(\overline{z})=\overline{Q_{{\bf C}_+}(z)}$,
such that $\Im Q_{{\bf C}_+}(z)=u(z)-u_{{\bf C}_+}(z)$ for $ z\in
\overline{{\bf C}_+}$. Therefore, if $ \alpha
>0$, we obtain $u(z)=\Im Q_{{\bf C}_+}(z)+u_{{\bf C}_+}(z)$ for
all $ z \in \overline{{\bf C}_+}$ and $\Im Q_{{\bf C}_+}(x)=0$ for
all $ x \in {\bf R}$. This completes the proof of Theorem.

\newpage
\ \thispagestyle{empty}
\newpage
\chapter{Integral Representations of Harmonic Functions in the Half Space}

\section{Introduction and Main Theorem}

\vspace{0.3cm}
  Let $\rho(R)\geq1$ is nondecreasing in $[0, +\infty)$ satisfying
$$
\varepsilon_0=\limsup_{R\rightarrow\infty}\frac{\rho'(R)R\log
R}{\rho(R)}<1. \eqno{(10.1.1)}
$$
For any real number  $ \alpha > 0 $, we denote by $(LU)_{\alpha }$
the space of all measurable  functions $f(x)$ in the  upper half
space $H$ which satisfy the following inequality:
$$
\int_{H}\frac{x_n|f(x)|dx}{1+|x|^{\rho(|x|)+n+\alpha +1}}
<\infty;\eqno{(10.1.2)}
$$
and $(LV)_{\alpha} $ the set of all measurable functions $g(x') $ in
${\bf R}^{n-1} $ which satisfy the following inequality:
$$
 \int_{\partial H}\frac{ |g(x')|dx'}{1+|x'|^{\rho(|x'|)+n+\alpha -1}} <\infty
.\eqno{(10.1.3)}
$$
We also denote by $(CH)_{\alpha }$ the set  of all continuous
functions $u(x)$ in the closed upper half space $\overline{H}$,
harmonic in the open upper half space $H$ with the positive part
$u^{+}(x)=\max\{u(x),0\}\in (LU)_{\alpha}$ and $ u^+(x')=u^+(x',0)
\in (LV)_{\alpha }$.

  The Poisson kernel for the upper half space $H $ is the
 function
$$
 P(x,y')=
\frac{2x_n}{\omega_n|x-y'|^n},
$$
where $ \omega_n = \frac{2\pi ^{\frac{n}{2}}}{\Gamma (\frac{n}{2})}$
is the area of the unit sphere in ${\bf R}^n$.

If $u(x)\leq 0  $ is harmonic in the open upper half space $H$,
continuous in the closed upper half space $\overline{H}$, then (see
\cite{HN}, \cite{SH} and \cite{ABR}) $ u\in (CH)_{\alpha}$ for each
$\alpha
> 0 $ and there exists a  constant $c\leq 0 $ such that
$$
u(x)=cx_n+  \int_{\partial H}P(x,y')u(y')dy' \eqno{(10.1.4)}
 $$
for all $ x \in  H $, the integral in (10.1.4) is absolutely
convergent. Motivated by this result, we will prove that if $u\in
(CH)_{\alpha }$, then $ u^{+}(x) \in (LU)_{\alpha }, u(x')\in
(LV)_{\alpha} $ and  a similar representation to (10.1.4) for the
function $ u \in (CH)_{\alpha }$ holds by modifying the Poisson
kernel $
 P_m(x,y')$.  It is well known (see
\cite{G}, \cite{HN} and \cite{SW}) that
  the Poisson kernel $
 P(x,y')$
is harmonic in $ x\in {\bf R}^n-\{y'\}$ and has a series expansion
in terms of the ultraspherical ( or Gegenbauer ) polynomials
 $ C^{\lambda}_{k}(t)\ (\lambda =\frac{n}{2})$.
The latter can be defined in terms of a generating function
$$
(1-2tr+r^2)^{-\lambda } = \sum _{k=0}^{\infty} C_{k}^{\lambda}(t)
r^k,
$$
where $|r|<1$, $ |t|\leq 1$ and $ \lambda > 0$. The coefficients $
C^{\lambda}_{k}(t) $ is called the ultraspherical ( or Gegenbauer )
polynomial of degree $ k $  associated with $ \lambda $, the
function $ C^{\lambda}_{k}(t) $ is a  polynomial of degree $ k $ in
$ t $ and satisfies the inequality (\cite{G}, p.82 and p.92)
$$
|C_k^{\lambda }(t)|\leq C_k^{\lambda }(1)=\frac{\Gamma (2\lambda
+k)}{\Gamma (2\lambda )\Gamma(k+1)}, \ \ |t|\leq 1 .
$$
Therefore,   a series expansion of the Poisson kernel $
 P(x,y')$ in terms of
the ultraspherical polynomials  $ C^{\lambda}_{k}(t)$ is
$$
 P(x,y')=\sum_{k=0}^{\infty}\frac{2x_n|x|^k}{\omega_n|y'|^{n+k}}C^{n/2}_{k}
 \left( \frac{x\cdot y'}{|x||y'|}\right ),
 $$
this series converges for $ |x|<|y'|$, each term is homogeneous in x
of degree $k+1$. Differentiating termwise in x gives
$$
\triangle_x P(x,y')=\sum_{k=0}^{\infty}\triangle_x \left
(\frac{2x_n|x|^k}{\omega_n|y'|^{n+k}}C^{n/2}_{k}
 \left( \frac{x\cdot y'}{|x||y'|}\right )\right ).
 $$
Each term $\triangle_x \left
(\frac{2x_n|x|^k}{\omega_n|y'|^{n+k}}C^{n/2}_{k}
 \left( \frac{x\cdot y'}{|x||y'|}\right )\right )$ is homogeneous in $ x$ of degree $
 k-1$, hence by the linear independence of homogenous
 functions, $\frac{x_n|x|^k}{\omega_n|y'|^{n+k}}C^{n/2}_{k}
 \left( \frac{x\cdot y'}{|x||y'|}\right )$ is harmonic on ${\bf
 R}^n$ for each $ k\ge 0 $.
 If $m\geq 0$ is an integer, we define a modified
 Poisson kernel of order $ m$ for $x\in H $ by
 $$
 P_m(x,y')=\left\{\begin{array}{ll}
 P(x,y') , &   \mbox{when }   |y'|\leq 1  ,\\
 P(x,y') - \sum_{k=0}^{m-1}\frac{2x_n|x|^k}{\omega_n|y'|^{n+k}}C^{n/2}_{k}
 \left( \frac{x\cdot y'}{|x||y'|}\right ),&
\mbox{when}\   |y'|> 1.
 \end{array}\right.\eqno{(10.1.5)}
  $$
The modified
 Poisson kernel $
 P_m(x,y')$ is harmonic in $ x\in H $.

  Up to now, a number of results about integral representations have
been achieved in \cite{DM}, \cite{ZM}, in this chapter, we will
establish the following theorem.

 \vspace{0.2cm}
 \noindent
{\bf Theorem 10.1.1} {\it  If  $ u\in (CH)_{\alpha}\ (\alpha > 0) $
, then the following properties hold:

\noindent {\rm (1)}\
$$
 \int_{\partial H}\frac{ |u(x')|}{1+|x'|^{\rho(|x'|)+n+\alpha -1}}dx' <\infty
;
$$

 \noindent
 {\rm (2)}\   the integral
$$
\int_{\partial H}P_{[\rho(|y'|)+\alpha]}(x,y')u(y')dy'
 $$
is absolutely convergent, it represents a harmonic function
$u_{H}(x)$ in $H$ and can be continuously extended to $\overline{H}
$ such that $u_{H}(y')=u(y')$ for $ y'\in
\partial H$;

\noindent
 {\rm (3)}\   There exists  a harmonic function $h(x)$  which
  vanishes on the boundary $\partial H$ such that $u(x)=h(x)+u_H(x)$ for all $ x \in \overline{H}$.
}

\section{ Main Lemma}

  In order to obtain the result, we need the lemma below:

\vspace{0.2cm}
 \noindent
{\bf Lemma 10.2.1 } {\it For any $ y'\in \partial H$ and $|x|>1,
x_n>0 $, the
 following inequalities
$$
 |P_m(x,y')|\leq \left\{\begin{array}{ll}
 \frac{(2^{m+n}+m2^{m}C^{n/2}_{m-1}
 \left( 1 \right ))|x|^{n+m-1}}
 {\omega_nx_n^{n-1}|y'|^{n+m-1}} , &   \mbox{when }
 1<|y'|\leq 2|x| ,  \\
\frac{2^{m+n+1}x_n|x|^{m}}{\omega_n|y'|^{n+m}}, & \mbox{when }
 |y'|> \max \{1,2|x|\}, \\
\frac{2}{\omega_nx_n^{n-1}}, & \mbox{when }
 |y'|\leq 1
 \end{array}\right.
$$
hold. }

\vspace{0.2cm}

 Proof:  When $ y'\in
\partial H , |y'|\leq 1 $, we have $ |x-y'|
\geq |x_n| $ and so
$$
|P_m(x,y')|\leq \frac{2}{\omega_nx_n^{n-1}};
$$
when $ y'\in
\partial H , 1<|y'|\leq 2|x|$, we also have $ |x-y'|
\geq |x_n| $ and so by (10.1.5)
\begin{eqnarray*}
|P_m(x,y')|
&\leq& \frac{2x_n}{\omega_n|x-y'|^n}+
\sum_{k=0}^{m-1}\frac{2x_n|x|^k}{\omega_n|y'|^{n+k}}C^{n/2}_{k}
 \left( 1 \right ) \\
&\leq&\frac{2}{\omega_nx_n^{n-1}} \left
(1+\frac{x_n^n|x|^k}{|y'|^{n+k}}C^{n/2}_{m-1}
 ( 1)\right ) \\
&\leq&  \frac{\big(2^{m+n}+m2^{m}C^{n/2}_{m-1}
 \left( 1 \right )\big)|x|^{n+m-1}}
 {\omega_nx_n^{n-1}|y'|^{n+m-1}};
\end{eqnarray*}
when $|y'|> \max \{1,2|x|\}$, we have by (10.1.5)
\begin{eqnarray*}
|P_m(x,y')|
&=&
\bigg|\sum_{k=m}^{\infty}\frac{2x_n|x|^k}{\omega_n|y'|^{n+k}}C^{n/2}_{k}
 \left( \frac{x\cdot y'}{|x||y'|}\right )\bigg|\\
&\leq& \frac{2}{\omega_n} \sum _{k=m}^{\infty
}\frac{x_n|x|^k}{|y'|^{n+k}}C^{n/2}_{k}
 ( 1) \\
&\leq& \frac{2^{m+n+1}x_n|x|^{m}}{\omega_n|y'|^{n+m}}.
\end{eqnarray*}
This proves the inequalities.

\section{ Proof of Theorem}

 If  $ u\in (CH)_{\alpha}\
(\alpha > 0) $ , suppose $ R>1$, then by the Carleman formula for
harmonic functions in the upper half space,
\begin{eqnarray*}
& &  \int_{\{x\in {\bf
R}^{n}:\;|x|=R,x_n>0\}}u(x)\frac{nx_n}{R^{n+1}}d\sigma(x)\\
&+ & \int_{\{x\in {\bf
R}^{n}:\;1<|x'|<R,x_n=0\}}u(x')\bigg(\frac{1}{|x'|^{n}}-\frac{1}{R^{n}}\bigg)dx'
=c_1+\frac{c_2}{R^{n}}, \\
\end{eqnarray*}
where
\begin{eqnarray*}
&c_1& =\int_{\{x\in {\bf
R}^{n}:\;|x|=1,x_n>0\}}\bigg[(n-1)x_nu(x)+x_n\frac{\partial
u(x)}{\partial n}\bigg]d\sigma(x),\\
&c_2&=\int_{\{x\in {\bf
R}^{n}:\;|x|=1,x_n>0\}}\bigg[x_nu(x)-x_n\frac{\partial
u(x)}{\partial
n}\bigg]d\sigma(x). \\
\end{eqnarray*}
Set
$$
m_{+}(R)=\int_{\{x\in {\bf
R}^{n}:\;|x|=R,x_n>0\}}u^{+}(x)\frac{nx_n}{R^{n+1}}d\sigma(x),
$$
$$
m_{-}(R)=\int_{\{x\in {\bf
R}^{n}:\;|x|=R,x_n>0\}}u^{-}(x)\frac{nx_n}{R^{n+1}}d\sigma(x),
$$
then
\begin{eqnarray*}
& & m_{-}(R)+\int_{\{x\in {\bf R}^{n}:\;1<|x'|<R,x_n=0\}}u^{-}(x')
\bigg(\frac{1}{|x'|^{n}}-\frac{1}{R^{n}}\bigg)dx' \\
&= &m_{+}(R)+\int_{\{x\in {\bf
R}^{n}:\;1<|x'|<R,x_n=0\}}u^{+}(x')\bigg(\frac{1}{|x'|^{n}}-\frac{1}{R^{n}}\bigg)dx'-
c_1-\frac{c_2}{R^{n}}, \hspace{3mm} (10.3.1)
\end{eqnarray*}
where $u^{+}(x)=\max\{u(x), 0\}, u^{-}(x)=(-u(x))^{+}$ and
$u(x)=u^{+}(x)-u^{-}(x)$.

Since $u\in (CH)_{\alpha}$, we obtain by (10.1.2)
$$
\int_1^\infty\frac{m_{+}(R)}{R^{\rho(R)+\alpha }}dR =
 n\int_{D}\frac{x_nu^{+}(x)}{|x|^{\rho(|x|)+n+\alpha +1}}dx <\infty ,\eqno{(10.3.2)}
$$
where $D=\{x\in H :\;|x|>1\}$.

By (10.1.3), we can also obtain
\begin{eqnarray*}
& & \int_1^\infty\frac{1}{R^{\rho(R)+\alpha }}\int_{\{x\in {\bf
R}^{n}:\;1<|x'|<R,x_n=0\}}u^{+}(x')\bigg(\frac{1}{|x'|^{n}}-\frac{1}{R^{n}}\bigg)dx'dR \\
& =& \int_{|x'|\geq
1}u^{+}(x')\int_{|x'|}^\infty\frac{1}{R^{\rho(R)+\alpha
}}\bigg(\frac{1}{|x'|^{n}}-\frac{1}{R^{n}}\bigg)dRdx' \\
&\leq &  \frac{n}{n+1}\int_{|x'|\geq 1}\frac{
u^{+}(x')}{|x'|^{\rho(|x'|)+n+\alpha -1}}dx'<\infty.
\hspace{50mm}(10.3.3)
\end{eqnarray*}
Similarly, we have
\begin{eqnarray*}
& & \int_1^\infty\frac{1}{R^{\rho(R)+\alpha/2 }}\int_{\{x\in {\bf
R}^{n}:\;1<|x'|<R,x_n=0\}}u^{-}(x')
\bigg(\frac{1}{|x'|^{n}}-\frac{1}{R^{n}}\bigg)dx'dR\\
& =& \int_{|x'|\geq
1}u^{-}(x')\int_{|x'|}^\infty\frac{1}{R^{\rho(R)+\alpha/2
}}\bigg(\frac{1}{|x'|^{n}}-\frac{1}{R^{n}}\bigg)dRdx'.\hspace{50mm}(10.3.4)
\end{eqnarray*}
So we have by (10.3.1), (10.3.2) and (10.3.3)
\begin{eqnarray*}
& & \int_1^\infty\frac{1}{R^{\rho(R)+\alpha/2 }}\int_{\{x\in {\bf
R}^{n}:\;1<|x'|<R,x_n=0\}}u^{-}(x')\bigg(\frac{1}{|x'|^{n}}-\frac{1}{R^{n}}\bigg)dx'dR \\
&\leq& \int_1^\infty\frac{1}{R^{\rho(R)+\alpha/2 }}
m_{+}(R) dR \\
& & +\int_1^\infty\frac{1}{R^{\rho(R)+\alpha/2 }} \bigg[\int_{\{x\in
{\bf
R}^{n}:\;1<|x'|<R,x_n=0\}}u^{+}(x')\bigg(\frac{1}{|x'|^{n}}-\frac{1}{R^{n}}\bigg)dx'\bigg] dR \\
& & -\int_1^\infty\frac{1}{R^{\rho(R)+\alpha/2 }} \bigg(
c_1+\frac{c_2}{R^{n}}\bigg) dR < \infty.\hspace{50mm}(10.3.5)
\end{eqnarray*}

  $\forall \alpha >0$, set
$$
I(\alpha)=\lim_{|x'|\rightarrow\infty}
\frac{\int_{|x'|}^\infty\frac{1}{R^{\rho(R)+\alpha/2
}}\bigg(\frac{1}{|x'|^{n}}-\frac{1}{R^{n}}\bigg)dR
}{|x'|^{-\big(\rho(|x'|)+n+\alpha -1\big)}},
$$
by the L'hospital's rule and (10.1.1), we have
$$
I(\alpha)=+\infty.
$$
Therefore, there exists $ \varepsilon_1>0$, such that for any $
|x'|\geq 1$,
$$
\int_{|x'|}^\infty\frac{1}{R^{\rho(R)+\alpha/2
}}\bigg(\frac{1}{|x'|^{n}}-\frac{1}{R^{n}}\bigg)dR \geq  \frac{
\varepsilon_1}{|x'|^{\rho(|x'|)+n+\alpha -1}}.
$$
Multiplying this by $u^{-}(x')$ and integrating with respect to
$x'$, we can obtain by (10.3.4) and (10.3.5)
\begin{eqnarray*}
& & \varepsilon_1\int_{|x'|\geq 1}\frac{
u^{-}(x')}{|x'|^{\rho(|x'|)+n+\alpha -1}}dx' \\
&\leq& \int_{|x'|\geq
1}u^{-}(x')\int_{|x'|}^\infty\frac{1}{R^{\rho(R)+\alpha/2
}}\bigg(\frac{1}{|x'|^{n}}-\frac{1}{R^{n}}\bigg)dRdx' < \infty.
\end{eqnarray*}
Thus
$$
\int_{|x'|\geq 1}\frac{ u^{-}(x')}{|x'|^{\rho(|x'|)+n+\alpha
-1}}dx'<\infty.
$$
by (10.3.3), we have
$$
\int_{|x'|\geq 1}\frac{ u^{+}(x')}{|x'|^{\rho(|x'|)+n+\alpha
-1}}dx'<\infty.
$$
Hence (1) holds.

 $\forall \alpha >0, R>1$, $\exists M(R)>0$, such that
for any $k>k_R=[2R]+1$, we have
$$
\frac{(2R)^{\rho(k+1)+\alpha+1 } }{k^{\alpha/2}}\leq M(R),
$$
so $\forall \alpha >0, R>1$, if $|x|\leq R, k>k_R=[2R]+1$, then
$|y'|\geq 2|x|$ and
\begin{eqnarray*}
& &
\sum_{k=k_R}^\infty\int_{k\leq|y'|<k+1}\frac{(2|x|)^{[\rho(|y'|)+\alpha]+1
} }{|y'|^{[\rho(|y'|)+\alpha]+n }}|u(y')|dy' \\
&\leq& \sum_{k=k_R}^\infty\frac{(2R)^{\rho(k+1)+\alpha+1 }
}{k^{\alpha/2}}\int_{k\leq|y'|<k+1}\frac{2|u(y')|}{1+|y'|^{\rho(|y'|)+\alpha/2+(n-1)
}}dy' \\
&\leq& 2M(R)\int_{|y'|\geq
k_R}\frac{|u(y')|}{1+|y'|^{\rho(|y'|)+\alpha/2+(n-1) }}dy'. \\
\end{eqnarray*}
So the integral is absolutely convergent.

  To verify the boundary behavior of $ u_H(x)$, fix a boundary point
$a'=(a_1,a_2,\cdots,\\a_{n-1})\in {\bf R}^{n-1}$, choose a large
$T>|a'|+1$, and write
\begin{eqnarray*}
u_H(x)
&=& \int_{|y'|\leq T}P(x,y')u(y')dy' \\
& &
-\sum_{k=0}^{[\rho(|y'|+\alpha)]-1}\frac{2x_n|x|^k}{\omega_n}\int_{1<
|y'|\leq T}\frac{1}{|y'|^{n+k}}C^{n/2}_{k}\left( \frac{x'\cdot
y'}{|x||y'|}\right)u(y')dy' \\
& & +\int_{|y'|> T}P_{[\rho(|y'|+\alpha)]}(x,y')u(y')dy'\\
& = & X(x)-Y(x)+Z(x). \\
\end{eqnarray*}

  Consider $ x\rightarrow a'$, the first term $X(x)$ approaches $
u(a')$ because it is the Poisson integral of $
u(y')\chi_{B(T)}(y')$, where $\chi_{B(T)}$ is the characteristic
function of the ball $B(T)=\{ y'\in {\bf R}^{n-1}:|y'|\leq T\}$; the
second term $Y(x)$ is a polynomial times $x_n$ and tends to $0$; and
the third term $Z(x)$is $O(x_n)$ and therefore also to $0$. So the
function $u_H(x) $ can be continuously extended to $\overline{H}$
such that $u_H(y')=u(y')$; consequently, $ u(x)-u_H(x) $ is harmonic
in $H$ and can be continuously extended to
 $\overline{H} $ with  $ 0$ in the boundary $ \partial H $
 of $ H $.
  The Schwarz reflection  principle
(\cite{ABR}, p.68 and \cite{GT}, p.28) applied to $u(x)-u_H(x)$
shows that there exists a harmonic  function $h(x)$ in ${\bf R}^n$
such that $h(x^*)=-h(x)=-(u(x)-u_H(x))$ for $ x\in \overline{H}$,
where $ x^*=(x', -x_n)$ is the reflection of $x$ in $ \partial H$.
Therefore, if $ \alpha >0$,  $ h(x)$ is a harmonic function  which
vanishes on the boundary $\partial H $ such that $u(x)=h(x)+u_H(x)$
for all $ x \in \overline{H} $. This completes the proof of Theorem.

\end{document}